\documentclass[10pt]{article}
\usepackage[english]{babel}

\usepackage{authblk}
\usepackage{color}
\usepackage{amsmath, amssymb, mathdots, accents, latexsym}
\usepackage[a4paper,left=3.5cm,right=4cm]{geometry}
\usepackage{ifthen}
\usepackage{amsthm}
\usepackage{caption,subcaption}
\usepackage{dsfont}
\usepackage{thmtools}
\usepackage{hyperref}
\usepackage{extarrows}

\usepackage{fancyhdr,enumerate}

\newtheorem{theorem}{Theorem}[section]
\newtheorem{lemma}[theorem]{Lemma}
\newtheorem{definition}[theorem]{Definition}
\newtheorem{remark}[theorem]{Remark}
\newtheorem{corollary}[theorem]{Corollary}
\newtheorem{example}[theorem]{Example}
\newtheorem{proposition}[theorem]{Proposition}

\newtheorem{openproblem}[theorem]{Open Problem}
\newtheorem{statement}[theorem]{Statement}
 
\usepackage{tikz, tikz-3dplot}
\usetikzlibrary{shapes.geometric, arrows}
\usepackage{pgfplotstable}
\usepackage{pgfplots}
\usepackage{pgffor,pgfmath}
\usepackage{booktabs}
\pgfplotsset{compat=1.14}

 \tikzstyle{startstop} = [rectangle, rounded corners, minimum width=1cm, minimum height=1cm,text centered, draw=black, fill=red!0]
\tikzstyle{io1} = [rectangle, trapezium left angle=80, trapezium right angle=100, minimum width=1cm, minimum height=1cm, text centered, draw=black, fill=blue!0]
\tikzstyle{io2} = [trapezium,  rounded corners, trapezium left angle=100, trapezium right angle=100, minimum width=1cm, minimum height=1cm, text centered, draw=black, fill=yellow!0]
\tikzstyle{process} = [rectangle, minimum width=1cm, minimum height=1cm, text centered, draw=black, fill=orange!0]
\tikzstyle{decision} = [circle, minimum width=1cm, minimum height=1cm, text centered, draw=black, fill=green!0]
\tikzstyle{decision2} = [ellipse, rounded corners=10mm, minimum width=2cm, minimum height=2cm, text centered, draw=black, fill=green!0]
\tikzstyle{arrow} = [thick,->,>=stealth]

\usepackage[sort]{natbib}
\usepackage{mathtools}
\usepackage{nameref,cleveref}
\usepackage{blindtext}
\usepackage{cancel}
\usepackage{autonum} 

\usepackage{titlesec}
\usepackage{enumitem}

\usepackage{booktabs}

\newcommand{\R}{\mathbb{R}}

\newcommand{\N}{\mathbb{N}}

\newcommand{\Gc}{\mathcal{G}}

\newcommand{\boundary}{B}

\newcommand{\Dc}{\mathcal{D}}

\newcommand{\Ind}{\mathds{1}}

\newcommand{\Fc}{\mathcal{F}}

\newcommand{\Nc}{\mathcal{N}}

\newcommand{\SNc}{\mathcal{SN}}
\newcommand{\WNc}{\mathcal{WN}}

\newcommand{\PNc}{\mathcal{PN}}

\newcommand{\Opc}[2][]
{%
  \ifthenelse{\equal{#1}{}}
  {\left\{#2\right\}}
  {\left\{#2\right\}_{#1}}
}

\newcommand{\locgrad}[2][]
{%
  \ifthenelse{\equal{#1}{}}
  {\Opc{\grad #2}}
  {\Opc[#1]{\grad #2}}
}

\newcommand{\length}{\mathrm{length}}

\newcommand{\internalnodes}{\nodeset\setminus\boundary}

\newcommand{\Hc}{\mathcal{H}}

\newcommand{\rayl}{\mathcal{R}}

\newcommand{\eigenfunction}{f}

\newcommand{\grad}{K}

\newcommand{\divg}{-\frac{1}{2}K^T}

\newcommand{\plap}{\Delta_p}

\newcommand{\inflap}{\Delta_{\infty}}

\newcommand{\lap}[1][]{
  \ifthenelse{\equal{#1}{}}
  {\Delta}
  {\Delta_{#1}}}

\newcommand{\edgeset}{E}

\newcommand{\nodeset}{V}

\newcommand{\nodemaxset}{V_{max}}

\newcommand{\edgelength}{\omega}

\newcommand{\edgeweight}[1][]{
  \ifthenelse{\equal{#1}{}}
  {\mu}
  {\mu_{#1}}}
\newcommand{\nodeweight}[1][]{
  \ifthenelse{\equal{#1}{}}
  {\nu}
  {\nu_{#1}}}

\newcommand{\genus}{\gamma}

\DeclareMathOperator*{\argmax}{arg\:max}

\newcommand{\La}{\Lambda}

\newcommand{\dist}{\mathrm{dist}}

\makeatletter
\def\wbar{\accentset{{\cc@style\underline{\mskip8mu}}}}

\makeatother

\renewcommand{\vec}[1]{\mbox{\boldmath \small $#1$}}

\newcommand{\diag}{\mathrm{diag}}

\newcommand{\supp}{\ensuremath{\mathrm{supp}}}

\newcommand{\matchingN}{\theta(\Gc)}
\newcommand{\betti}{\beta(\Gc)}

\DeclareRobustCommand{\rchi}{{\mathpalette\irchi\relax}}
\newcommand{\irchi}[2]{\raisebox{\depth}{$#1\chi$}}

\newcommand{\morse}[1][]{
  \ifthenelse{\equal{#1}{}}
  {\mathcal{MI}}
  {\mathcal{MI}_{#1}}}

\newcommand{\PD}[2][]{
  \ifthenelse{\equal{#1}{}}
  {\textcolor{red}{\textbf{PD: #2}\marginpar{\textcolor{red}{PD}}}}
  {\textcolor{red}{\textbf{PD: #2}\marginpar{\textcolor{red}{PD: #1}}}}
  }

  \newcommand{\DZ}[2][]{
  \ifthenelse{\equal{#1}{}}
  {\textcolor{purple}{\textbf{DZ: #2}\marginpar{\textcolor{purple}{DZ}}}}
  {\textcolor{purple}{\textbf{DZ: #2}\marginpar{\textcolor{purple}{DZ: #1}}}}
  }

\title{Nonlinear spectral graph theory}

\author[1]{Piero Deidda \thanks{deidda@mathc.rwth-aachen.de}}
\author[2,3]{ Francesco Tudisco \thanks{f.tudisco@ed.ac.uk}}
\author[4]{Dong Zhang \thanks{dongzhang@math.pku.edu.cn}}
\affil[1]{\small Department of Mathematics, RWTH Aachen, Pontdriesch 10-12,
52062 Aachen, Germany}
\affil[2]{\small School of Mathematics and Maxwell Institute for Mathematical Sciences, University of Edinburgh,  Peter Guthrie Tait Road, EH9 3FD, Edinburgh,
UK.}
\affil[3]{\small Gran Sasso Science Institute, Viale F. Crispi, 7, 67100 L'Aquila, Italy.}
\affil[4]{\small LMAM and School of Mathematical Sciences, Peking University, Beijing 100871, China.}
\date{}

\begin{document}

\maketitle

\begin{abstract}
Nonlinear spectral graph theory is an extension of the traditional (linear) spectral graph theory and studies relationships between spectral properties of nonlinear operators defined on a graph and topological properties of the graph itself. Many of these relationships get tighter when going from the linear to the nonlinear case. 
In this manuscript, we discuss the spectral theory of the graph $p$-Laplacian operator. In particular, we report links between the $p$-Laplacian spectrum and higher-order Cheeger (or isoperimetric) constants, sphere packing constants, independence, and matching numbers of the graph.  
The main aim of this paper is to present a complete and self-contained introduction to the problem, accompanied by a discussion of the main results and the proof of new results that fill some gaps in the theory. The majority of the new results are devoted to the study of the graph $\infty$-Laplacian spectrum and the information that it yields about the packing radii, the independence numbers, and the matching number of the graph. This is accompanied by a novel discussion about the nodal domains induced by the $\infty$-eigenfunctions. There are also new results about the variational spectrum of the $p$-Laplacian, the regularity of the $p$-Laplacian spectrum varying $p$, and the relations between the $1$-Laplacian spectrum and new Cheeger constants.   

\bigskip 


\bigskip 

\end{abstract}

\tableofcontents

\section{Introduction}

\subsection{Motivation, nonlinear Rayleigh quotients, and the graph $p$-Laplacian}
spectral graph theory has a long history. It studies the relationships between the topology of a graph and the spectrum of linear operators associated with it, such as the adjacency or Laplacian matrices~\cite{Stadler_book-Laplacian-eigenvectors,chung1997spectral,spielman2019spectralandalgebraic,trevisan2017lecture}. These matrices encode structural properties of the graph, which are reflected in their spectral characteristics. Recently, it has been observed that replacing linear operators with nonlinear ones can enhance results and address limitations in several applications. Examples include semisupervised learning~\cite{slepcev2019analysis,tudisco2021nonlinear}, signal processing and variational filtering~\cite{rudin1992nonlinear,osher2005iterative,Elmoataz1}, and spectral graph theory itself~\cite{Zhang_1-Lap_cheeger_cut,Bhuler,hein2010inverse,Tudisco1}. At the same time, nonlinearities make both the analysis and the algorithmic aspects significantly more challenging.

This manuscript focuses on the emerging field of nonlinear spectral graph theory. Nonlinear eigenvalues and eigenfunctions are often defined as critical values and points of a Rayleigh quotient $\rayl(f)$, expressed as the ratio of two homogeneous functionals. A classical example is the Rayleigh quotient associated with a matrix $A$ in $\ell_p$ and $\ell_q$ spaces:
\[
\rayl(f) = \frac{\|Af\|_p}{\|f\|_q}.
\]
Nonlinear eigenpairs arising from such Rayleigh quotients have diverse applications. In variational filtering, for instance, the term $\|Af\|_p$ acts as a regularization functional~\cite{burger2016spectral,BungertNonlineardecomp,gilboa2018nonlinear,Burger23}. Regularised inverse problems, such as
\[
t \mapsto \min_f \|f - g\|_2^2 + t J(f),
\]
where $J(f)$ is a positive, homogeneous functional, have been shown to yield spectral decompositions of input signals. When $J(f) = \|Af\|_1$ for a matrix $A$ satisfying certain conditions (e.g., diagonal dominance of $AA^*$), this decomposition corresponds to the eigenfunctions of $J$ in $\ell_2$, i.e., the critical points of the Rayleigh quotient:
\[
\rayl(f) = \frac{J(f)}{\|f\|_2}.
\]

In spectral graph theory, applications of nonlinear Rayleigh quotients include data clustering and community detection~\cite{Bhuler,chang2016spectrum,hein2010inverse,tudisco2018community,tudisco2022nonlinear}, graph partitioning~\cite{deidda2024_inf_eigenproblem}, and core-periphery detection~\cite{rombach2014core,Tudisco_core_periphery}. For the latter, $A$ is often defined as a tensor such that, for a function $f$ on the graph's nodes, $Af$ is the matrix whose rows capture the values of $f$ on the endpoints of each edge. The corresponding Rayleigh quotient is then:
\[
\rayl(f) = \frac{\|Af\|_{p,1}}{\|f\|_q} =  \frac 1 {\|f\|_q} \sup_{g\neq 0}\frac{\|(Af)g\|_p}{\|g\|_1},
\]
where $p$ is such that $p \to \infty$. The maximizer of this Rayleigh quotient, which solves a nonlinear eigenvalue problem, has been proposed as an effective centrality measure~\cite{Tudisco_core_periphery}. This approach distinguishes core nodes from peripheral ones and has been extended to hypergraphs and multiplex networks~\cite{Tudisco_core_periphery_hypergraphs,bergermann2023nonlinear}.
Finally, applications to data clustering and partitioning arise when $A$ is taken as the directed incidence matrix of the graph, denoted here by $\grad$, and $f$ is a function on the graph nodes. In this setting, $p$ and $q$ are usually assumed equal, with the most significant applications arising at the extremal values $p = 1$ and $p = \infty$. Specifically, for $p = 1$, each eigenvalue corresponds to the isoperimetric constant of a subgraph, while for $p = \infty$, each eigenvalue corresponds to the distance between certain node pairs in the graph.

The same properties hold in the infinite dimensional setting, i.e. when $f$ belongs to the Sobolev space $W^{2,p}(\Omega)$ for some open subset $\Omega\subset \R^N$, and we consider the $p$-Rayleigh quotient induced by the gradient operator, $\nabla f$~\cite{Lind2,Lind3,Lind4,Kawohl2003,parini2010second}. As we will elaborate later, the directed incidence matrix, when applied to a function $f$, provides a discrete approximation of the directional derivatives of $f$ along the graph's edges, and thus a discrete analog of the gradient in the continuous setting.

The nonlinear eigenproblem associated with the critical points of the $p$-Rayleigh quotient
\[
\rayl_p(f) = \frac{\|\grad f\|_p}{\|f\|_p},
\]
is known as the $p$-Laplacian eigenvalue problem, which is a central topic of this manuscript. This nomenclature is justified by the fact that, for $p = 2$, the problem reduces to the classical eigenvalue problem of the graph Laplacian matrix. The operator obtained by differentiating $\|\grad f\|_p^p/p$ with respect to $f$ is called the graph $p$-Laplacian operator, denoted $\plap(f)$. This operator serves as a graph counterpart to the differential $p$-Laplacian operator in the infinite dimensional setting, which is defined as the differential operator in the Euler--Lagrange equations of the functional $\|\nabla f\|_p^p/p$.

The graph $p$-Laplacian operator finds applications not only in spectral graph theory but also in image processing~\cite{Elmoataz1} and semisupervised learning~\cite{slepcev2019analysis,calder2018game,flores2022analysis}. 

Additionally, for a given pair of nodes $\{u,v\}$, minimizing $\|\grad f\|_p$ over functions satisfying $f(u) - f(v) = 1$ defines a nonlinear analogue of the classical graph resistance, known as the $p$-resistance $r_p(u,v)$~\cite{herbster2009predicting,Alagmir_2011_p-resistances}. The $p$-resistance is a distance metric on the graph~\cite{herbster2010triangle} and is tied to the graph's topological properties, much like the $p$-Laplacian spectrum. Specifically, as $p \to 1$ or $p \to \infty$, the $p$-resistance of an edge $(u,v)$ reduces to the $(u,v)$-cut of the graph and the shortest path distance between $u$ and $v$, respectively~\cite{Alagmir_2011_p-resistances}. Furthermore, the $p$-resistance is directly connected to semisupervised learning problems, as the function minimizing $\|\grad f\|_p$ under the constraints $f(u) = 0$ and $f(v) = 1$ grows on the unlabeled nodes in proportion to the $p$-resistance~\cite{saito23effective_p_resistance}.

Extensions of the $p$-Laplacian eigenvalue problem and its applications to clustering and partitioning include investigations on hypergraphs~\cite{Mulas2022_plap_hyper,saito2018hypergraph,saito23,Fazeny23,Fazeny24,Li}, simplicial complexes~\cite{jostzhang24+,Schaub20}, and signed graphs~\cite{ge2023new,zhang2023pLap_noddom}. These extensions also encompass the study of the signless $p$-Laplacian operator, which is derived by using the undirected incidence matrix~\cite{chang2016signless,Mulas2022_plap_hyper}. Such generalizations provide valuable tools for analyzing a wide range of network structures and their associated spectral properties.

\subsection{Spectrum, variational eigenvalues, and geometric limits}
Looking into the graph $p$-Laplacian eigenvalue problem in more detail, we note that letting $\partial_f$ indicate the derivative operator, the critical point condition $\partial_f\rayl_{p}(f)=0$ yields the equation: 
\begin{equation}\label{intro_p-Lap-eig_eq}
    \plap f=\lambda |f|^{p-2}f\,,
\end{equation}
where absolute value and powers are taken point-wise, 
and we say that $(f,\lambda)$ is a $p$-Laplacian eigenpair if it solves \eqref{intro_p-Lap-eig_eq}. Despite its natural definition, several questions arise when moving from the linear case of the Laplacian matrix to the nonlinear case with $p$ different from $2$. 

A first issue concerns the cardinality of the spectrum. In the linear case, the spectrum of the Laplacian is well-defined and well-understood: it is countable in the infinite-dimensional setting and finite, equal to the space's dimension, in the graph setting. For the $p$-Laplacian with $p \neq 2$, however, the countability and finiteness of the spectrum remain open problems. In particular, counterexamples in the graph setting show that the number of eigenvalues of the $p$-Laplacian can exceed the dimension of the space~\cite{Amghibech1,zhang2021homological,DEIDDA2023_nod_dom}.

Another complication arises from the behaviour of the eigenfunctions. In the linear case, the eigenfunctions form an orthogonal basis, and the algebraic and geometric multiplicities of each eigenvalue coincide; that is, the multiplicity of any eigenvalue equals the dimension of its eigenspace. These properties are fundamental for signal decomposition or approximation using eigenfunctions as frequency components. Unfortunately, these features are lost in the nonlinear $p$-Laplacian case. Specifically:
\begin{itemize}[leftmargin=*,noitemsep,topsep=0pt]
    \item The cardinality of the spectrum is unknown, making the concept of algebraic multiplicity not directly well-defined.
    \item Nonlinearity leads to eigenfunctions that are generally non-orthogonal.
    \item Multiple eigenfunctions can exist for the same eigenvalue, but they are not necessarily infinitely many~\cite{Amghibech1}. Consequently, the notion of geometric multiplicity does not directly transfer over to the nonlinear setting.
\end{itemize}
These challenges will be discussed in detail in \Cref{subsec:p-lap_spectrum}.

One common approach to address these difficulties is through the definition of variational eigenvalues. This involves constructing a family of eigenvalues whose cardinality aligns with the dimension of the space in the graph setting or is countable in the infinite-dimensional setting. These variational eigenvalues are typically defined via a $\min\max$ principle, analogous to the Fisher--Courant characterization of linear eigenvalues. However, the specific choice of families of sets over which the $\min\max$ is performed can yield different families of variational eigenvalues. Whether these families produce equivalent eigenvalues remains an open question in many cases.

The most standard family of variational eigenvalues is defined using Lusternik--Schnirelmann theory and the concept of Krasnoselskii genus~\cite{struwe,Ghoussoub,Fucik}. In this approach, a generalized notion of dimension, the genus, is applied to symmetric subsets of the function space. The $k$-th variational eigenvalue is then defined as:
\[
\lambda_k(\plap) = \min_{\mathrm{genus}(A) \geq k} \max_{f \in A} \rayl_p^p(f),
\]
where $\mathrm{genus}(A)$ refers to the Krasnoselskii genus of $A$ (see Section \ref{sec:Krasnoselskii}), and it can be shown that $\lambda_k(\plap)$ is a critical value of 
$\rayl_p^p$ for any $k$.

An important advantage of this definition is the ability to introduce a notion of algebraic multiplicity, defined as the number of times an eigenvalue appears in the variational sequence. Furthermore, if a variational eigenvalue $\lambda$ has multiplicity $k$, then there exists a subset of eigenfunctions associated with $\lambda$ that has genus greater than $k$. This subset can be interpreted as a generalized eigenspace, and it contains at least $k$ linearly independent eigenfunctions~\cite{struwe}. 

This framework provides a structured way to study the $p$-Laplacian spectrum and its properties, addressing some of the challenges introduced by the nonlinearity of the operator. We refer to \Cref{subsec:variational_spectrum} for a deeper discussion of the variational spectrum and its implications.

Somewhat surprisingly, the introduction of variational eigenvalues not only recovers but also strengthens the relationships between certain geometric properties of the domain and the spectrum of the Laplacian. A prominent example of this is the Cheeger inequality, which becomes an exact equality as $p$ approaches $1$. Specifically, considering the possibility of Dirichlet boundary conditions, the limits of the first two variational eigenvalues of the $p$-Laplacian operator as $p \to 1$ are equal to the first and second Cheeger constants of the graph~\cite{Bhuler,chang2016spectrum,Hua}.
Moreover, when considering the limit of higher variational eigenvalues of the graph $p$-Laplacian as $p \to 1$, it can be shown that:
\[
\lim_{p \to 1} \lambda_k(\plap) \leq h_k(\Gc),
\]
where $h_k(\Gc)$ denotes the higher-order Cheeger constant of index $k$~\cite{Tudisco1}.
For further details, we refer to \Cref{Sec:1-Lapl_spectrum}.

Analogously, when studying the $\infty$-limit of the variational $p$-Laplacian eigenvalues, we find that the limits of the first two variational eigenvalues of the $p$-Laplacian are equal to the reciprocals of the first two spherical packing radii of the graph. These radii correspond to the maximal radii that allow the inscription of one and two disjoint balls in the interior of the graph, respectively~\cite{deidda2024_inf_eigenproblem}. For further details, we refer to \Cref{Sec:inf_eigenproblem}.
Additionally, for the higher variational eigenvalues in the $\infty$-limit, it holds that:
\[
\lim_{p \to \infty} \lambda_k(\plap) \leq \frac{1}{R_k(\Gc)},
\]
where $R_k(\Gc)$ is the packing radius of order $k$ for the graph~\cite{grove1995new}. The packing radius $R_k(\Gc)$ represents the maximal radius that permits the inscription of $k$ disjoint balls within the graph's interior. For analogous results in the infinite dimensional setting, we refer to~\cite{Lind2,Lind3,parini2010second,Kawohl2003}.

When studying the relationships between $p$-Laplacian eigenvalues and higher-order Cheeger constants or packing radii of the domain, the analysis of the nodal domains of eigenfunctions becomes crucial. A nodal domain induced by an eigenfunction $f$ is defined as one of the maximal subdomains where $f$ is strictly positive or strictly negative.

The significance of nodal domains is twofold. First, consider an eigenpair $(f, \lambda)$ on a domain $\Omega$, which can be a graph or a bounded open domain in $\R^n$, with homogeneous Dirichlet boundary conditions. For any nodal domain $A$ induced by $f$, it is straightforward to observe that the restriction of $f$ to $A$ is itself an eigenfunction for $\lambda$ on $A$, again satisfying homogeneous Dirichlet boundary conditions. This implies that the nodal domains induced by an eigenfunction $f$ are subdomains of $\Omega$ that inherit the same eigenvalue $\lambda$ as $f$ on $\Omega$. Using this property, if we denote by $\Nc(f)$ the number of nodal domains induced by a function $f$, we can establish a lower bound for the limit of the $\plap$ eigenvalues:
\begin{equation}\label{intro_eq_geometrical_relations}
h_{\Nc(f_1)}(\Gc) \leq \lim_{p \to 1} \lambda_k(\plap) \leq h_k(\Gc)\:,\qquad 
1/R_{\Nc(f_{\infty})(\Gc)} \leq \lim_{p \to \infty} \lambda_k(\plap) \leq 1/R_k(\Gc),
\end{equation}
where $f_1$ and $f_{\infty}$ are appropriate limits of the $k$-th variational eigenfunctions of the $p$-Laplacian as $p \to 1$ and $p \to \infty$, respectively~\cite{ZhangNodalDO, Tudisco1, deidda2023PhdThesis, deidda2024_inf_eigenproblem, KELLER201680}.

In addition, the number of nodal domains induced by an eigenfunction is strongly connected to the frequency (or index) of its corresponding eigenvalue. This relationship was first observed for the Laplacian operator in Sturm's oscillation theorem for strings and later generalized by Courant to higher-dimensional domains~\cite{Courant}. For the linear graph Laplacian, it has been shown that the nodal count of any eigenfunction $f_k$ corresponding to $\lambda_k(\Delta)$ satisfies the following inequality:
\begin{equation}\label{intro_nodal_count_bound}
  k + m - 1 - \betti - z \leq \Nc(f_k) \leq k + m - 1,
\end{equation}
where $\betti$ is the total number of independent loops in the graph, $z$ is the number of zeros of $f_k$, and $m$ is the multiplicity of $\lambda_k(\Delta)$~\cite{Davies, Duval, Xu, Stadler_book-Laplacian-eigenvectors, Berkolaiko1, Berkolaiko2}.
These inequalities have been extended to the nonlinear graph $p$-Laplacian operator in~\cite{Tudisco1, DEIDDA2023_nod_dom, zhang2023pLap_noddom}. Notably, similar inequalities hold not only for variational eigenfunctions but also for nonvariational ones. The bounds depend solely on the position of the corresponding eigenvalue with respect to the variational spectrum, specifically satisfying $\lambda_k(\plap) < \lambda \leq \lambda_h(\plap)$. In particular, we note that the bounds for the nodal domain count, when combined with the inequalities in \eqref{intro_eq_geometrical_relations}, yield relationships between the $p$-Laplacian eigenvalues and the geometric quantities $R_k(\Gc)$ and $h_k(\Gc)$. We will discuss this topic in greater detail in \Cref{sec:nodal_domains}.

Considering the above discussion and in view of \eqref{intro_eq_geometrical_relations}, it is natural to focus further attention on the investigation of the $1$-Laplacian and the $\infty$-Laplacian eigenvalue problems. However, the definition of eigenpairs as critical values and points does not directly carry over to the cases $p = 1$ and $p = \infty$, due to the lack of differentiability of the corresponding Rayleigh quotients. Various approaches have been proposed to address this difficulty, which we will review in \Cref{subsec:infinity_1_Notation}.

In~\cite{Lind2,Lind3}, the authors study $\infty$-eigenpairs defined as solutions to the $\infty$-limit eigenvalue equation. In~\cite{hein2010inverse}, while investigating the $1$-Laplacian eigenvalue problem, the authors propose a generalized eigenfunction as a Clarke critical point of the Rayleigh quotient $\rayl_p$, that is, a function $f$ such that:
\begin{equation}\label{intro_Hein}
0 \in \partial^{Cl} \rayl_p(f),
\end{equation}
where $\partial^{Cl} \rayl_p(f)$ denotes the Clarke subgradient of the locally Lipschitz function $\rayl_p$~\cite{clarke1990optimization}.

Finally, in the context of the graph $1$-eigenvalue problem, the authors of~\cite{chang2016spectrum,chang2021nonsmooth} define a generalized eigenfunction $f$ as a function satisfying:
\begin{equation}\label{intro_CHang}
0 \in \partial \|\grad f\|_p \cap \bigcup_{\lambda \geq 0} \lambda \partial \|f\|_p,
\end{equation}
where $\partial \|\grad f\|_p$ and $\partial \|f\|_p$ represent the subgradients of the convex functions $f \mapsto \|\grad f\|_p$ and $f \mapsto \|f\|_p$, respectively~\cite{rockafellar2015convex,ekeland1999convex}. 

It is worth noting that the formulation in \eqref{intro_CHang} is the weakest among these definitions. Specifically, it follows from the properties of the Clarke subdifferential that any solution to \eqref{intro_Hein} also satisfies \eqref{intro_CHang}. Similarly, any solution of the $\infty$-limit eigenvalue equation proposed in~\cite{Lind2,Lind3} has been shown to satisfy \eqref{intro_CHang}~\cite{deidda2024_inf_eigenproblem}. However, there exist examples of functions that satisfy \eqref{intro_CHang} but not \eqref{intro_Hein}~\cite{zhang2021homological}, as well as functions that solve \eqref{intro_CHang} but not the $\infty$-limit eigenvalue equation~\cite{deidda2024_inf_eigenproblem}.

Notably, the authors of~\cite{chang2021nonsmooth} have demonstrated that the definition of variational eigenvalues can be extended to the degenerate cases $p = 1$ and $p = \infty$, where eigenvalues are understood as critical values in the sense of \eqref{intro_CHang}.

These results clearly open the door to the investigation of direct relationships between the $1$- and $\infty$-Laplacian variational eigenvalues and the topology of the graph, without the need to consider limits as in \eqref{intro_eq_geometrical_relations}. 

The $1$-Laplacian spectrum has been thoroughly studied in~\cite{ZhangNodalDO, ZHANG_top_mult, Zhang_1-Lap_cheeger_cut}. Notably, it can be shown that any $1$-Laplacian eigenvalue corresponds to the isoperimetric constant of the subgraphs induced by the nodal domains of the associated eigenfunction. Furthermore, inequalities analogous to \eqref{intro_eq_geometrical_relations} directly relate the $1$-Laplacian eigenpairs and their nodal counts to the Cheeger constants of the graph:
\begin{equation}
    h_{\Nc(f_k)}(\Gc) \leq \Lambda_k(\Delta_1) \leq h_k(\Gc).
\end{equation}
We explore the $1$-Laplacian spectrum in greater detail in \Cref{Sec:1-Lapl_spectrum}.
In that section, beyond discussing the Cheeger constants of the graph, we recall results from~\cite{ZHANG_top_mult} that relate the multiplicity of the largest $1$-Laplacian eigenvalue to the independence number of the graph. Moreover, we extend a result from~\cite{zhang2021homological} to relate the $1$-Laplacian variational eigenvalues to novel Cheeger-like constants defined via the homogeneous Dirichlet eigenvalue problem.
Finally, in \Cref{subsection_1-Lap_nodal_domains}, we review results on the nodal count of $1$-Laplacian eigenpairs. However, unlike the case of the linear Laplacian, there exists no lower bound on the nodal count in terms of the index, as in \eqref{intro_nodal_count_bound}; see also~\cite{ZhangNodalDO}.

In contrast, the $\infty$-Laplacian variational spectrum has been only preliminarily explored in~\cite{deidda2024_inf_eigenproblem}. Therefore, \Cref{Sec:inf_eigenproblem} aims to address this gap. After proving the continuity of the $p$-Laplacian variational spectrum as $p \to \infty$ in \Cref{Sec:regularity_of_p_lap_spectrum}, we establish an $\infty$-analogue of \eqref{intro_eq_geometrical_relations}:
\begin{equation}
    \frac{1}{R_{\PNc(f_k)}(\Gc)} \leq \Lambda_k(\Delta_\infty) \leq \frac{1}{R_k(\Gc)},
\end{equation}
where only certain nodal domains, referred to as ``perfect nodal domains'', are considered for the lower bound. Perfect nodal domains are those where the oscillations of the eigenfunction have the same amplitude.
We provide a comprehensive study of the perfect nodal domains of $\infty$-eigenfunctions in \Cref{subsec_infinity_lap_nodal_dom}, where we also examine the nodal count of viscosity eigenpairs, i.e. solutions to the $\infty$-limit eigenvalue equation introduced in~\cite{deidda2024_inf_eigenproblem}.

In \Cref{Sec:inf_eigenproblem}, beyond establishing the connection between the $\infty$-Laplacian variational spectrum and the sphere packing problem of the graph, we investigate relationships between the variational spectrum and the independence numbers of the graph. We prove that the number of perfect nodal domains of $\infty$-eigenfunctions can be bounded in terms of an independence number defined with respect to the corresponding eigenvalue. Additionally, we show that the matching number of the graph can be related to the multiplicity of the largest $\infty$-eigenvalue.

Finally, we discuss a minimal $k$-partition problem based on the homogeneous Dirichlet $\infty$-eigenproblem and relate this problem to the sphere packing problem of the graph. This analysis highlights the interplay between spectral properties, graph topology, and optimization problems in the $\infty$-Laplacian framework.

Overall, the primary goal of this manuscript is to review the state of the art in the spectral theory of the graph $p$-Laplacian operator, while integrating several novel contributions that address important gaps in the existing theory. Where appropriate, we enrich the discussions with illustrative examples and heuristic insights, paying special attention to building connections between the infinite-dimensional and graph settings, as well as between the linear and nonlinear frameworks. Furthermore, with the aim of fostering further research in this field, we highlight several open problems that we believe warrant future investigation.

\subsection{Summary of Main Contributions}
The main contributions of this work are twofold:

On the one hand, we provide a comprehensive review of the state of the art in spectral theory for nonlinear graph Laplacian operators. We offer a unifying perspective on the connections between graph properties and nonlinear eigenfunctions, shedding light on this fascinating research area with numerous modern applications. Alongside reviewing and synthesizing existing results, we present alternative proofs for several known results, using different arguments that, in our opinion, provide greater clarity and precision.

The key results summarized from the literature include the relationships among three types of $\infty$-Laplacian eigenpairs (Theorem \ref{theorem:viscosity_eigenpairs_are_generalized_eigenpairs}); results about the regularity of the $p$-Laplacian spectrum \Cref{lemma_upp_semicontinuity_plap_spectrum} \Cref{thm:homological-eigen}; the duality between node $p$-Laplacian and edge $q$-Laplacian (Theorem \ref{Thm:duality}); nodal counts for $p$-Laplacian eigenfunctions (Theorems \ref{Theorem_1st_eigen_characterization}, \ref{thm:nodal-count}, \ref{Thm_nodal_domains_decomposition_1-Lap_eigenfunctions}, \ref{Thm:nodal-count-for_1_lap}); relations between $1$-eigenvalues and Cheeger constants (Theorems \ref{thm:1-eigenpair_and_cheeger_constants} \ref{ThM:p-eigenpairs_and_cheeger_constants}) and between $\infty$-eigenvalues and packing radii (\Cref{Cor:packingradii_inequality});  a Sandwich theorem for estimating the multiplicity of the largest eigenvalue of the normalized $1$-Laplacian (Theorem \ref{Thm_1_lap_mult_independence_number}).

On the other hand, we significantly extend the existing theory with several new results, including, a new bound for the variational index of a nonlinear eigenvalue in terms of its linear index (\Cref{Comparison_theorem}),
novel relation between homological eigenvalues and Clarke critical values (Proposition \ref{prop:homology-Clarke}); new monotonicity inequalities on variational eigenvalues of the $p$-Laplacian (Lemma \ref{lemma:p-monotonic}); weak nodal counts for the second eigenvalue of the $1$-Laplacian (Theorem \ref{thm:second-eigen-1}); new relations between the variational  $\infty$-eigenpairs and the packing radii of the graph (\Cref{thm:k-inequality}), new bounds for packing radii in terms of the spectral minimal partition costs of the $\infty$-Laplacian (Theorem~\ref{Thm_spectral_min_partitions_2}). 

    More specifically, the main theoretical novelties include:

    \begin{itemize}
        \item An extension of the theoretical understanding of the graph $\infty$-Laplacian spectrum, particularly in relation to the packing radii (Theorems \ref{thm:k-inequality}, \ref{Thm_spectral_min_partitions}), independence numbers (Propositions \ref{pro:independence}, \ref{prop_k-Indipendece-perfect_nodal_domains-inequalities}, \ref{prop_k-Indipendece-perfect_nodal_domains-inequalities_unweighted_case}), and the matching number of the graph (Proposition \ref{pro:matching-2}).
        
        \item A novel detailed study of the nodal domain theory of the $\infty$-Laplacian (Theorems \ref{Thm_nodal_count_viscosity_eigenfunctions}, \ref{thm:relate-3bounds}, \ref{Thm:lower_bound_perfect_nodal_count_limit_eigenpairs}), introducing the concept of perfect nodal domains to address difficulties arising from the high oscillation frequencies of $\infty$-Laplacian eigenfunctions (Theorem \ref{thm:infty_nodal-count}).
        
        \item The introduction of the notion of Dirichlet Cheeger constants and establish new inequalities linking higher-order Cheeger constants, Dirichlet Cheeger constants, and $1$-Laplacian eigenvalues (Theorem~\ref{thm:Dirichlet_Cheeger_constants_1-Lap_eigenvalues}).   
    \end{itemize}

\subsection{Related relevant topics}

While aiming to provide a comprehensive overview and introduction to modern nonlinear spectral graph theory, we have had to leave several key topics unexplored due to space constraints. This subsection offers a concise summary of topics that have received significant attention in recent years but are not covered in this paper. We briefly discuss these topics and point to the main references and criticisms. Readers interested only in the topics discussed in this paper may skip this section without any loss of continuity or clarity.

\paragraph{From the graph setting to the infinite-dimensional setting.}
We mentioned earlier that results in the graph setting are closely related to analogous results in the infinite-dimensional setting. However, precise connections between the two settings and, in particular, whether key results carry over in the limit as the number of nodes in the graph grows to infinity, have only been marginally investigated. This is an active area of research that deserves further study.
Convergence results for the discrete $p$-Laplacian spectrum toward the spectrum of the differential $p$-Laplacian operator on grids and geometric random graphs would be of significant impact. To the best of our knowledge, this topic has been investigated only for the linear Laplacian operator in \cite{trillos2018variational, trillos2020error}. Moreover, results on the convergence of viscosity $\infty$-eigenpairs on grid graphs toward viscosity solutions of the limit eigenvalue equation on a bounded open domain have been discussed in \cite{bozorgnia2024infinity}. In \cite{Burger23, roith2023continuum} it has also been proved that $\Gamma$-convergence of some graph functional, e.g. $\|\grad f\|_p^p$, toward some functional on an infinite dimensional space, e.g. $\|\nabla f\|_p^p$, can be used to study convergence of the corresponding ground states, i.e. eigenfunctions of minimal energy. 
Additionally, in the framework of semi-supervised learning, different authors have studied the convergence of the solutions of a $p$-Poisson equation with a fixed source term from the graph setting to the infinite-dimensional one \cite{calder2018game,calder2019consistency, slepcev2019analysis, roith2023continuum}.

\paragraph{Other nonlinear eigenvalue problems on graphs.}
The graph $p$-Laplacian is a central object in nonlinear spectral graph theory. However, many other nonlinear eigenproblems that encode geometric information are not discussed in this manuscript. See e.g. the signless $p$-Laplacian \cite{ZhangNodalDO, chang2016signless}, the functionals studied in the core periphery literature \cite{Tudisco_core_periphery}, and extensions of the $p$-Laplacian to hypergraphs \cite{Mulas2022_plap_hyper} and signed graphs \cite{zhang2023pLap_noddom, ge2023new}. Most of these nonlinear eigenproblems have been only marginally investigated so far and do not enjoy the same level of understanding as the $p$-Laplacian on graphs. Their deeper investigation could positively impact the further development of the field.

\paragraph{Computational aspects in nonlinear spectral graph theory.}
A significant challenge associated with nonlinearity lies in computational aspects. While nonlinear spectral graph theory often yields sharper insights than its linear counterpart, the numerical investigation of the relevant nonlinear eigenproblems is considerably more intricate than in the linear case. Consequently, this has somewhat restricted the practical application of nonlinear spectral approaches. Nevertheless, there is growing interest in the numerical study of these problems, and we encourage the community to actively engage in this crucial matter.

From a computational standpoint, it is natural to distinguish between the extreme eigenpair (corresponding to the smallest nontrivial eigenvalue) and the rest of the spectrum.

The smallest nontrivial eigenpair of the $p$-Laplacian operator corresponds to the first variational eigenpair when the graph has a boundary or the operator is shifted by a regularization term, and to the second variational one otherwise. For the computation of these eigenpairs, different discrete and continuous methods have been developed, including nonlinear power methods and nonlinear optimization strategies \cite{hein2010inverse, Hein_Beyond, Bozorgnia_power_method, BOZORGNIA_approximation, laubmann2025duality, tudisco2018community, tudisco2021nonlinear,cristofari2020total, Jarlebring_an_inverse}, gradient-like flows for the $p$-Rayleigh quotient \cite{gilboa2018nonlinear, gilboa2021iterative, Gilboa_rayleight_quot_minimization} whose steady states are eigenvectors, and gradient flows admitting asymptotic profiles that are eigenvectors \cite{BUNGERT_gradient_flows_and_nonlinear_power_method, HYND_approximation_of_the_least_Rayleight_quotient}. These latter methods guarantee a decrease in the $p$-Rayleigh quotient at each iteration or throughout the flow and are thus mainly aimed at computing eigenpairs of minimal energy. In particular, some of these methods provide convergence guarantees to the first variational eigenvalue when the starting point is strictly positive. Additionally, some formal interpretations of the nonlinear power methods and their variants in terms of discretized gradient flows have been discussed in \cite{BUNGERT_gradient_flows_and_nonlinear_power_method, Hein_Beyond}. A further algorithm is provided by the self-consistent field iteration and its variants \cite{BAI_SCF, upadhyaya2021self}, which find the minimum eigenpair as the limit of a sequence of linear eigenpairs.

The computation of higher-order eigenpairs is an even more challenging and less studied problem. In \cite{Yao07}, the authors propose a scheme for computing a sequence of $N$ eigenpairs using a local min-max method. The approach relies on the heuristic assumption that a $p$-Laplacian eigenvalue can be characterized as a local maximum of the $p$-Rayleigh quotient on a linear subspace spanned by eigenvectors corresponding to lower frequencies. 

The authors of \cite{deidda2024_spec_energy} introduce novel gradient-based numerical methods for computing $p$-Laplacian eigenpairs for any $p\in(2,\infty)$. The methods are obtained by reformulating the graph $p$-Laplacian eigenproblem into a constrained weighted Laplacian eigenvalue problem and representing $p$-Laplacian eigenpairs as saddle points of a family of spectral energy functions. Notably, the gradient flows presented in \cite{deidda2024_spec_energy} allow the computation of $p$-Laplacian eigenpairs as limits of eigenpairs of generalized linear eigenvalue problems that evolve over time, while the index of the linear eigenpair remains constant throughout the flow. The resulting discretization of this flow is conceptually similar to the self-consistent field iteration \cite{BAI_SCF, upadhyaya2021self}. However, to our knowledge, a precise comparison of these two methods has not been discussed.

In \cite{Luo10}, the authors propose the idea of $p$-orthogonal eigenfunctions and show that $p$-Laplacian eigenfunctions are $p$-orthogonal up to a second-order Taylor expansion. Using these results, the authors of \cite{Lanza24} establish a projected gradient descent method, and an alternating direction method, to find approximated $p$-Laplacian eigenpairs by solving a constrained minimization under $p$-orthogonality constraints. 

In \cite{berkolaiko2025eigenvalues}, the authors propose a Hellmann--Feynman type perturbation theory for the $p$-Laplacian. This theory suggests that the eigenpairs of the $p$-Laplacian on a generic graph can be computed by perturbing the $p$-Laplacian spectrum of a spanning subtree. This is particularly intriguing because the $p$-Laplacian spectrum on a tree is known to be finite and has the same dimension as the graph \cite{DEIDDA2023_nod_dom}. Moreover, \cite{DEIDDA2023_nod_dom} provides a constructive, albeit computationally expensive, method to compute the entire $p$-Laplacian spectrum on a tree.

Finally, in \cite{ge2025computing}, the authors establish an equivalence between the $p$-Laplacian eigenvalue problem and a tensor eigenvalue problem, for even values of $p$. Consequently, they explore how algorithms designed for tensor eigenvalue problems can be utilized to compute the $p$-eigenpairs when $p$ is even.

We conclude by emphasizing that, in many applications, one is interested in computing specific variational eigenpairs, rather than just any eigenpair. While desirable, we observe that this might be particularly challenging. For instance, the Cheeger constants of a graph are known to be NP-complete, so, by the equivalence between variational eigenpairs and Cheeger constants, the computational complexity of the $p$-Laplacian variational eigenvalues as $p$ approaches $1$ should be similarly challenging. However, studying the properties of any computable eigenpairs is also crucial and may allow us to establish their proximity to a determinate variational one. For example, both the methods presented in \cite{deidda2024_spec_energy, Yao07} compute $p$-Laplacian eigenpairs such that the eigenfunction is a critical point of the $p$-Rayleigh quotient of a given Morse index. 

\paragraph{Nonlinear spectral decompositions.}
In the linear case, the eigenfunctions of the Laplacian form an orthonormal basis of the underlying space, which has significant applications in denoising and dimensionality reduction methods such as principal component analysis. This is because high-frequency components are typically associated with noise and fine-scale details. Therefore, projecting a signal onto the subspace spanned by low-frequency eigenfunctions often captures most of its relevant features. This raises the question of whether the $p$-Laplacian eigenpairs also have similar applications. As we will discuss in \cref{sec:nodal_domains}, for $p$-Laplacian eigenpairs, the higher the eigenvalue, the higher the frequency of the eigenfunction. The difference lies in the type of ``wave'': eigenfunctions become more piecewise constant as $p$ approaches $1$ and more peaked as $p$ approaches $\infty$. However, for $p \neq 2$, we lose the most important property used in the linear case, namely that the eigenfunctions form an orthonormal basis. While they still span the entire space, as we will discuss in \cref{sec:p_Lap_eigenvalue}, they are not orthogonal and can be more than the dimension of the space.

A further problem is that the nonlinearity of the operator does not guarantee that any point in the linear span of a set of low-frequency eigenfunctions has only low frequencies. Assume, for example, that we want to solve $\plap u=f$ but only have $\tilde{f}$, which is a noisy measurement of $f$. Then, unlike the linear case, if, instead of $f$, we consider the projection of $\tilde{f}$ onto the space $\pi$ spanned by $m$ low-frequency eigenfunctions of $\Delta_p$, we have no guarantee that $u$ remains in $\pi$. This leads to the question of whether there exists a subspace (not necessarily linear), e.g., a manifold embedded in $\R^N$, that is invariant under the $p$-Laplacian operator and has bounded frequency, i.e., such that the $p$-Rayleigh quotient is smaller than a certain eigenvalue. 

An alternative approach to this problem is given by the ``nonlinear spectral decompositions'' of a signal $f$ with respect to a convex functional $J(\cdot)$, defined on a Hilbert space $\mathcal{X}$. This is a map from $f$ to a vector valued Radon measure $\Phi_s$ such that (i) $f$ can be reconstructed in terms of $\Phi_s$, i.e. $f=\int_0^\infty d\Phi_s$; and (ii) if $f$ is a nonlinear eigenvector associated to $J$, i.e. $\lambda f\in \partial J(f)$ where $\partial J$ denotes the subgradient of $J$, then $\Phi_s=f \delta_\lambda(s)$. Note that condition (ii) is the same as asking that $(f,\lambda)$ is a generalized critical pair (see \cref{Def:generalized_p-eigenpair}) of the Rayleigh quotient $J(f)/\|f\|_\mathcal{X}$, where $\|\cdot\|_\mathcal{X}$ is the norm induced by the scalar product.
Interestingly, nonlinear spectral decompositions can be defined in terms of gradient flow solutions with $f$ as the starting point, and solutions of regularized inverse problems \cite{gilboa2018nonlinear, burger2016spectral, BungertNonlineardecomp, Gilboa_total_variation_spectral_framework, COHEN_introducing_p_lap_spectra}. Additionally, if the convex functional $J(\cdot)$ can be written as $\|Af\|_1$ and the matrix $A$ satisfies some additional properties, then the spectral decomposition produces exactly a decomposition of the signal into finitely many nonlinear eigenvectors associated with the functional $J(\cdot)$ \cite{burger2016spectral, BungertNonlineardecomp}. These results also apply to the $p$-Laplacian by considering the functional $\|\nabla f\|_p$ and, in particular, the Total Variation by taking $p=1$. However, the eigenpairs considered in these spectral decompositions are critical pairs of the Rayleigh quotient $\|\nabla f\|_p/\|f\|_2$ and thus are not the same as the ones we consider in this manuscript, which are critical pairs of $\|\nabla f\|_p/\|f\|_p$. The extension of these results to Banach spaces would fill this gap and allow us to understand how the eigenpairs investigated in this manuscript can be used in spectral decompositions; however, we are not aware of results in this direction. {Furthermore, the connection between the eigenpairs corresponding to $\|\nabla f\|_p/\|f\|_2$ and the geometry of the graph is not yet clear.}

\paragraph{Domain characterization by the the full $p$-Laplacian spectrum}
A further natural question is whether the eigenpairs of the $p$-Laplacian operator contain sufficient information to characterize the underlying geometry. It is well known that for the standard graph Laplacian, the full set of eigenpairs determines the Laplacian matrix and hence the connectivity structure and edge weights of the graph. An analogous phenomenon holds in the continuous setting, where the complete spectral decomposition of the Laplace-Beltrami operator determines the operator itself, which in turn encodes the Riemannian metric of the underlying manifold. Thus, in both the discrete and continuous linear settings, the full spectral structure provides, in principle, a means of recovering the underlying geometry of the domain. Whether an analogous reconstruction principle holds for the nonlinear $p$-Laplacian is, to the best of our knowledge, an open question. More specifically, it is natural to ask whether all, or even only part, of the $p$-Laplacian eigenpairs can determine the underlying graph or manifold. 

\section{Graph Setting: Definitions and Notation}\label{Sec:intro-sec_graph_setting}
We devote this section to establishing the notation used throughout this work.

We consider an undirected graph denoted by the triple $\Gc := (\nodeset, \edgeset, \edgelength)$, where $\nodeset$ is the finite set of nodes (or vertices) of the graph, and $\edgeset \subset \nodeset \times \nodeset$ is the set of edges. Since the graph is undirected, the edge set $\edgeset$ is symmetric, meaning that if $(u, v) \in \edgeset$, then $(v, u) \in \edgeset$ as well. Finally, $\edgelength : \edgeset \to \R$ is a weight function on the edges such that $\edgelength(u, v) = \edgelength(v, u)$. For convenience, we will often write $\edgelength_{uv}$ to denote $\edgelength(u, v)$.

The value of $\edgelength(u, v)$ can be interpreted as the reciprocal of the length of the edge $(u, v)$. This interpretation aligns with the common notion that the edge weight $\edgelength(u, v)$ reflects the robustness of the relationship between nodes $u$ and $v$: the stronger the relationship, the closer the nodes are considered to be. According to this interpretation, we define the distance between two nodes $u$ and $v$ as the length of the shortest path joining them:
\begin{equation}
d(u, v) = \min_{\Gamma \in \mathrm{path}_{u,v}} \length(\Gamma),
\end{equation}
where $\mathrm{path}_{u,v}$ denotes the set of paths joining $u$ to $v$:
\begin{equation}
\mathrm{path}_{u,v} = \big\{ \Gamma = \{u_i\}_{i=1}^n \;\big|\; u_1 = u, \; u_n = v, \; (u_i, u_{i+1}) \in \edgeset, \; n \; \text{arbitrary} \big\}.
\end{equation}
The length of a path $\Gamma = \{u_i\}_{i=1}^n$ is defined as:
\begin{equation}
\length(\Gamma) = \sum_{i=1}^{n-1} \frac{1}{\edgelength(u_i, u_{i+1})}.
\end{equation}

Next, we introduce the object of our study: the graph $p$-Laplacian operator. This graph operator is defined in a manner consistent with its differential counterpart, defined in the infinite dimensional setting, and expressed as:
\[
\plap f = -\mathrm{div}(|\nabla f|^{p-2} \nabla f).
\]

Throughout the manuscript, when referring to the infinite dimensional setting, we refer to a bounded open domain $\Omega\subset \R^n$ and the associated function space $W^{1,p}(\Omega)$.

Let $\Hc(\nodeset)$ and $\Hc(\edgeset)$ denote the Hilbert spaces of functions on the nodes and edges of the graph, respectively, which are isomorphic to $\R^{|\nodeset|}$ and $\R^{|\edgeset|}$. To fix the notation and avoid confusion, we use lowercase letters (e.g., $f$) to denote functions in $\Hc(\nodeset)$ and uppercase letters (e.g., $G$) to denote functions in $\Hc(\edgeset)$.

We now recall the definition of the directed weighted incidence matrix of $\Gc$, denoted by $K \in \R^{|\edgeset| \times |\nodeset|}$:
\begin{equation}
    K\big((u_1, u_2), u_3\big) := 
    \begin{cases}
        \edgelength_{u_1u_2}, & \text{if } u_3 = u_2, \\
        -\edgelength_{u_1u_2}, & \text{if } u_3 = u_1, \\
        0, & \text{otherwise}.
    \end{cases}
\end{equation}
Observe that, based on the interpretation of $\edgelength_{uv}$ as $1/\length(u, v)$, for any function $f \in \Hc(\nodeset)$, we have:
\begin{equation}
    Kf(u, v) = \edgelength_{uv} \big(f(v) - f(u)\big) = \frac{f(v) - f(u)}{\length(u, v)}.
\end{equation}
The equation above shows that $Kf(u, v)$ provides a discrete approximation of the derivative along the edge $(u, v)$ of a continuous function that takes values $f(u)$ and $f(v)$ at the endpoints of the edge. For this reason, $K$ can be interpreted as an operator from $\Hc(\nodeset)$ to $\Hc(\edgeset)$, serving as a graph analogue of the differential gradient operator, as it measures the slope of a function $f$ along the edges. In particular, note that since the graph is undirected, $\grad f(u, v) = -\grad f(v, u)$.

Next, we introduce the divergence operator. On graphs, the absence of a boundary is usually interpreted as analogous to homogeneous Neumann boundary conditions. According to the classical divergence theorem, if $K$ is viewed as the graph counterpart of the differential gradient operator, then $-1/2 \, K^T$ can be understood as a graph divergence operator. This interpretation aligns with the following identity:
\begin{equation}
    -\langle f, -\frac{1}{2}K^T G \rangle_{\Hc(\nodeset)} = \langle \grad f, G \rangle_{\Hc(\edgeset)}.
\end{equation}

Where the Hilbert spaces $\mathcal{H}(\nodeset)$ and $\mathcal{H}(\edgeset)$ have been equipped with the  scalar products $
\langle f, g \rangle_{\Hc(\nodeset)} = \sum_{u \in \nodeset} f(u) g(u)$ and $\langle F, G \rangle_{\Hc(\edgeset)} = \frac{1}{2} \sum_{(u, v) \in \edgeset} F(u, v) G(u, v)$ with the factor $1/2$ included in the scalar product on $\Hc(\edgeset)$ to account for the undirected nature of the graph.

In particular, we have 
\begin{equation}
\begin{aligned}
\divg:\;\Hc(\edgeset)&\longrightarrow \Hc(\nodeset)\\
G\;&\longrightarrow \divg\,G(u)=\frac{1}{2}\sum_{v\sim u}\edgelength_{uv}\big(G(u,v)-G(v,u)\big)\,,
\end{aligned}
\end{equation}
where $\{v \;|\; v \sim u\}$ represents the set of nodes connected to $u$ by an edge, i.e., such that $(u, v) \in \edgeset$.

Having provided a graph interpretation of the differential gradient and divergence operators, we can now show that both the classical graph Laplacian operator ($p = 2$) and the more general $p$-Laplacian operator ($p \in (1, \infty)$), as defined below, align with their differential counterparts defined in the infinite dimensional setting. 

\begin{definition}
\label{def:graph_p-Lap_operator}
The graph $p$-Laplacian operator is defined as:
\begin{equation}
\plap f(u) = \frac{1}{2}K^T\big(|\grad f|^{p-2} \odot \grad f\big)(u) = \sum_{v \sim u} \edgelength_{uv} |\grad f(v, u)|^{p-2} \grad f(v, u),
\end{equation}
where $|\grad f|^{p-2}$ is understood entrywise, and $\odot$ denotes the entrywise product (we will omit this symbol in subsequent expressions).
\end{definition}

When $p = 2$, it can be verified that this definition corresponds to the classical definition of the graph Laplacian $\lap[2]$ in terms of the adjacency matrix. Specifically, consider $A$, the weighted adjacency matrix of the graph (where the weights are squared), defined as $A_{uv} = \edgelength_{uv}^2 \Ind_{\edgeset}\big((u, v)\big)$, where $\Ind$ denotes the indicator function. A direct computation shows that
\begin{equation}
\lap[2] = \mathrm{diag}(A \mathbf{1}) - A,
\end{equation}
where $\mathbf{1}$ denotes the vector with all entries equal to $1$.
Moreover, for $p \neq 2$, using the expression of $K$, we have:
\begin{equation}
\plap f(u) = \sum_{v \sim u} \edgelength_{uv}^p |f(u) - f(v)|^{p-2} \big(f(u) - f(v)\big),
\end{equation}
which corresponds to the popular definition of the graph $p$-Laplacian found in e.g.~\cite{Tudisco1, Amghibech1, Amghibech2}.
We emphasize that the graph and differential versions of the $p$-Laplacian are not only similar in form but are also consistent from an approximation point of view. Indeed, when we consider a sequence of geometric graphs that converge, in some sense, to a bounded open domain $\Omega\subset\R^n$, the graph $p$-Laplacian operators converge to the differential $p$-Laplacian of that domain \cite{slepcev2019analysis, GarciaSlepcev16, Garcia20, Garcia16}.

It is worth mentioning that the graph $p$-Laplacian operator can be defined, consistently with the differential $p$-Laplacian operator, in alternative ways to \cref{def:graph_p-Lap_operator}. A prominent alternative is given by the game theoretic $p$-Laplacian operator \cite{calder2018game, manfredi2015nonlinear, Elmoataz1}. This operator is a linear combination of the linear Laplacian operator and the $\infty$-Laplacian (we will introduce the last later in \cref{def_Inf_lap_operator}). The game theoretic $p$-Laplacian has found significative applications in image processing \cite{Elmoataz1} and semi-supervised learning \cite{calder2018game} and admits interesting interpretations from a game-theory point of view \cite{Peres, manfredi2010asymptotic}. However its relevance from a spectral point of view appears to be limited. Indeed, differently from the operator in \cref{def:graph_p-Lap_operator}, which can be derived as $\partial_f\|Kf\|_p^p$, the graph game theoretic $p$-Laplacian does not admit a variational formulation. Indeed, as we will see in the next section, the variational formulation is the crucial property that allows us to introduce the notion of 
$p$-Laplacian eigenpairs.

We devote the remainder of this section to discuss the boundary case. While graphs typically do not have a boundary, one can be artificially introduced, and we will see that nonlinear spectral graph theory in the presence of a boundary exhibits some intriguing properties that are of interest for graph applications.

The boundary of a graph is defined as a subset of its nodes, $\boundary \subset \nodeset$. Homogeneous Dirichlet boundary conditions can then be imposed by considering the subspace of node functions:
\begin{equation}
\mathcal{H}_0(\nodeset) := \{f : \nodeset \to \R \;|\; f(u) = 0 \;\forall u \in \boundary\}.
\end{equation}
If we consider the set of internal nodes $\nodeset \setminus \boundary$ and the function space defined on it:
\begin{equation}
\mathcal{H}(\internalnodes) := \{f : \internalnodes \to \R\},
\end{equation}
then $\mathcal{H}(\internalnodes)$ and $\mathcal{H}_0(\nodeset)$ are clearly isomorphic through the trivial embedding map. 

Furthermore, identifying a function $f \in \mathcal{H}(\internalnodes)$ with its lifting $\tilde{f} \in \mathcal{H}_0(\nodeset)$, it is straightforward to observe that $\grad \tilde{f} = \tilde{\grad} f$, where $\tilde{\grad}$ is the submatrix of $K$ obtained by removing the columns corresponding to the boundary nodes. As a result, in the boundary case, the graph $p$-Laplacian, as an operator from $\mathcal{H}(\internalnodes)$ to itself, can be expressed as:
\begin{equation}\label{Def_boundary_p-Lap}
\begin{aligned}
    \tilde{\Delta}_p f(u) &= \tilde{\grad}^T\big(|\tilde{\grad} f|^{p-2} \odot \tilde{\grad} f\big)(u)
    \\
    &= \sum_{\substack{v \sim u \\ v \in \internalnodes}} \edgelength_{uv}^p |f(u) - f(v)|^{p-2} \big(f(u) - f(v)\big)
    + \sum_{\substack{v \sim u \\ v \in \boundary}} \edgelength_{uv}^p |f(u)|^{p-2} f(u).
\end{aligned}
\end{equation}

It is straightforward to observe that, given $f \in \mathcal{H}(\internalnodes)$ and its corresponding lifting $\tilde{f} \in \mathcal{H}_0(\nodeset)$, we have $\tilde{\Delta}_p f(u) = \plap(\tilde{f})(u)$ for any internal node $u \in \internalnodes$. For this reason, in later discussions about the boundary case, we will implicitly assume that functions belong to $\mathcal{H}(\internalnodes)$, and we will omit the tilde superscript on the operators $\grad$ and $\plap$. 

We emphasize that the $p$-Laplacian operator defined in the boundary case \eqref{Def_boundary_p-Lap} matches the definition of the generalized $p$-Laplacian operator introduced in~\cite{DEIDDA2023_nod_dom} and has also been studied in~\cite{Hua, PARK2011}. 
In particular, note that the non-boundary case is equivalent to setting $\boundary = \emptyset$. Therefore, unless otherwise stated, we allow our graph to have a boundary in this work.

Moreover, throughout the paper we will denote by $N=|\internalnodes|$ the dimension of $\mathcal{H}(\internalnodes)$.

\section{The $p$-Laplacian eigenvalue problem}\label{sec:p_Lap_eigenvalue}
The $p$-Laplacian eigenvalue problem can be introduced by mimicking the infinite dimensional setting. Recalling the interpretation of the directed incidence matrix as a graph gradient operator, we define the Rayleigh quotient:
\begin{equation}
\rayl_{\nodeweight,p}(f) = \frac{\|\grad f\|_p}{\|f\|_{\nodeweight,p}} = \frac{\Big(\frac{1}{2} \sum_{(u, v) \in \edgeset} |\grad f(u, v)|^p\Big)^{\frac{1}{p}}}{\Big(\sum_{u \in \internalnodes}\nodeweight_u |f(u)|^p\Big)^{\frac{1}{p}}} \qquad f\in \mathcal{H}(\internalnodes)
\end{equation}
where a strictly positive measure $\nodeweight : \internalnodes \to \R_+$ has been included in the space $\mathcal{H}(\internalnodes)$, inducing the weighted $p$-norm above.
For $p \in (1, \infty)$, the critical point equation of $\rayl_p(f)$, up to scaling factors, is given by:
\begin{equation}\label{Def_p-lap_eigenequation}
\plap f(u) = \big(\rayl_{\nodeweight,p}(f)\big)^p \nodeweight_u |f(u)|^{p-2} f(u), \qquad \forall u \in \internalnodes.
\end{equation} 
This leads to the following definition:
\begin{definition}\label{Def_p-eigenpair_by_eq}
A pair $(f, \lambda)$ is called a $p$-Laplacian eigenpair if:
\[
\plap f(u) = \lambda \nodeweight_u |f(u)|^{p-2} f(u) \quad \forall u \in \internalnodes.
\]
\end{definition}
Multiplying \eqref{Def_p-lap_eigenequation} by $f(u)$ and summing over all vertices $u \in \internalnodes$ reveals that, if $\lambda$ is an eigenvalue corresponding to the eigenfunction $f$, it necessarily holds that $\lambda = \big(\rayl_{\nodeweight,p}(f)\big)^p$. 
Given an eigenpair $(f, \lambda)$, we refer to $\Lambda = \rayl_{\nodeweight,p}(f) = \lambda^{\frac{1}{p}}$ as a normalized eigenvalue. Throughout this manuscript, we adopt the notation $\lambda$ to refer to eigenvalues, while the capital $\Lambda$ always refers to a normalized eigenvalue. Normalized eigenvalues are particularly relevant in the study of the limit case $p$ going to $\infty$. Indeed, as we will discuss in more detail in \cref{subsec:infinity_1_Notation} and \cref{Sec:regularity_of_p_lap_spectrum}, the set of normalized $p$-Laplacian eigenvalues, varying $p$, is bounded and thus one can study accumulation points of the normalized eigenvalues in the limit $p$ going to $\infty$.

%
\begin{remark}\label{weighted_p_lap_eigenvalue_problem}
In more generality, we mention that also the Hilbert space $\Hc(\edgeset)$ could be endowed with a measure, $\edgeweight : \edgeset \to \R$, inducing the following norm:
\begin{equation} 
\|G\|_{\edgeweight, p}^p = \frac{1}{2} \sum_{(u,v) \in \edgeset} \edgeweight_{uv} |G(u, v)|^p.
\end{equation}
This norm as also $\|f\|_{\nodeweight,p}$ can be interpreted as graph counterparts of the integrals $\int |\nabla f|^p \, d\edgeweight$ and $\int |f|^p \, d\nodeweight$, respectively. Notably, unlike $\edgelength$, which represents a length and is included in the graph counterpart of the gradient, the weights $\edgeweight$ and $\nodeweight$ are not raised to the power $p$. This distinction aligns with the idea of ``discretization'' of the aforementioned integrals.
In this weighted setting, differentiating the Rayleigh quotient 
\begin{equation}\label{eq_p_rayleigh_quotient}
    \rayl_{p, \edgeweight, \nodeweight}^p(f) = \frac{ \|\grad f\|_{p, \edgeweight}^p }{ \|f\|_{p, \nodeweight}^p}= \frac{\sum\limits_{(u,v)\in E} \edgeweight_{uv} \edgelength_{uv}^p|f(u)-f(v)|^p}{2 \sum\limits_{v\in \internalnodes} \nodeweight_v |f(v)|^p},
\end{equation}
    yields the (weighted) $p$-Laplacian eigenvalue equation:
\begin{equation}
\lap[p, \edgeweight] f(u) = \frac{1}{2}K^T\big(\edgeweight \odot |\grad f|^{p-2} \grad f\big)(u) = \lambda \nodeweight_u |f(u)|^{p-2} f(u), \quad \forall u \in \internalnodes.
\end{equation}
\end{remark}
Such a generalized problem has been studied in certain cases~\cite{Tudisco1, DEIDDA2023_nod_dom}. 
To avoid overly heavy notation, we will not explicitly consider measures on the edge space. 
These can be effectively incorporated into the edge length parameter $\edgelength$, requiring only minor adjustments to the results. 
On the other hand, unless otherwise specified, we always allow for the presence of an everywhere nonzero measure (or weight) $\nodeweight$ on $\mathcal{H}(\internalnodes)$. 
In particular, unless relevant to the discussion, we will typically omit the subscript $\nodeweight$ in the notation $\|f\|_{\nodeweight, p}$ and $\rayl_{\nodeweight,p}$. However, we always assume that the norm $\|f\|_p$, as also the Rayleigh quotient $\rayl_p$, are weighted with respect to some everywhere nonzero measure $\nodeweight$.
\begin{definition}[Unweighted graph]\label{Def:Unweighted_graph}When referring to an \textbf{unweighted graph}, we mean that $\edgelength_{uv} = 1$ for all $(u, v) \in \edgeset$ and $\nodeweight_u = 1$ for all $u \in \internalnodes$.
\end{definition}
Finally, note that when measures $\edgeweight$ and $\nodeweight$ are present on the edges and nodes, they naturally induce weighted inner products:
\begin{equation}
\langle F, G \rangle_{\edgeweight} = \sum_{(u, v) \in \edgeset} \edgeweight_{uv} F(u, v) G(u, v), \qquad
\langle f, g \rangle_{\nodeweight} = \sum_{u \in \internalnodes} \nodeweight_u f(u) g(u).
\end{equation}
These inner products can be useful when studying generalized derivatives and dual formulations.

\begin{remark}\label{remark_discrete_to_continuous_anistropic_plap}
We point out here that the graph operator and the differential operator differ significantly in the type of diffusion they induce. Indeed, by looking at \cref{def:graph_p-Lap_operator}, we observe that the diffusion coefficient $|Kf(u,v)|^{p-2}$ depends on the direction that we consider. Differently, in the differentiable $p$-Laplacian operator $\plap^{\text{diff}}f:=-\text{div}(|\nabla f|^{p-2}\nabla f)$, the diffusion coefficient $|\nabla f|^{p-2}$ depends on the point but not on the direction. As a consequence, it is easy to prove e.g. that knowing the graph $p$-Laplacian spectrum on a path graph $P_{n_i}$ of order $n_i$, then the $p$-Laplacian spectrum of the grid graph $P_{n_1}\square \cdots\square  P_{n_j}$, obtained as the product of $j$ path graphs of different order, includes the sum of $p$-eigenvalues on $P_{n_i}$, i.e.,  $$\mathrm{spec}(\Delta_p(P_{n_1}\square \cdots\square  P_{n_j})) \supset \left\{\sum_{i=1}^j\lambda(\Delta_p(P_{n_i}))\right\},$$
where $\lambda(\Delta_p(P_{n_i}))$ varies among the $p$-Laplacian eigenvalues on the line graph $P_{n_i}$. Because of the different nature of the graph operator and the differential one, the same result is not known to hold in the infinite-dimensional setting, where it seems unlikely to hold. The only case we know it holds is the linear one. Indeed, it is well known that $\Delta_2$ on the one-dimensional interval $[0,1]$ admits the eigenfunction $f_j(t)=\sin (j \pi t)$ with corresponding eigenvalue $\lambda_j=(\pi j)^2$. Then, on the rectangle $[0,1]\times
\ldots \times [0,1]\subset \mathbb{R}^n$, $\Delta_2$ can be proved to have eigenvalues
$\pi^2\left(j_1^2+\ldots + j_n^2\right)$ for positive integers $j_1,\ldots ,j_n$ with corresponding eigenfunction $f(x)=f_1(x_1)\cdots f_n(x_n)=\sin (j_1 \pi x_1)\cdots\sin (j_1 \pi x_n)$. 
Nonetheless, we highlight that various authors have studied variations of the standard differential $p$-Laplacian operator in the infinite-dimensional setting, which account for different diffusion coefficients in different directions. These operators are commonly referred to as anisotropic $p$-Laplacian operators or orthotropic $p$-Laplacian operators. To the best of our knowledge, the consistency of the different differential 
$p$-Laplacian operators with respect to various graph-to-continuum limit notions has not yet been investigated.
\end{remark}

%
%
%
%
%
\subsection{$p$-Laplacian Spectrum count: criticisms and open problems}\label{subsec:p-lap_spectrum}

In this subsection we highlight that, unlike the linear case $p = 2$, where $\Delta_2$ is a symmetric positive semidefinite matrix, for a generic $p$ the $p$-Laplacian operator $\plap$ exhibits several key differences. In particular:
\begin{itemize}
    \item The number of $p$-Laplacian eigenpairs can exceed the dimension of the space, i.e., $|\internalnodes|$.
    \item The eigenpairs are generally not orthogonal.
    \item There is no clear notion of eigenspace or multiplicity.
\end{itemize}
These issues have been explored in \cite{DEIDDA2023_nod_dom, zhang2021homological, Amghibech1}.

For instance, in \cite{DEIDDA2023_nod_dom}, the authors provide an example of a graph with 4 nodes that has at least 5 distinct eigenvalues. In \cite{Amghibech1}, the author computes all the eigenvalues of $\plap$ for unweighted complete graphs and shows that their number is given by:
\[
1 + \#\{\{i, j\} : i + j \leq N, \; i, j \in \mathbb{Z}_+\} = 1 + \lfloor N / 2 \rfloor (N - \lfloor N / 2 \rfloor) = \lfloor N^2 / 4 \rfloor + 1,
\]
where $N$ is the cardinality of the node set.

In \cite{zhang2021homological}, the author introduces the concept of $p$-Laplacian homological eigenvalues.    This class of eigenvalues is shown to include the variational eigenvalues, defined as the min-max of the Rayleigh quotient $\rayl_p$ over subspaces of fixed genus (a generalization of the notion of dimension) in the range from $1$ to $|\nodeset|$. Since each variational eigenvalue corresponds to a given index between $1$ and $N$, the number of these eigenvalues is equal to the dimension of the space and so of the graph. We will discuss these eigenvalues extensively in \cref{subsec:variational_spectrum}. However, importantly, the study in \cite{zhang2021homological} demonstrates that, for certain graphs and $p$ values close to $1$, there exist homological eigenvalues that are not variational, thereby proving the existence of more eigenvalues than the dimension of the space.

\begin{figure}[t]
    \centering
    \includegraphics[width=0.32\textwidth]{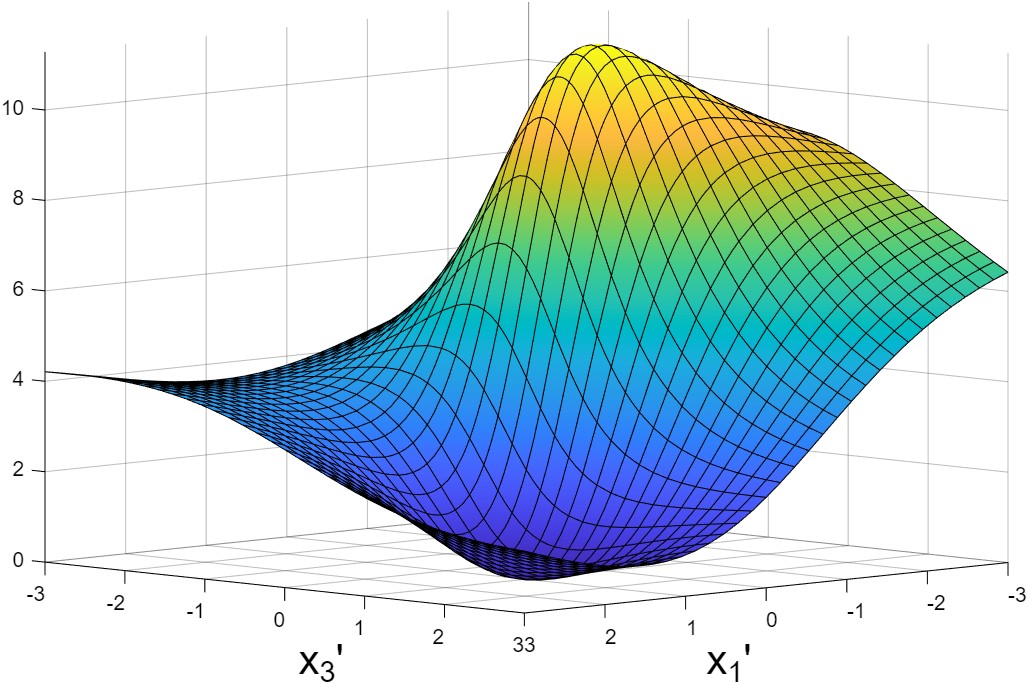}
    \includegraphics[width=0.32\textwidth]{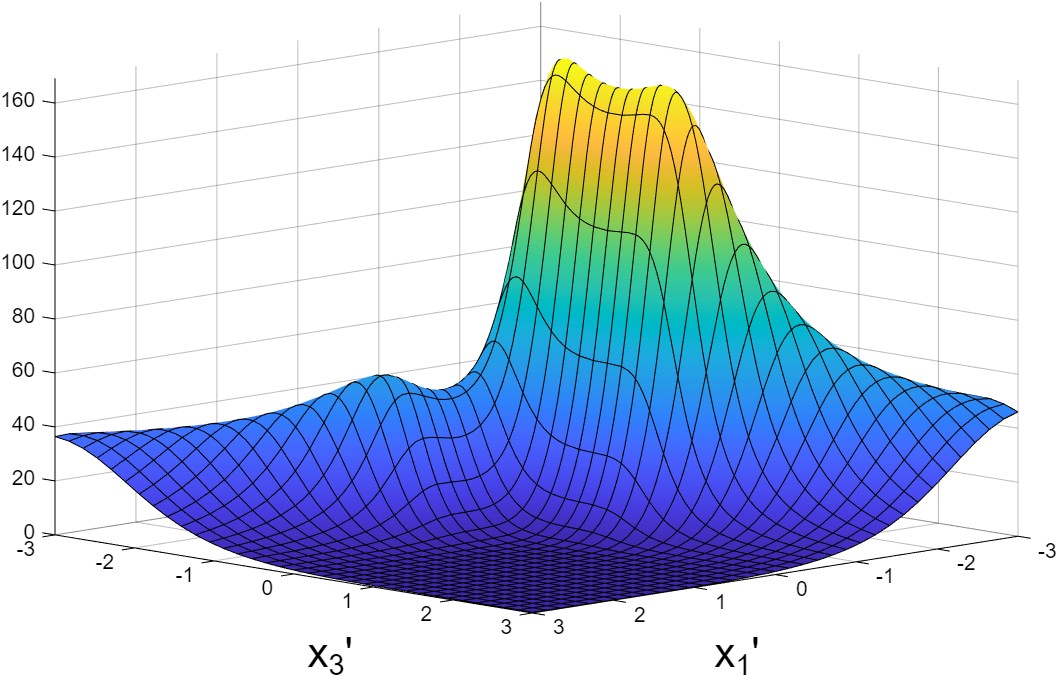}
    \includegraphics[width=0.32\textwidth]{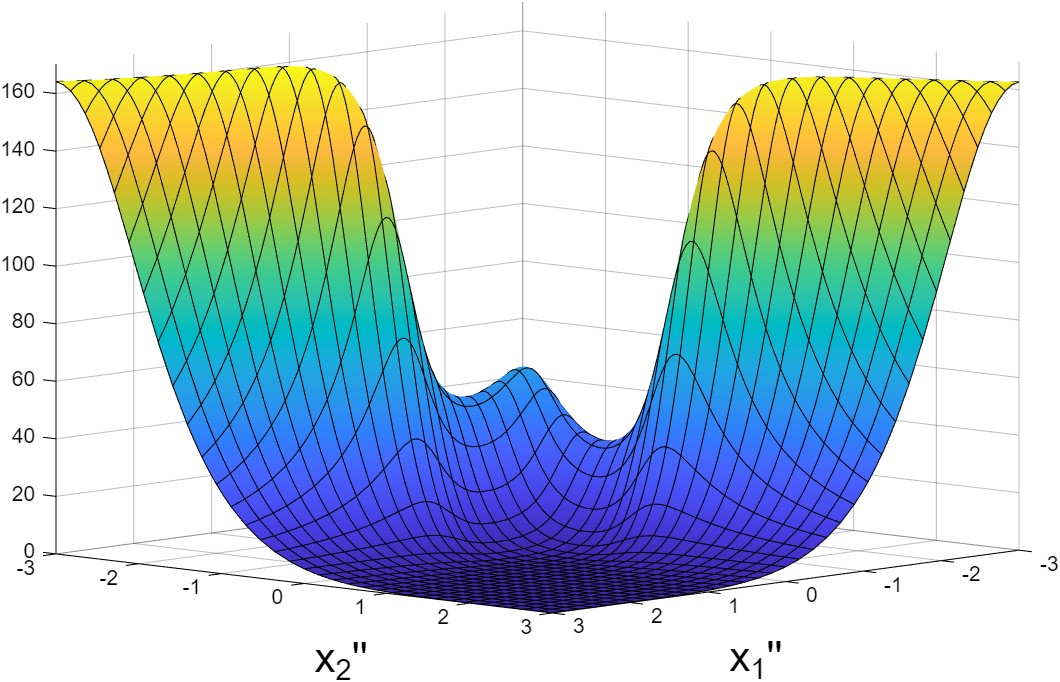}
    \caption{
    We consider a complete graph on three vertices, with edge weights $1^{1/p}$, $1.5^{1/p}$, and $5^{1/p}$. The corresponding Rayleigh quotient:
    $\rayl_p^p(x) = \big(5|x_1 - x_2|^p + |x_1 - x_3|^p + 1.5|x_2 - x_3|^p\big)/\big(|x_1|^p + |x_2|^p + |x_3|^p\big)$
    is plotted from left to right for $p = 2$ in the projective space $x_2 \neq 0$ as a function of $x_1'=x_1/x_2$ and $x_3'=x_1/x_2$, for $p = 6$ in the projective space $x_2 \neq 0$ as a function of $x_1'=x_1/x_2$ and $x_3'=x_1/x_2$, and for $p = 6$ in the projective space $x_3 \neq 0$ as a function of $x_1''=x_1/x_3$ and $x_2''=x_2/x_3$. In the linear case ($p = 2$), we observe the existence of only three critical points: one maximum, one minimum, and one saddle point. Conversely, in the nonlinear case ($p = 6$), we observe three saddle points, three maxima, and one minimum at $x_1 = x_2 = x_3$, yielding a total of seven distinct eigenvalues.}
    \label{fig:example-triangle-graph}
\end{figure}

In the examples above, the $p$-Laplacian operator is assumed to be unweighted. To provide an intuitive illustration of this issue, in Figure~\ref{fig:example-triangle-graph} we present a graphical plot of the $p$-Rayleigh quotient for a weighted complete graph of three nodes without boundary. We observe that, when considering weighted graphs with different weight functions, the number of distinct eigenvalues can further increase.

Note, however, that if $p > N - 1 \geq 1$, where $N$ is the number of vertices, then $\rayl_p(\cdot)$ is $C^{N-1}$-smooth. In particular, by Sard's Theorem, the set of critical values of $\rayl_p|_{S_2}$ has measure zero, where $S_2=\{f\,|\,\|f\|_2=1\}$ and since $\rayl_p$ is zero homogeneous it is sufficient to study the critical points $f$ with $\|f\|_2=1$. Consequently, the spectrum of the $p$-Laplacian has measure zero for any $p > N - 1$, and since the spectrum must be topologically closed and compact, it is indeed nowhere dense. 
Indeed, Sard's Theorem states that if $f: M \to M'$ is at least $C^r$-smooth with $r > \min\{0, m - m'\}$, then the set of critical values of $f$ has measure zero in $M'$, where $M$ and $M'$ are smooth manifolds of dimensions $m$ and $m'$, respectively.
In the case of the $p$-Laplacian we let $M$ be the smooth sphere in $\R^N$ (thus $m=N-1$), and let $M'=\R$  (thus $m'=1$). 
Then we know that $\rayl_p$ is  $C^{\lfloor p\rfloor}$-smooth and thus it is at least $C^{N-1}$-smooth when $p>N-1$.

The finiteness of the $p$-Laplacian spectrum on general graphs, as well as the existence of upper bounds for the cardinality of the spectrum, remains a significant open problem.
\begin{openproblem}\label{open_problem_spectrum_finiteness}
Is the graph $p$-Laplacian spectrum always finite? 
\end{openproblem}
{\it Update: We notice that a recent preprint \cite{Colbrook}, posted while this manuscript was in review, addressed and confirmed the above open problem about the finiteness of the $p$-Laplacian spectrum for any value of $p$ and provided an upper bound for its cardinality.}

This problem has potentially far-reaching implications. For instance, in \cite{zhang2021homological}, it is shown that isolated homological eigenvalues are continuous in $p$. However, since the $p$-Laplacian spectrum is always bounded in the graph setting, the existence or non-existence of isolated eigenvalues is closely tied to the problem of the finiteness of the $p$-Laplacian spectrum.

So far, all the results and examples suggest a positive answer to the open problem \ref{open_problem_spectrum_finiteness}. Indeed
structured examples, such as complete graphs~\cite{Amghibech1} and trees~\cite{DEIDDA2023_nod_dom, Tudisco1, zhang2023pLap_noddom} exhibit finite spectra. In addition, reformulating the graph $p$-Laplacian eigenvalue problem in terms of a tensor eigenvalue problem \cite[Proposition 7.1]{ge2023new}, the authors of \cite{ge2025computing} further proved the finiteness of the $p$-Laplacian spectrum under generic condition for any even value of $p$. Last we note that for any rational number $p > 1$, we can write $p - 1 = k / m$, where $k, m \geq 1$ are integers. Using the notation $[a]^p = |a|^p \mathrm{sign}(a)$, we define $[g(u)] = [f(u)]^{1/m}$ for all $u \in \nodeset$. In terms of the new variable $g$, the $p$-Laplacian eigenproblem becomes:
\[
\sum_{v \sim u} \edgelength^p_{uv} \big[[g(u)]^m - [g(v)]^m\big]^{k/m} = \lambda [g(u)]^k.
\]
Introducing a further variable $h: \edgeset \to \R$, defined by $[h(u, v)]^m = [g(u)]^m - [g(v)]^m$ for all $(u, v) \in \edgeset$, the Rayleigh quotient $\rayl_p^p(f)$ can be rewritten as:
\begin{equation}
\rayl_p^p(f) = \frac{\sum_{(u, v) \in \edgeset} \edgelength^p_{uv} |h(u, v)|^{m + k}}{\sum_{u \in \nodeset} |g(u)|^{m + k}}.
\end{equation}
Recalling that if $(f, \lambda)$ is a $p$-eigenpair, then $\lambda = \rayl_p^p(f)$, the $p$-Laplacian eigenvalue equation becomes:
\begin{equation}\label{rational_p_lap_eigenequation}
\begin{cases}
\displaystyle{\Big(\sum_{u \in \nodeset} |g(u)|^{m + k}\Big) \sum_{v \sim u} \edgelength^p_{uv} [h(u, v)]^k = \Big(\sum_{(u, v) \in \edgeset} |h(u, v)|^{m + k}\Big) [g(u)]^k} & \forall u \in \nodeset, \\[0.8em]
[h(u, v)]^m = [g(u)]^m - [g(v)]^m & \forall (u, v) \in \edgeset.
\end{cases}
\end{equation}

Now note that, up to the position of a function $g$ in $\R^{|\nodeset|}$, for any two nodes $u$ and $v$, $|g(u) - g(v)|$ can be written as $\pm (g(u) - g(v))$, and analogously for $|g(u)| = \pm g(u)$. The regions of $\R^{|\nodeset|}$ where the signs of $\big(g(u) - g(v)\big)$ and $g(u)$ are fixed for all $(u, v) \in \edgeset$ and $u \in \nodeset$ are defined by:
\begin{equation}
g(u_{i_1}) \geq g(u_{i_2}) \geq \dots \geq g(u_{i_N}),
\end{equation}
and by the sign of the nodes:
\begin{equation}
g(u_{i_1}), \dots, g(u_{i_k}) \geq 0, \quad g(u_{i_{k+1}}), \dots, g(u_{i_N}) \leq 0.
\end{equation}
Within each such region, since all signs are prescribed, the eigenequation \eqref{rational_p_lap_eigenequation} reduces to a system of $|\nodeset|$ $(m + 2k)$-homogeneous polynomial equations and $|\edgeset|$ $m$-homogeneous polynomial equations. Because of the homogeneity of the system of equations, the problem of finding its solutions can be restricted,  without loss of generality, to the projective space $\mathbb{P}^{|\nodeset| + |\edgeset| - 1}$. By Bézout's Theorem, assuming the polynomial equations are in general position, there are at most $(m + 2k)^{|\nodeset|} k^{|\edgeset|}$ solutions in $\mathbb{P}^{|\nodeset| + |\edgeset| - 1}$ for each $\pm$-assignment. Since there are $2^{|\nodeset| + |\edgeset|}$ such regions, the total number of solutions is bounded by:
\[
2^{|\nodeset| + |\edgeset|} (m + 2k)^{|\nodeset|} k^{|\edgeset|}.
\]

Thus, for any rational $p \geq 1$, the $p$-Laplacian operator likely has a finite number of eigenvalues under generic conditions.  This discussion leads us to conjecture the finiteness of the graph $p$-Laplacian spectrum and the following open problem 
\begin{openproblem}
    If the $p$-Laplacian spectrum is finite, how many eigenvalues can it have?
\end{openproblem}

\subsection{The variational spectrum}\label{subsec:variational_spectrum}
Despite the criticisms highlighted so far, classical results from the calculus of variations guarantee that it is always possible to characterize a set of ``variational'' eigenvalues, whose number, counted with their multiplicity, is equal to the dimension of the space, $|\internalnodes|$. 
To define these eigenpairs, we observe that, due to the homogeneity of the Rayleigh quotient $\rayl_p$, the study of its critical points can be restricted to the $p$-sphere:
\[
S_p := \{f \in \Hc(\internalnodes) \;|\; \|f\|_p = 1\}.
\]
The proof of the existence of variational eigenvalues is a direct consequence of the deformation lemma presented below. This celebrated result, along with its direct consequence in Theorem \ref{Critical_point_theorem}, represents a particular case of more general and classical results (see, e.g., \cite{palais1970critical, struwe, Ghoussoub, Fucik}). For this reason, we do not include the proof.

Defining $
\rayl_p^c := \{f\in S_p \;|\; \rayl_p(f) \le c\}$, we say that $c$ is a critical value of $\rayl_p$ if there exists some $f$ such that $\nabla_f \rayl_p(f) = 0$ and $ \rayl_p(f) = c$. Any value $c$ that is not a critical value is referred to as a regular value of $\rayl_p$. We have
\begin{lemma}[Deformation Lemma]\label{Deformation_lemma_Intro}
Assume $c$ is a regular value of $\rayl_p$. Then, there exist $\epsilon > 0$ and a continuous family of deformations $\phi \in C([0,1] \times S_p, S_p)$ such that:
\begin{itemize}[noitemsep,topsep=0pt,leftmargin=*]
    \item $\phi(t, f) = -\phi(t, -f) \quad \forall (t, f)$,
    \item $\phi(1, \rayl_p^{c+\epsilon}) \subset \rayl_p^{c-\epsilon}$,
    \item $\phi(0, f) = f$.
\end{itemize}
\end{lemma}

As a direct consequence of the Deformation Lemma, we can provide sufficient conditions for the existence of critical values of $\rayl_p$, (see e.g. \cite{struwe}, where the result is discussed for general functionals, for the proof)

\begin{theorem}\label{Critical_point_theorem}
Let $\mathcal{F}$ be a family of subsets of $S_p$ such that, for any regular value $c \in \R$ of $\rayl_p$, there exist $\epsilon > 0$ and a continuous deformation of the domain $\phi : [0, 1] \times S_p \to S_p$ satisfying:
\begin{equation}
\begin{cases}
\phi(0, \cdot) = id_{S_p}(\cdot), \\
\phi(1, \rayl_p^{c+\epsilon}) \subset \rayl_p^{c-\epsilon}, \\
\phi(t, A) \in \mathcal{F}, \;\forall A \in \mathcal{F}, \; \forall t \in [0, 1].
\end{cases}
\end{equation}
Then, the value:
\[
\Lambda := \underset{A \in \mathcal{F}}{\inf} \; \underset{f \in A}{\sup} \rayl_p(f),
\]
is a critical value of $\rayl_p$, i.e., the $p$-th root of an eigenvalue of $\plap$.
\end{theorem}

{We would like to point out that the last \cref{Critical_point_theorem} mentioned above provides a constructive definition for certain critical values of $\rayl_p$, but it does not necessarily provide a constructive method for studying the critical points of $\rayl_p$.} Indeed, even if $A\in \mathcal{F}$ and $f\in A$ realize the min-max in \cref{Critical_point_theorem}, $f$ could not be a critical point. This is, e.g., the case when $A$ is the whole sublevel set $\rayl_p^{\Lambda}$.

Based on the above theorem, it is now possible to introduce the variational eigenpairs of the $p$-Laplacian operator.
Theorem \ref{Critical_point_theorem} indicates that to define the variational eigenvalues, we must construct $N=|\internalnodes|$ families, $\{\mathcal{F}_k\}$, of subsets of $S_p$ that are stable under deformations. That is, if $A \in \mathcal{F}_k$ and $\phi$ is a deformation, then $\phi(A) \in \mathcal{F}_k$.

To understand how this works, let us recall the Fisher-Courant $\min\max$ characterization of the eigenvalues of a symmetric matrix (e.g., the graph Laplacian). Specifically:
\[
\lambda_k(\Delta_2) = \min_{\mathrm{dim}(A) \geq k} \max_{f \in A \setminus \{0\}} \frac{\langle \Delta_2 f, f \rangle}{\langle f, f \rangle} = \min_{\mathrm{dim}(A) \geq k} \max_{f \in A \setminus \{0\}} \rayl_2(f)^2.
\]

In the $p$-Laplacian setting, the definition of variational eigenvalues extends this $\min\max$ principle by considering families of subsets of $\R^{N}$ that are more general than linear subspaces of fixed dimension and are closed under the deformations described in \Cref{Critical_point_theorem}.

We explore different ways to construct such families of subsets, which allow us to generalize the notion of eigenvalues from the linear case ($p=2$) to the $p$-Laplacian operator ($ p> 2$). In the following, we discuss the most standard class of variational eigenvalues, which is given by the Krasnoselskii eigenvalues. We postpone to \Cref{appendix_variational_spectrum} the discussion of different other classes of variational eigenvalues that have been introduced and studied in the literature.

\subsubsection{The Krasnoselskii spectrum}\label{sec:Krasnoselskii}

The most commonly studied family of variational eigenvalues is the Krasnoselskii spectrum. Its definition is based on the idea of using a generalized notion of dimension, known as the Krasnoselskii genus, which is related to the Lyusternik–Schnirelmann category of a space~\cite{struwe, Ghoussoub, Fucik}. 
To construct this family of eigenvalues, we first note that since we are interested in studying critical points of $\rayl_p$, which is an even functional, it suffices to generalize the notion of dimension to symmetric subsets.
Thus, we introduce the family $\mathcal{A}$ of subsets of $\R^N$ that are symmetric and closed, i.e.:
$$\mathcal{A}=\lbrace A \subseteq \R^N | \; A\;\: \mathrm{closed}\,,\;\: A=-A\rbrace\,.$$
Then, we observe that in the case $A$ is a linear subspace of dimension $k$, $A\setminus\{0\}$ can be mapped with continuity onto the unit sphere of dimension $k-1$, 
\begin{equation}
    S^{k-1}_2:=\{x\in \R^k|\;\|x\|_2=1\}.
\end{equation}
This notion can be generalized by defining, for any $A\in\mathcal{A}$, the Krasnoselskii genus of $A$: 
\[
\gamma(A)=\begin{cases}
0 \qquad  &\text{if}\; A=\emptyset\\
\inf\lbrace k\in\mathbb{N}\;\;| & \exists\,\psi\in C(A,S^{k-1}_2)\;\, s.t. \;\, \psi(f)=-\psi(-f)\rbrace  \\
+\infty \qquad &\text{if $\nexists \; k$ as above}
\end{cases}\,.
\]
Note that, if $\gamma(A)\geq k$ and $\phi\in C(\R^N,\R^N)$ is a symmetric map, then $\gamma(\phi(A))\geq \gamma(A)$.
Hence, the families $\mathcal{F}_k(S_p)=\lbrace A\in \mathcal{A} \,|A\subset S_p,\;\gamma(A)\geq k \rbrace$, where $S_p=\{f\,|\,\|f\|_p=1\}$,  satisfy the hypotheses of Theorem \ref{Critical_point_theorem}.
Thus, we can define the Krasnoselskii variational eigenvalues of $\plap$ as
\begin{equation}\label{Krasnoselskii_variationa_eigenvalues_Intro}
\Lambda_k=\lambda_k^{\frac{1}{p}}=\inf_{A\in\mathcal{F}_k}\sup_{f\in A}\mathcal{R}_{p}(f)\,.
\end{equation}
Since we are working in a finite-dimensional space, it is also possible to prove that the above $\inf\sup$ formulation is actually a $\min\max$.
The advantage of defining eigenvalues in this manner is twofold. First, it ensures that the number of selected eigenvalues matches the dimension of the space. Second, due to the properties of the Krasnoselskii genus, we recover certain multiplicity results recalled below (e.g., Lemma \ref{lemma_genus_and_multiplicity} and Proposition \ref{prop.5.3_struwe} ) from the linear case (see also~\cite{DEIDDA2023_nod_dom, ZhangNodalDO}). 
Indeed, since the variational eigenvalues $\{\lambda_k\}$ form an increasing sequence by definition, it is possible to introduce a first notion of variational multiplicity (or simply the multiplicity of a variational eigenvalue) as follows:
%
\begin{definition}\label{DEf_variational_multiplicty}
Let $\lambda _{k}$ be a variational eigenvalue. If $\lambda _{k}$ appears $m$ times in the sequence of the variational eigenvalues
\begin{equation}
\lambda _{1}\leq \lambda _{2}\leq \dots \leq \lambda _{k-1}<\lambda _{k}
=\dots = \lambda _{k+m-1}<\lambda _{k+r}\leq \dots \leq \lambda _{N}
\,,
\nonumber
\end{equation}
we say that $\lambda _{k}$ has (variational) multiplicity $m$ and write
$\mathrm{mult}(\lambda _{k})=m$.
\end{definition}

Clearly, the notion of variational multiplicity applies only to the variational spectrum. However, for a generic eigenvalue $\lambda $, it is possible to introduce a notion of ``geometric'' or $\gamma$-multiplicity by means of the Krasnoselskii
genus of the corresponding eigenspace:
%
\begin{definition}\label{DEf_gamma-mult}
Assume $\lambda $ is a $p$-Laplacian eigenvalue, if
\begin{equation}
\gamma \Big(\{f\in S_{p} \,:\,\plap(f)=\lambda |f|^{p-2}f
\}\Big)=m,
\nonumber
\end{equation}
we say that $\lambda $ has $\gamma $-multiplicity $m$ and we write
$\gamma \text{-}\mathrm{mult}(\lambda )=m$.
\end{definition}
Within this setting, we say that an eigenvalue $\lambda$ is simple if there exists a unique eigenfunction (up to a $-1$ multiplicative factor) $f\in  S_{p}$ associated to $\lambda$. Note that, if $\lambda $ is a simple eigenvalue, then necessarily $\gamma \text{-}\mathrm{mult}(\lambda )=1$ and, if
$\lambda $ is variational then also $\mathrm{mult}(\lambda )=1$.
Vice-versa, an eigenvalue $\lambda$ with $\gamma \text{-}\mathrm{mult}(\lambda )=1$ is
not necessarily simple. See, for instance, Example \ref{Ex:p-Lap_eingenspace_varying_p_on_complete_graphs}.
Observe that the definitions of variational multiplicity and $\gamma$-multiplicity recall the concepts of algebraic and geometric multiplicity of a linear eigenvalue, respectively. However, while in the linear case of a self-adjoint operator, like $\Delta_2$, algebraic and geometric multiplicities coincide, this is not the case of the $p$-Laplacian operator, where the $\gamma$-multiplicity can be larger than the variational multiplicity. The proof of the next lemma follows directly
from Lemma 5.6 and Proposition 5.3, see also Chapter II of \cite{struwe}:
\begin{lemma}\label{lemma_genus_and_multiplicity}
If $\lambda $ is a variational eigenvalue, then
\begin{equation}
\gamma \text{-}\mathrm{mult}(\lambda )\geq \mathrm{mult}(\lambda )\, .%
\end{equation}
\end{lemma}
However, as a corollary of the proposition 
that we report below and the multiplicity inequality above, we have that any variational
eigenvalue $\lambda$ is associated to, at least,
$\mathrm{mult}(\lambda )$ orthogonal eigenfunctions, similarly to the linear case.

\begin{proposition}\label{prop.5.3_struwe} (See Proposition 5.3 Chapter II \cite{struwe}). Suppose $A\subset H$ is a compact symmetric subset of a Hilbert space $H$ and suppose $\gamma(A)=m<\infty\,$. Then $A$ contains at least $m$ mutually orthogonal vectors $\{v_i\}_{i=1}^m$.
\end{proposition}

Equality or difference of the $\gamma$-multiplicity and variational multiplicity is a delicate subject. Indeed for a generic $p\in(1,\infty)$ it holds the following open problem
\begin{openproblem}\label{open_prob_gamma_mult_vs_var_mult}
    Does there exists a graph $\Gc$ and a $p$-Laplacian variational eigenvalue $\lambda$, for some $p\in(1,\infty)$, such that 
    $$\gamma\text{-mult}(\lambda)>\text{mult}(\lambda)?$$
\end{openproblem}
At this point in the paper, we have not yet formally introduced the eigenvalues and variational eigenvalues in the extremal cases $ p=1$ and $ p=\infty$. We will introduce them in \cref{subsec:infinity_1_Notation}. However, we can anticipate that \cref{open_prob_gamma_mult_vs_var_mult} is formulated only for $p\in (1,\infty)$, not by chance. Indeed, it can be proved that in the case $p=1,\infty$ there exist variational eigenvalues whose $\gamma$-multiplicity is strictly larger than the variational one, see e.g. \cref{exam:path}. 
On the other hand, in \cite{DEIDDA2023_nod_dom} it has been proved that on trees all the eigenvalues are variational and their $\gamma$-multiplicity and variational multiplicity are the same.
\begin{theorem}[Theorem 3.7 in \cite{DEIDDA2023_nod_dom}]\label{variational_eigenvalues_on_trees}
Let $\Gc$ be a tree; then $\Delta_p$ for $p\in(1,\infty)$ admits only variational eigenvalues and any such eigenvalue $\lambda$ is such that 
$\gamma\text{-mult}(\lambda)=\text{mult}(\lambda).$
\end{theorem}

We conclude by stating the following technical lemma from~\cite{Solimini}, which provides an upper bound for the 
Krasnoselskii genus of the intersection of different subsets. This result will be useful in later derivations.
\begin{lemma}{\cite[Prop.\ 4.4]{Solimini}}\label{Lemma_Krasnoselskii_intersection}
  Let $X$ be a Banach space and $\mathcal A$ the class of the closed symmetric subsets of $X$. Given $A\in \mathcal{A}$,
    consider a Karsnoselskii test map $\varphi: A\to \mathbb R^k$ continuous and such that  $\varphi(x)=-\varphi(-x)$, where
    $k<\gamma(A)$. Then,
  $\gamma(\varphi^{-1}(0))\geq \gamma(A)-k$.
\end{lemma}
\begin{remark}[The second variational eigenvalue]\label{characterization_2nd_variational_eigenvalue}
We remark that the second variational eigenvalue admits a further characterization. It is known from, e.g.~\cite{Tudisco1,DEIDDA2023_nod_dom}, that $\lambda_2$ is the smallest eigenvalue strictly greater than $\lambda_1$, meaning that no eigenvalues exist in the interval $(\lambda_1, \lambda_2)$. This leads to the following characterization via the mountain pass theorem:
\begin{align}\label{Eq_lambda_2_on_curves}
\lambda_2(\Delta_p) &= \inf\limits_{\substack{\text{curve } c: [-1,1] \to \mathbb{R}^N \setminus \{0\},\\
c(\pm 1) = \pm f_1}} \max\limits_{f \in c([-1,1])} \mathcal{R}_p^p(f) 
= \inf\limits_{\substack{S \text{ is a centrally symmetric} \\ \text{simple closed curve in } \R^N}} \max\limits_{f \in S} \mathcal{R}_p^p(f),
\end{align}
where $f_1$ is an eigenfunction corresponding to the first eigenvalue of $\Delta_p$, and we identify $C(\internalnodes)$ with $\R^N$. The final equality follows from the fact that any symmetric subset of genus $2$ must be connected and, therefore, contains a centrally symmetric simple closed curve, which in turn is a subset of genus $2$. See also~\cite{Cuesta99} for a discussion of this problem in the infinite dimensional setting.
Moreover, in the absence of a potential term (i.e., when there is no boundary), we have the additional equality~\cite{zhang2021discrete, Bhuler}:
\begin{equation}
\lambda_2(\Delta_p) = \min\limits_{f \text{ nonconstant}} \max\limits_{t \in \R} \mathcal{R}_p^p(f+t)
= \min\limits_{f \perp 1} \frac{\sum_{(u,v) \in E} |\grad f(u,v)|^{p}}{2 \min\limits_{t \in \R} \|f - t\|_p^p}.
\end{equation}

Interestingly, this characterization of the second variational eigenvalue allows one to investigate the spectral gap of the graph $p$-Laplacian operator, i.e., the distance between the second variational eigenvalue from zero on a graph without boundary. A preliminary investigation of this topic has been conducted in \cite{Berkolaiko_2017} for unweighted graphs. In that work, the authors provide bounds on the second variational eigenvalue of the $p$-Laplacian in terms of the edge connectivity, denoted by $\eta$, both for discrete and quantum graphs. We recall that the edge connectivity of $\Gc$ is the minimum number of edges that must be removed to disconnect the graph. In more details, \cite{Berkolaiko_2017} show that the second variational $p$-Laplacian eigenvalue on an unweighted graph $\Gc$ can be bounded from below as $\lambda_2(\Delta_p)\geq \eta {\lambda}_2^{\text{PATH}}(\Delta_p)$, where ${\lambda}_2^{\text{PATH}}(\Delta_p)$ denotes the second variational eigenvalue of the graph $p$-Laplacian on a path graph having the same number of nodes as $\Gc$. Furthermore,  the authors prove that equality holds if and only if $\Gc$ is an $\eta$-regular pumpkin chain, i.e., a path graph where each node is connected to the following one by exactly $\eta$ edges.
\end{remark}

Several other families of variational eigenvalues can be introduced by modifying the classes of subsets considered in the $\min\max$ characterization, as we discuss in appendix \cref{appendix_variational_spectrum}.





\subsection{Degenerate cases $p=1,\infty$}\label{subsec:infinity_1_Notation}

A separate section is necessary to discuss the $p$-Laplacian eigenvalue problem in the two extremal cases $p=1$ and $p=\infty$.
Indeed, in these cases the $p$-Laplacian operator cannot be defined as in \cref{def:graph_p-Lap_operator}; however, the Rayleigh quotients $\rayl_1(f)$ and $\rayl_{\infty}(f)$ are still well defined, even though they are no longer differentiable.
This raises the problem of how to define the $1$- and the $\infty$-eigenpairs.
The answer to this problem is not unique, and several approaches have been proposed in the literature.
Here, we discuss an approach that was initially introduced for the case $p=1$ \cite{chang2016spectrum, hein2010inverse}, but has recently been extended also to the case $p=\infty$ \cite{Bungert1,bungert2021eigenvalue, deidda2024_inf_eigenproblem}.
The idea is to define a generalized notion of critical points for the Rayleigh quotients $\rayl_1(f)$ and $\rayl_{\infty}(f)$.
Of course, different notions of generalized critical points can be introduced. However, not all of them are equally amenable to analysis. Here we focus on a standard approach, which in practice is also the most tractable. Nevertheless, in the remainder of the section we briefly comment on other possible definitions of generalized critical points, see \cref{Clarke_eigenpairs}.

Note that in the $p$-Rayleigh quotient both the numerator $\|\grad f\|_p$ and the denominator $\|f\|_p$ are convex functions of $f$. Thus, even when they are not differentiable, they admit a subgradient.
We briefly recall the definition of the set-valued subgradient $\partial$ of a convex function $\Psi:\R^N\rightarrow\R$ at a point $f_0$:
\begin{equation}\
\partial \Psi(f_0)=\big\{\xi\,|\;\Psi(g)-\Psi(f_0)\geq\, \langle\xi,g-f_0\rangle\;\forall g\in \R^N\big\}\,.
\end{equation}
We also recall that this is a consistent generalization of the notion of gradient, in the sense that if the function $\Psi$ is differentiable at the point $f_0$, then $\partial\Psi(f_0)$ coincides with the usual gradient in $\R^N$.

Using subgradients, we can therefore generalize the notion of $p$-eigenpairs in a meaningful way for any $p\in[1,\infty]$.

\begin{definition}\label{Def:generalized_p-eigenpair}
    $f$ is a $p$-eigenfunction and $\Lambda$ is a normalized $p$-eigenvalue if and only if
    $$\emptyset \neq\partial\|\grad f\|_p\cap\Lambda\partial\|f\|_p\,.$$
\end{definition}
In particular, it is not difficult to observe that if $p\in(1,\infty)$ and $(f, \Lambda)$ satisfy the condition in \Cref{Def:generalized_p-eigenpair}, then $(f,\lambda)$ with $\lambda=\Lambda^p=\rayl_p^p(f)$ satisfies the eigenvalue equation in \Cref{Def_p-eigenpair_by_eq}. Thus, the notions of $p$-eigenfunction and normalized $p$-eigenvalue introduced in \Cref{Def:generalized_p-eigenpair} extend the definition given in \Cref{Def_p-eigenpair_by_eq} to the full range $p \in [1,\infty]$.
In particular, with a slight abuse of notation, for the two extremal cases $p=1$ and $p=\infty$ we call $\Lambda$ in \Cref{Def:generalized_p-eigenpair} the eigenvalue corresponding to $f$, omitting the attribute ``normalized''. This is motivated by the fact that for $p=1$ there is no distinction between eigenvalues and normalized eigenvalues, while for $p=\infty$ only normalized eigenvalues can be defined.

%
%

Next, we recall some useful results about subgradients that allow us to characterize explicitly the $p$-Laplacian eigenvalue problem also in the limit cases $p=1$ and $p=\infty$. From \cite{burger2016spectral}, we recall a characterization of the subgradient of a convex and absolutely $1$-homogeneous function $\Psi$ on $\R^N$, where by absolutely $1$-homogeneous we mean $\Psi(c f)=|c|\Psi(f)$ for all $f\in \R^N$ and $c\in \R$.

\begin{lemma}[see \cite{burger2016spectral}]\label{Lemma_norm_subgradient}
Let $\Psi$ be a convex and absolutely $1$-homogeneous function on $(\R^N, \langle\cdot, \cdot \rangle$) and let $f_0\in \R^N$. Then 
    $$\partial\Psi(f_0)=\{\xi\;|\; \Psi(g)\geq \langle\xi,g\rangle\;\forall\,	 g,\;\: \Psi(f_0)= \langle\xi,f_0\rangle  \}\,.$$
\end{lemma}
As a first corollary of \cref{Lemma_norm_subgradient} we note that, also in the extremal cases $p=1,\infty$, the eigenvalue matches the value of the Rayleigh quotient computed in the eigenfunction
\begin{corollary}
    Let $(f,\Lambda)$ be a $p$-eigenpair as in \cref{Def:generalized_p-eigenpair}, then $$\Lambda=\rayl_p(f).$$
\end{corollary}
As a second corollary of \cref{Lemma_norm_subgradient}, and by invoking the characterization of the subdifferential of a composition of a convex function with a linear transformation (see Theorem~23.9 in~\cite{rockafellar2015convex}), we obtain a characterization of the subgradient of the function \( f \mapsto \|A f\|_{\mu,p} \), where we are allowing for a strictly positive measure $\mu$ in $\R^M$, with $M$ equal to the dimension of $Af$, see also \cref{weighted_p_lap_eigenvalue_problem}.
\begin{corollary}\label{Thm_subgradient_chain_rule}
Let $\Phi(f)=\|Af\|_{\edgeweight,p}$, where $A$ is a linear transformation from $(\R^N, \langle\cdot,\cdot\rangle_{\nodeweight})$ to $(\R^M,\langle\cdot,\cdot\rangle_{\edgeweight})$, with $\nodeweight$ and $\edgeweight$ entrywise strictly positive, then 
$$\partial \Phi(f)=\Big\{\nodeweight^{-1}\odot\big( A^T \big(\edgeweight\odot \Xi\big)\big) \Big| \; \|\Xi\|_{\edgeweight,p}\geq \langle\Xi,g\rangle_{\edgeweight}\;\forall\,	 g \in \R^N,\;\: \|f_0\|_{\edgeweight,p}= \langle \Xi , f_0\rangle_{\edgeweight} \Big\},$$
where $\nodeweight^{-1}$ is the entrywise inverse of $\nodeweight$.
\end{corollary}

For completeness, we included a proof of the last result in \cref{Appendix:a}.

To complete the discussion about the two nonsmooth cases we characterize the sets $\partial\|f\|_{\nodeweight,1}$, $\partial\|\grad f\|_{1}$, $\partial\|f\|_{\infty}$ and $\partial\|\grad f\|_{\infty}$ where we are considering all these functions as functions from $\mathcal{H}(\internalnodes)=\big(\R^{N}, \langle\cdot,\cdot\rangle_{\nodeweight}\big)$ to $\R$ and we remember that throughout the paper we are assuming $\edgeweight_{uv}=1$ for any $(u,v)\in \edgeset$ (see \cref{weighted_p_lap_eigenvalue_problem}).
The subgradients of the $1$-norms can be calculated from \Cref{Lemma_norm_subgradient} and \Cref{Thm_subgradient_chain_rule}, yielding the following formulas, where $\xi$ and $\Xi$ denote functions in $\Hc(\internalnodes)$ and $\Hc(\edgeset)$, respectively. 

\begin{equation}\label{1_subgradient}
\begin{aligned}
\partial\|f\|_{\nodeweight,1} &=\Big\{\xi\,\Big|\;\xi(u)\in\,\mathrm{sign} \big(f(u)\big) \;\forall u\in \internalnodes\Big\}\\
\partial\|\grad f\|_1 &=\Big\{\nodeweight^{-1}\odot\big(K^T\,\Xi\big)\,\Big|\;\Xi(u,v)\in\mathrm{sign} \big(\grad f(u,v)\big)\;\forall (u,v)\in \edgeset\Big\}
\end{aligned}
\end{equation}
where $\mathrm{sign}(x)$ is the set valued function, $\mathrm{sign}(x)=\begin{cases} 1\quad &\text{if }\; x>0\\[-0.33em]
[-1,1]\quad &\text{if }\; x=0\\[-0.33em]
-1\quad &\text{if }\; x<0\end{cases}\,.$

An analogous computation yields the following formulas for the subgradients of the $\infty$-norm 
\begin{equation}\label{inf_subgradient}
\begin{aligned}
\partial\|f\|_{\infty} &=\left\{ \xi\,\Bigg|\begin{array}{lr}\| \xi\|_{\nodeweight,1}=1,\;\:\xi(u)=0\;\:\text{if}\,\,|f(u)|<\|f\|_{\infty}\\[0.5em] \xi(u)f(u)=\big|\xi(u)f(u)\big|\end{array}\right\}\\[0.5em]
\partial\|\grad f\|_{\infty} &=\left\{ \nodeweight^{-1}\odot K^T\,\Xi\,\Bigg|\begin{array}{lr}\|\Xi\|_1=1/2,\;\Xi(u,v)=0\;\:\text{if}\,\,|\grad f(u,v)|<\|\grad f\|_{\infty}\\[0.5em] \Xi(u,v)\grad f(u,v)=\big|\Xi(u,v)\grad f(u,v)\big|\end{array}\right\}
\end{aligned}
\end{equation}
where we recall that $\|\Xi\|_1=1/2 \sum_{(u,v)\in\edgeset}|\Xi(u,v)|$ and $\nodeweight_u>0$ for any node $u$. 
Note that, as long as the density $\nodeweight$ is everywhere nonzero, by definition $\|f\|_{\infty,\nodeweight}=\max_{u\in\internalnodes}|f(u)|$ and it does not depend on $\nodeweight$. For this reason, in the following, when considering the $\infty$-eigenvalue problem, we will always consider just the case $\nodeweight_u=1$ for any node $u$, without loss of generality.

\begin{remark}
    We want to specify that even if we have formally defined the $1$- and $\infty$-Laplacian eigenpairs, we have not introduced a formal definition of the $1$- and $\infty$-Laplacian operators. Different authors \cite{chang2016spectrum, hein2010inverse} introduce the $1$-Laplacian operator as the set-valued operator $f\rightarrow \partial\|\grad f\|_1$. However, we note that adopting the same choice for $p=\infty$ would lead to several inconsistencies. The first one is that when considering $1<p<\infty$, the definition of the $p$-Laplacian operator given in \cref{def:graph_p-Lap_operator} coincides with the expression $\partial \|\grad f\|_p^p/p$. However, such expression is not well defined in the case $p=\infty$, where we should rather consider $\partial\|\grad f\|_\infty$. The second inconsistency is that there exists an $\infty$-Laplacian operator deeply studied both in the infinite-dimensional and graph setting \cite{Lind2, deidda2024_inf_eigenproblem, elmoataz2011infinity} that is obtained by looking at the limit of the $p$-Laplace equation for $p$ diverging to $\infty$. We will discuss such an operator in \cref{viscosity_eigenpairs}, but we can anticipate that its definition is different from $\partial\|\grad f\|_\infty$.
    For these reasons, we choose not to formally define the $p$-Laplacian operator as $\partial\|\grad f\|_p$ either for $p=1$ or $p=\infty$. Consequently, we also point out that with a small abuse of notation, we use $\Lambda(\Delta_1)$ and $\Lambda(\inflap)$ to denote an eigenvalue that satisfies \cref{DEF:limit_eigenpairs} for $p=1$ and $p=\infty$ respectively. 
    \end{remark}

Finally we note that the generalization of the $p$-Laplacian eigenvalue problem to the degenerate cases $p=1$ and $p=\infty$ presented so far is not the only possible one. In the following subsections, we review alternative extensions of the $p$-Laplacian eigenproblem to the two degenerate cases $p=1$ and $p=\infty$.

\subsubsection{Clarke $p$-eigenpairs and the variational eigenvalues for $p=1,\infty$}\label{Clarke_eigenpairs}
An alternative approach to defining eigenpairs is based on Clarke's generalized derivative. In fact, in the degenerate cases, the term \emph{eigenfunction} is often used to denote Clarke generalized critical points of the $1$- and $\infty$-Rayleigh quotients~\cite{hein2010inverse, zhang2021homological}. Since the Rayleigh quotients are locally Lipschitz functions, their Clarke subgradient $\partial^{Cl} \rayl_{1/\infty}(f)$ is well defined~\cite{clarke1990optimization}. In this framework, $f$ is said to be an eigenfunction if
\begin{equation}\label{Clarke_eigenfunction}
    0 \in \partial^{Cl} \rayl_{1/\infty}(f)\,
\end{equation}
where
\begin{equation}
    \partial^{Cl} \rayl_{1/\infty}(f)=\Big\{\xi \in \R^N \;\Big|\; \langle \xi, v \rangle \leq \lim_{g\rightarrow f}\lim_{h \downarrow 0} \frac{\rayl_{1/\infty}(g+hv)-\rayl_{1/\infty}(g)}{h} \,,\;\forall v \in \R^N\Big\}.
\end{equation}

In particular, if $f$ is a Clarke eigenfunction, the corresponding eigenvalue $\lambda$ is given by the value of the $1,\infty$-Rayleigh quotient evaluated at $f$.
Note that the definition of the $1,\infty$-eigenpairs in terms of Clarke generalized critical points is stronger than the definition of generalized eigenpairs introduced in \Cref{Def:generalized_p-eigenpair}. Indeed, from Proposition~2.3.14 of~\cite{clarke1990optimization}, we have
\begin{equation}
    \partial^{Cl} \rayl_p(f) \subseteq \frac{\|f\|_p \partial\|\grad f\|_p - \|\grad f\|_p \partial \|f\|_p}{\|f\|_p^2}.
\end{equation}
Thus, whenever $f$ satisfies $0 \in \partial^{Cl} \rayl_{1/\infty}(f)$, it also satisfies the condition in \Cref{Def:generalized_p-eigenpair}. On the other hand, in~\cite{zhang2021homological} it is shown that not every solution of \Cref{Def:generalized_p-eigenpair} is a critical point in the Clarke sense.

Despite this distinction, exploiting the explicit expressions for the subgradients in \eqref{1_subgradient} and \eqref{inf_subgradient} significantly simplifies the investigation of the eigenvalue equation in \Cref{Def:generalized_p-eigenpair}, when compared to the Clarke eigenfunction formulation in \eqref{Clarke_eigenfunction}.

We conclude this section by highlighting that, by means of Clarke subgradients, it can be proved that the min--max principle in Eq.~\eqref{Krasnoselskii_variationa_eigenvalues_Intro} characterizes eigenvalues as generalized critical values of \Cref{Def:generalized_p-eigenpair}, thereby allowing the definition of variational eigenvalues also for $p=1$ and $p=\infty$.
This fact follows from Theorems~6.1 and~6.4, and Theorems~5.1 and~5.8 of~\cite{chang2021nonsmooth}. More precisely, we say that $c$ is a critical value of $\rayl_p$ for $p=1,\infty$ if there exists $f$ such that $0\in \partial^{Cl} \rayl_p(f)$ and $\rayl_p(f)=c$, where $\partial^{Cl}$ denotes the Clarke subgradient. Any value $c$ that is not a critical value is referred to as a regular value of the quotient $\rayl_p$. With these definitions, it can be observed that the deformation theory for sublevel sets of $\rayl_p$ at regular values, including \Cref{Deformation_lemma_Intro} and \Cref{Critical_point_theorem}, still holds in the degenerate cases; see~\cite{chang2021nonsmooth} and the references therein.

\subsubsection{Viscosity $p$-eigenpairs}\label{viscosity_eigenpairs}
A further approach to characterizing eigenpairs in the degenerate cases $p=1$ and $p=\infty$ consists in studying the limiting points of the $p$-eigenpairs as $p$ approaches $1$ or $\infty$, rather than considering a critical point equation for the $1$- and $\infty$-Rayleigh quotients. This approach has been explored in both the infinite-dimensional and graph settings, though primarily for $p=\infty$ (see~\cite{deidda2024_inf_eigenproblem, Lind2, Lind3, Esposito}).

In particular, for any node $u$, it is possible to introduce the local gradient $\locgrad{f}(u)$ of a function $f$ at $u$:
\begin{equation}
  \locgrad{f}(u):=\{\grad f (u,v)\;|\;v\sim u \}\,,
\end{equation}
which represents the set of ``derivatives'' of $f$ along the edges emanating from node $u$.
Consequently, following \cite{Elmoataz2, Elmoataz1, elmoataz2011infinity}, the $\infty$-Laplacian operator can be introduced as
\begin{equation}\label{def_Inf_lap_operator}
    \inflap(f)(u) = \Big(\| \Opc{\grad f}^{-}(u) \|_{\infty} - \| \Opc{\grad f}^{+}(u) \|_{\infty} \Big),
\end{equation}
where the operators $x^+ := \max\{ x, 0 \}$ and $x^- := \max\{ -x, 0 \}$ are applied entrywise to the local gradient vector of $f$ at $u$.

In~\cite{deidda2024_inf_eigenproblem} it has been shown that any limit of $p$-Laplacian eigenpairs satisfies a specific set of equations.
\begin{theorem}\label{thm:limitin_infty_eigenvalue_theorem}
Let $(f_{p_j},\lambda_{p_j})$ be a sequence of $p$-Laplacian eigenpairs, and assume that
\[
\lim_{j\rightarrow \infty}(f_{p_j},\lambda_{p_j}^{\frac{1}{p_j}})=(f,\Lambda).
\]
Then $(f,\Lambda)$ satisfies the following set of equations:
\begin{equation}\label{limitin_inf_eigenvalue_eq}
0=
\begin{cases}
\min\lbrace\|\Opc{\grad f}(u)\|_{\infty}-\Lambda f(u)\;,\; \inflap f(u)\rbrace, & \text{if } u\in\internalnodes \text{ and } f(u)>0,\\
\inflap f(u), & \text{if } u\in\internalnodes \text{ and } f(u)=0,\\
\max\lbrace -\|\Opc{\grad f}(u)\|_{\infty}-\Lambda f(u) \;,\; \inflap f(u)\rbrace, & \text{if } u\in\internalnodes \text{ and } f(u)<0.
\end{cases}
\end{equation}
\end{theorem}

In the infinite-dimensional setting, equations analogous to those in \Cref{thm:limitin_infty_eigenvalue_theorem} must be interpreted in the viscosity sense. This motivates the following definition.
\begin{definition}\label{DEF:viscosity_eigenpair}
A pair $(f,\Lambda)$ is called a \textbf{viscosity} $\infty$-eigenpair if it satisfies the set of equations in \Cref{thm:limitin_infty_eigenvalue_theorem}.
\end{definition}

Moreover, we formally introduce the notion of limit eigenpairs.
\begin{definition}\label{DEF:limit_eigenpairs}
A pair $(f,\Lambda)$ is called a \textbf{limit} $\infty$-eigenpair if there exists a sequence $\{(f_{p_j},\Lambda_{p_j})\}_{j}$ of $p_j$-Laplacian eigenpairs such that
\begin{equation}
    \lim_{j\rightarrow\infty}(f_{p_j},\Lambda_{p_j})=(f,\Lambda).
\end{equation}
\end{definition}

The following inclusion relations hold between the different types of eigenpairs (see~\cite{deidda2024_inf_eigenproblem}).
\begin{theorem}\label{theorem:viscosity_eigenpairs_are_generalized_eigenpairs}
If $(f,\Lambda)$ is a limit eigenpair, then it is also a viscosity eigenpair. If $(f,\Lambda)$ is a viscosity $\infty$-eigenpair, then it is also a generalized $\infty$-eigenpair according to \Cref{Def:generalized_p-eigenpair}.
Moreover, there exist viscosity $\infty$-eigenpairs that are not limit eigenpairs, as well as generalized eigenpairs that are not viscosity eigenpairs.
\end{theorem}

\begin{openproblem}
\Cref{thm:limitin_infty_eigenvalue_theorem} shows that any limit of $p$-Laplacian eigenpairs as $p\to\infty$, say $(f,\Lambda)$, satisfies a limit eigenvalue equation. The solutions of this equation are referred to as viscosity $\infty$-eigenpairs, and it is known that these $\infty$-eigenpairs form a subset of the generalized $\infty$-eigenpairs defined in \Cref{Def:generalized_p-eigenpair}.
The existence or nonexistence of viscosity $1$-eigenpairs defined via a $1$-limit eigenvalue equation sharper than the generalized eigenvalue equation given in terms of subgradients (see \Cref{Def:generalized_p-eigenpair}) is currently an open problem.
\end{openproblem}

\subsection{Homological eigenvalues}
In this section, we introduce a different class of eigenvalues by considering variations in the homological groups of the sublevel sets of the $p$-Rayleigh quotient.

These eigenvalues are defined according to the notion of homological critical values introduced in~\cite{cohen2005stability} and formalized in~\cite{govc2016definition}. The case of the $p$-Laplacian operator has been studied in depth in~\cite{zhang2021homological}. We begin with the definition of homological regular and critical values of $\rayl_p$.
\begin{definition}
    A real number $\tau$ is a homological regular value of $\rayl_p$ if there exists $\epsilon>0$ such that for any $t_1<t_2$ with $t_1,t_2\in (\tau-\epsilon,\tau+\epsilon)$, the inclusion
    $$\big\{f\,\big|\;\text{s.t. } \rayl_p(f)\leq t_1\big\} \quad  \xhookrightarrow{\qquad} \quad \big\{f\,\big|\;\text{s.t. } \rayl_p(f)\leq t_2\big\}$$
    induces isomorphisms on all homology groups. A real number $\tau$ that is not a homological regular value of $\rayl_p$ is said to be a homological critical value.
\end{definition}

In particular, we have the following definition of $p$-Laplacian homological eigenvalues.
\begin{definition}
    We say that $\Lambda$ is a $p$-Laplacian homological eigenvalue if $\Lambda$ is a homological critical value of $\rayl_p$.
\end{definition}

According to the following result and~\cite{govc2016definition}, any $p$-Laplacian homological critical value is also a critical value in the Clarke sense.
\begin{proposition}\label{prop:homology-Clarke}
    Any homological critical value of $\rayl_p$, for $p\in[1,\infty]$, is a generalized eigenvalue in the Clarke sense.
\end{proposition}
\begin{proof}
Suppose, by contradiction, that $c$ is a homological critical value but not a critical value in the Clarke sense. Since Clarke critical values form a closed subset of $\R$, there exists $\epsilon_0>0$ such that there is no Clarke critical value in the interval $(c-\epsilon_0,c+\epsilon_0)$. Since $\rayl_p$ is zero-homogeneous, we may restrict $\rayl_p$ to the unit $\ell^2$-sphere.
For convenience, we use $\{\mathcal{R}_p\le t\}$ to denote the sublevel sets $\{f\in \mathcal{H}(\internalnodes):\|f\|_2=1 \text{ and } \mathcal{R}_p(f)\le t\}$.
According to Theorem~3.1 and Remarks~3.3 and~3.4 in~\cite{chang1981variational}, there exists an odd homeomorphism $\eta:S_2\to S_2$ induced by a decreasing flow such that $\eta\{\mathcal{R}_p\le c+\epsilon\}\subset \{\mathcal{R}_p\le c-\epsilon\}$, where $S_2$ denotes the unit $\ell^2$-sphere.
Now, for any $t_1,t_2$ satisfying $c-\epsilon<t_1<t_2<c+\epsilon$, consider the inclusion
\begin{equation}
\eta\{\mathcal{R}_p\le t_1\}\mathop{\hookrightarrow}\limits^{i_1} \eta\{\mathcal{R}_p\le t_2\}\mathop{\hookrightarrow}\limits^{i_2} \{\mathcal{R}_p\le t_1 \} \mathop{\hookrightarrow}\limits^{i_3}\{\mathcal{R}_p\le t_2 \}.
\end{equation}
Since $\eta$ is a homeomorphism, both $i_2\circ i_1$ and $i_3\circ i_2$ induce isomorphisms on all homology groups. By Lemma~3.1 in~\cite{govc2016definition}, the inclusions $i_1$, $i_2$, and $i_3$ also induce isomorphisms on all homology groups. In particular, the inclusion $\{\mathcal{R}_p\le t_1 \} \hookrightarrow \{\mathcal{R}_p\le t_2 \}$ induces isomorphisms on all homology groups. Therefore, $c$ is a homological regular value of $\rayl_p$, contradicting the assumption.
\end{proof}

We have thus shown that any homological eigenvalue of the $p$-Laplacian, for any $p\in [1,\infty]$, is also an eigenvalue in the Clarke sense. In turn, by the discussion in \Cref{Clarke_eigenpairs}, this implies that any homological eigenvalue is a generalized eigenvalue as in \Cref{Def:generalized_p-eigenpair}, but not necessarily a variational one.

The class of homological eigenvalues introduced in this section is of particular interest when studying the regularity of the $p$-Laplacian spectrum; we discuss this aspect in greater detail in \Cref{Sec:regularity_of_p_lap_spectrum}. For $p\in(1,\infty)$, it remains an open problem whether eigenvalues that are not homological exist.

\begin{openproblem}\label{openproblem_homological_eigenvalues}
    Do the homological eigenvalues exhaust the $p$-Laplacian spectrum for $p\in (1,\infty)$?
\end{openproblem}

In this regard, we note that in~\cite{zhang2021homological} it has been proved that for $p=1$ there may exist eigenvalues that are not homological, explaining why \Cref{openproblem_homological_eigenvalues} concerns only the case $p\in (1,\infty)$.

We include illustrative examples related to the themes discussed in this section in \Cref{Sec:1-Lapl_spectrum}; see \cref{exam:complete} and \cref{exam:path}. Using different results presented in the subsequent sections, we illustrate situations in which: (i) the variational multiplicity and $\gamma$-multiplicity differ; (ii) the homological spectrum does not exhaust the spectrum in the $1$-Laplacian case; (iii) the number of orthogonal eigenvectors associated with a given eigenvalue can exceed the $\gamma$-multiplicity; (iv) eigenvalues exceed the dimension of the space; and (v) non-variational eigenvalues exist.


\section{Regularity of the $p$-Laplacian spectrum}\label{Sec:regularity_of_p_lap_spectrum}
In this section we summarize a series of results on the regularity of the $p$-Laplacian spectrum as $p$ varies.
We start by recalling that, since the limit of $p$-eigenpairs as $p\to p^*$ yields a $p^*$-eigenpair, the $p$-Laplacian spectrum is always upper semicontinuous.
To this end, for any $p\in[1,\infty]$, we denote the $p$-Laplacian spectrum by $\mathrm{Spec}(\plap):=\{(f,\Lambda) :\text{ $p$-Laplacian eigenpair}\}$.

\begin{lemma}[Lemma 2.1 and Proposition 2.3  \cite{zhang2021homological}]\label{lemma_upp_semicontinuity_plap_spectrum}
    Let $\mathrm{Spec}(\plap)$ denote the $p$-Laplacian spectrum. Then the set valued map that associates $p$
 to $\mathrm{Spec}(\plap)$ is upper semicontinuous, i.e. for any $p\geq 1$ and any $\epsilon>0$ there exists $\delta>0$ such that for any $p'\in(p-\delta,p+\delta)$ with $p'\geq 1$
 $$\mathrm{Spec}(\Delta_{p'})\subset \bigcup_{(f,\Lambda)\in \mathrm{Spec}(\plap)}B_{\epsilon}(f)\times\big(\Lambda-\epsilon, \Lambda+\epsilon \big).$$
\end{lemma}

Moreover, in the following lemma, which is a new contribution, we discuss the regularity of the variational eigenvalues.
It is worth noting that our result essentially generalizes a monotonicity inequality by Dru\c{t}u and Mackay (see Lemma 2.7 in \cite{dructu2019random}).

\begin{lemma}\label{lemma:p-monotonic}
Let $\Gc$ be a graph with a measure $\nodeweight$ on the nodes. Then, for any $k$, the variational eigenvalue $\lambda_k(\Delta_p)$ is locally Lipschitz continuous with respect to $p$, and $\lim_{p\to+\infty}\lambda_k(\Delta_p)^{\frac1p}=\Lambda_k(\Delta_{\infty})$. In particular:
\begin{itemize}
    \item the function 
$p\mapsto \big(|E|/2\big)^{-\frac1p}\lambda_k(\Delta_p)^{\frac1p}$ 
is increasing on $
[1,+\infty)$,
\item 
the function 
$p\mapsto|\nodeweight|^{\frac{1}{p}}\lambda_k(\Delta_p)^{\frac1p}$ is decreasing on $[1,+\infty)$,
\end{itemize}
where $|\nodeweight|=\sum_{u\in\internalnodes} |\nodeweight_u|$ is the total variation of the measure $\nodeweight$ defined on the nodes.
\end{lemma}
\begin{proof}
For $1\le p\le q$, from H\"older's inequality there hold:
\begin{equation}
\|f\|_q\leq \|f\|_p\leq |\nodeweight|^{\frac{1}{p}-\frac{1}{q}}\|f\|_q \qquad \text{and} \qquad 
\|Kf\|_q\leq \|Kf\|_p\leq \Big(\frac{|\edgeset|}{2}\Big)^{\frac{1}{p}-\frac{1}{q}}\|Kf\|_q. 
\end{equation}
Thus,
\begin{equation}
|\nodeweight|^{\frac1q-\frac1p} \rayl_q(f)\le \rayl_p(f)\le \big(\frac{|E|}{2}\big)^{\frac1p-\frac1q}\rayl_q(f).
\end{equation}
In particular, the two functions
\begin{equation}
p\longmapsto |\nodeweight|^{\frac1p}\; \rayl_p(f) \qquad \text{and} \qquad p\longmapsto\Big(\frac{|E|}{2}\Big)^{-\frac1p} \rayl_p(f)
\end{equation}
are monotonically decreasing and increasing on $[0,+\infty)$, respectively.
Hence, by definition of variational eigenvalues, $|\nodeweight|^{\frac1p}\lambda_k^{\frac{1}{p}}(\Delta_p)$ is monotonically decreasing with respect to $p$, while $\big(|E|/2\big)^{-\frac1p}\lambda_k^{\frac{1}{p}}(\Delta_p)$ is monotonically increasing with respect to $p$.
In addition, taking $q=\infty$, we have
\begin{equation}
\frac{|E|^{-\frac1p}}{2}\lambda_k(\Delta_p)^{\frac1p}\le \Lambda_k(\Delta_{\infty})\le |\nodeweight|^{\frac1p}\lambda_k(\Delta_p)^{\frac1p}
\end{equation}
which proves 
$\lim_{p\to+\infty}\lambda_k(\Delta_p)^{\frac1p}=\Lambda_k(\Delta_{\infty})$.
\end{proof}

In particular, we have the following consequence of the previous lemma.
\begin{corollary}
The $\infty$-variational eigenvalues are limit (and thus also viscosity) eigenvalues.
\end{corollary}

Another monotonicity inequality regarding the regularity of variational eigenvalues is established in \cite{zhang2021homological}, which claims the following.
\begin{lemma}
Let $\Gc$ be a graph with a measure $\nodeweight$ on the nodes. Denote by $t_G:=\max_{v\in V}
\deg (v)/\nodeweight_v$, where $\deg(v)$ is the number of vertices that are adjacent to $v$.
Then, for any $1\le p < q<\infty$,
we have
$$\big(2\max_{\{u,v\}\in E}w_{uv}\big)^{-p}\lambda_k(\Delta_{p})\ge \big(2\max_{\{u,v\}\in E}w_{uv}\big)^{-q}\lambda_k(\Delta_{q})$$
and
$$p\left(\frac{2\lambda_k(\Delta_{p})}{t_G}\right)^{\frac{1}{p}}\leq q\left(\frac{2\lambda_k(\Delta_{q})}{t_G}\right)^{\frac{1}{q}}$$
that is,
the function
$$p\mapsto\big(2\max_{\{u,v\}\in E}w_{uv}\big)^{-p}\lambda_k(\Delta_p)$$
is decreasing on $[1,+\infty)$, while the function $$p\mapsto p\Big(\frac{2}{t_G}\lambda_k(\Delta_p)\Big)^{\frac1p}$$
is increasing on $[1,+\infty)$.
\end{lemma}

\begin{example}\label{Ex:p-Lap_eingenspace_varying_p_on_complete_graphs}
We illustrate how the spectrum of the $p$-Laplacian varies with $p$ in the simple case of an unweighted complete graph on $3$ nodes, $K_3$. In this case, the eigenvalues of the $p$-Laplacian are:
\begin{equation}
 \lambda_0(\plap)=0, \quad \lambda_a(\Delta_p)=1+2^{p-1}, \quad \lambda_b(\Delta_p)=(1+2^{q-1})^{p-1}  \qquad \forall\; 1<p<+\infty,
\end{equation}
where $q$ is the H\"older conjugate of $p$. Then these eigenvalues degenerate to the two eigenvalues $\Lambda_0(\plap)=0$, $\Lambda_{a}(\plap)=\Lambda_{b}(\plap)=2$ when $p\in\{1,+\infty\}$.
Moreover, when $p\in (1,2)\cup (2,+\infty)$, the two positive eigenvalues $\lambda_a(\plap)$ and $\lambda_b(\plap)$ are distinct and both of them have multiplicity $1$. However, when $p=1$ or $p=2$ or $p=+\infty$, the multiplicity of the only positive eigenvalue of $\Delta_p$ is two.
Next, it is possible to observe that for any $p\ne 1,2,\infty$, $X_{0}(\plap)$, $X_a(\plap)$, and $X_b(\plap)$ defined below are the eigenspaces corresponding to $\lambda_0(\plap)$, $\lambda_a(\Delta_p)$, and $\lambda_b(\Delta_p)$, respectively.
\begin{equation}\label{eq.eigenspaces_triangle}
    \begin{aligned}
        X_0(\Delta_p)&=\mathrm{span}(\vec 1)\setminus\vec0,\\
        X_a(\Delta_p)&=\bigcup_{i\ne j}\mathrm{span}(f_{i,j})\setminus\vec0,\quad  \text{where}\quad  f_{i,j}(v)=\begin{cases}
        1 \quad &\text{ if } v=i,\\
        -1 \quad &\text{ if } v=j,\\
        0 \quad &\text{ if } v\not\in \{i,j\}.
        \end{cases}\\
         X_b(\Delta_p)&=\bigcup_{i=1}^3\mathrm{span}(g_{i})\setminus\vec 0,\quad \text{where}\quad g_i(v)=\begin{cases}
             1 \quad &\text{ if } v=i,\\
             -2^{-\frac{1}{p-1}} \quad &\text{ if } v\ne i.
         \end{cases}
    \end{aligned}
\end{equation}
\begin{figure}[ht]
    \begin{center}
\begin{minipage}{.48\textwidth} 
\begin{tikzpicture}[scale=0.42]
{\draw[color=black,->] (0,0,0) -- (4.5,0,0) node[anchor=north east]{$x_2$};}
{\draw[color=black,->] (0,0,0) -- (0,4.5,0) node[anchor=north west]{$x_3$};}
{\draw[color=black,->] (0,0,0) -- (0,0, 4.5) node[anchor=south]{$x_1$};}
\filldraw[red,fill opacity=0.10,line width=0,draw =red!35] (0,0,0)--(0,3+1,0)--(-3-1,3+1,0)--(-3-1,0,0);
\filldraw[red,fill opacity=0.20,line width=0,draw =red!35] (0,0,0)--(-3-1,0,0)--(-3-1,0,4)--(0,0,3+1);
\filldraw[red,fill opacity=0.30,line width=0,draw =red!35] (0,0,0)--(0,0,3+1)--(0,-4,3+1)--(0,-3-1,0);
\filldraw[red,fill opacity=0.10,line width=0,draw =red!35] (0,0,0)--(0,-3-1,0)--(4,-3-1,0)--(3+1,0,0);
\filldraw[red,fill opacity=0.20,line width=0,draw =red!35] (0,0,0)--(3+1,0,0)--(3+1,0,-4)--(0,0,-3-1);
\filldraw[red,fill opacity=0.30,line width=0,draw =red!35] (0,0,0)--(0,0,-3-1)--(0,4,-3-1)--(0,3+1,0);
\draw[color=blue,thick,->] (0,0,0)  -- (3,3,3)node[anchor=south west]{\footnotesize $(1,1,1)$};
\draw[color=red,very thin,dashed]  (4,-4,0) --(-4,4,0); 
\draw[color=red,very thin,dashed]  (4,0,-4) --(-4,0,4); 
\draw[color=red,very thin,dashed]  (0,4,-4) --(0,-4,4);
\node (a) at (4,4.6) {\scriptsize \color{white!50!black}$p=1$};
\end{tikzpicture}~~
\end{minipage}
\begin{minipage}{.48\textwidth}    
\begin{tikzpicture}[scale=0.51]
{\draw[color=black,->] (0,0,0) -- (4.5,0,0) node[anchor=north east]{$x_2$};}
{\draw[color=black,->] (0,0,0) -- (0,4.5,0) node[anchor=north west]{$x_3$};}
{\draw[color=black,->] (0,0,0) -- (0,0, 4.5) node[anchor=south]{$x_1$};}
\filldraw[red,fill opacity=0.10,line width=0,draw =red!15]  (3,-3,0) -- (3,0,-3)-- (0,3,-3)--(-3,3,0)--(-3,0,3)--(0,-3,3)--(3,-3,0);
\draw[color=red,very thin,dashed]  (3,-3,0) --(-3,3,0); 
\draw[color=red,very thin,dashed]  (3,0,-3) --(-3,0,3); 
\draw[color=red,very thin,dashed]  (0,3,-3) --(0,-3,3); 
\draw[color=blue,thick,->] (0,0,0)  -- (3,3,3)node[anchor=south west]{\footnotesize $(1,1,1)$};
\node (a) at (4,4.6) {\scriptsize \color{white!50!black}$p=2$};
\end{tikzpicture}
\end{minipage}
\begin{minipage}{.48\textwidth}            
\begin{tikzpicture}[scale=0.51]
{\draw[color=black,->] (0,0,0) -- (4.5,0,0) node[anchor=north east]{$x_2$};}
{\draw[color=black,->] (0,0,0) -- (0,4.5,0) node[anchor=north west]{$x_3$};}
{\draw[color=black,->] (0,0,0) -- (0,0, 4.5) node[anchor=south]{$x_1$};}
\draw[color=red,very thin,dashed]  (3,-3,0) --(-3,3,0); 
\draw[color=red,very thin,dashed]  (3,0,-3) --(-3,0,3); 
\draw[color=red,very thin,dashed]  (0,3,-3) --(0,-3,3);
\draw[color=purple]  (2,-2.82,-2.82) --(-2,2.82,2.82); 
\draw[color=purple]  (-2.82,2,-2.82) --(2.82,-2,2.82); 
\draw[color=purple]  (-2.82,-2.82,2) --(2.82,2.82,-2) ; 
\draw[color=blue,thick,->] (0,0,0)  -- (3,3,3)node[anchor=south west]{\footnotesize $(1,1,1)$};
\node (a) at (4,4.6) {\scriptsize \color{white!50!black}$p=3$};
\end{tikzpicture}
\end{minipage}
\begin{minipage}{.48\textwidth}
\begin{tikzpicture}[scale=0.42]
{\draw[color=black,->] (0,0,0) -- (4.5,0,0) node[anchor=north east]{$x_2$};}
{\draw[color=black,->] (0,0,0) -- (0,4.5,0) node[anchor=north west]{$x_3$};}
{\draw[color=black,->] (0,0,0) -- (0,0, 4.5) node[anchor=south]{$x_1$};}
\filldraw[red,fill opacity=0.20,line width=0,draw =red!35] (0,0,0)--(-3,3,3)--(-3,0,3)--(-3,-3,3);
\filldraw[red,fill opacity=0.30,line width=0,draw =red!35] (0,0,0)--(-3,3,-3)--(-3,3,0)--(-3,3,3);
\filldraw[red,fill opacity=0.30,line width=0,draw =red!35] (0,0,0)--(3,-3,3)--(3,-3,0)--(3,-3,-3);
\filldraw[red,fill opacity=0.10,line width=0,draw =red!35] (0,0,0)--(-3,-3,3)--(0,-3,3)--(3,-3,3);
\filldraw[red,fill opacity=0.10,line width=0,draw =red!35] (0,0,0)--(3,3,-3)--(0,3,-3)--(-3,3,-3);
\filldraw[red,fill opacity=0.20,line width=0,draw =red!35] (0,0,0)--(3,-3,-3)--(3,0,-3)--(3,3,-3);
\draw[color=blue,thick,->] (0,0,0)  -- (3,3,3)node[anchor=south west]{\footnotesize $(1,1,1)$};
\draw[color=red,very thin,dashed]  (4,-4,0) --(-4,4,0); 
\draw[color=red,very thin,dashed]  (4,0,-4) --(-4,0,4); 
\draw[color=red,very thin,dashed]  (0,4,-4) --(0,-4,4); 
\node (a) at (4,4.6) {\scriptsize \color{white!50!black}$p=\infty$};
\end{tikzpicture}
\end{minipage}
    \end{center}
    \caption{From left to right and from top to bottom, we report the eigenspaces relative to the $p$-Laplacian eigenvalues for $p$ taking values $1,2,3$ and $\infty$ in Example \ref{Ex:p-Lap_eingenspace_varying_p_on_complete_graphs}. Dashed lines denote the vectors $f_{i,j}$ varying $i$ and $j$ in $\{1,2,3\}$. Similarly, for $p=3$, the solid red lines denote the vectors $g_i$ and $-g_i$ for $i=1,2,3$. Finally the red regions denote the eigenspaces relative to $\Lambda_{a,b}(\plap)$; different opacities correspond to different supporting planes.}
    \label{fig:p-eigenspaces_triangle}
\end{figure}
\end{example}
Differently, for $p=1,2,\infty$, the eigenspace corresponding to $\Lambda_{a}(\plap)$ or $\Lambda_{b}(\plap)$ is given by
 \begin{equation}
 X(\Delta_p):=\bigcup\limits_{i\ne j}\left(\mathrm{cone}(f_{i,j},g_i)\cup \mathrm{cone}(f_{i,j},-g_j)\right)\setminus\mathbf{0},
 \end{equation}
where $f_{i,j}$ and $g_i$ are defined as above in \eqref{eq.eigenspaces_triangle} and $\text{cone}(f,g):=\{t_1f+t_2g\,|\,t_1\ge 0,t_2\ge 0\}$. Observe that, just as in the linear case with respect to matrix perturbations, the multiplicity may change when varying $p$. Thus, if the $\gamma$-multiplicity increases for a certain value $p^*$, there are eigenfunctions that are not limits of $p$-eigenfunctions for $p$ close to, but different from, $p^*$.
A graphical representation of the $p$-Laplacian eigenspaces on $K_3$ is reported in \Cref{fig:p-eigenspaces_triangle}.
We observe that $X(\Delta_1)$ and $X(\Delta_\infty)$ can be written as the union of six angular planar regions that intersect orthogonally along $X_a(\plap)$.

We conclude this section by discussing the regularity of the homological eigenvalues. Indeed, this class of eigenvalues deserves particular attention in the study of the continuity of the spectrum. Intuitively, since the $p$-Rayleigh quotient varies continuously with $p$, the sublevel sets are locally isomorphic in $p$. In particular, a change in the homology of the sublevel sets of the $p$-Rayleigh quotient yields a change in the homology of the sublevel sets of the $q$-Rayleigh quotient for $q$ close to $p$, heuristically supporting the continuity of the homological eigenvalues. Formally, we report the following theorem from \cite{zhang2021homological}.
\begin{theorem}\label{thm:homological-eigen}
 Let $p\geq 1$ and assume $\Lambda(\plap)$ is an isolated homological eigenvalue of $\plap$, then for any $\epsilon$, there exists $\delta$ such that for any $q\in (p-\delta,p+\delta)$, $\Delta_q$ admits a homological eigenvalue $\Lambda(\Delta_q)\in \big(\Lambda(\plap)-\epsilon,\Lambda(\plap)+\epsilon\big)$.
\end{theorem}

An interesting problem is the differentiability of the $p$-Laplacian spectrum with respect to perturbations of the graph weights. In this regard, in the recent work \cite{berkolaiko2025eigenvalues}, the authors introduce the notion of regular eigenpairs and show that they are differentiable as functions of the edge weights. They call an eigenpair $(f,\lambda)$ regular if
\begin{equation}
\text{dim Ker}\Big(\text{Hess}_f L(\lambda,f)\Big)=1
\end{equation}
where $L(\lambda,f)=\|\grad f\|_p^p-\lambda\big( \|f\|_{\nodeweight,p}^p -1 \big)$ is assumed to be twice differentiable in $f$ and $\text{Hess}_f$ denotes the Hessian matrix obtained by differentiating with respect to the function $f$.
The manuscript \cite{berkolaiko2025eigenvalues} also provides expressions for the derivatives of the eigenpair. These expressions generalize the Hellmann--Feynman formulas holding for the spectrum of linear self-adjoint operators.

%
\section{Duality}\label{Sec:Duality}
In this section we discuss a duality result between the $p$-Laplacian eigenproblem on the nodes considered so far and a $q$-``Laplacian'' eigenproblem defined on the edges, where $q$ is the conjugate of $p$. This result extends a well-known correspondence in the linear case. Indeed, if we write the Laplacian matrix as $\Delta=\grad^T\grad$, it is immediate to note that the nonzero eigenvalues of $\Delta$ coincide with those of $\Delta^E=\grad \grad^T$. In particular, there is a one-to-one correspondence between the nonzero eigenvalues of an operator acting on node functions and those of an operator acting on edge functions.
Some recent results extend this equivalence to the nonlinear case \cite{zhang2021discrete, tudisco2022nonlinear, laubmann2025duality}. We refer to \cite{bungert2021eigenvalue, laubmann2025duality} for a corresponding investigation in the infinite-dimensional setting.

By the $p$-Laplacian eigenvalue problem on the nodes we mean the generalized one introduced in \Cref{Def:generalized_p-eigenpair}. By the $q$-``Laplacian'' eigenvalue problem on the edges, instead, we refer to the eigenvalue problem given by the critical point equation of the $q$-Rayleigh quotient $\rayl_q^{\edgeset}$ defined by
\begin{equation}
    \rayl_q^{\edgeset}(G):=\frac{\|\hat{\grad}^T G\|_{\nodeweight,q}}{\|G\|_q}\,,
\end{equation}
where $G$ varies in $\mathcal{H}(\edgeset)=\{G:\edgeset\rightarrow\R\}$ and $\hat{\grad}^T$ is a scaled transpose incidence matrix of the graph, i.e. $\hat{\grad}^T=\mathrm{diag}(\nodeweight^{-1})\grad^T$.
A critical pair $(G,\eta)$ of $\rayl_q^{\edgeset}$ is again considered in the generalized sense of \Cref{Def:generalized_p-eigenpair}.
In particular, we note that $(G,\eta)$ is an edge $q$-eigenpair for $q\in(1,\infty)$ if it satisfies the nonlinear eigenvalue equation:
\begin{equation}\label{eig_eq_on_edges}
  \left({\hat{\grad}}
    \left|{\hat{\grad}}^TG\right|^{q-2}
    {\hat{\grad}}^TG\right)(u,v)
  =\eta |G(u,v)|^{q-2}G(u,v)\qquad \forall (u,v)\in\edgeset
\end{equation}
More in detail, in~\cite{zhang2021discrete, tudisco2022nonlinear} the authors show that the nonzero critical values and points of $\rayl_q^{\edgeset}$ correspond to the nonzero critical values and points of $\rayl_p$, where $p$ is the conjugate of $q$, with corresponding statements also for variational eigenvalues and multiplicities.

The proof is based on the following lemma from \cite{zhang2021discrete}. Since we adapted the proof to the weighted case, for completeness we report it in \cref{appendix_duality}.
\begin{lemma}\label{Lemma_duality}
Assume we have $\mathcal{H}(\internalnodes)=(\R^N,\langle\cdot,\cdot\rangle_{\nodeweight})$, $\mathcal{H}(\edgeset)=(\R^M, \langle\cdot,\cdot\rangle_{\edgeweight})$, a linear function, say $Ax$ with $A\in \R^{M\times N}$, from the first to the second Hilbert space, and two convex $1$-homogeneous and positive functionals $\Phi:(\R^M, \langle\cdot,\cdot\rangle_{\edgeweight})\rightarrow \R^+$ and $\Psi:(\R^N,\langle\cdot,\cdot\rangle_{\nodeweight})\rightarrow \R^+$ with $\Phi(G)=0$ iff $G=0$ and $\Psi(f)=0$ iff $f=0$.
Then $(f,\Lambda)$ is a critical pair of
$$\rayl(f):=\frac{\Phi(A f)}{\Psi(f)}$$
with $\Lambda\neq 0$ and $\Xi\in \partial_{G=Af}\Phi(G)$, $\xi \in \partial_f \Psi(f)$ that satisfy $A^T(\edgeweight \odot\Xi)=\Lambda \nodeweight\odot f$,
if and only if $(\Xi,\Lambda)$ is a critical pair of
$$\rayl^*(\Theta)=\frac{\Psi^*(\nodeweight^{-1}\odot(A^T(\edgeweight\odot \Theta))}{\Phi^*(\Theta)},$$
where $\Phi^*(\Theta)=\sup_{G} \langle \Theta, G\rangle_{\edgeweight}/\Phi(G)$ and $\Psi^*(\eta)=\sup_{f} \langle \eta, f\rangle_\nu/\Psi(f)$.
\end{lemma}
We highlight that \Cref{Lemma_duality} shows a correspondence between eigenpairs of a primal problem and those of a dual one, where the dual problem is expressed in terms of norm-dual-type functions. However, analogous results also hold for different dual problems obtained, e.g., by means of Legendre--Fenchel dual functions \cite{tudisco2022nonlinear, laubmann2025duality}.

In the remainder of this section we use $\Delta_p^\nodeset$ and $\Delta_q^\edgeset$ to denote the usual $p$-Laplacian operator acting on the space of node functions and the $q$-``Laplacian'' operator, obtained by differentiating (or subdifferentiating) $\Psi(G)=\|\hat{\grad}^\top G\|_q$, acting on the space of edge functions, respectively.

\begin{theorem}\label{Thm:duality}
Suppose $(f,\Lambda)$ is a $p$-eigenpair on the nodes with $p\in[1,\infty]$, $\Lambda \ne 0$, and let $\Xi\in \partial\big|_{G=\grad f} \|G\|_p$ be such that
$$\hat{\grad}^\top \Xi \in  \Lambda \partial \|f\|_{\nodeweight,p}.$$
Then $(\Xi,\Lambda)$ is a $q$-eigenpair on the edges.
Analogously, if $(G,\Lambda)$ is a $q$-eigenpair on the edges with $q\in[1,\infty]$, $\Lambda\neq 0$, and $\xi\in \partial\big|_{g=\hat{\grad}^\top G} \|g\|_{\nodeweight,q}$ is such that
$$\grad \xi \in  \Lambda \partial \|G\|_q,$$
then $(\xi,\Lambda)$ is a $p$-eigenpair on the nodes.
Moreover the eigenvalues satisfy:
\begin{itemize}
    \item $\gamma\text{-}\mathrm{mult}_{\Delta_p^\nodeset}(\Lambda)=\gamma\text{-}\mathrm{mult}_{\Delta_q^\edgeset}(\Lambda)$,
    \item $\Lambda_{d_1+k}(\plap^\nodeset)=\Lambda_{d_2+k}(\Delta_q^\edgeset)$ for any $k=0,1,\cdots,N$,
    \item $\mathrm{mult}_{\Delta_p^{\nodeset}}(\Lambda)=\mathrm{mult}_{\Delta_q^{\edgeset}}(\Lambda)$,
\end{itemize}
where $d_1=\text{dimension}\big(\mathrm{Ker}(\grad)\big)$, $d_2=\text{dimension}\big(\mathrm{Ker}(\grad^\top)\big)$, and the subscripts $\Delta_p^{\nodeset}$ and $\Delta_q^{\edgeset}$ in $\gamma\text{-}\mathrm{mult}_{\Delta_p^{\nodeset}}(\cdot)$, $\gamma\text{-}\mathrm{mult}_{\Delta_q^{\edgeset}}(\cdot)$, $\mathrm{mult}_{\Delta_p^{\nodeset}}(\cdot)$, and $\mathrm{mult}_{\Delta_q^{\edgeset}}(\cdot)$ distinguish between the multiplicity and $\gamma$-multiplicity of eigenvalues of the operators $\Delta_p^\nodeset$ and $\Delta_q^\edgeset$, respectively.
\end{theorem}

Note that the above theorem shows that the normalized eigenvalues of the $p$- and $q$-Laplacians on the nodes and edges coincide. Moreover, given an eigenvalue, the corresponding eigenfunctions are given by the subgradients that satisfy the generalized eigenproblem in \Cref{Def:generalized_p-eigenpair}.
Moreover, in the particular case $p,q\in(1,\infty)$, \Cref{Thm:duality} states that if $(f,\lambda)$ is an eigenpair of $\plap$ with $\lambda\neq 0$, then $(|\grad f|^{p-2}\grad f,\lambda^{\frac{q}{p}})$ is a $q$-eigenpair on the edges. Vice versa, if $(G,\eta)$ is a $q$-eigenpair on the edges with $\eta\neq 0$, then $\big(\eta^{\frac{p}{q}},|{\hat{\grad}}^T G|^{q-2}{\hat{\grad}}^TG\big)$ is a $\plap$-eigenpair.

\begin{remark}\label{Remark_edge-q-Laplacian_operator}
We spend a few words on the edge $q$-``Laplacian'' operator. We first point out that the linear edge Laplacian has much potential to handle the leaderless consensus problem; see \cite{Zelazo2007} for details. Thus, an appropriate nonlinear version of the edge Laplacian, such as the edge $q$-``Laplacian'' operator, is natural and could be useful in particular tasks \cite{Mulas2022_plap_hyper, Fazeny, Fazeny24}.
On the other hand, we note that, unlike the node $p$-Laplacian operator, which has a ``natural'' interpretation as a discretized differential operator, the same does not hold for the edge version. Nevertheless, there are different reasons why we use this name for the operator in \eqref{eig_eq_on_edges}. The first reason is the clear algebraic analogy with the node operator in \Cref{def:graph_p-Lap_operator}. The second reason is that it matches the definition of a weighted $q$-Laplacian on a dual (oriented) hypergraph $\mathcal{HG}$ of $\Gc$ \cite{MulasZhang,Reff}. To define $\mathcal{HG}$, we let any edge $e=(u,v)$ of $\Gc$ become a node $u_e$ of $\mathcal{HG}$, and we define an oriented hyperedge $h_u$ for any node $u$ of $\Gc$. We recall that an oriented hyperedge $h$ is given by a set $h^+$ of input vertices (from which, heuristically, the flow enters) and a set $h^-$ of output vertices (from which the flow goes out) \cite{shi1992signed}. In the particular case of $h_u$, we define $h_u^-$ as the set of all nodes $u_{(u,v)}$ varying $v\in \Gc$ with $v$ connected to $u$, and $h_u^+$ as the set of all nodes $u_{(v,u)}$ varying $v\in \Gc$ with $v$ connected to $u$. Then, for a function $G\in \mathcal{H}(\edgeset)$ defined on the edges we have
\begin{equation}
    \begin{aligned}
    \grad^\top G(u)&=\sum_{(v,u)\in \edgeset} \edgeweight_{uv} G(v,u)- \sum_{(u,v)\in \edgeset} \edgeweight_{uv} G(u,v)\\
    &=\sum_{v_e\in h_u^+} \edgeweight_e G(v_e)- \sum_{v_e\in h_u^-} \edgeweight_e G(v_e)=\grad^\top G(h_u)
    \end{aligned}
\end{equation}
where we have identified the function $G\in \mathcal{H}(\edgeset)$ with the corresponding function defined on $\mathcal{HG}$ and given by $G(v_e)=G(e)$, and analogously for the functions $\grad^\top G$ defined on the nodes.
In particular, it follows the following expression for the $q$-Laplacian defined on the edges, where we set $e=(u,v)$:
\begin{equation}
\begin{aligned}
    \grad(|\grad^\top G|^{q-2})(u,v)=\edgeweight_{uv} |\grad^\top G(v)-\grad^\top G(u)|^{q-2}(\grad^\top G(v)-\grad^\top G(u))=&\\
    =\edgeweight_{uv} |\grad^\top G(h_v)-\grad^\top G(h_u)|^{q-2}(\grad^\top G(h_v)-\grad^\top G(h_u))=\grad(|\grad^\top G|^{q-2}&)(u_e).
\end{aligned}
\end{equation}
Note that
\begin{equation}
\begin{aligned}
    (\grad^\top G(h_v)-\grad^\top G(h_u))=& \sum_{v_e\in h_v^+} \edgeweight_e G(v_e)+ \sum_{v_e\in h_u^-} \edgeweight_e G(v_e)-\\
    &- \sum_{v_e\in h_v^-} \edgeweight_e G(v_e)-\sum_{v_e\in h_u^+} \edgeweight_e G(v_e)
\end{aligned}
\end{equation}
and that $h_v^+=\{(w,v)|\;w\sim v\}$, $h_u^-=\{(u,w)|\;w\sim u\}$ is the set of nodes of the hypergraph connected to $(u,v)$ by some hyperedge ($h_v$ or $h_u$) and co-oriented with $(u,v)$ in $h_v$ and $h_u$ respectively (two nodes in a hyperedge $h_w$ are co-oriented if both belong to $h_w^+$ or $h_w^-$). Analogously, $h_v^-=\{(v,w)|\;w\sim v\}$, $h_u^+=\{(w,u)|\;w\sim u\}$ is the set of nodes of $\mathcal{HG}$ connected to $(u,v)$ by some hyperedge ($h_v$ or $h_u$) and anti-oriented with $(u,v)$ in $h_v$ and $h_u$ respectively. Finally, noting that $h_u$ and $h_v$ are the only two hyperedges incident to the node $u_e=(u,v)$, this proves that up to constants the $q$-Laplacian operator defined on the edges matches the definition given in \cite{tudisco2022nonlinear}. 
This perspective stems from the 
dual hypergraph \cite{Berge,Robeva}.
\end{remark}

\begin{figure}[t]
\begin{center}
  \begin{tikzpicture}[scale=1.2]
 \draw[thick,black]  (30:1)-- (90:1)--(150:1)--(210:1)--(270:1)--(-30:1)--(30:1);
  \node (0) at (30:1) {\color{red} $\bullet$};
  \node (0) at (-30:1) {\color{red}$\bullet$};
  \node (0) at (270:1) {\color{red}$\bullet$};
  \node (0) at (210:1) {\color{red}$\bullet$};
  \node (0) at (150:1) {\color{red}$\bullet$};
  \node (0) at (90:1) {\color{red}$\bullet$};
  \node(0) at (90:0.01) {$C_6$};
 \end{tikzpicture}
 \begin{tikzpicture}[scale=1]
  \node (0) at (0,0) {};
  \node (0) at (0,1) {$\leftarrow$ dual to each other $\rightarrow$};
 \end{tikzpicture}
  \begin{tikzpicture}[scale=1.2]
 \draw[thick,red]  (60:1)--(120:1)--(180:1)--(240:1)--(-60:1)--(0:1)--(60:1);
   \node (0) at (60:1) {\color{black}$\bullet$};
   \node (0) at (0:1) {\color{black}$\bullet$};
   \node (0) at (-60:1) {\color{black}$\bullet$};
   \node (0) at (240:1) {\color{black}$\bullet$};
    \node (0) at (180:1) {\color{black}$\bullet$};
  \node (0) at (120:1) {\color{black}$\bullet$};
   \node(0) at (90:0.01) {$\mathcal{H}C_6$};
 \end{tikzpicture}    
\end{center}
    \caption{A cycle graph on $N$ nodes is dual to itself (here $N=6$).}
    \label{fig:cycle}
\end{figure}

\begin{figure}[t]
\begin{center}
$P_6$: 
\begin{tikzpicture}[scale=0.7]
\draw (1,0)--(2,0)--(3,0)--(4,0)--(5,0)--(6,0);
\node (1) at (1,0) {\color{red}$\bullet$}; 
\node (2) at (2,0) {\color{red}$\bullet$}; 
\node (2) at (3,0) {\color{red}$\bullet$}; 
\node (2) at (4,0) {\color{red}$\bullet$}; 
\node (2) at (5,0) {\color{red}$\bullet$}; 
\node (2) at (6,0) {\color{red}$\bullet$}; 
\end{tikzpicture}    
~~~~~ $\mathcal{H}P_6$: \begin{tikzpicture}[scale=0.7]
\draw[red] (1.1,0)--(2,0)--(3,0)--(4,0)--(5,0)--(6,0)--(6.9,0);
\node (2) at (2,0) {$\bullet$}; 
\node (2) at (3,0) {$\bullet$}; 
\node (2) at (4,0) {$\bullet$}; 
\node (2) at (5,0) {$\bullet$}; 
\node (2) at (6,0) {$\bullet$}; 
\end{tikzpicture}   

\vspace{0.3cm}

$P_5^D$: 
\begin{tikzpicture}[scale=0.7]
\draw[red] (1.1,0)--(2,0)--(3,0)--(4,0)--(5,0)--(6,0)--(6.9,0);
\node (1) at (1,0) {$\circ$}; 
\node (2) at (2,0) {$\bullet$}; 
\node (2) at (3,0) {$\bullet$}; 
\node (2) at (4,0) {$\bullet$}; 
\node (2) at (5,0) {$\bullet$}; 
\node (2) at (6,0) {$\bullet$}; 
\node (2) at (7,0) {$\circ$}; 
\end{tikzpicture}    
\end{center}
    \caption{$P_N$ and $P_{N-1}^D$ path graphs for $N=6$, with $\bullet$ representing  interior vertices  and $\circ$ denoting boundary vertices.}
    \label{fig:path}
\end{figure}

\begin{example}[dual of $C_N$]
Let $\mathcal{G}=C_N$ be the cycle graph on $N$ vertices. Then $\mathcal{H}C_N$ is also a cycle graph on $N$ vertices, and thus we write $\mathcal{H}C_N\cong C_N$. We refer to Figure~\ref{fig:cycle} for an illustration in which the black edges in $\mathcal{G}$ correspond to the black vertices in $\mathcal{HG}$, and the red vertices in $\mathcal{G}$ correspond to the red edges in $\mathcal{HG}$. This will be used in Example~\ref{ex:odd-cycle}.
\end{example}

\begin{example}[dual of $P_N$]
Let $P_N$ be the path graph on $N$ vertices, and let $P_{N-1}^D$ be the path graph with $(N-1)$ interior vertices and two boundary vertices. Then $\mathcal{H}P_N$ is the hypergraph with $N-1$ vertices and $N$ hyperedges (see Figure~\ref{fig:path}). We point out that a canonical way to identify $\mathcal{H}P_N$ with $P_{N-1}^D$ is to disregard the two boundary vertices of $P_{N-1}^D$. Hence, we write $\mathcal{H}P_N\approx P_{N-1}^D$, and \Cref{Thm:duality} yields a relationship between the eigenvalues of the $p$-Laplacian on path graphs (unsupervised learning) and the $q$-Laplacian on path graphs under homogeneous Dirichlet boundary conditions (semi-supervised learning), where $q$ is the H\"older conjugate of $p$.
\end{example}


\section{Nodal Domains of the graph $p$-Laplacian eigenfunctions}\label{sec:nodal_domains}
The nodal domains induced by a function $f$ are the maximal subdomains where $f$ has constant sign.
The study of nodal domains of the linear Laplacian modes and their count has a long history and many important applications, both in the infinite-dimensional and the graph settings \cite{daneshgar2010isoperimetric, daneshgar2012nodal, Band, Blum, Courant}. In particular, for the linear Laplacian operator it is well known that the number of nodal domains induced by the modes is strongly related to their frequencies \cite{Band, Berkolaiko1, Berkolaiko2, Biy1, Biy2, Davies, Fiedler, Xu, Oren}.
The same problem has also been investigated for the nonlinear graph $p$-Laplacian. In this case, the nodal domain count of the eigenfunctions is strongly related to the position of the corresponding eigenvalues with respect to the variational spectrum \cite{Tudisco1,DEIDDA2023_nod_dom, zhang2023pLap_noddom, ZhangNodalDO}. As we will discuss in the next sections, such results have relevant applications in the study of graph partitioning and clustering (see, e.g., \cite{Tudisco1, deidda2023PhdThesis, ZhangNodalDO}).
We start with the formal definition of nodal domains.
\begin{definition}[Nodal domains]\label{Intro_def:nodal-domains}
  Given a graph $\Gc$ and a function
  $f:\internalnodes\rightarrow \mathbb{R}$, a subset of vertices
  $A\subseteq \internalnodes$ is a strong nodal domain induced by $f$ if the subgraph $\Gc_A\subset\Gc$ with vertices in $A$ is a maximal connected subgraph of $\Gc$ where $f$ is nonzero and has constant sign. Similarly, $A\subseteq \internalnodes$ is a weak nodal domain induced by $f$ if the subgraph $\Gc_A\subset\Gc$ with vertices in $A$ is a maximal connected subgraph of $\Gc$ where $f$ is (not strictly) positive or negative.
  We denote by $\SNc(f)$ and $\WNc(f)$ the number of strong and weak nodal domains induced by a function $f$.
\end{definition}
Note that, by definition, the inequality $\WNc(f)\leq \SNc(f)$ always holds, with equality whenever $f$ is everywhere different from zero.
Moreover, as we rigorously prove in \cref{lemma_nodal_domain_boundary_eigenvalue}, the behavior of a $p$-Laplacian eigenfunction $f$ on a nodal domain (either weak or strong) $A$ fully characterizes the corresponding eigenvalue $\lambda$. Before stating the lemma, we introduce the notation $E(A,A^c)$ for the set of edges connecting the nodal domain $A$ with its complement:
\begin{equation}\label{nod_domain_boundary}
E(A,A^c)=\{(u,v)\in\edgeset\;s.t.\; u\in A,\;v\in A^c\;\text{or}\; u\in A^c,\;v\in A\}\,,
\end{equation}
where $A^c=\nodeset\setminus A$ (i.e. it includes the boundary if the boundary is non empty).
\begin{lemma}\label{lemma_nodal_domain_boundary_eigenvalue}
Let $(f,\Lambda)$ be a $p$-Laplacian eigenpair, $p>1$, $A$ a nodal domain induced by $f$ and $E(A,A^c)$ the boundary of $A$. Then
$$\lambda=\frac{\|\grad f|_{E(A,A^c)}\|_{\edgelength,p-1}^{p-1}}{\|f|_A\|_{\nodeweight,p-1}^{p-1}}$$
where $\grad f|_{\partial A}(u,v)=\Ind_{\edgeset(A,A^c)}(u,v)\grad f(u,v)$ and $f|_A(u)=\Ind_A(u) f(u)$.
\end{lemma}
\begin{proof}
The proof follows by summing the eigenvalue equations over all nodes $u\in A$:
\begin{equation}
    \sum_{u\in A}\plap f(u)=\sum_{u\in A}\lambda\nodeweight_u|f(u)|^{p-2}f(u).
\end{equation}
Since $A$ is a nodal domain, we may assume without loss of generality that $f(u)\geq 0$ for any $u\in A$, which means that the sum on the right-hand side equals the $p-1$-``norm'' of $f|_A$. Similarly, on the left-hand side we have
\begin{equation}\label{eq1_p-p-1eigenvalue}
\begin{aligned}
     \sum_{u\in A}\sum_{v\sim u}\edgelength_{uv}|\grad f(u,v)|^{p-2}\grad f(u,v)=& \sum_{u\in A}\sum_{\substack{v\sim u\\v\in A}}\edgelength_{uv}|\grad f(u,v)|^{p-2}\grad f(u,v)+\\
     &+\sum_{u\in A}\sum_{\substack{v\sim u\\v\in A^c}}\edgelength_{uv}\big(\grad f(v,u)\big)^{p-1},
\end{aligned}
\end{equation}
where we have used that $Kf(v,u)$ is positive for any $u\in A$ and $v\in A^c$. Finally, it is enough to observe that the first sum on the left-hand side of \eqref{eq1_p-p-1eigenvalue} is equal to zero by definition of $\grad$: if both $u,v\in A$, then $Kf(u,v)=-Kf(v,u)$, while the second sum equals $\|\grad f|_{E(A,A^c)}\|_{\edgelength,p-1}^{p-1}$.
\end{proof}
Note that in \cref{lemma_nodal_domain_boundary_eigenvalue}, with a small abuse of notation, we wrote $\|f\|_{\nodeweight,p-1}^{p-1}=\sum_u \nodeweight_u|f(u)|^{p-1}$ and $\|\grad f\|_{\edgelength,p-1}^{p-1}=(1/2)\sum_{(u,v)}  \edgelength_{uv} |\grad f(u,v)|^{p-1}$ also for $p\in(1,2)$, although these functionals do not define norms in this range.
\Cref{lemma_nodal_domain_boundary_eigenvalue} is particularly relevant because it shows that the nodal domains associated with an eigenfunction induce a form of balanced partition of the graph. This fact, as we will see in more detail in \cref{Sec:1-Lapl_spectrum} and \cref{Sec:inf_eigenproblem}, allows us to establish relations between certain ``optimal'' partitioning problems of the graph and the $p$-Laplacian variational spectrum.
Note also that the right-hand-side expression for $\lambda$ in \Cref{lemma_nodal_domain_boundary_eigenvalue} is homogeneous. Thus, larger eigenvalues are expected to correspond to smaller nodal domains and higher nodal counts.
Classical results relate the number of nodal domains induced by an eigenfunction of the linear Laplacian to the corresponding frequency \cite{Berkolaiko1, Berkolaiko2, Xu}.
In the remainder of this chapter we recall that analogous results extend to the nonlinear $p$-Laplacian setting for values of $p\in (1,\infty)$.
We start by discussing the nodal count of the first variational eigenpair, which admits a Perron--Frobenius characterization (see \cite{DEIDDA2023_nod_dom, Hua}).
\begin{theorem}\label{Theorem_1st_eigen_characterization}
Let $\lambda _{1}$ be the first variational $p$-Laplacian eigenvalue on a connected
graph $\Gc$. Then
\begin{enumerate}
\item $\lambda _{1}$ is simple and the corresponding eigenfunction
$f_{1}$ (up to a multiplicative constant) is strictly positive, i.e., $f_{1}(u)>0$ for all $u\in \internalnodes$;
\item if $g$ is an eigenfunction corresponding to an eigenvalue
$\lambda $ and $g(u)>0$ for all $u\in \internalnodes$, then
$\lambda =\lambda _{1}$.
\end{enumerate}
\end{theorem}
Note that the above theorem is trivial in the non-boundary case in which the matrix $\grad$ has a non-empty kernel given by the nodewise constant function. The boundary case, instead, is more technical, and we refer to \cite{DEIDDA2023_nod_dom, Hua} for a proof.
As a consequence of \cref{Theorem_1st_eigen_characterization}, any eigenfunction that is not the first variational one admits at least two nodal domains. 
More in detail, the number of nodal domains generated by $p$-Laplacian eigenfunctions can be bounded both from above and below in terms of the position of the corresponding eigenvalue with respect to the variational spectrum; see \cite{Tudisco1,DEIDDA2023_nod_dom, zhang2023pLap_noddom} for details. We summarize these results in \cref{thm:nodal-count}.

In the theorem we denote by $\betti$ the number of independent cycles of the graph. This number is also referred to as the first Betti number of the graph from an algebraic topology perspective. In particular, $\betti$ equals the dimension of the kernel of the transpose of the graph incidence matrix. Indeed, each loop in the graph can be represented by an indicator function in this kernel, and conversely each element in the kernel can be written as a linear combination of indicator functions of loops. In particular, $\betti=0$ means that each connected component of the graph is a tree.

\begin{theorem}\label{thm:nodal-count}
  Suppose that $\Gc$ is a connected graph with
$\betti =|\edgeset|-|\internalnodes|+1$ independent cycles, $1<p<\infty$ and $\lambda_1<\lambda_2\leq\dots\leq\lambda_N$ are the variational eigenvalues of $\plap$ and let $f$ be an eigenfunction of $\plap$ with eigenvalue $\lambda$. Then
  \begin{itemize}
\item If $\lambda < \lambda_{k+1}$
$$ \SNc(f) \leq k\,;$$
\item If $\lambda=\lambda_k$ 
$$ \WNc(f)\leq k.$$
\end{itemize}

Moreover for a function $f:\internalnodes\to \mathbb R$, let $l(f)$ be the number
of independent cycles in $\Gc$ where $f$ has constant sign and
$\lbrace v_{i} \rbrace _{i=1}^{z(f)}$ the nodes such
that $f(v_{i})=0$,
with $z(f)$ being the number of such nodes. Let
$\Gc'=\Gc\setminus \big\{\lbrace v_{i} \rbrace _{i=1}^{z(f)}\cup \boundary\big\}$
be the graph obtained by removing from $\Gc$ all the nodes where
$f$ is zero, as well as all the edges connected to those nodes. Let
$c(f)$ be the number of connected components of
$\Gc'$ and
$\beta(\Gc')=|\edgeset'|-|\nodeset'|+c(f)$ the number
of independent cycles of the graph $\Gc'$. Then:
\begin{itemize}
    \item
If $f$ is an eigenfunction with eigenvalue
$\lambda $ such that $\lambda >\lambda _{k}$, then $f$ induces strictly
more than $k-\beta +l(f)-z(f)$ nodal domains. Precisely, it holds
$$
    \SNc(f)\geq k-\beta(\Gc')+l(f)-z(f)+c(f).
$$
\item
If $f$ is an eigenfunction corresponding to the variational
eigenvalue $\lambda _{k}>\lambda _{k-1}$ with
$\mathrm{mult}(\lambda_{k})=m$, then:
\begin{equation}
\SNc(f)\geq \max \left\{ \begin{array}{lr}\Big(k-\beta(\Gc')+l(f)-z(f)+(c
(f)-1)\Big)\\ \Big(k-\beta(\Gc')+l(f)-z(f)+(m-1)\Big) \end{array}
\right\}
\end{equation}
\end{itemize}
\end{theorem}
Observing that generally $l(f)\geq 0$, $c(f)\geq 1$, and $\beta(\Gc')\leq \betti$ (where $\betti$ is the number of independent cycles of $\Gc$), from the previous theorem we can conclude that if $f$ is an eigenfunction of $\plap$ with eigenvalue $\lambda$ such that $\lambda>\lambda_{k}$, then
\begin{equation}\label{lower_bound_1_weak}
        \SNc(f)\geq k-\betti-z(f)+1\,,
\end{equation}
and moreover if $f$ is an eigenfunction of $\plap$ corresponding to the variational eigenvalue $\lambda_k$ with $\lambda_k>\lambda_{k-1}$, then
\begin{equation}\label{lower_bound_2_weak}
        \SNc(f)\geq k-\betti-z(f)+\mathrm{mult}(\lambda_k)-1\,,
\end{equation}
where $z(f)$ is the number of nodes where $f$ is zero.
Despite the possible theoretical interest of the quantities $l(f)$, $c(f)$, and $\beta(\Gc')$ in \Cref{thm:nodal-count}, these are usually unknown and depend on the eigenfunction. In contrast, in the weaker formulations \eqref{lower_bound_1_weak} and \eqref{lower_bound_2_weak} the only term that depends on the eigenfunction is the number of zeros $z(f)$. 
It is worth noting that the bounds reported in \cref{thm:nodal-count} depend on the number $\betti$ of independent cycles in the graph, and they become sharper as the graph becomes more tree-like.
In particular, it can be proved that on trees the nodal count reflects exactly the eigenvalue index \cite{DEIDDA2023_nod_dom}. We recall that on such graphs the variational spectrum exhausts the spectrum; see \cref{variational_eigenvalues_on_trees}.

\begin{theorem}[Theorem 3.8 in \cite{DEIDDA2023_nod_dom}]\label{thm_nodal_count_on-trees}
Let $\Gc$ be a tree and assume $f$ is an everywhere nonzero $p$-Laplacian eigenfunction, with $p\in(1,\infty)$. Assume also that $f$ is associated to the eigenvalue $\lambda_k(\plap)$, then $f_k$ induces exactly $k$ strong nodal domains.
\end{theorem}

The corresponding studies for the $1$- and $\infty$-Laplacians are more complicated, mainly because the $1$- and $\infty$-eigenfunctions are not locally uniquely defined.


\subsection{$1$-Laplacian Nodal Domains}\label{subsection_1-Lap_nodal_domains}

In the case $p=1$, the relevance of studying nodal domains becomes particularly evident. Indeed, there is an exact correspondence between topological properties of the nodal domains and the eigenvalues.

To see this, recall that if $f$ is a $1$-Laplacian eigenfunction associated to the eigenvalue $\Lambda$, then by \Cref{Def:generalized_p-eigenpair} there exist $\xi\in\partial \|f\|_1$ and $\Xi\in\partial \|\grad f\|_1$ such that
\begin{equation}\label{EQ:1_Lap_eigenvalue_eq}
-\mathrm{div}\, \Xi(u)=\frac{1}{2}\sum_{v\sim u} \edgelength_{uv} \Big(\Xi(v,u)-\Xi(u,v)\Big)=\Lambda  \nodeweight_u\xi(u)\qquad \forall u\in \internalnodes.
\end{equation}

Now consider a strong nodal domain $A$ induced by $f$, i.e. a maximal connected component where $f$ is strictly positive or strictly negative, and assume without loss of generality that $f>0$ on $A$. The following result establishes the equality between the eigenvalue and the normalized cut $c(A)$ of $A$, defined as the ratio between the total weight of the edges leaving $A$ and its ``mass'':
\begin{equation}\label{cut}
c(A):=\frac{\sum\limits_{\substack{u\in A,\;
v\in A^c, v\sim u }}\edgelength_{uv}}{\nodeweight(A)},
\end{equation}
where $\nodeweight(A)=\sum_{u\in A}\nodeweight_u$.
\begin{lemma}\label{lemma:1-eigenvalues_are_isoperimetric_constants}
Let $(f,\Lambda)$ be a $1$-eigenpair and assume $A$ is a strong nodal domain induced by $f$. Then
$$
\Lambda=c(A)\,.
$$
\end{lemma}
\begin{proof}
Simply sum the eigenvalue equation \eqref{EQ:1_Lap_eigenvalue_eq} over $A$ and use the characterization of the subgradients \eqref{1_subgradient}, i.e.
$\xi(u)=1\;\forall \,u\in A$ and $\Xi(v,u)=-\Xi(u,v)=1\; \forall u\in A,\; v\in A^c$ such that $(u,v)\in \edgeset$.
\end{proof}
In particular, for the $1$-Laplacian operator it is not difficult to see that any admissible nodal domain decomposition induces piecewise constant eigenfunctions. Indeed, assume $(f,\Lambda)$ is a $1$-eigenfunction and $A$ is one of the strong nodal domains induced by $f$. We know that there exist $\Xi\in \partial\|\grad f\|_1$ and $\xi\in\partial\|f\|_1$ such that $K^T \Xi=\Lambda \nodeweight\odot \xi$. It is easy to check that the subgradients $\Xi$ and $\xi$ of $\grad f$ and $f$ also satisfy $\Xi\in \partial\|\grad \chi_A\|_1$ and $\xi\in\partial\|\chi_A\|_1$, where $\chi_A$ denotes the characteristic function of the set $A$, taking value $1$ on $A$ and zero elsewhere. Thus, $\chi_A$ is also an eigenvector corresponding to the eigenvalue $\Lambda$.

In more detail, we have the following result from \cite{ZhangNodalDO}.
\begin{theorem}\label{Thm_nodal_domains_decomposition_1-Lap_eigenfunctions} Let $(f,\Lambda)$ be a $1$-Laplacian eigenpair and let $\{A_i^+\}$, $\{A_i^-\}$ be the strong positive and negative nodal domains induced by $f$, respectively. Then there exists a connected subset  $\Omega_{f,\Lambda}$ of $S_1$ defined as 
$$\Omega_{f,\Lambda}:=\left\{\left.t\sum_i \alpha_i \frac{\chi_{A_i^+}}{\|\chi_{A_i^+}\|_1}-(1-t)\sum_j \beta_j \frac{\chi_{A_i^-}}{\|\chi_{A_i^-}\|_1}\,\right|\,0\le t\le 1,\,\sum\alpha_i=\sum\beta_j=1\right\}$$
such that any element in $\Omega_{f,\Lambda}$ is a $1$-eigenfunction of the eigenvalue $\Lambda$. 
\end{theorem}

In particular, it is possible to prove the following nonlinear version of the Courant nodal domain theorem for the $1$-Laplacian eigenfunctions \cite{Tudisco1,ZhangNodalDO}:

\begin{theorem}\label{Thm:nodal-count-for_1_lap}
    Let $(f,\Lambda)$ be a $1$-Laplacian eigenpair with $\Lambda<\Lambda_k(\Delta_1)$; then $f$ induces at most $k-1$ strong nodal domains.
\end{theorem}
We note that, unlike the case $p>1$, when $p=1$ both the upper bound for the weak nodal count and the lower bound for the strong nodal count from \Cref{thm:nodal-count} are no longer valid. The failure of the lower bound is easily seen as a consequence of \Cref{Thm_nodal_domains_decomposition_1-Lap_eigenfunctions}, while we refer to Example 10 of \cite{ZhangNodalDO} for a counterexample to the upper bound for the weak nodal count.

If we consider special eigenpairs, the upper bound might be improved (e.g., Theorem 4 in \cite{zhang2023pLap_noddom} shows that for any eigenfunction $f$ with minimal support, $\mathcal{SN}(f)= 1$).

Note also that if $\Lambda_2(\Delta_1)$ has multiplicity $1$, then by \Cref{Thm:nodal-count-for_1_lap} any associated eigenfunction has at most two nodal domains (weak or strong). However, when the multiplicity is higher, the nodal count could in principle grow with the multiplicity. The next theorem studies the weak nodal count.

\begin{theorem}\label{thm:second-eigen-1}
Let $G$ be a connected graph. 
Then, any eigenfunction $f$ corresponding to $\Lambda_2(\Delta_1)$ has exactly one positive weak nodal domain or exactly one negative weak nodal domain. Moreover, if $f$ has at least two positive strong nodal domains and two negative strong nodal domains, then the number of weak nodal domains of $f$ is exactly two.
\end{theorem}
\begin{proof}
Assume, without loss of generality, that $\{f> 0\}\ne\emptyset$. If $\{f< 0\}=\emptyset$, then $f$ has exactly one positive weak nodal domain. When $A:=\{f< 0\}\ne \emptyset$, then the property of $1$-Laplacian eigenfunctions implies that the cut $(A,A^c)$ is a Cheeger cut of $G$ (see \cref{thm:1-eigenpair_and_cheeger_constants}). By Theorem 2.11 in \cite{Zhang_1-Lap_cheeger_cut}, one of $A$ and $A^c$ is connected (this property is folklore, and the formulation we used here was first proposed by Zhang in the conference paper \cite{oberwolfach15}). Therefore, either $\{f\ge0\}$ or $\{f< 0\}$ induces a connected subgraph of $G$, which implies that $f$ has either exactly one positive weak nodal domain or exactly one negative weak nodal domain. 

Now assume that $f$ has two positive strong nodal domains and two negative strong nodal domains. Then both $\{f>0\}$ and $\{f<0\}$ induce disconnected subgraphs. Similar to the proof in the previous paragraph, both $\{f\ge0\}$ and $\{f\le0\}$ induce connected subgraphs, which implies that $f$ has exactly one positive weak nodal domain and exactly one negative weak nodal domain.
\end{proof}

We illustrate the discussion in this section on the nodal domain count of $1$-Laplacian eigenpairs in \cref{ex_1_nodal_domain}, \cref{ex_3_nodal_domains} and \cref{ex_4_nodal_domains}. In particular, in \cref{ex_1_nodal_domain} we show that, unlike the case $p>1$, there exist $1$-Laplacian eigenfunctions corresponding to eigenvalues larger than the first variational one that induce a unique strong nodal domain. In \cref{ex_3_nodal_domains} we show that the statement of \cref{thm:second-eigen-1} does not hold if the number of strong positive (or negative) nodal domains is smaller than $2$. Finally, \cref{ex_4_nodal_domains} shows a situation where \cref{thm:second-eigen-1} applies.


\subsection{$\infty$-Laplacian Nodal Domains}\label{subsec_infinity_lap_nodal_dom}

This section presents a novel discussion of the nodal domain count for $\infty$-Laplacian eigenpairs. 
In general, the nodal domain theorems for the graph $p$-Laplacian with $1<p<\infty$ (see \Cref{thm:nodal-count}) fail when we consider $p=\infty$. 
The reason for the failure of the nodal count upper bounds in the case $p=\infty$ is that the $\infty$-eigenequation is unable to capture small oscillations. Indeed, unlike the case $p=1$ (\Cref{lemma:1-eigenvalues_are_isoperimetric_constants}), the eigenvalue does not encode a property of the nodal domains.
To see this, consider an eigenfunction $f$ and a nodal domain $A$ where $|f|< \|f\|_{\infty}$ and $|\grad f|< \|\grad f\|_{\infty}$. From the characterization \eqref{inf_subgradient}, we see that any subgradient $\xi\in \partial\|f\|_{\infty}$ and $\Xi\in \partial\|\grad f\|_{\infty}$ is zero on $A$. Hence, we can modify $f$ to $\tilde{f}$ on $A$, and as long as we maintain $|\tilde{f}|< \|f\|_{\infty}$ and $|\grad \tilde{f}|< \|\grad f\|_{\infty}$, the function $\tilde{f}$ remains an eigenfunction with the same eigenvalue. In particular, we can choose $\tilde{f}$ so that it induces arbitrarily many oscillations.
We analyze this phenomenon in the next example.

\begin{example}\label{example:nodal_domains_inf_eigenpairs}
Let $\Gc_n$ be the graph in \Cref{fig.counterexample_inf_nodal_count} with node set
\begin{equation}
V_n:=\{v,x_1,\cdots,x_n,y_1,\cdots,y_n,z_1,\cdots,z_n\},
\end{equation}
and edge set
\begin{equation}
E_n:=\{\{v,x_1\},\{v,y_1\},\{v,z_1\},\{x_i,x_{i+1}\},\{y_i,y_{i+1}\},\{z_i,z_{i+1}\}:i=1,\cdots,n-1\},
\end{equation}
where $n\ge 2$, and we have no weights on the nodes or on the edges.
\begin{figure}[t]
  \centering
  \begin{tikzpicture}[inner sep=1.5mm, scale=.5, thick]

    \node (1) at (-1,0) [circle,draw] {$v$};
    \node (11) at (3,3) [circle,draw] {$x_1$};
    \node (21) at (3,0) [circle,draw] {$y_1$};
    \node (31) at (3,-3) [circle,draw] {$z_1$};
    
    \node (12) at (6,3) [circle,draw] {$x_2$};
    \node (22) at (6,0) [circle,draw] {$y_2$};
    \node (32) at (6,-3) [circle,draw] {$z_2$};

    \node (13) at (9,3) [circle,draw] {$x_3$};
    \node (23) at (9,0) [circle,draw] {$y_3$};
    \node (33) at (9,-3) [circle,draw] {$z_3$};
    
    \draw [-] (1.north) -- (11.west);
    \draw [-] (1.east) -- (21.west);
    \draw [-] (1.south) -- (31.west);
    \draw [-] (11.east) -- (12.west);
    \draw [-] (21.east) -- (22.west);
    \draw [-] (31.east) -- (32.west);
    \draw [-] (12.east) -- (13.west);
    \draw [-] (22.east) -- (23.west);
    \draw [-] (32.east) -- (33.west);
    
     \end{tikzpicture}
     \caption{Example graph where the standard nodal domain count bounds fail for the $p=\infty$ case, see \cref{example:nodal_domains_inf_eigenpairs} for the discussion.}\label{fig.counterexample_inf_nodal_count}
\end{figure}
Consider $f$ defined by
\begin{equation}
f(v)=0, \quad f(x_i)=\frac in,\quad f(y_i)=-\frac in \quad \text{and}\quad f(z_i)=(-1)^i\frac {1}{2n}.
\end{equation}
Then, by \Cref{thm:k-inequality} (which we will recap in the next section), $\Lambda_2(\inflap)=1/n$ is the second variational $\infty$-eigenvalue, and $f$ is an eigenfunction corresponding to $\Lambda_2(\inflap)$.
Clearly, the number of strong nodal domains of $f$ is $2+n$, and the number of weak nodal domains of $f$ is $1+n$. Next, we show that the multiplicity of $\Lambda_2(\inflap)$ is $2$.
In fact, by \Cref{prop:Infinity_eigenfunctions_characterization}, the set $X$ of eigenfunctions corresponding to $\Lambda_2$ in the unit $l^\infty$-sphere is $X=X_{xy}\cup X_{yx}\cup X_{xz}\cup X_{zx}\cup X_{yz}\cup X_{zy}$, where
\begin{equation}
X_{xy}=\left\{f \; \Bigg|\; \begin{array}{lr} f(x_i)= i/n=-f(y_i)\,,
\\f(v)=0\,,
\\|f(z_1)|\leq 1/n,\;|f(z_i)-f(z_{i+1})|\leq 1/n\end{array}\right\}
\end{equation}
and the other five sets $X_{yx}$, $X_{xz}$, $X_{zx}$, $X_{yz}$, $X_{zy}$ are defined in a similar manner. We let $\Phi:X\to \R^3\setminus\{0\}$ be defined as $\Phi(f)=(f(x_n),f(y_n),f(z_n))$. Then $\Phi$ is an odd continuous map, and it is not difficult to see that $\Phi(X)\cong \mathbb{S}^1$, which implies that the genus of $X$ is $2$ and the variational multiplicity of $\Lambda_2(\inflap)$ is at most $2$; see \Cref{lemma_genus_and_multiplicity}. In consequence, since the strong nodal domain count exceeds $3$ and the weak nodal domain count exceeds $2$, the usual Courant nodal domain theorem fails if we use the ordinary strong or weak nodal domains.
\end{example}

To partially overcome the failure of the nodal domain theorems in the case $p=\infty$, in the remainder of this section we restrict attention to a particular class of nodal domains, which we call \textit{perfect} ones. A perfect nodal domain is one where the eigenfunction reaches its maximal oscillatory amplitude. Relying on these nodal domains, we show that it is possible to recover some Courant-like nodal domain theorems.
We start with the formal definition.

\begin{definition}[perfect nodal domain]\label{def:perfect-nodal}
A subset $U\subset \internalnodes$ is called a positive (resp., negative) perfect nodal domain of a function $f:\internalnodes\to\R$, if $U$ is a nodal domain and $U$ contains a vertex $v$ satisfying $|f(v)|=\|f\|_\infty$. We use $\PNc(f)$ to denote the number of perfect nodal domains induced by a function $f$.
\end{definition}

Fortunately, by adding some natural conditions, we can establish the following Courant nodal domain theorem for the graph $\infty$-Laplacian.

\begin{theorem}\label{thm:infty_nodal-count}
Let $f$ be an eigenfunction corresponding to an $\infty$-eigenvalue $\Lambda$ with $\Lambda< \Lambda_{k+1}(\inflap)$. Then we have:
$$\PNc(f)\le k.$$
\end{theorem}
\begin{proof}
 Suppose that $f$ has $c$ perfect nodal domains, denoted by $V_1,\ldots,V_{c}$.
Consider a linear function-space $X$ spanned by $f|_{V_1},\ldots,f|_{V_{c}}$, where the restriction $f|_{V_i}$ is defined by
\begin{equation}
    f|_{V_i}(v)=\begin{cases}f(v),&\text{ if }v\in V_i,\\ 0,&\text{ if } v\not\in V_i.\end{cases}
\end{equation}
Since $V_1,\ldots,V_{c}$ are pairwise disjoint, $\dim X=c$.
We prove that $\rayl_\infty(g)\leq \Lambda_k(\inflap)$ for any $g\in X$. Note that for any $g\in X\setminus 0$, there exists $(t_1,\ldots,t_{c})\ne\mathbf{0}$ such that $g=\sum_{i=1}^{c} t_i f|_{V_i}$. Then, since the $V_i$ are disjoint,
\begin{align}
\|g\|_{\infty}&=\|\sum_{i=1}^{c}t_if|_{V_i}\|_{\infty}=\max\limits_{v\in \internalnodes}\sum_{i=1}^{c}|t_i|\cdot|f|_{V_i}(v)|=\max\limits_{i}\max\limits_{v\in V_i}\sum_{i=1}^{c}|t_i|\cdot|f|_{V_i}(v)|
\\&=\max\limits_{i}\max\limits_{v\in V_i}|t_i|\cdot|f|_{V_i}(v)|=\max\limits_{i}|t_i|\cdot\|f|_{V_i}\|_\infty=\max\limits_{i=1,\ldots,c}|t_i|\cdot\|f\|_\infty.
\end{align}
Moreover, for any $(v,v')\in E$
\begin{equation}\label{eq:equa-1}
\begin{aligned}
|g(v)-g(v')|&=\left| \sum_{i=1}^{c}t_i (f|_{V_i}(v)-f|_{V_i}(v'))\right|\leq\max\limits_{i=1,\ldots,c}|t_i|\cdot\sum_{i=1}^{c} \Big|f|_{V_i}(v)-f|_{V_i}(v')\Big|
\\&
\leq\max\limits_{i=1,\ldots,c}|t_i|\cdot |f(v)-f(v')|
\end{aligned}
\end{equation}
where the second inequality in \eqref{eq:equa-1} is derived by the nodal domain decomposition. Indeed, if $v$ and $v'$ are in the same nodal domain $V_i$, then $f(v)=f_{V_i}(v)$, $f(v')=f_{V_i}(v')$ and $f_{V_j}(v)=f_{V_j}(v')=0$ for any $j\neq i$, proving the equality. Otherwise, if $v\in V_i$ and $v'\not\in V_i$ (or vice versa), then by the nodal domain definition, necessarily $f(v)f(v')\leq 0$, yielding
\begin{equation}
    \sum_{i=1}^{c} \big|f|_{V_i}(v)-f|_{V_i}(v')\big|=\big|f|_{V_i}(v)\big| +\big|\sum_{j\neq i} f|_{V_j}(v')\big|\leq \big|f(v)\big| +\big| f(v')\big| =\big|f(v)-f(v')\big|.
\end{equation}

In particular, for any $g\in X$,
$$\mathcal{R}_\infty(g)\le \frac{\max\limits_{i=1,\ldots,c}|t_i|\cdot\max\limits_{(v,v')\in E}\edgelength_{vv'}|f(v)-f(v')|}{\max\limits_{i=1,\ldots,c}|t_i|\cdot\|f\|_\infty}=\rayl_{\infty}(f)=\Lambda.$$
By the definition of the variational eigenvalues, we have
\begin{equation}
\Lambda_{c}(\inflap)=\inf\limits_{X'\in \mathcal{F}_c(S_{\infty})}\sup\limits_{g'\in X'}
\mathcal{R}_\infty(g')\le \sup\limits_{g\in X\cap S_{\infty}}
\mathcal{R}_\infty(g)< \Lambda_{k+1}(\inflap)
\end{equation}
which implies $c\le k$.
\end{proof}


\subsubsection{Nodal domains of limit and viscosity $\infty$-eigenpairs}
We devote this section to develop the nodal domain theory for the viscosity $\infty$-eigenpairs and the $\infty$-eigenpairs obtained as limits of $p$-eigenpairs.

We start by studying the regularity of the nodal domains when varying $p$. In particular, we consider $\WNc(f)$, $\SNc(f)$ and $\PNc(f)$ as functionals from the node functions space $\mathcal{H}(\internalnodes)$ to $\N$, associating to any function the corresponding number of weak, strong and perfect nodal domains, respectively. Then, in the next lemma we discuss the regularity of these discrete functions.

\begin{lemma}\label{lemma:nodal-domain-limit}
Consider the functions $\WNc, \SNc, \PNc$ from $\mathcal{H}(\internalnodes)$ to $\N$, then: 
\begin{itemize}
    \item The function $\WNc(f)$ is lower semicontinuous, i.e. given a convergent sequence of functions $\{f_k\}$ with $\lim\limits_{k\to+\infty} f_k=f$ then $$\liminf_{k\rightarrow \infty}\mathcal{WN}(f_k)\ge\mathcal{WN}(f).$$
In fact, $\mathcal{WN}(\tilde{f})\ge\mathcal{WN}(f)$ for any $\tilde{f}$ that is sufficiently close to $f$.
    \item The function $\PNc(f)$ is upper semicontinuous, i.e. given a convergent sequence of functions $\{f_k\}$ with $\lim\limits_{k\to+\infty} f_k=f$ then $$\limsup_{k\rightarrow \infty}\mathcal{PN}(f_k)\le\mathcal{PN}(f).$$
In fact, for any $\tilde{f}$ sufficiently close to $f$, we have $\mathcal{PN}(\tilde{f})\le\mathcal{PN}(f)$.
    \item The function $\SNc(f)$ is neither lower nor upper semicontinuous. 
\end{itemize}
\end{lemma}

\begin{proof}
Start with the lower semicontinuity of $\WNc$ and assume by contradiction that there exist $f:\internalnodes\to\R$ and $f_k\to f$ such that $\liminf_{k\to\infty} \WNc(f_k)< \WNc(f)$.
 Note that for sufficiently large $k$, 
 \begin{equation}
 \{v:f(v)>0\}\subset \{v:f_k(v)>0\}\quad  \text{and}\quad  \{v:f(v)<0\}\subset \{v:f_k(v)<0\}.
 \end{equation}
 Let $U^+_1,\cdots,U^+_{m_+}$ be the positive weak nodal domains of $f$, and let $U^-_1,\cdots,U^-_{m_-}$ be the negative weak nodal domains of $f$. Let $v^+_i\in U^+_i$ and $v^+_j\in U^+_j$ with $i\ne j$ and $f(v^+_i),f(v^+_j)>0$. Then, for any path $v_0=v^+_i\sim v_1\sim \cdots \sim v_l=v^+_j$, there exists a vertex $v_m$ in the path such that $f(v_m)<0$. Thus we have, 
 \begin{equation}
     f_k(v_m)<0<\min\{f_k(v^+_i),f_k(v^+_j)\},
 \end{equation}
 meaning that $v^+_i$ and $v^+_j$ lie in different positive weak nodal domains of $f_k$. In consequence, the vertices $v^+_1,\cdots,v^+_{m_+}$, $v^-_1,\cdots,v^-_{m_-}$, 
are in some weak nodal domains of $f_k$ that are different from each other. Thus, $\WNc(f)=m_++m_-\le \WNc(f_k)$. 

To prove the second part assume that $A_1,\dots, A_k$ are the $\delta$-height strong nodal domains of $f$ i.e. the nodal domains with $\|f|_{A_i}\|_{\infty}\geq \delta$ and denote by $\SNc_\delta(f)$ their count.
Define 
\begin{equation}
    \nodeset_{\delta}(f):=\{v\in\internalnodes:\; |f(v)|\geq \delta\}\subset \cup_{i=1}^k A_i.
\end{equation} 
Next, let 
\begin{equation}
    \epsilon_1=\min_{v\in \nodeset\setminus \nodeset_{\delta}(f)} \delta-|f(v)| \quad \text{and}\quad  \epsilon_2=\min_{v\in\cup A_i}|f(v)|,
\end{equation} then for any $f_k$ such that $\|f_k-f\|_{\infty}<\min \{\epsilon_1/2,\epsilon_2\}$
it is easily proved that $\nodeset_{\delta}(f_k)\subset \nodeset_{\delta}(f)$ and $f_k(v)f(v)>0$ for any $v\in\cup A_i$. The last properties imply that for any $\delta$-height nodal domain of $f_k$, say $A_j^k$, there exists some $\delta$-height nodal domain of $f$, say $A_{i_j}$, such that $A_{i_j}\subset A_j^k$. This, in turn, implies that $\SNc_{\delta}(f_k)\leq \SNc_{\delta}(f)$ which means $\limsup_{f_k \rightarrow f} \SNc_{\delta}(f_k)\leq \SNc_{\delta}(f)$, i.e. upper semicontinuity of the $\delta$-height
 strong nodal domain count.
 To conclude about the perfect nodal domains, just consider 
 \begin{equation}
     \epsilon=\max\{|f(v)| \text{ with } v \text{ s.t. } |f(v)|<\|f\|_{\infty}\} \quad \text{and}\quad \delta=\|f\|_{\infty}-\epsilon/2,
 \end{equation}
 then $\SNc_{\delta}(f)=\PNc(f)$.
 Moreover, if we consider $f_k$ sufficiently close to $f$, it is easy to observe that any perfect nodal domain of $f_k$ is also a $\delta$-height nodal domain of $f_k$, i.e. $\SNc_{\delta}(f_k)\geq \PNc(f_k)$. Thus 
 \begin{equation}
     \limsup_{f_k\rightarrow f} \PNc(f_k)\leq \limsup_{f_k\rightarrow f} \SNc_{\delta}(f_k)\leq \SNc_{\delta}(f)=\PNc(f)\,.
 \end{equation}

Consider $f\ge0$ with $\SNc(f)\ge2$ and let $f_k=f+1/k>0$, then $\SNc(f_k)=\WNc(f_k)=\WNc(f)=1<\SNc(f)$. This means that $\SNc(f)$ is not lower semicontinuous.

Consider $A\neq \internalnodes$ and $f=1_A$ with $\SNc(f)=1$,  and let $f_k=f-1/2k$. Then,  $\SNc(f_k)=\WNc(f_k)\ge2>1=\WNc(f)=\SNc(f)$, which proves that $\SNc(f)$ is not upper semicontinuous.
\end{proof}
Next we need another technical result about the growth rate of viscosity solutions on any nodal domain. 
To this end, let $A\subset \internalnodes$ and denote by $\chi^\nodeset_A$ the corresponding characteristic function on the nodes and by $\chi^\edgeset_A$ the characteristic function of the edges having at least one ending point in $A$. Then, for any function $f\in \mathcal{H}(\internalnodes)$ we write $f|_A(u):=f(u) \cdot\chi^\nodeset_A(u)$ and $K|_A f(u,v):=Kf(u,v)\cdot\chi^\edgeset_A(u,v)$.

\begin{lemma}\label{Lemma:grow_rate_viscosity_solution}
 Let $(f,\Lambda)$ be an $\infty$-viscosity eigenpair and let $A$ be a strong nodal domain of $f$. Then 
 $$\|\grad|_Af \|_{\infty}\leq \Lambda \|f|_A\|_{\infty}.$$
\end{lemma}
\begin{proof}
To prove it, we have to show that
\begin{equation}\label{eq_1_Lip_viscosity_eigenpairs}
\|\grad|_Af \|_{\infty}=\max\{\edgelength_{vv'}|f(v)-f(v')|\,:\; v\in A,\, v'\sim v\}\leq \Lambda\|f|_{A}\|_\infty.
\end{equation}
Note that the vertices in $A:=\{v_1,\cdots,v_c\}$ can be ordered  based on the magnitude of $f$ as $|f(v_1)|\ge |f(v_2)|\ge\cdots\ge |f(v_c)|>0$. Next, we prove by induction on the index $i$ induced by this ordering that $\|\{\grad f(v_i)\}\|_{\infty}\leq \Lambda \|f|_A\|_{\infty}$ for all $i$, which is  equivalent to \eqref{eq_1_Lip_viscosity_eigenpairs}.
Since $f$ is a viscosity eigenfunction corresponding to  $\Lambda$, if $u=v_1$, $f(u)=\|f|_{A}\|_\infty$ and we have
\begin{equation}
   \|\{Kf(u)\}\|_{\infty}= \max_{v'\sim u}\;\edgelength_{vv'}|f(u)-f(v')|=\Lambda\|f|_{A}\|_\infty.
\end{equation}

Now assume the thesis to be true for any $i < k$ and prove it for $i=k$. Observe that $|f(v_{k})|\leq \|f|_{A}\|_\infty$ and so either  
\begin{equation}\label{eq_2_Lip_viscosity_eigenpairs}
\|\{Kf(v_{k})\}\|_{\infty}=
\max_{v'\sim v_{k}}\,\edgelength_{v_k v'}|f(v_k)-f(v')|=\Lambda |f(v_k)|<\Lambda\|f|_{A}\|_\infty, 
\end{equation}
or 
\begin{equation}\label{eq_3_Lip_viscosity_eigenpairs}
\begin{aligned}
     \max_{v'\sim v_k}\,\edgelength_{v_k v'}(f(v_k)-f(v'))_+&=\|\{Kf(v_k)^-\}\|_{\infty}\\
     &=\|\{Kf(v_k)^+\}\|_{\infty}=\max_{v'\sim v_k}\edgelength_{v_k v'}(f(v')-f(v_k))_+.
\end{aligned}
\end{equation}
If \eqref{eq_2_Lip_viscosity_eigenpairs} holds we already have the thesis so consider the case where \eqref{eq_3_Lip_viscosity_eigenpairs} holds.  In this case, there exists $v_i\in A$ with $|f(v_i)|> |f(u)|$, meaning $i<k$, and $\|\{Kf(v_i)\}\|_{\infty}\geq \|\{Kf(v_k)\}\|_{\infty}$. In particular, by induction, we have
\begin{equation}
\|\{Kf(v_k)\}\|_{\infty}\leq \|\{Kf(v_i)\}\|_{\infty}\leq \Lambda\|f|_A\|_{\infty}, 
\end{equation}
which concludes the proof.
\end{proof}
\begin{theorem}\label{Thm_nodal_count_viscosity_eigenfunctions}
    Let $f$ be an  eigenfunction  corresponding to an $\infty$-eigenvalue $\Lambda$ with $\Lambda<\Lambda_{k+1}(\inflap)$. The following inequalities are satisfied:

\begin{itemize}
\item If  $f$ is a viscosity  eigenfunction, then $ \SNc(f)\leq k.$ 
\item If  $f$ is a limit eigenfunction and $$\Lambda_{k-1}(\inflap)<\Lambda=\Lambda_k(\inflap)=\Lambda_{k+m-1}(\inflap)<\Lambda_{k+m}(\inflap),$$
then $\WNc(f_k)\le k$ and $\SNc(f_k)\le k+m-1$.
\end{itemize}

\end{theorem}

\begin{proof}

 Suppose that $f$ is a viscosity eigenfunction corresponding to  $\Lambda$, and $f$ has $m$ strong nodal domains which are denoted by $V_1,\ldots,V_{m}$.   
Similarly to the proof of \Cref{thm:infty_nodal-count}, we have
$\|\sum_{i=1}^{m} t_i f|_{V_i}\|_{\infty}=\max_i |t_i|\|f|_{V_i}\|_{\infty}$. Moreover, for any $(v,v')\in\edgeset$, assuming that $v\in V_{i_0}$ and $v\in V_{i_1}$  (the case $v$ and/or $v'\not\in \cup V_i$ can be studied analogously)
\begin{equation}
\begin{aligned}
|\grad f(v,v')|&=   \edgelength_{vv'}\Big|\sum_{i=1}^{m} t_i \Big( f|_{V_i}(v)-f|_{V_i}(v')\Big)\Big|= \edgelength_{vv'} \Big| t_{i_0} f(v)- t_{i_1} f(v')\Big|\\
    &\leq \max_{i}\{|t_{i_0}|,|t_{i_1}|\}\,\edgelength_{vv'}|f(v)-f(v')|\leq \max_{i}|t_i|\,\max_{u\in V_i, w\sim u }\edgelength_{uw}|f(u)-f(w)|\\
    &=\max_i |t_i|\:\|K|_{V_i}\|_{\infty} 
\end{aligned}
\end{equation}
where the first inequality follows noting that either $i_0=i_1$, i.e. $v$ and $v'$ are in the same nodal domain or $f(v)f(v')<0$, when $v$ and $v'$ are connected but belong to different nodal domains.
As a direct consequence of \Cref{Lemma:grow_rate_viscosity_solution}, we have the following upper bound
\begin{equation}
\begin{aligned}
    \mathcal{R}_\infty\left(\sum_{i=1}^{m} t_i f|_{V_i}\right)
    &\le \frac{\max_i |t_i|\:\|K|_{V_i}\|_{\infty} }{\max_i|t_i|\cdot\|f|_{V_i}\|_\infty}\leq \max\limits_{i=1,\ldots,m}\frac{\|K|_{V_i}\|_{\infty}}{\|f|_{V_i}\|_\infty}\leq \Lambda.
\end{aligned}
\end{equation}
Thus, the same argument from the proof of \Cref{thm:infty_nodal-count} yields $\SNc(f)\leq k$ .

About the second point, the upper bound for the number of strong nodal domains is a direct consequence of point 1 and \Cref{theorem:viscosity_eigenpairs_are_generalized_eigenpairs}.
To bound the weak nodal count observe that by Lemma \ref{lemma:p-monotonic}, $\lim_{p\to+\infty}\lambda_j(\Delta_p)^{\frac1p}=\Lambda_j$ for any $j$. Thus, for sufficiently
large $p>1$, 
\begin{equation}
    \lambda_{k-1}(\Delta_p)<\lambda_{k}(\Delta_p)\le\cdots\le \lambda_{k+m-1}(\Delta_p)<\lambda_{k+m}(\Delta_p),
\end{equation}
which implies that the multiplicity of $\lambda_{k}(\Delta_p)$ is at most $m$. Then, by Theorems 3.3 and 3.5 in \cite{Tudisco1} (see also \Cref{thm:nodal-count}), $\WNc(f_{k,p})\leq k$ and $\SNc(f_{k,p})\le k+m-1$ for any eigenfunction $f_{k,p}$ corresponding to $\lambda_k(\Delta_p)$. 
Finally, the inequality $\WNc(f_k)\le k$ is a direct consequence of Lemma \ref{lemma:nodal-domain-limit}.
\end{proof}

We remark that, as a side consequence of the last proof, if $(f,\Lambda)$ is a limit eigenpair with $(f_{p_j},\Lambda_{p_j})$ converging to $(f,\Lambda)$, then the bounds on the nodal domains in \Cref{Thm_nodal_count_viscosity_eigenfunctions} transfer to the eigenfunctions $f_{p_j}$ for sufficiently large $p_j$. 
We introduce the notation $X_{\Lambda}$ to denote the set of $\infty$-eigenfunctions corresponding to an eigenvalue $\Lambda$ as in \Cref{Def:generalized_p-eigenpair}. Moreover for $\Lambda$ that is a viscosity eigenvalue, we define $X_{\Lambda}^{v}\subset X_\Lambda$  the set of all viscosity eigenfunctions corresponding to $\Lambda$. If $\Lambda$ is not a viscosity eigenvalue we simply set $X_{\Lambda}^v=\emptyset$. The next theorems relate the maximum number of perfect nodal domains associated to a viscosity, or limit, eigenvalue to the number of strong nodal domains.
Before presenting the Theorem, we recall a technical result from \cite{Tudisco1}.
\begin{lemma}\label{Lemma:p-Ryleigh_quotient_nodal_domain_decomposition_space}
    Assume $(f,\lambda)$ is a $p$-Laplacian eigenpair with $p\in[1,\infty)$. Assume that  $A_1,\dots,A_n$ are the weak or strong nodal domains of $f$ and $F=\text{span}\{f|_{A_i}\}$, then 
    $$\sup_{g\in F}\rayl_p^p(g)\leq \lambda.$$ 
\end{lemma}

\begin{theorem}\label{thm:relate-3bounds}
Let $\Lambda$ be an $\infty$-eigenvalue. 
Then 
$$\max\limits_{f\in X_\Lambda}\PNc(f)\ge \max\limits_{f\in X_\Lambda^v}\SNc(f).$$
\end{theorem}

\begin{proof} Let $f$ be a viscosity  eigenfunction   corresponding to $\Lambda$, then from \Cref{theorem:viscosity_eigenpairs_are_generalized_eigenpairs} it is also a general $\infty$-eigenfunction meaning that there exist $\Xi\in \partial\|\grad f\|_{\infty}$ and $\xi\in \partial \|f\|_{\infty}$ such that 
\begin{equation}
    \grad^T \Xi=\Lambda \xi.
\end{equation}
Now suppose that $f$ has $m$ strong nodal domains, which are denoted by $V_1,\cdots,V_m$.
Let 
$$\hat{f}=\sum_{i=1}^m\frac{\|f\|_\infty}{\|f|_{V_i}\|_\infty}f|_{V_i}.$$
Then, $\|\hat{f}\|_{\infty}=\|f\|_{\infty}$ and $\|\hat{f}|_{V_i}\|_\infty=\|f\|_\infty$ for all $i=1,\dots,m$.
Moreover, $\|K\hat{f}\|_{\infty}= \|K f\|_{\infty}$. In order to prove it, assume  that $(u,v)\in \edgeset$, with $u\in V_i$, satisfies $|K\hat{f}(u,v)|=\|K\hat{f}\|_{\infty}$, and assume $v\in V_j$ with $i\neq j$. Then, using $f(u)f(v)<0$, \Cref{Lemma:grow_rate_viscosity_solution} and the expression $\Lambda=\|\grad f\|_{\infty}/\|f\|_{\infty}$, we have the following set of inequalities
\begin{equation}
\begin{aligned}
    \|K\hat{f}\|_{\infty}&=\edgelength_{uv}\bigg|\frac{\|f\|_{\infty}}{\|f|_{v_i}\|_{\infty}}f(u)-\frac{\|f\|_{\infty}}{\|f|_{v_j}\|_{\infty}}f(v) \bigg|\\
    &\leq  \max \bigg\{\frac{\|f\|_{\infty}}{\|f|_{v_i}\|_{\infty}}, \frac{\|f\|_{\infty}}{\|f|_{v_j}\|_{\infty}}\bigg\} \edgelength_{uv}\big| f(u)-f(v)\big|\\
    &\leq \max \bigg\{\frac{\|f\|_{\infty}}{\|f|_{v_i}\|_{\infty}}, \frac{\|f\|_{\infty}}{\|f|_{v_j}\|_{\infty}}\bigg\} \min\Big\{\|K|_{V_i}f\|_{\infty}, \|K|_{V_j}f\|_{\infty}\Big\}\\
    &\leq \Lambda\frac{\|f\|_{\infty}}{\min\big\{\|f|_{v_i}\|_{\infty},\|f|_{v_j}\|_{\infty}\big\}}\min\Big\{\|f|_{V_i}\|_{\infty}, \|f|_{V_j}\|_{\infty}\Big\}
    =\|\grad f\|_{\infty},
\end{aligned}
\end{equation}
the cases $v_j\not\in \cup V_j$ and $v_j\in V_i$
can be dealt analogously. 
Moreover, from the limit eigenvalue equation in \Cref{thm:limitin_infty_eigenvalue_theorem}, if $|f(u)|=\|f\|_{\infty}$, there exists $v\sim u$ such that $|Kf(u,v)|=\|Kf\|_{\infty}$ and then necessarily also $|K\hat{f}(u,v)|=\|Kf\|_{\infty}$.
In particular, it follows that $\Xi\in \partial\|\grad \hat{f}\|_{\infty}$ and $\xi\in \partial\|\hat{f}\|_{\infty}$, proving that also $\hat{f}$ is an $\infty$-eigenfunction of eigenvalue $\Lambda$.
By the definition of (perfect) nodal domains, $\PNc(\hat{f})=\SNc(f)$, and we can conclude $\max_{f\in X_\Lambda}\PNc(f)\ge \max_{f\in X_\Lambda^G}\SNc(f)$. 
\end{proof}

\begin{theorem}\label{Thm:lower_bound_perfect_nodal_count_limit_eigenpairs}
 Let $\Lambda$ be an $\infty$-limit eigenvalue and assume that $\{\Lambda_{p_i}\}$ is a sequence of $p_i$-Laplacian normalized eigenvalues converging to $\Lambda$. Then
 $$\max\limits_{f\in X_{\Lambda}}\PNc(f)\geq \limsup\limits_{p_i\to+\infty}\max\limits_{f\in X_{\Lambda_{p_i}}}\SNc(f).$$
\end{theorem}
\begin{proof}
Note that, up to subsequences, we can consider $f_{p_i}\in X_{\Lambda_{p_i}}$ that converges to some $f\in X_{\Lambda}$
with $\|f_{p_i}\|_{\infty}=1$, $\|f\|_{\infty}=1$ and
\begin{equation}
\SNc(f_{p_i})=\max\limits_{f\in X_{\Lambda_{p_i}}}\SNc(f_{p_i})=\limsup\limits
    _{p_i\to+\infty}\max\limits_{f\in X_{\Lambda_{p_i}}}\SNc(f).
\end{equation}
Then, we construct a sequence $\{\hat{f}_{p_i}\}$, where any function $\hat{f}_{p_i}$ is obtained by stretching  $f_{p_i}$ up to $\|f_{p_i}\|_{\infty}$ on any nodal domain. Literarily, for a generic function $g$ we define
\begin{equation}
\hat{g}=\sum_{i=1}^m\frac{\|g\|_\infty}{\|g|_{V_i}\|_\infty}g|_{V_i},
\end{equation}
where $\{V_i\}$ are the strong nodal domains induced by $g$.
Note that, since the $\{V_i\}$ are strong nodal domains, for any $(u,v)\in \edgeset$ 
\begin{equation}\label{eq1:thm_limisup_nodal_count}
   |Kg(u,v)| \leq |K\hat{g}(u,v)|.
\end{equation}
Next, since $\|\hat{f}_{p_i}\|_{\infty}=1$ for any $i$, the sequence $\{\hat{f}_{p_i}\}$ is bounded and thus, eventually considering a subsequence, we can assume $\hat{f}_{p_i}\to f'$ for some function $f'$ with $\|f'\|_{\infty}=1$. 
According to \Cref{Lemma:p-Ryleigh_quotient_nodal_domain_decomposition_space}
 we know $\rayl_{p_i}(\hat{f}_{p_i})\leq \mathcal{R}_{p_i}(f_{p_i})=\Lambda_{p_i}$. So, for $p_i$ that goes to $\infty$, we obtain
 \begin{equation}
 \mathcal{R}_\infty(f')\leq \mathcal{R}_\infty(f)=\Lambda.
 \end{equation}
 On the other hand, since from \eqref{eq1:thm_limisup_nodal_count} $|Kf_{p_1}(u,v)|\leq |K\hat{f}_{p_1}(u,v)|$, passing to the limit we have $|Kf(u,v)|\leq |Kf'(u,v)|$ for any edge, i.e. $\|Kf\|_{\infty}\leq \|Kf'\|_{\infty}$. In addition, $\|f\|_{\infty}=\|f'\|_{\infty}$, i.e. 
 \begin{equation}
 \mathcal{R}_\infty(f')\geq\mathcal{R}_\infty(f)=\Lambda,
 \end{equation}
yielding the equality. As a direct consequence, if $(u,v)$ is such that $|Kf(u,v)|=\|Kf\|_{\infty}$, the same holds for $Kf'(u,v)$. So
\begin{equation}\label{eq2:thm_limisup_nodal_count}
\begin{aligned}
    \Big\{u\in \internalnodes\,\Big|\; |f(u)|=\|f\|_{\infty}\Big\}&\subseteq \Big\{u\in \internalnodes\,\Big|\; |f'(u)|=\|f'\|_{\infty}\Big\}\\
    \Big\{(u,v)\in \edgeset\,\Big|\; |Kf(u,v)|=\|Kf\|_{\infty}\Big\}&\subseteq \Big\{(u,v)\in \edgeset\,\Big|\; |Kf'(u,v)|=\|Kf'\|_{\infty}\Big\}.
    \end{aligned}
\end{equation}
Now recall that, since $f$ is a viscosity eigenfunction relative to the eigenvalue $\Lambda$, it is also an $\infty$-eigenpair as in \cref{Def:generalized_p-eigenpair}. So there exist $\Xi\in \partial\|Kf\|_{\infty}$ and $\xi\in \partial\|f\|_{\infty}$
that satisfy the eigenvalue equation in \Cref{Def:generalized_p-eigenpair}. As a direct consequence of \eqref{eq2:thm_limisup_nodal_count} and \eqref{inf_subgradient}, it is easy to prove that $\Xi\in \partial\|Kf'\|_{\infty}$ and $\xi\in \partial\|f'\|_{\infty}$, proving that $f'$ is an eigenfunction corresponding to $\Lambda$. Finally, since $f'\in X_\Lambda$, from Theorem \ref{lemma:nodal-domain-limit} it follows that 
\begin{equation}
    \max\limits_{f\in X_{\Lambda}}\PNc(f)\geq \PNc(f')\ge \PNc(\hat{f}_{p_i})=\SNc(\hat{f}_{p_i})=\SNc(f_{p_i})=\limsup\limits_{p_i\to+\infty}\max\limits_{f\in X_{\Lambda_{p_i}}}\SNc(f),
\end{equation} 
i.e. the thesis.
\end{proof}

Recall from \Cref{lemma:p-monotonic} that the normalized $p$-Laplacian variational eigenvalues converge toward the $\infty$-Laplacian variational eigenvalues, i.e. $ \lambda_k(\plap)^{\frac{1}{p}}\rightarrow \Lambda_k$, so we can formulate the following corollary of \cref{Thm:lower_bound_perfect_nodal_count_limit_eigenpairs}.

\begin{corollary}\label{Cor:lower_bound_perfect_nodal_count_variational_limit_eigenpairs}
For any $k=1,\dots,N$
     $$\max\limits_{f\in X_{\Lambda_k}}\PNc(f)\ge \limsup\limits
    _{p\to+\infty}\max\limits_{f_{p}\in X_{\lambda_k(\Delta_p)}(\Delta_p)}\SNc(f_{p}),$$
    where $X_{\lambda_k}(\Delta_p)$ denotes the set of eigenfunctions corresponding to the $k$-th variational eigenvalue $\lambda_k(\Delta_p)$ of the $p$-Laplacian.
\end{corollary}

In particular, we remark that by combining the last \Cref{Thm:lower_bound_perfect_nodal_count_limit_eigenpairs} and \Cref{Cor:lower_bound_perfect_nodal_count_variational_limit_eigenpairs} with the lower bounds in \Cref{thm:nodal-count}, we can give a lower bound for the maximal numbers of perfect nodal domains of eigenfunctions corresponding to $\Lambda$.


\subsubsection{Zeros and decomposition of the eigenfunction support}

In this section, we study the decomposition of the support of an eigenfunction given by its zeros. To this end, given a nonzero function $f$, for each connected component $U$ of the support of $f$, 
\begin{equation}
    \mathrm{supp}(f):=\{u\in \internalnodes\,|\;f(u)\neq0\},    
\end{equation}
we say that $U$ is an \textbf{equal-component} of $f$ if $\mathcal{R}_\infty(f|_U)=\mathcal{R}_\infty(f)$,  a \textbf{super-oscillation component} if  $\mathcal{R}_\infty(f|_U)>\mathcal{R}_\infty(f)$ or a \textbf{sub-oscillation component} if $\mathcal{R}_\infty(f|_U)<\mathcal{R}_\infty(f)$.
Next, let $U_1,\cdots,U_m$ be the connected components of $\mathrm{supp}(f)$, and let $\{1.\cdots,l\}=I_e\cup I_+\cup I_-$ be the unique index decomposition such that $\mathcal{U}_e(f):=\{U_i\big|\,i\in I_e\}$, $\mathcal{U}_+(f):=\{U_i\big|\,i\in I_+\}$, and $\mathcal{U}_-(f):=\{U_i\big|\,i\in I_-\}$ are the collections of all the equal-components, super-oscillation components, and sub-oscillation components of $f$, respectively. 
The next Lemma proves that viscosity $\infty$-eigenfunctions admit only equal components.

  \begin{lemma}\label{lemma:equal-component}
Any $\infty$-eigenfunction has an equal-component. Moreover, if $f$ is a viscosity eigenfunction, any connected component of the support of $f$ is an equal-component. 
  \end{lemma}

\begin{proof}
From the assumptions and the definition of nodal domains, there exists some strong nodal domain $U'\subset A$ such that $v\in U'$ with $|f(v)|=\|f\|_{\infty}$. In particular, it is not difficult to check that $f|_U$ is an eigenfunction  corresponding to  $\Lambda$, which means that $\mathcal{R}_\infty(f|_U)=\mathcal{R}_\infty(f)$.
If $f$ is further assumed to be a viscosity eigenfunction, we can simply let $U$ be any connected component of the support of $f$. In this case, since $f|_U$ solves the limit eigenvalue equation in \Cref{thm:limitin_infty_eigenvalue_theorem} on $U$, from \Cref{theorem:viscosity_eigenpairs_are_generalized_eigenpairs} it is an eigenfunction corresponding to  $\Lambda$ on $U$, which implies that $\mathcal{R}_\infty(f|_U)=\mathcal{R}_\infty(f)$.
\end{proof}

Finally, in the next proposition we show that the number of super-oscillation components and sub-oscillation components of an $\infty$-eigenfunction can be bounded in terms of the position of the corresponding eigenvalue with respect to the variational spectrum and vice versa.

\begin{proposition}\label{pro:connected-components}
 Let  $(f,\Lambda)$ be an $\infty$-eigenpair. 
 \begin{itemize}
     \item If $\Lambda<\Lambda_k(\inflap)$, we have $|\mathcal{U}_e(f)|+|\mathcal{U}_-(f)|\le k-1$; in addition if $|\mathcal{U}_-(f)|\ge k$, then $\Lambda>\Lambda_k(\inflap)$. 
     \item If $\Lambda>\Lambda_k(\inflap)$, we have $|\mathcal{U}_e(f)|+|\mathcal{U}_+(f)|\le N-k$; in addition if $|\mathcal{U}_+(f)|\ge N-k+1$, then $\Lambda<\Lambda_k(\inflap)$. 
 \end{itemize}   
\end{proposition}

\begin{proof}
Let $U_1,\cdots,U_m$ be the connected components of the support of $f$. 
For any $(t_1,\cdots,t_m)\ne 0$, 
\begin{align}
\mathcal{R}_\infty(\sum_{i=1}^{m} t_i f|_{U_i})&= \frac{\max\limits_{i=1,\ldots,m}|t_i|\cdot\max\limits_{\{v,v'\}\in E}|f|_{U_i}(v)-f|_{U_i}(v')|}{\max\limits_{i=1,\ldots,m}|t_i|\cdot\|f|_{U_i}\|_\infty}
\\&\le \max\limits_{i\text{ with }t_i\ne0}\frac{\max\limits_{\{v,v'\}\in E}|f|_{U_i}(v)-f|_{U_i}(v')|}{\|f|_{U_i}\|_\infty}=\max\limits_{i\text{ with }t_i\ne0}\mathcal{R}_\infty(f|_{U_i}),
\end{align}
and similarly,
\begin{equation}
\mathcal{R}_\infty(\sum_{i=1}^{m} t_i f|_{U_i})\ge \min\limits_{i\text{ with }t_i\ne0}\mathcal{R}_\infty(f|_{U_i}).
\end{equation}
Let $\{1.\cdots,m\}=I_e\cup I_+\cup I_-$ be the  index decomposition, and let $\{U_1,\cdots,U_m\}=\mathcal{U}_e(f)\cup \mathcal{U}_+(f)\cup \mathcal{U}_-(f)$ be the corresponding set decomposition. By Proposition \ref{lemma:equal-component}, $\mathcal{U}_e(f)\ne\varnothing$. 
So consider the space of functions $\mathcal{X}_-=\mathrm{span}\{f|_{U}:U\in \mathcal{U}_e(f)\cup \mathcal{U}_-(f)\}$ and $\mathcal{X}_+=\mathrm{span}\{f|_{U}:U\in \mathcal{U}_e(f)\cup \mathcal{U}_+(f)\}$, and observe that for any $g\in \mathcal{X}_-\setminus\{0\}$,
\begin{equation}
\mathcal{R}_\infty(g)\le \max_{U\in \mathcal{U}_e(f)\,\cup\, \mathcal{U}_-(f)}\mathcal{R}_\infty(f|_U)= \mathcal{R}_\infty(f)=\Lambda.
\end{equation}
We can now prove the first thesis in the statement of the proposition. Assume by contradiction that $|\mathcal{U}_e(f)|+|\mathcal{U}_-(f)|\ge k$. From 
 $\gamma(\mathbb{S}^{N-1}\cap \mathcal{X}_-)=\dim \mathcal{X}_-=|\mathcal{U}_e(f)|+|\mathcal{U}_-(f)|\ge k$, we have that 
 \begin{equation}
 \Lambda_{k}(\inflap)\le
\sup\limits_{g\in \mathbb{S}^{N-1}\cap \mathcal{X}_-}\mathcal{R}_\infty(g)\le\Lambda <\Lambda_{k}(\inflap)
 \end{equation}
which yields the contradiction. 

Similarly, if $|\mathcal{U}_-(f)|\ge k$, then we obtain
\begin{equation}
    \Lambda_{k}(\inflap)\le
\sup\limits_{g\in \mathbb{S}^{N-1}\cap \mathcal{X}_{<0}}\mathcal{R}_\infty(g)\le \max_{U\in  \mathcal{U}_-(f)}\mathcal{R}_\infty(f|_U)< \mathcal{R}_\infty(f)=\Lambda  
\end{equation}
where  $\mathcal{X}_{<0}=\mathrm{span}\{f|_{U}:U\in  \mathcal{U}_-(f)\}$. We have proved also the first statement. 

Next we prove the second statement. Let $g\in \mathcal{X}_+\setminus\{0\}$, $\mathcal{R}_\infty(g)\ge\Lambda$ and suppose by contradiction that $\dim \mathcal{X}_+=|\mathcal{U}_e(f)|+|\mathcal{U}_+(f)|\ge N-k+1$. 
By the intersection property of  Krasnoselskii's $\mathbb{Z}_2$-genus \Cref{Lemma_Krasnoselskii_intersection}, for any centrally symmetric subset $S\subset\R^{N}$ with $\gamma(S)\ge k$,  $S\cap (\mathcal{X}_+\setminus0)\ne\varnothing$. 
Therefore,
\begin{equation}
    \Lambda_{k}(\inflap)=
\inf\limits_{ \gamma(S)\ge k}\sup\limits_{g\in S}\mathcal{R}_\infty(g) \ge \inf\limits_{ \gamma(S)\ge k}\inf\limits_{g\in S\cap (\mathcal{X}_+\setminus0)}\mathcal{R}_\infty(g) 
\ge
\inf\limits_{g\in \mathcal{X}_+\setminus0}\mathcal{R}_\infty(g)\ge\Lambda
\end{equation}
which is a contradiction. 
Analogously, if $|\mathcal{U}_+(f)|\le N-k+1$, then we obtain
\begin{equation}
\Lambda_{k}(\inflap)=\inf\limits_{ \gamma(S)\ge k}\sup\limits_{g\in S}\mathcal{R}_\infty(g)\ge \inf\limits_{g\in \mathcal{X}_{>0}\setminus0}\mathcal{R}_\infty(g)\ge \min_{U\in  \mathcal{U}_+(f)}\mathcal{R}_\infty(f|_U)>\Lambda
\end{equation}
where  $\mathcal{X}_{>0}=\mathrm{span}\{f|_{U}:U\in  \mathcal{U}_+(f)\}$. 
\end{proof}

We point out that as long as $\Lambda_{k+1}(\inflap)>\Lambda_k(\inflap)>\Lambda_{k-1}(\inflap)$, since $|\mathcal{U}_e(f)|\geq 1$, the additional statements in the thesis of \Cref{pro:connected-components} are direct consequences of points $1$ and $2$. However, if $\Lambda_k$ has higher multiplicity, then these become of independent interest. 
As a corollary of \Cref{lemma:equal-component} and \Cref{pro:connected-components} we have the following result.

\begin{corollary}\label{cor:connected-components}If $f$ is a viscosity $\infty$-eigenfunction corresponding to $\Lambda$, with $\Lambda_h(\inflap)<\Lambda<\Lambda_k(\inflap)$ then the number of connected components of the support of $f$ is less than or equal to $ \min\{k-1,N-h\}$.
\end{corollary}


\section{Graph $1$-Laplacian spectrum}\label{Sec:1-Lapl_spectrum}

In this section, we discuss the properties of the $1$-Laplacian spectrum and its relevance to graph clustering. In particular, we discuss tight relations between the $1$-Laplacian eigenpairs and the Cheeger constants of the graph. We start by recalling the geometric characterization of the $1$-Laplacian eigenvalues from \Cref{lemma:1-eigenvalues_are_isoperimetric_constants}.

To this end, given a subset of nodes $A\subset \internalnodes$, we recall the definition of the normalized cut of $A$ in \cref{cut}:
\begin{equation}
c(A)=\frac{\|\edgelength\big(E(A,A^c)\big)\|_1}{\nodeweight(A)}=\frac{\sum_{(u,v)\in E(A,A^c)}\edgelength_{uv}}{2\sum_{u\in A}\nodeweight_u}\,,
\end{equation}
where $\nodeweight(A)=\sum_{u\in A}\nodeweight_u$ is the measure of the set $A$, $E(A,A^c)$ denotes the set of edges connecting $A$ to its complement $A^c$, and $\edgelength\big(E(A,A^c)\big)$ denotes the sum of the inverse lengths of the edges in $E(A, A^c)$.
Recall from \Cref{sec:nodal_domains} and \Cref{lemma:1-eigenvalues_are_isoperimetric_constants} that any nodal domain $A$ induced by an eigenfunction $f$ satisfies $c(A)=\Lambda$.
Thus, any $1$-Laplacian eigenvalue corresponds exactly to the normalized cut of some subgraph of $\Gc$. As a direct consequence, the spectrum is always finite, and we have the following result about eigenvalue multiplicities.

\begin{proposition}\label{Prop_general_simplicity_of_1eigenvalues}
For almost all edge lengths $\edgelength$ and node weights $\nodeweight$, any $1$-Laplacian eigenvalue is simple and the corresponding eigenfunction induces a unique strong nodal domain.
\end{proposition}
\begin{proof}
Note that the equation $c(A)=c(A')$ for two different subsets $A,A'\subset \internalnodes$ corresponds to an equation in the node weights and edge lengths. In particular, a solution of $c(A)=c(A')$ lies on a codimension-$1$ hypersurface in $\R^{|\internalnodes|+|\edgeset|}$. There are at most $2^N-1$ distinct values of $c(A)$, and thus at most $(2^N-1)(2^N-2)/2$ different equations of the form $c(A)=c(A')$.
It follows that the set of weights for which $c(A)=c(A')$ holds for some $A\ne A'$ lies in a union of finitely many codimension-$1$ hypersurfaces, which has Lebesgue measure zero. Thus, for almost all lengths and weights, $c(A)\ne c(A')$ whenever $A\ne A'$.

Now consider a $1$-Laplacian eigenpair $(f,\Lambda)$ and let $A:=\{f\ge t\geq 0 \}$ (or $A:=\{f\le -t\leq 0\}$) be a super-level set (or a sub-level set) such that $A\ne\varnothing$. We observe that $(\chi_{A},\Lambda)$ is also an eigenpair, where $\chi_{A}$ denotes the characteristic function of the set $A$, taking value $1$ on $A$ and $0$ elsewhere. Indeed, we know that there exist $\Xi\in \partial\|\grad f\|_1$ and $\xi\in\partial\|f\|_1$ such that $K^T \Xi=\Lambda \nodeweight\odot \xi$, and it is easy to check that the subgradients $\Xi$ and $\xi$ also satisfy $\Xi\in \partial\|\grad \chi_A\|_1$ and $\xi\in\partial\|\chi_A\|_1$. In particular, any superlevel or sublevel set of $f$ satisfies $\Lambda=c(A)$. Since we have proved that for almost all edge lengths and node weights $c(A')\ne c(A)$ whenever $A'\ne A$, the family $\big\{\{f\ge t\},\{f\le -t\}:t\geq 0\big\}\setminus\{\varnothing\}$ must consist of a single element. This implies that $f=t\chi_A$ for some $t\ne 0$ and $A\subset \internalnodes$, proving that $\Lambda$ is simple.
\end{proof}

\begin{remark}
Additionally, we recall from \cite{zhang2021homological} that the $1$-Laplacian variational eigenvalues admit an equivalent definition that makes their geometric content more explicit. Indeed, consider $\mathcal{P}_2(V)$, the family of pairs of disjoint and nonempty subsets of nodes,
$$\mathcal{P}_2(V):=\big\{(A,B)|\;A,B\subset \internalnodes,\;A\cap B=\varnothing, \; A\cup B\neq \varnothing\big\}.$$
Clearly, $\mathcal{P}_2(V)$ is a partially ordered set with respect to inclusion, i.e. $(A,B) \preceq (A’,B’)$ iff $A\subseteq A'$ and $B\subseteq B'$. We say that $\mathcal{C}\subset\mathcal{P}_2(V)$ is an inclusion chain if for any $(A,B),(A’,B’)\in\mathcal{C}$, either $(A,B) \preceq (A’,B’)$ or $(A',B') \preceq (A,B)$. Then, for each subfamily $\mathcal{A}\subseteq\mathcal{P}_2(V)$, we define the following set of functions induced by the inclusion chains contained in $\mathcal{A}$:
$$\mathbf{K}(\mathcal{A}):=\bigcup_{\substack{\mathcal{C}\subset\mathcal{A}\\ \text{ inclusion chain }}}\mathrm{conv}\big\{\chi_A-\chi_B:(A,B)\in \mathcal{C}\big\},$$
where $\chi_A$ and $\chi_B$ denote the characteristic functions of the sets $A$ and $B$, respectively. Finally, following \cite{zhang2021homological}, by considering the subfamilies of $\mathcal{P}_2(\nodeset)$ that induce sets of functions with genus at least $k$,
$\mathcal{S}_k:=\{\mathcal{A}\subset \mathcal{P}_2(V):\gamma(\mathbf{K}(\mathcal{A}))\ge k\}$,
we can equivalently define the $k$-th variational eigenvalue of the graph $1$-Laplacian as
\begin{equation}\label{eq:Lambda_k}
\Lambda_k(\Delta_1)=\min\limits_{\mathcal{A}\in \mathcal{S}_k}\max\limits_{(A,B)\in\mathcal{A}}\frac{w(E(A,A^c))+w(E(B,B^c))}{2\nu(A\cup B)} 
\end{equation}

Thus, in principle, all variational eigenvalues can be obtained by examining families of set-pairs of $V$. In particular, $\Lambda_k(\Delta_1)$ can be computed exactly in a finite number of steps; however, this observation appears to be of purely theoretical significance.
\end{remark}

\subsection{Relations with Cheeger constants}
As anticipated, the study of the $1$-Laplacian spectrum is closely related to the study of different families of higher-order isoperimetric constants (or Cheeger constants) of the graph. Indeed, the bounds that relate the eigenvalues of the linear Laplacian operator and the isoperimetric constants of the domain~\cite{Cheeger} become sharper when we consider the nonlinear $p$-Laplacian with $p<2$.

We recall that Cheeger constants are typically used in data analysis applications, as they provide information about the number and the quality of the clusters of a graph. Next, we recall their definition.
Given an integer $k$, we consider all possible families of $k$ nonempty and disjoint subsets of $\internalnodes$:
\begin{equation}
\mathcal{D}_k(\Gc)=\big\{A_1,\dots,A_k\subset\internalnodes\,|\;A_i\neq\emptyset,\;A_i\cap A_j=\emptyset\;\:\forall\,i,j\big\}\,.
\end{equation}
\begin{definition}[Cheeger constants]\label{DEf:Cheeger_constant}
The $k$-th Cheeger constant (or isoperimetric constant), see~\cite{lee2014multiway, lawler1988bounds, daneshgar2010isoperimetric}, is defined as
$$
h_k(\Gc):=\min_{\{A_1,\dots, A_k\}\in\mathcal{D}_k(\Gc)}\;\max_{i=1,\dots,k} c(A_i).
$$
\end{definition}

Observe that a ``small'' value of $h_k(\Gc)$ means that there exist $k$ subsets of nodes that are simultaneously fairly massive (i.e., $\nodeweight(A_i)$ is sufficiently large for all $i$) and weakly connected to the rest of the graph. Indeed, having small $\edgelength\big(E(A_i,A_i^c)\big)$ for each $i$ heuristically indicates that the subsets $\{A_i\}$ are far from each other. This is precisely the structure expected of $k$ clusters of nodes.
The constant $h_k(\Gc)$ can thus be regarded as an indicator of how well the graph can be clustered into $k$ subgraphs, with the corresponding family of subsets providing approximate clusters.

It is a well known fact (see, e.g.,~\cite{chung1997spectral, trevisan2017lecture}) that the second eigenvalue $\lambda_2(\Delta_2)$ of the linear normalized graph Laplacian satisfies
\begin{equation}
\frac{h_2^2(\Gc)}{2}\leq \lambda_2(\Delta_2)\leq 2h_2(\Gc).
\end{equation}
Further results relate higher eigenvalues to higher Cheeger constants~\cite{trevisan2017lecture,lee2014multiway,daneshgar2012nodal}. Interestingly, these relationships become sharper when we consider the $p$-Laplacian and let $p$ approach $1$. Indeed, given a subset $A\subset\internalnodes$ and its characteristic function $\chi_A$, we have
\begin{equation}\label{eq-cheeger_constant_and_characteristic_functn}
\rayl_1(\chi_A)=\frac{\sum_{(u,v)\in\edgeset}\edgelength_{uv}|\chi_A(u)-\chi_A(v)|}{2\sum_{v\in\internalnodes}\nodeweight_u|\chi_A(u)|}=c(A)\,.
\end{equation}
Moreover, we recall from \Cref{subsection_1-Lap_nodal_domains} that if $A$ is a strong nodal domain of an eigenfunction $f$, then $c(A)=\Lambda$, where $\Lambda$ is the corresponding eigenvalue.
These remarks suggest direct connections between Cheeger constants and $1$-Laplacian eigenpairs, which we formalize in the remainder of this section.
In particular, we consider the general case of the $1$-Laplacian eigenproblem with Dirichlet boundary conditions, which has been partially investigated in~\cite{Hua}. Moreover, we refer to~\cite{Kawohl2003, parini2010second, BobkovParini18} for an investigation of the same problem on bounded open domains in $\R^n$. Finally, we refer to~\cite{chang2016spectrum, Bhuler, ZhangNodalDO} for the $1$-Laplacian spectrum on graphs without boundary conditions; it is worth mentioning that many results from the non-boundary case extend to the boundary case with minimal effort. It is also interesting to recall that the Cheeger cuts of a sequence of graphs approximate the Cheeger cut of an open bounded domain in an appropriate sense~\cite{GarciaSlepcev16}.

The next theorem establishes connections between higher Cheeger constants and the $1$-Laplacian spectrum~\cite{chang2016spectrum, Bhuler}. Since an explicit proof of this theorem in the boundary case does not appear in the literature, for completeness we provide a proof in the appendix.
\begin{theorem}\label{thm:1-eigenpair_and_cheeger_constants}
Let $(f_k,\Lambda_k)$ be the $1$-Laplacian variational eigenpairs and let $h_k(\Gc)$ be the Cheeger constants of the graph. Then
$$h_{\SNc(f_k)}(\Gc)\leq\Lambda_k(\Delta_1)\leq h_k(\Gc).$$
Moreover, when $k=1,2$ it holds the equality $h_k(\Gc)=\Lambda_k(\Delta_1)$.
\end{theorem}
%
Below we report a technical remark that is useful in proving several results about the $1$-Laplacian spectrum.
\begin{remark}\label{remark_cheeger_constants_and_support}
Given a function $f$, we denote its support by $\supp(f):=\{u\in\internalnodes \,|\;\text{s.t.}\; f(u)\neq 0\}$. Let $\Gc(\supp(f))$ be the graph induced by $\supp(f)$, obtained by setting $\nodeset\setminus \supp(f)$ as the boundary, and let $\mathcal{H}_0(\Gc(\supp(f)))$ be the set of functions defined on $\Gc(\supp(f))$. From the equality $\Lambda_1(\Delta_1)=h_1(\Gc)$ in \Cref{thm:1-eigenpair_and_cheeger_constants}, we have
\begin{equation}
\rayl_1(f)\geq \min_{g\in \mathcal{H}_0(\Gc(\supp(f)))}\rayl_1(g)= \Lambda_1(\Gc(\supp(f)))=h_1\big(\Gc(\supp(f)\big).
\end{equation}
In particular, if $f_m$ is a minimizer of $\rayl_1$ in $\mathcal{H}_0(\Gc(\supp(f)))$, then $f_m$ induces one nodal domain $A=\supp(f_m)$ such that
\begin{equation}
c(A)=\Lambda_1(\Gc(\supp(f)))=\rayl_1(1_{A})=\rayl_1(f_m)\leq \rayl_1(f).
\end{equation}
In particular, there exists some $A\subset \supp(f)$ such that $c(A)\leq \rayl_1(f)$. The same conclusion was proved in~\cite{Tudisco1} using a different argument.
\end{remark}
The lower bound in \Cref{thm:1-eigenpair_and_cheeger_constants} can be strengthened by considering sub-eigenfunctions. In particular, we say that a function $f$ is a sub-eigenfunction corresponding to the $1$-eigenvalue $\Lambda$ if $\rayl_1(f|_A)\le \Lambda$ for any nodal domain $A$.
\begin{proposition}\label{prop:sub-eigen}
Suppose $f$ is a sub-eigenfunction corresponding to $\Lambda_k(\Delta_1)$. Then
$$h_{\mathcal{SN}(f)}(\Gc)\le \Lambda_k(\Delta_1). $$
\end{proposition}
\begin{proof}
Suppose $f$ is a sub-eigenfunction corresponding to $\Lambda_k(\Delta_1)$. Then, for any nodal domain $A_i$ of $f$, $\rayl_1(f|_{A_i})\le \Lambda_k(\Delta_1)$.
Moreover, by \Cref{remark_cheeger_constants_and_support}, there exists a nonempty subset $\tilde{A}_i\subset A_i$ such that $\rayl_1(\chi_{\tilde{A}_i})\le \rayl_1(f|_{A_i})$. Thus,
\begin{equation}
    h_{\mathcal{SN}(f)}(\Gc)\le \max_i c(\tilde{A}_i)=\max_i \rayl_1(\chi_{\tilde{A}_i})\le \max_i \rayl_1(f|_{A_i})\le \Lambda_k(\Delta_1).    
\end{equation}
\end{proof}

The last result has been generalized to $p$-Laplacian eigenpairs with $p>1$ in~\cite{daneshgar2012nodal,Tudisco1, KELLER201680}, where the authors relate $p$-eigenpairs to Cheeger constants for any $p>1$.%

Below we report an adapted version of Theorem 5.1 in~\cite{Tudisco1}, where we take into account the power $p$ on the edge lengths. For a discussion of these bounds also in the context of infinite graphs, we refer to~\cite{KELLER201680}.
\begin{theorem}\label{ThM:p-eigenpairs_and_cheeger_constants}
Let $(f_k(\plap),\lambda_k(\plap))$ be the $k$-th variational eigenpair of the $p$-Laplacian, $p>1$. Then
$$\min_{(u,v)\in\edgeset}\edgelength_{uv}^{p-1}\frac{2^{p-1}}{\tau(\Gc)^{p-1}}\frac{\big(h_{\SNc(f)}(\Gc)\big)^p}{p^{p}}\leq \lambda_k(\plap)\leq \max_{(u,v)\in\edgeset}\edgelength_{uv}^{p-1} 2^{p-1}h_k(\Gc)\,,$$
where $\tau(\Gc)=\max_{u\in\internalnodes}\sum_{v\sim u}\edgelength_{uv}/\nodeweight_u$.
\end{theorem}
Looking at the last two theorems, we note that whenever we have a variational eigenfunction whose nodal domain count reflects the corresponding frequency, letting $p$ approach $1$ yields an eigenvalue that reproduces exactly a higher-order Cheeger constant.
In addition, we observe the following corollary of the last theorem.
\begin{corollary}\label{corollary_1-lap_cheeger_lower_bound_limsup}
    Let $\Nc(k)=\limsup_{p\rightarrow 1} h_{\SNc(f_k(\plap))}(\Gc)$. Then
    $$h_{\Nc(k)}(\Gc)\leq\Lambda_k(\Delta_1)\leq h_k(\Gc) \qquad \forall k=1,\dots, N.$$
\end{corollary}
\begin{proof}
The proof is an immediate consequence of the continuity of the variational eigenvalues (\Cref{lemma:p-monotonic}) and \Cref{ThM:p-eigenpairs_and_cheeger_constants}.
\end{proof}

\begin{remark}\label{remark_variational_eigenvalues_on_trees_cheeger_constants}
    We recall that on trees the variational spectrum of the $p$-Laplacian (for $p\in(1,\infty)$\,) exhausts the spectrum, see \cref{variational_eigenvalues_on_trees}, and whenever an eigenfunction $f$ associated to the $k$-th variational eigenvalue is everywhere different from zero, it induces exactly $k$ nodal domains, see \cref{thm_nodal_count_on-trees}. Considering the discussion in \cite{DEIDDA2023_nod_dom}, it is possible to observe that the condition of an eigenfunction being everywhere nonzero is generic in the weights $\{\edgelength_{uv}\}_{(u,v)\in \edgeset}$. Thus, by continuity of the variational spectrum with respect to the weights and of the Cheeger constants, it is possible to observe that on trees the variational eigenvalues are exactly equal to the Cheeger constants:
    \begin{equation}
        \text{if $\Gc$ is a tree, then  } \Lambda_k(\Delta_1)=h_k(\Gc) \quad \forall k=1,\dots, N
    \end{equation}
    We refer to a joint paper currently in progress~\cite{ge2025} for a formal proof and to~\cite{MengZhang25} for a combinatorial proof.
\end{remark}

Note that, whenever $\SNc(f_k(\Delta_1))\leq \limsup_{p\rightarrow 1} \SNc(f_k(\plap))$, the lower bound provided in \Cref{corollary_1-lap_cheeger_lower_bound_limsup} is stronger than the lower bound provided by \Cref{thm:1-eigenpair_and_cheeger_constants}. By the third item in Lemma~\ref{lemma:nodal-domain-limit}, we do not know whether $\SNc(f_k(\Delta_1))\leq \limsup_{p\rightarrow 1} \SNc(f_k(\plap))$ holds. However, for $f_k(\plap)$ with $p>1$, we have lower bounds and, on some special graphs such as trees, we know the nodal count exactly. In contrast, $f_k(\Delta_1)$ generically has a unique nodal domain (\cref{Prop_general_simplicity_of_1eigenvalues}). Hence, in a generic situation, we expect the lower bound in \Cref{thm:1-eigenpair_and_cheeger_constants} to be trivial, while the lower bound in \Cref{corollary_1-lap_cheeger_lower_bound_limsup} might not be.

We conclude this section with illustrative examples in which we discuss, in practice, several properties of the $p$-Laplacian spectrum that we have discussed theoretically so far.

\begin{example}\label{exam:path}
Let $P_7$ denote the path graph on seven vertices. Since any $1$-Laplacian eigenvalue is equal to the normalized cut of some subgraph \cref{lemma:1-eigenvalues_are_isoperimetric_constants}, we can state that the $1$-Laplacian eigenvalues, listed in increasing order, are: 
\begin{equation}
0,\;\frac13,\;\frac12,\;\frac23,\;1,\;2,    
\end{equation}
where we have $\gamma\text{-}\mathrm{mult}(2)=3$, $\gamma\text{-}\mathrm{mult}(1)=2$, and
\[\gamma\text{-}\mathrm{mult}\big(0\big)=\gamma\text{-}\mathrm{mult}\left(\frac13\right)=\gamma\text{-}\mathrm{mult}\left(\frac12\right)=\gamma\text{-}\mathrm{mult}\left(\frac23\right)=1.\]
To see it, we recall from \cref{lemma:1-eigenvalues_are_isoperimetric_constants} that if $f$ is a $1$-Laplacian eigenfunction associated with the eigenvalue $\Lambda$ and $A$ is a nodal domain induced by $f$, then $\Lambda=\sum_{u\in A, v\in A^c}\edgelength_{uv}/\nodeweight(A)$. Thus, considering, e.g. $\lambda=1/3$, we can consider the characteristic functions
\begin{equation}
    1_{\{1,2,3\}} \quad \text{and} \quad 1_{\{5,6,7\}},
\end{equation}
where $1_{\{1,2,3\}}$ indicates the function taking value $1$ on $\{1,2,3\}$ and 0 otherwise and analogously $1_{\{5,6,7\}}$.
Below we prove that they are both eigenfunctions.
Indeed, consider the two functions
\begin{equation}\label{eq:indicator_gradient}
\begin{aligned}
    \xi(v)&=\begin{cases}
    1 \quad \text{if } v=1,2,3\\
    -1 \quad \text{if } v=4\\
    0 \quad \text{otherwise}
\end{cases}\\
\Xi(u,v)=-\Xi(v,u)&=\begin{cases}
    -1 \quad \text{if } (u,v)=(3,4)\\
    -2/3 \quad \text{if } (u,v)=(2,3),(4,5),(5,6)\\
    -1/3 \quad \text{if } (u,v)=(1,2)\\
    0 \quad \text{otherwise}.
\end{cases}
\end{aligned}
\end{equation}
It is not difficult to show that $\xi\in \partial\|1_{\{1,2,3\}}\|_1$,  $\Xi\in\partial\|\grad 1_{\{1,2,3\}}\|_1$, and they satisfy the generalized eigenvalue equation in \Cref{Def:generalized_p-eigenpair}. Thus $1_{\{1,2,3\}}$ is an eigenfunction corresponding to the eigenvalue $\lambda=1/3$ and, by symmetry, $1_{\{5,6,7\}}$ has the same property. Similarly we can verify that any nonzero combination $c_1 1_{\{1,2,3\}}-c_2 1_{\{4,5,6\}}$ with both $c_1,c_2\geq 0$ is an eigenfunction corresponding to $1/3$. In particular, we can prove that the set
\begin{equation}
X_{1/3}:=\pm\{c_1 1_{1,2,3}-c_2 1_{5,6,7} \text{ with both } c_1,c_2\geq 0\}\setminus\{0\}
\end{equation}
exhaust the entire set of eigenfunctions corresponding to $1/3$. Indeed, by the nodal domain property from \cref{lemma:1-eigenvalues_are_isoperimetric_constants} recalled above, any eigenfunction corresponding to the eigenvalue $\Lambda=1/3$ does not admit any nodal domain other than $\{1,2,3\}$ and $\{5,6,7\}$. Moreover, if an eigenfunction $f$ is positive on $\{1,2,3\}$, it necessarily takes the value $0$ in $\{4\}$, and so the subgradient of $\|Kf\|_1$ restricted to $\{1,2,3,4\}$ has to be necessarily as in \eqref{eq:indicator_gradient}, implying that $f$ is constant on $\{1,2,3\}$. Finally, we can observe that if we take $f$ that is positive both on $\{1,2,3\}$ and $\{4,5,6\}$ then the subgradient eigenvalue equation in $\{4\}$ takes the form 
\begin{equation}
   -2=\Xi(3,4)+\Xi(5,4)=\lambda \xi(4)\in\frac{1}{3}[-1,1]
\end{equation}
which has no solution. In particular, we observe that $X_{1/3}\cap S_2$ is the disconnected union of two arcs of circle and thus $\gamma(X_{1/3}\cap S_2)=1$.

A similar discussion shows that 
\begin{equation}
    X_{1/2}=\{\pm c_1 1_{1,2}\mp c_2 1_{6,7} \text{ with both } c_1, c_2\geq 0\}\setminus\{\mathbf{0}\}, \qquad X_{2/3}=\{c 1_{3,4,5} \;|\, c\in \R\setminus\{0\} \}.
\end{equation}
Moreover, for $X_{1}$ we have:
\begin{equation}
\begin{aligned}
    X_{1}=&\{c_1 1_1+c_2 1_7  \;|\, c_1,c_2\in \R\}\;\cup \;\pm\{c_1 1_{2,3}-c_2 1_{4,5} \text{ with } c_1, c_2\geq 0\}\;\cup\\
    &\pm\{c_1 1_{2,3}-c_2 1_{5,6} \text{ with  } c_1, c_2\geq 0\}\,\cup\,\pm\{c_1 1_{3,4}-c_2 1_{5,6} \text{ with } c_1, c_2\geq 0\},
\end{aligned}
\end{equation}
thus $X_{1}\cap S_2$ is the union of a circle and some disconnected arcs of circle and it has genus $2$. 
Finally, for $X_{2}$ we have 
\begin{equation}
    X_2=\{c_21_2+c_31_3+c_41_4+c_51_5+c_61_6 \;\text{ where }\: c_ic_{i+1}\leq 0\}\setminus\{\mathbf{0}\}.
\end{equation}
It is possible to see that $\gamma(X_2)= 3$ by inspecting its structure. Indeed, the $3$-dimensional linear subspace of $X_2$ spanned by $\{\chi_{\{2\}},\chi_{\{4\}},\chi_{\{6\}}\}$, where $\chi_{x}$ denotes the characteristic function of the node $x$, is included in $X_2\cup\{0\}$ meaning that $\gamma(X_2)\ge 3$. 
In addition, an odd continuous function from $X_2$ to $S_2$ that yields $\genus(X_2)\leq 3$ is given by:
\begin{equation}
f(x)=\frac{c_2-c_3,c_4-c_5,c_6}{\|(c_2-c_3,c_4-c_5,c_6)\|_2} \qquad \forall x\in X_2.
\end{equation}
Thus, we have shown that $\gamma(X_2)= 3$.
On the other hand, as it has been discussed in \cite{ MengZhang25} and \cref{remark_variational_eigenvalues_on_trees_cheeger_constants}, the variational eigenvalues of the $1$-Laplacian on trees are equal to the higher order Cheeger constants $\Lambda_k(\Delta_1)=h_k(\Gc)$.

From this characterization, we can deduce the following list of Krasnoselkii variational eigenvalues counted with their variational multiplicity:
\begin{equation}
\begin{aligned}
\Lambda_1(\Delta_1)=&\, 0,\quad \Lambda_2(\Delta_1)=\frac13,\quad \Lambda_3(\Delta_1)=\frac23,\quad \Lambda_4(\Delta_1)=1,\\[.2em]
&\Lambda_5(\Delta_1)=2,\quad \Lambda_6(\Delta_1)=2,\quad \Lambda_7(\Delta_1)=2. 
\end{aligned}
\end{equation}
From this example we can practically observe some of the critical situations that we have discussed theoretically in the previous sections. Precisely:
\begin{itemize}
\item {\bf The $\gamma$-multiplicity and the variational multiplicity can be different.}

Indeed the variational multiplicity of the eigenvalue $1$ is 1, but its $\gamma$-multiplicity is $2$, thus generally $\gamma\text{-}\mathrm{mult}(\lambda)\geq \mathrm{mult}(\lambda)$ for a variational eigenvalue $\lambda$ (see \Cref{lemma_genus_and_multiplicity})
\item
{\bf
The variational spectrum and the homological spectrum generally  do not exhaust the spectrum.}

Indeed $1/2$ is neither a variational eigenvalue nor a homological one. To prove it, recall that there exist finitely many $1$-eigenvalues, so if $1/2$ were a homological eigenvalue, then it must be isolated, and we can apply Theorem \ref{thm:homological-eigen} to derive that for $p>1$ sufficiently close to $1$, $\Delta_p$ has a homological eigenvalue that is close to $1/2$. However, all the eigenvalues of $\Delta_p$ on a tree are variational, see \cref{variational_eigenvalues_on_trees} and their limit for $p$ converging to $1$ is given by the $1$-Laplacian variational eigenvalues (\cref{lemma:p-monotonic}), which do not include $1/2$. This proves that $1/2$ is a non-variational and non-homological eigenvalue.

\item {\bf The number of orthogonal eigenvectors associated with an eigenvalue can be larger than its $\gamma$-multiplicity, \Cref{prop.5.3_struwe}.}

We can observe this by looking at $X_{1/3}$.
\item {\bf The eigenvalues counted with their $\gamma$-multiplicities may exceed the dimension of the space.}

In this case we have $\sum \gamma\text{-}\mathrm{mult}(\lambda_i)=9>7=N$.
\end{itemize}
\end{example}

Next, we discuss several examples concerning the nodal count of $1$-Laplacian eigenfunctions. In all these examples, we strongly rely on the geometrical properties of the $1$-Laplacian eigenpairs discussed in \cref{thm:1-eigenpair_and_cheeger_constants}.

\begin{example}\label{ex_1_nodal_domain}
Consider the connected graph illustrated in \Cref{fig:example_nodal_}, with $\nodeweight_u=\deg(u)$ and $\edgeweight_{uv}=1$ for all nodes and edges. 
As a consequence of \Cref{thm:1-eigenpair_and_cheeger_constants}, one can verify that the second eigenvalue satisfies $\Lambda_2(\Delta_1)=1/3$. Moreover, by \Cref{lemma:1-eigenvalues_are_isoperimetric_constants} and \Cref{Thm_nodal_domains_decomposition_1-Lap_eigenfunctions}, any eigenfunction $f$ corresponding to the eigenvalue $1/3$ can be shown to have the form $f=c\,1_{\{4,5\}}$ for some $c\neq 0$. 
This implies that any eigenfunction corresponding to $\Lambda_2(\Delta_1)$ has exactly one strong nodal domain. However, it is known that any eigenfunction corresponding to $\lambda_2(\Delta_p)$ for $p>1$ has two strong nodal domains.
\begin{figure}[t]
    \centering
        \begin{tikzpicture}[scale=1]
    
    \node (1) at  (-1,1) [circle,draw] {$1$};
    \node (2) at  (-1,-1) [circle,draw] {$2$};
    \node (3) at  (0,0) [circle,draw] {$3$};
    \node (4) at  (1,0) [circle,draw] {$4$};
    \node (5) at  (2,0) [circle,draw] {$5$};
    
    \draw (3)--(1)--(2)--(3)--(4)--(5);
    
    \end{tikzpicture}
    \caption{Graph discussed in \cref{ex_1_nodal_domain} where any $1$-Laplacian eigenfunction associated to the second variational eigenvalue $\Lambda_2(\Delta_1)$ induces a unique strong nodal domain.}
    \label{fig:example_nodal_}
\end{figure}
\end{example}

The next examples illustrate the discussion in \Cref{thm:second-eigen-1}.

\begin{example}\label{ex_3_nodal_domains}
Consider the unweighted graph with $12$ vertices shown in \Cref{fig:4-weak-nodal-2nd_eigenfunction}, and let $f$ be the function taking value $1$ on the vertices labeled `$+$’ and value $-1$ on the vertices labeled `$-$’. It can be proved that $f$ is a $1$-Laplacian eigenfunction corresponding to the second variational eigenvalue, whose variational multiplicity equals $2$. 
Indeed, by studying the Cheeger constants of the graph in \Cref{fig:4-weak-nodal-2nd_eigenfunction}, we deduce from \Cref{thm:1-eigenpair_and_cheeger_constants} that $\Lambda_2(\Delta_1)=\Lambda_3(\Delta_1)=1/3$, which yields $\mathrm{mult}(1/3)\ge 2$. Moreover, by inspecting all characteristic functions $\chi_A$, we know from \Cref{Thm_nodal_domains_decomposition_1-Lap_eigenfunctions} that any nodal domain of any eigenfunction corresponding to $\Lambda_2(\Delta_1)$ must coincide with one of the three nodal domains displayed in \Cref{fig:4-weak-nodal-2nd_eigenfunction}. 

Furthermore, by considering the generalized eigenvalue equation in \Cref{Def:generalized_p-eigenpair} together with the characterization of the subgradients of $\|\grad f\|_1$ and $\|f\|_1$ in \eqref{1_subgradient}, one can check that any such eigenfunction must be constant on each of the nodal domains shown in \Cref{fig:4-weak-nodal-2nd_eigenfunction}. Finally, we observe that the set of nonzero functions that are linear combinations of the characteristic functions of the three nodal domains has genus $3$, and that at least one such linear combination (for instance, the sum of the three characteristic functions) is not an eigenfunction corresponding to $\Lambda_2(\Delta_1)$. Hence,
$\gamma\text{-}\mathrm{mult}(\Lambda_2(\Delta_1))
=\gamma\big(X_{\Lambda_2(\Delta_1)}\big)\le 3-1=2.
$ 
Since $\gamma\text{-}\mathrm{mult}(\lambda)\ge \mathrm{mult}(\lambda)$, we conclude that
$\gamma\text{-}\mathrm{mult}(\Lambda_2(\Delta_1))\ge \mathrm{mult}(\Lambda_2(\Delta_1))=2$.
This example shows that the conclusion of \Cref{thm:second-eigen-1} does not hold when either the number of positive or the number of negative weak nodal domains is smaller than $2$.

\begin{figure}[t] \centering
\begin{tikzpicture}[scale=1.2]
\node (1) at (-1,0) {$\mathbf{-}$};
\node (2) at (1/2,1) {$\mathbf{-}$};
\node (4) at (1/2,-1) {$\mathbf{-}$};
\node (3) at (1,0) {$\mathbf{-}$};
\node (12) at (-1/2,1) {$\mathbf{-}$};
\node (11) at (-1/2,-1) {$\mathbf{-}$};
\node (5) at (-2,0) {$\mathbf{+}$};
\node (6) at (2,0) {$\mathbf{+}$};
\node (7) at (-4,1) {$\mathbf{+}$};
\node (8) at (-4,-1) {$\mathbf{+}$};
\node (9) at (4,1) {$\mathbf{+}$};
\node (10) at (4,-1) {$\mathbf{+}$};
\draw (1) to (2);
\draw (1) to (3);
\draw (1) to (4);
\draw (1) to (5);
\draw (2) to (3);
\draw (2) to (4);
\draw (3) to (4);
\draw (3) to (6);
\draw (5) to (7);
\draw (5) to (8);
\draw (7) to (8);
\draw (6) to (9);
\draw (10) to (9);
\draw (6) to (10);
\draw (11) to (1);
\draw (11) to (2);
\draw (11) to (3);
\draw (11) to (4);
\draw (12) to (1);
\draw (12) to (2);
\draw (12) to (3);
\draw (12) to (4);
\draw (12) to (11);
\draw (1) circle(0.23);
\draw (2) circle(0.23);
\draw (3) circle(0.23);
\draw (4) circle(0.23);
\draw (5) circle(0.23);
\draw (6) circle(0.23);
\draw (7) circle(0.23);
\draw (8) circle(0.23);
\draw (9) circle(0.23);
\draw (10) circle(0.23);
\draw (11) circle(0.23);
\draw (12) circle(0.23);
\end{tikzpicture}
    \caption{An example graph discussed in \cref{ex_3_nodal_domains} proving that if a $1$-Laplacian eigenfunction has either the number of strong or weak nodal domains smaller than 2, then it can induce more than 2 weak nodal domains. In particular, it proves that the hypotheses of \Cref{thm:second-eigen-1} are necessary.}
    \label{fig:4-weak-nodal-2nd_eigenfunction}
\end{figure}

\end{example}
Differently, the next example shows a situation where \Cref{thm:second-eigen-1} applies.

\begin{example}\label{ex_4_nodal_domains}
Consider the unweighted graph with $13$ vertices shown in \Cref{fig:4-strong_nodal_domains-2nd_eigenfunction}, and let $f$ be the function taking value $1$ on the vertices labeled `$+$’ and value $-1$ on the vertices labeled `$-$’. It can be shown that $f$ is a $1$-Laplacian eigenfunction corresponding to the second variational eigenvalue, whose variational multiplicity equals $3$. 
Indeed, by studying the Cheeger constants of the graph in \Cref{fig:4-strong_nodal_domains-2nd_eigenfunction}, we deduce from \Cref{thm:1-eigenpair_and_cheeger_constants} that
$
\Lambda_2(\Delta_1)=\Lambda_3(\Delta_1)=\Lambda_4(\Delta_1)=\tfrac13,
$ 
which yields $\mathrm{mult}(1/3)\ge 3$. By checking all characteristic functions $\chi_A$, we infer from \Cref{Thm_nodal_domains_decomposition_1-Lap_eigenfunctions} that any nodal domain of any eigenfunction corresponding to $\Lambda_2(\Delta_1)$ must coincide with one of the four nodal domains displayed in \Cref{fig:4-strong_nodal_domains-2nd_eigenfunction}. 

Moreover, by considering the generalized eigenvalue equation in \Cref{Def:generalized_p-eigenpair} together with the characterization of the subgradients of $\|\grad f\|_1$ and $\|f\|_1$ in \eqref{1_subgradient}, one can verify that any such eigenfunction must be constant on each of the nodal domains shown in \Cref{fig:4-strong_nodal_domains-2nd_eigenfunction}. Finally, we observe that the set of nonzero functions obtained as linear combinations of the characteristic functions of the four nodal domains has genus $4$, and that at least one such linear combination (for instance, the sum of the four characteristic functions) is not an eigenfunction corresponding to $\Lambda_2(\Delta_1)$. Hence,
\[
\gamma\text{-}\mathrm{mult}(\Lambda_2(\Delta_1))
=\gamma\big(X_{\Lambda_2(\Delta_1)}\big)\le 4-1=3.
\]
Since $\gamma\text{-}\mathrm{mult}(\lambda)\ge \mathrm{mult}(\lambda)$, we conclude that
\[
\gamma\text{-}\mathrm{mult}(\Lambda_2(\Delta_1))\ge \mathrm{mult}(\Lambda_2(\Delta_1))=3.
\]

\begin{figure}[t]
    \centering
\begin{tikzpicture}[scale=1.4]
\node (1) at (0,1) {$+$};
\node (2) at (1,0) {$+$};
\node (3) at (1,1) {$+$};
\node (4) at (2,1.5) {$0$};
\node (5) at (3,1) {$-$};
\node (6) at (4,1) {$-$};
\node (7) at (3,0) {$-$};
\node (8) at (1,2) {$-$};
\node (9) at (1,3) {$-$};
\node (10) at (0,2) {$-$};
\node (11) at (3,2) {$+$};
\node (12) at (3,3) {$+$};
\node (13) at (4,2) {$+$};
\draw[thick] (1) to (3);
\draw[thick] (2) to (3);
\draw[thick] (3) to (4);
\draw[thick] (4) to (5);
\draw[thick] (5) to (6);
\draw[thick] (5) to (7);
\draw[thick] (4) to (8);
\draw[thick] (8) to (9);
\draw[thick] (8) to (10);
\draw[thick] (4) to (11);
\draw[thick] (11) to (12);
\draw[thick] (11) to (13);
\draw (1) circle(0.2);
\draw (2) circle(0.2);
\draw (3) circle(0.2);
\draw (4) circle(0.2);
\draw (5) circle(0.2);
\draw (6) circle(0.2);
\draw (7) circle(0.2);
\draw (8) circle(0.2);
\draw (9) circle(0.2);
\draw (10) circle(0.2);
\draw (11) circle(0.2);
\draw (12) circle(0.2);
\draw (13) circle(0.2);
\end{tikzpicture}
    \caption{An example graph discussed in \cref{ex_4_nodal_domains} where the second $1$-Laplacian eigenvalue has variational multiplicity $3$ and the hypotheses of \Cref{thm:second-eigen-1} are satisfied.}
    \label{fig:4-strong_nodal_domains-2nd_eigenfunction}
\end{figure}
\end{example}

The next example is devoted to a preliminary discussion on how the spectrum of the $p$-Laplacian varies with $p$. 

\begin{example}\label{exam:complete} 
Consider the complete graph of order $N$ with uniform weights. In \cite{Amghibech1}, it was shown that for $p \not\in \{1,2,\infty\}$, the number of eigenvalues of $\Delta_p$ is  
\(
[N/2](N-[N/2]) + 1,
\)
and every eigenvalue has $\gamma$-multiplicity equal to $1$. However, the proof provided in \cite{Amghibech1} contained certain flaws. For completeness, we provide a corrected proof in the appendix (see \Cref{thm:complete}).

Notably, for $N \geq 7$, there exist approximately $O(N^2)$ distinct non-variational eigenvalues of the $p$-Laplacian. In contrast, it has been established in \cite{zhang2021homological} that all $1$-Laplacian eigenvalues of a complete graph with $N$ vertices are variational. Specifically, it can be verified that the sum of the $\gamma$-multiplicities of all eigenvalues of $\Delta_1$ on the complete graph $K_N$ equals $N$, implying that all eigenvalues are variational.

Similarly, as we will discuss in \Cref{Sec:inf_eigenproblem}, when $p = \infty$, there do not exist non-variational eigenvalues, since the only nonzero eigenvalue is given by $2$.
\end{example}

Combining the results of the previous two examples, \Cref{exam:path} and \Cref{exam:complete}, we summarize in \Cref{tab:variational?} a classification scheme indicating whether non-variational eigenvalues exist for the complete graph and the path graph (as a specific case of a tree) for different choices of $p$. 
For completeness, we preliminarily include results on the $\infty$-Laplacian, which will be further discussed in \Cref{Sec:inf_eigenproblem}.
Further examples and discussions can be found in \cite{zhang2021homological}. 

\begin{table}[t]
\resizebox{\textwidth}{!}{
\begin{tabular}{cccc}
\toprule
 & $\Delta_p$ ($p\in (1,\infty)\setminus\{2\}$) & $\Delta_1$ & $\Delta_{\infty}$ \\ 
\midrule 
complete graph, e.g., $K_7$ &  {Yes} \cite{Amghibech1}  & No \cite{zhang2021homological} & No, \Cref{Ex_inf_eigenvalues_complete_graph} \\ \hline
tree graph, e.g., $P_7$ &  {No} \cite{DEIDDA2023_nod_dom} & Yes \cite{zhang2021homological} &  Yes, \Cref{ex:path-nonvariational}\\
\bottomrule
\end{tabular}
}
\caption{\textbf{Is there a non-variational eigenvalue?} The table answers the question for different values of $p$ and two choices of special graph structures.}\label{tab:variational?}
\end{table}

\subsection{Relations with Dirichlet Cheeger constants}

In this section, we recall some combinatorial lower bounds for the variational eigenvalues of graph $1$-Laplacian that have been introduced in \cite{zhang2021homological}. And later, as a partial consequence of these lower bounds, we introduce a discussion about the Dirichlet Cheeger constants of the graph. In particular, similarly to \Cref{thm:1-eigenpair_and_cheeger_constants}, in \Cref{thm:Dirichlet_Cheeger_constants_1-Lap_eigenvalues} we establish new upper and lower bounds for the $1$-Laplacian variational eigenvalues in terms of the Dirichlet Cheeger constants.
Let $P:=\{A_1,\cdots,A_k\}\in\mathcal{D}_k(\Gc)$ be a family of subsets of the nodes of the graph such that each $A_i$ is either a \textbf{singleton} or induces a \textbf{triangle}, where we say that $A$ induces a triangle if $A=\{u_1,u_2,u_3\}$ and $w_{u_iu_j}>0$ for any $i\ne j$. We call such $P$ an \textbf{ST subpartition}. 
We denote by $\mathcal{SD}(\Gc):= \cup_{k=1}^N\mathcal{D}_k(\Gc)$ the collection of all these ST subpartitions.  
For any  $P=\{A_1,\cdots,A_k\}\in \mathcal{SD}(\Gc)$,  let $$h_*(P):=
\min\limits_{A\subset \cup_{i=1}^kA_i:|A\cap A_i|\le 1}c(A)$$
and  let $$\ell(P):=\sum_{i=1}^k \ell(A_i),\;\text{ where }\ell(A_i)=\begin{cases}1,&\text{ if }A_i\text{ is a singleton},\\
2,&\text{ if } {A_i}\text{ induces a triangle}.
\end{cases}$$
We define $\mathcal{SD}_{k}(\Gc)=\{P\in \mathcal{SD}(\Gc):\ell(P)=k\}$. Then any triangle $A_i=\{u_1,u_2,u_3\}$ can be equipped with a space hexagon defined as $\text{HEX}[u_1,u_2,u_3]=[\chi_{u_1},-\chi_{u_2}]\cup [\chi_{u_2},-\chi_{u_3}]\cup [\chi_{u_3},-\chi_{u_1}]\cup[-\chi_{u_1},\chi_{u_2}]\cup [-\chi_{u_2},\chi_{u_3}]\cup [-\chi_{u_3},\chi_{u_1}]$, where $\chi_u$ is the indicator function of the node $u$ that takes value $1$ on $u$ and zero elsewhere. Such hexagon $\text{HEX}[u_1,u_2,u_3]$ is odd-homeomorphic to $\mathbb{S}^1$, meaning that there exists a function $\Phi:A\to \mathbb{S}^{1}$ such that $\Phi$ is a homeomorphism, and $\Phi(-x)=-\Phi(x)$. Thus, $\text{HEX}[u_1,u_2,u_3]$ has genus $2$. Note that any function from $\text{HEX}[u_1,u_2,u_3]$ is such that its nodal domains are composed exactly by one vertex. As a consequence, it is possible to prove the following Theorem from \cite{zhang2021homological}.

\begin{theorem}[\cite{zhang2021homological}]\label{lemma:minmax-eigen-lower} 
Let $\Lambda_k(\Delta_1)$ be the $k$-th $1$-Laplacian variational eigenvalue on a graph $\Gc$. Then, for any $k=1,\cdots,N$ 
$$\Lambda_{k}(\Delta_1)\ge \max_{P\in \mathcal{SD}_{N-k+1}}h_*(P).$$ 
{Moreover, if the graph $\Gc$ is a tree, the equality  $\Lambda_{k}(\Delta_1)= \max\limits_{P\in \mathcal{SD}_{N-k+1}}h_*(P)$ holds for all $k$.}
\end{theorem}

Next, we introduce the Dirichlet Cheeger constants of the graph.
\begin{definition}[Dirichlet Cheeger constant]\label{DEF:dirichlet_cheeger_constants}Given a graph $\Gc$ and a nonempty subset $A\subset \internalnodes$, the Dirichlet Cheeger constant of $\Gc$ on $A$ is defined as $$h_1\big(\Gc(A)\big)=\min\limits_{B\subseteq A}c(B).$$
\end{definition}
\begin{remark}
We remark that $h_1\big(\Gc(A)\big)$ is indeed the Cheeger constant of the graph $\Gc(A)$ induced by $A$ and with boundary $\nodeset\setminus A$. In particular, from \Cref{thm:1-eigenpair_and_cheeger_constants}, we have that $h_1\big(\Gc(A)\big)$ coincides with the first $1$-Laplacian eigenvalue of the graph $\Gc(A)$, which we denote by $\Lambda_1\big(\Delta_1,\Gc(A)\big)$. Namely,
\begin{equation}
    h_1\big(\Gc(A)\big)=\Lambda_1\big(\Delta_1,\Gc(A)\big).
\end{equation}
\end{remark}

Next, we define $\mathcal{P}_k$ as the family of all subsets of $\internalnodes$ containing at least $k$ vertices:
\begin{equation}
    \mathcal{P}_k:=\{A\subset \internalnodes\,\big|\;|A|\geq k\},
\end{equation}
and accordingly define the \textbf{$k$-th Dirichlet Cheeger constant} of the graph as
\begin{equation}
    H_k(\Gc)=\max_{A\in \mathcal{P}_{N-k+1}}h_1\big(\Gc(A)\big).
\end{equation}
It is immediate to observe that if $k_1>k_2$, then $H_{k_1}(\Gc)\geq H_{k_2}(\Gc)$.
Moreover, as we prove in the next theorem, the variational $1$-Laplacian eigenvalues can be bounded in terms of the Dirichlet Cheeger constants.
To this end, note that, as a direct consequence of Theorem~\ref{lemma:minmax-eigen-lower}, if we take all ST subpartitions $P:=(A_1,\dots,A_k)$ in $\mathcal{SD}_k$ to be composed of singletons, then
\[
h_*(P)=h_1\big(\Gc(\cup_{i=1}^k A_i)\big).
\]
Now, for any $A:=\{v_1,\dots,v_{N-i+1}\}\in \mathcal{P}_{N-k+1}$ realizing $H_k(\Gc)$ (i.e., $H_k(\Gc)=h_1\big(\Gc(A)\big)$), we have $|A|=N-i+1$ with $i\le k$, and thus
$P:=(\{v_1\},\dots,\{v_{N-i+1}\})\in \mathcal{SD}_{N-i+1}$.
Then, as a consequence of Theorem~\ref{lemma:minmax-eigen-lower} and the identity
$h_*(P)=h_1\big(\Gc(A)\big)$, we obtain
\[
\Lambda_{k}(\Delta_1)\ge \Lambda_{i}(\Delta_1)\ge h_*(P)
= h_1\big(\Gc(A)\big)=H_k(\Gc).
\]
We summarize this lower bound, together with a corresponding upper bound, in the next theorem.

\begin{theorem}\label{thm:Dirichlet_Cheeger_constants_1-Lap_eigenvalues}
Let $\Lambda_k(\Delta_1)$ be the $k$-th variational eigenvalue of the $1$-Laplacian and assume $f_k$ to be a corresponding eigenfunction. Let $V^{+}:=\{v\in\internalnodes\text{ s.t. } f_k(v)>0\}$ and $m:=N-|V^+|$, where $N=|\internalnodes|$, then 
$$
H_{k}(\Gc)\leq \Lambda_k(\Delta_1)\leq H_{m+1}(\Gc) .  
$$
\end{theorem}
\begin{proof}
We start from the lower bound. For any subset $A\subset \internalnodes$ with $|A|=N-k+1$, let $\mathcal{H}_0\big(\Gc(A)\big)$ be the set of functions supported on $\Gc(A)$, i.e $\mathcal{H}_0\big(\Gc(A)\big)=\{f:\internalnodes\to \R: \mathrm{supp}(f)\subset A\}$. Note that $\dim \mathcal{H}_0\big(\Gc(A)\big)=N-k+1$. 
So, for any subset $S\subset \R^N$ with $\gamma(S)\ge k$, the intersection property of genus, see \cref{Lemma_Krasnoselskii_intersection}, gives $S\cap \mathcal{H}_0\big(\Gc(A)\big)\setminus \{0\}\ne\varnothing$, and hence
\begin{equation}
\Lambda_k(\Delta_1)=\inf_{\gamma(S)\ge k}\sup_{f\in S}\mathcal{R}_1(f)\ge \inf_{f\in \mathcal{H}_0(\Gc(A))}\mathcal{R}_1(f).
\end{equation}
Moreover, from \Cref{remark_cheeger_constants_and_support}, for any $f\in \mathcal{H}_0\big(\Gc(A)\big)$, there exists $B\subset A$ such that $\mathcal{R}_1(f)\ge c(B)$, which yields
\begin{equation}
 \inf_{f\in \mathcal{H}_0(\Gc(A))}\mathcal{R}_1(f) \ge \min\limits_{B\subset A}c(B)=h_1\big(\Gc(A)\big).
\end{equation}
In consequence, 
$\Lambda_k(\Delta_1)\ge h(A)$, and the proof of the lower bound is concluded. 

To prove the opposite inequality note that if $f_k$ is an eigenfunction relative to $\Lambda_k(\Delta_1)$, there exist $\Xi\in \partial_{\grad f_k}\|\grad f_k\|_1$ and $\xi\in\partial_{f_k}\|f_k\|_1$ such that 
\begin{equation}
    \grad^T\Xi=\Lambda_k(\Delta_1) \nodeweight\odot \xi.
\end{equation}
Now summing the eigenvalue equation over the nodes in $V^+$ and using the characterization of the subgradients in \eqref{1_subgradient}, we observe that all the contributions of the edges connecting two nodes both in $V^+$ are deleted and we obtain 
\begin{equation}
 c(V^+)=\Lambda_k(\Delta_1).   
\end{equation}
Similarly if we take any $B\subset V^+$, since $|\Xi(u,v)|\leq 1$  for any edge, summing the eigenvalue equation over $B$, we obtain 
\begin{equation}
    c(B)\geq \frac{\sum_{u\in B,\, v\in B^c}\edgelength_{uv}\Xi(v,u)}{\sum_{u\in B}\nodeweight_u}=\Lambda_k(\Delta_1).
\end{equation} 
In particular,
$h_1\big(\Gc(V^+)\big)=\min_{B\subset V^+}c(B)=\Lambda_k(\Delta_1),$
which finally yields the upper bound
$\Lambda_k(\Gc)\leq h_1\big(\Gc(V^+)\big)\leq H_{m+1}(\Gc).$
\end{proof}

\begin{corollary}
    Assume that $\Lambda_{k-1}(\Delta_1)<\Lambda_k(\Delta_1)$. Then any eigenfunction corresponding to $\Lambda_k$ has no more than $N-k+1$ positive (or negative) entries. 
\end{corollary}
\begin{proof}
    If by absurd there was an eigenfunction $f_k$ such that $f_k^+$ was supported on more than $N-k+1$ nodes, then from \Cref{thm:Dirichlet_Cheeger_constants_1-Lap_eigenvalues}, $m<N-(N-k+1)=k-1$ (i.e., $m+1\le k-1$):
    \begin{equation}
        H_{k-1}(\Gc)\leq \Lambda_{k-1}(\Delta_1)< \Lambda_{k}(\Delta_1)\leq H_{m+1}(\Gc)\leq H_{k-1}(\Gc)
    \end{equation}
which is absurd.
\end{proof}

Combining the $k$-way Cheeger inequality $h_k(\Gc)\ge \Lambda_k(\Delta_1)$ and Proposition  \ref{thm:Dirichlet_Cheeger_constants_1-Lap_eigenvalues}, we also have 
\begin{equation}\label{eq_cheeger_dirichlet_inequality}
H_{k}(\Gc) = \max_{|A|=N-k+1}h_1\big(\Gc(A)\big) \leq h_k(\Gc).    
\end{equation}

It can be directly verified that \eqref{eq_cheeger_dirichlet_inequality} is actually an equality for path graphs and star graphs, and in \cite{MengZhang25}, it is shown that \eqref{eq_cheeger_dirichlet_inequality} shrinks to an equality for tree graphs. 

\begin{theorem}
Let $\Gc$ be a tree graph. Then, the $k$-way Cheeger constants introduced in \Cref{DEf:Cheeger_constant} are equal to the $k$-way Dirichlet Cheeger constants introduced in \Cref{DEF:dirichlet_cheeger_constants}, i.e.
\begin{equation}\label{eq:conj-equal}
h_k(\Gc)= H_k(\Gc),\;\; \forall k.     
\end{equation}    
\end{theorem} 

Note that we can show that $h_k(\Gc)$ (resp., $H_k(\Gc)$) is a lower bound (resp., an upper bound) of the min-max and max-min eigenvalues 
$\Lambda_k(\Delta_1)$, $\Lambda_k^D(\Delta_1)$, $\Lambda_k^Y(\Delta_1)$, $\tilde{\Lambda}_k(\Delta_1)$, $\tilde{\Lambda}_k^D(\Delta_1)$, $\tilde{\Lambda}_k^Y(\Delta_1)$ defined in \cref{appendix_variational_spectrum}. It follows from the equality \eqref{eq:conj-equal} that all the 8 quantities coincide on trees, which strengthen the known identity $\Lambda_k(\Delta_1)=h_k(\Gc)$, and also gives evidence to  Open Problem \ref{eq:conj-minmax}.

\begin{remark}
About the max min eigenvalue, note that, for any family of sets $\mathcal{F}_k^*$ used to define the three types of min-max variational eigenvalues in \cref{appendix_variational_spectrum}, if $A\in\mathcal{P}_{N-k+1}$, the set $X_A=\{f:\internalnodes\to \R: \mathrm{supp}(f)\subset A\}\cap \mathbb{S}^{N-1}\in \mathcal{F}_{N-k+1}^*$. Thus the max-min eigenvalue satisfies
\begin{equation}
\begin{aligned}
      \tilde{\Lambda}_k^{*}(\Delta_1):=&\sup_{B\in\Fc_{N-k+1}^*}\inf_{f\in B}\rayl_1(f)\ge  \sup_{A\in \mathcal{P}_{N-k+1}}\inf_{f\in X_A}\rayl_1(f)=\\ &=\max_{A\in \mathcal{P}_{N-k+1}} h_1\big(\Gc(A)\big)= H_k(\Gc).
\end{aligned}
\end{equation}
On the other hand, given a linear subspace $Y=\mathrm{span}(\chi_{A_1},\cdots,\chi_{A_k})$ with the family $\{A_1,\cdots,A_k\}\in \mathcal{D}_k(\Gc)$, and given any $B\in\Fc_{N-k+1}^*$, then from the intersection property in \cref{Lemma_Krasnoselskii_intersection} we have $Y\cap B\ne\varnothing$. Thus
\begin{equation}
\begin{aligned}
\tilde{\Lambda}_k^{*}(\Delta_1):&=\sup_{B\in\Fc_{N-k+1}^*}\inf_{f\in B}\rayl_1(f)\le \tilde{\Lambda}_k^{*}(\Delta_1)=\sup_{B\in\Fc_{N-k+1}^*}\inf_{f\in B\cap Y}\rayl_1(f) \le
\\ & \le \sup_{f\in Y}\rayl_1(f)  = \max_{i=1,\cdots,k}c(A_i)
\end{aligned}
\end{equation}
and due to the arbitrariness  of $(A_1,\cdots,A_k)\in \mathcal{D}_k(\Gc)$, we have
$\tilde{\Lambda}_k^{*}(\Delta_1)\le h_k(\Gc)$.
\end{remark}

\subsection{Relations with independence numbers}\label{Subsec:1_lap_indep_numb}
In this section, we recall a result from \cite{ZHANG_top_mult} that relates the multiplicity of the largest eigenvalue of the graph $1$-Laplacian to the independence number of the graph. \textbf{In this section we assume the graph not to have a boundary and we assume the length of any edge to be uniform, $\edgelength_{uv}=1$ for all $(u,v)\in \edgeset$}. The independence number of a graph is the size of the largest set of vertices such that any two distinct vertices in the set are at distance $2$.  Formally:
\begin{definition}\label{2_indip_number}
The independence number of a graph $\Gc$ is defined as:
$$\alpha_2(\Gc)=\max  \{k:\exists 
v_1,\cdots,v_k\in \nodeset \text{ with }\mathrm{dist}(v_i,v_j)\ge 2,\forall i\ne j \}.$$
A set $\{v_1,\cdots, v_{\alpha_2(\Gc)}\}$ as above is defined as a maximum independent set.
\end{definition}

Then we have the following result for the normalized $1$-Laplacian from \cite{ZHANG_top_mult} 
\begin{theorem}\label{Thm_1_lap_mult_independence_number}
 Let $\Gc$ be a graph with uniform weights on the edges $\edgelength_{uv}=1$ for all $(u,v)\in\edgeset$ and measure $\nodeweight_u=\deg(u)$ for any node $u\in \nodeset$. 
 Assume $\Lambda_N(\Delta_1)$ to be the largest eigenvalue of the $1$-Laplacian, then
 $$\alpha_2(\Gc) \leq \mathrm{mult}\big(\Lambda_N(\Delta_1)\big)\leq 2\alpha_2(\Gc).$$
 In addition 
$$\gamma\text{-}\mathrm{mult}\big(\Lambda_N(\Delta_1)\big)\geq \alpha^*(\Gc)\geq \alpha_2(\Gc),$$
 where $\alpha^*(\Gc):=\max\{b+2q$ such that there are $b$ nodes and $q$ third order circles in $\Gc$ that are pairwise non-adjacent$\}$. In particular, the equality $\alpha_2(\mathcal{G}) =\mathrm{mult}\big(\Lambda_N(\Delta_1)\big)$ holds for any triangle-free perfect graph $\mathcal{G}$. Similarly, the equality $\mathrm{mult}\big(\Lambda_N(\Delta_1)\big)=2\alpha_2(\mathcal{G})$ holds for the triangle graph, and $\gamma\text{-}\mathrm{mult}\big(\Lambda_N(\Delta_1)\big)=\alpha^*(\mathcal{G})$ holds for complete graphs and triangle-free perfect graphs. 
\end{theorem}

We recall that in \cref{Thm_1_lap_mult_independence_number} a simple graph $\Gc$ is called a \emph{triangle-free perfect graph} if $\Gc$ contains no cycles of length three, and for every induced subgraph of $\Gc$, the clique covering number equals the chromatic number \cite{godsil01}.

\begin{remark}
The author of \cite{ZHANG_top_mult} gives  another upper bound for the variational multiplicity of the maximum eigenvalue of the $1$-Laplacian, $\mathrm{mult}\big(\Lambda_N(\Delta_1)\big)$, using a variant of  clique covering number. 
The clique covering number $\kappa(\Gc)$ is the smallest number of cliques that cover $\Gc$, where a clique is a complete subgraph. We have $$\mathrm{mult}\big(\Lambda_N(\Delta_1)\big)\le \kappa^*(\Gc)\le  2\kappa(\Gc)$$
where $\kappa^*(\Gc):=\min\{b+2q:\exists b\text{ edges and }q\text{ cliques of order larger than }3 \text{ that cover }\nodeset\}$, equivalently, $\kappa^*(\Gc)$ stands for the minimum possible $b+2q$ such that  there are $V_1,\cdots,V_b$ (each of them is a clique of cardinality 1 or 2), and $V_{b+1},\cdots,V_{b+q}$ (each of them is a clique of cardinality at least 3) satisfying $V_1\cup \cdots \cup V_{b+q}=\nodeset$. 
\end{remark}


\section{Graph $\infty$-Laplacian spectrum}\label{Sec:inf_eigenproblem}

Similarly to the $1$-Laplacian case, the $\infty$-Laplacian spectrum is also closely linked to certain topological properties of the graph. It has been observed in \cite{deidda2024_inf_eigenproblem, Lind2, Lind3, bungert2021eigenvalue} that there are connections between the radii of balls inscribed in the graph (or in a bounded Euclidean domain) and the $\infty$-eigenvalues. 
In \Cref{Sec:intro-sec_graph_setting}, we introduced a notion of distance between nodes of the graph defined in terms of the length of the shortest path. Accordingly, given a node $u$ and a real number $R$, we say that the ball of radius $R$ centered at $u$,
\[
B_R(u)=\{v\in \nodeset \,|\; d(u,v)< R\},
\]
is inscribed in the graph if the distance of $u$ from the boundary is greater than $R$, $d_{\boundary}(u)\geq R$.
Next, we recall a characterization of the $\infty$-eigenpairs from \cite{deidda2023PhdThesis}, where we adopt the notation
$\nodemaxset(f)=\{v\in \internalnodes\,|\;|f(v)|=\|f\|_{\infty}\}$.

\begin{proposition}\label{prop:Infinity_eigenfunctions_characterization}
$(f,\Lambda)$ is an $\infty$-eigenpair with $f$ not constant if and only if there exists a shortest path $\Gamma=\{(u_i,u_{i+1})\}_{i=0}^{n-1}$ such that
\begin{enumerate}
    \item $u_0, u_n\in \nodemaxset(f)\cup \boundary.$
    \item $Kf(u_i,u_{i+1})=\|\grad f\|_{\infty}\:\;\forall i=1,\dots,n.$
    \item $\Lambda=\Big(|f(u_0)|+|f(u_n)|\Big)\Big/\Big(\length(\Gamma)\;\|f\|_{\infty}\Big)$    
\end{enumerate}
where $f(u)=0$ if $u\in \boundary$.
In particular, whenever 1.2.3 are satisfied and $u_i\not\in \boundary$ with $i=0$ or $i=n$:  
    $$\frac{1}{\Lambda}=\min\bigg\{\min_{\{v|f(v)=-f(u_i)\}}\frac{d(u_i,v)}{2},d_{\boundary}(u_0)\bigg\}\,.$$
\end{proposition}

The last proposition proves that if $u_0\in \internalnodes$, the ball $B_{1/\Lambda}(u_0)$ is inscribed in the graph. Moreover $B_{1/\Lambda}(u_0)\cap B_{1/\Lambda}(v)=\emptyset$ for any $v$ such that $f(v)=-f(u_0)$, which is a non-emptyset condition whenever also $u_n\in \internalnodes$.
\begin{lemma}\label{Lemma_inf_eigenvector}
    Given any node $u\in \internalnodes$ and any pair of nodes $v_1,v_2\in \internalnodes$ with $d(v_1,v_2)/2\leq \min\{d_{\boundary}(v_1),d_{\boundary}(v_2)\}$ we can associate them the $\infty$-eigenpairs 
    $\big(f_{u}^{\Lambda^*},\Lambda^*\big)$ and $\big(f_{v_1}^{\Lambda^{**}}-f_{v_2}^{\Lambda^{**}},\Lambda^{**}\big)$, respectively; where 
    $$\Lambda^*=\frac{1}{\dist_{\boundary}(u)}, \qquad \Lambda^{**}=\frac{2}{\dist(v_1,v_2)} \qquad\text{and}\qquad 
 f_{w}(v)=\max\{1-\Lambda d(v,w),0\}.$$
\end{lemma}
\begin{proof}
     It is not difficult to check that the functions $f_u$ and $f_{v_1}-f_{v_2}$ with the corresponding eigenvalues $\Lambda_1$ and $\Lambda_2$ satisfy the conditions of \Cref{prop:Infinity_eigenfunctions_characterization}.    
\end{proof}

In particular, we have the following result 
 \begin{proposition}\label{Prop:gamma_mult_inf_eigenvalues}
     Let $\Lambda$ be an $\infty$-eigenvalue then  $\gamma\text{-}\mathrm{mult}(\Lambda)\geq C$ where 
     \begin{equation}
         \begin{aligned}
            C = \max k+l \text{ such that }\;\exists 
         v_1,\dots, v_{2k},\; u_{1},\dots, u_{l} \in \internalnodes \text{ with } \\
         \begin{cases} \dist(u,v) \geq 2/\Lambda \qquad &\forall u,v\in \{v_i\}_{i=1}^{2k}\cup\{u_i\}_{i=1}^l\\
           \dist_{\boundary}(v_i) \geq 1/\Lambda\qquad &\forall i=1,\dots,2k\\
           \dist(v_{2i-1}, v_{2i})=2/\Lambda \qquad &\forall i=1,\dots,k\\
           \dist_{\boundary}(u_{i})=1/\Lambda \qquad &\forall i=1,\dots,l
          \end{cases}
         \end{aligned}
     \end{equation}
      
 \end{proposition}
 \begin{proof}
      Suppose $v_1,\cdots,v_{2k},\; u_{1},\dots, u_{l}\in \internalnodes$ satisfy the hypothesis. Then, we shall check that for any $(t_1,\cdots,t_k, \bar{t}_1,\dots, \bar{t}_l)\ne 0$, $\Big(\sum_{i=1}^k t_i(f_{v_{2i-1}}-f_{v_{2i}})+\sum_{i=1}^l \bar{t}_i f_{u_i},\;\Lambda \Big)$ satisfies the conditions of \Cref{prop:Infinity_eigenfunctions_characterization}, where any function $f_{w}$ is defined as in \Cref{Lemma_inf_eigenvector}.
      
       In fact, we may assume without loss of generality that $0\neq|t_1|=\max\{\{t_i,\bar{t}_j\}$ varying $i$ and $j$. Hence, $f$ is nonzero because $\|f\|_\infty=|t_1|$ and $f$ satisfies the conditions in \Cref{prop:Infinity_eigenfunctions_characterization} on the  shortest path joining $v_1$ to $v_2$. We have proved that  any nonzero $f\in \pi=\mathrm{span}(f_{l,v_{1}}-f_{l,v_{2}},\ldots,f_{l,v_{2k-1}}-f_{l,v_{2k}}, f_{u_1},\dots, f_{u_l})$ is an  eigenfunction corresponding to the eigenvalue  $\Lambda$. Finally the genus of $\pi\setminus \{0\}$ is $k+l$, which means that the  multiplicity of the eigenvalue $\Lambda$ is at least $k+l$.
 \end{proof}

As a consequence of \cref{Prop:gamma_mult_inf_eigenvalues} we can characterize all of the $\infty$-Laplacian eigenvalues on a generic unweighted graph without boundary.

\begin{example}\label{EX:inf_noweights_spectrum}
    Assume $\Gc$ is an unweighted graph without boundary, i.e. $\edgelength_{uv}=1$ for all $(u,v)\in\edgeset$ and $\boundary=\emptyset$. Then the set of the $\infty$-eigenvalues is $\{0,2/l:l=1,\cdots,\mathrm{diam}(\Gc)\}$. 
    Indeed, from \cref{prop:Infinity_eigenfunctions_characterization}, any eigenvalue $\Lambda$ can be written as $\Lambda=2/l$, where $l$ is the distance between some nodes of the graph. Vice versa, given two nodes whose reciprocal distance is $l$, based on \cref{Lemma_inf_eigenvector} we can define an $\infty$-eigenvector corresponding to the eigenvalue $\Lambda=2/l$. 
    Moreover $f$ is an eigenfunction corresponding to $2/l$, if and only if $f$ satisfies the  following $l$-\textit{Shortest Path condition} ($\mathrm{SP}_l$ for brevity): 
    there exists a shortest path $v_0,v_1,\cdots,v_\ell$ s.t. 
$f(v_i)=(1-\frac{2i}{\ell})\|f\|_\infty$,  $i=0,\ldots,\ell$, and   $|f(v')-f(u')|\le \frac{2}{\ell}\|f\|_\infty$, $\forall\{v',u'\}\in E$.
Finally, the  $\gamma$-multiplicity of the eigenvalue $2/l$ is at least 
\begin{equation}
\max\left\{k\Big|\begin{array}{lr}
    \exists v_1,\cdots,v_{2k}\textrm{ s.t. }\mathrm{dist}(v_i,v_j)\ge l\text{ for }i\ne j\\
    \mathrm{dist}(v_{2i-1},v_{2i})=l\text{ for }i=1,\ldots,k\end{array}\right\}.
\end{equation}
\end{example}


\subsection{Relations with sphere packing}

We have seen that any $\infty$-eigenvalue basically represents the radius of some ball inscribed in the graph \cref{prop:Infinity_eigenfunctions_characterization}. In this section we use this information to establish relationships between the $\infty$-variational eigenvalues and the packing radii of the graph. We start from the definition of $k$-th packing radius, which is the maximal radius such that we can inscribe $k$ disjoint balls in the domain.
\begin{definition}\label{def:packing_radius}
    The $k$-th packing radius, $R_k(\Gc)$ of the graph $\Gc$ is defined as 
$$R_k(\Gc)=\max \{R:\exists 
v_1,\cdots,v_k\in \internalnodes \text{ with }\dist(v_i,v_j)\ge 2R,\; d_{\boundary}(v_i)\geq R,\; \forall i\neq j\}$$
\end{definition}

The following result relates the packing radii to the $\infty$-variational eigenvalues. This is a generalization of other results proved in \cite{deidda2024_inf_eigenproblem} for the graph setting and \cite{Lind2,Lind3} for the infinite dimensional setting.
\begin{theorem}\label{thm:k-inequality}
Given a connected graph $\Gc$, we have the equalities  $\Lambda_1(\inflap)=1/R_1(\Gc)$ and $\Lambda_2(\inflap)=1/R_2(\Gc)$. Moreover, for any $k=1,2,\cdots,N$,
$$\frac{1}{R_{\PNc(f_k)}(\Gc)}\le \Lambda_k(\inflap)\le \frac{1}{R_k(\Gc)}$$ 
 where $f_k$ is any eigenfunction corresponding to $\Lambda_k(\inflap)$ and $\PNc(f_k)$ is the number of perfect nodal domains induced by $f_k$. 
\end{theorem}

\begin{proof}
    The upper bound is easily proved taking $\{V_i\}_{i=1}^k$ that realizes $R_k(\Gc)$, the functions $f_{v_i}$ as in \Cref{Lemma_inf_eigenvector} and $\pi=\text{span}\{f_{v_i}\}$. Indeed one can easily observe that $\rayl_{\infty}(f)\leq 1/R_k(\Gc)$ for any $f\in \pi$. Thus, since $\gamma(\pi\cap S_\infty)=k$, by the definition of the variational eigenvalues, we have the upper bound in the thesis. We refer to \cite{deidda2024_inf_eigenproblem} for all the details as also for the equality $\Lambda_1(\inflap)=1/R_1(\Gc)$, $\Lambda_2=1/R_2(\Gc)$ (see also \cref{remark:
equality_first_variational_inf_eigenvalues_packing_radii} for the equality).

    To prove the lower bound proceed as in \Cref{thm:infty_nodal-count}. Take $\{V_i\}_{i=1}^{\PNc(f_k)}$ the perfect nodal domain induced by $f_k$ and for any perfect nodal domain $V_i$ consider $f_k|_{V_i}$ the function equal to $f_k$ on $V_i$ and equal to zero elsewhere. One can prove that $\max_{g\in \pi'}\rayl_{\infty}(g)\geq \Lambda_k(\inflap)$, where $\pi':=\mathrm{span}\{f_k|_{V_i}\}$ (see the proof of \Cref{thm:infty_nodal-count} for the computations). Thus, since $1/R_{\PNc(f_k)}(\Gc)\leq \max_{g\in \pi'}\rayl_{\infty}(g)$, we get the lower bound and conclude the proof.
\end{proof}

The result in \cref{thm:k-inequality} can be 

\begin{corollary}\label{Cor:packingradii_inequality}
    Under the same hypotheses of \Cref{thm:k-inequality}, we also have
$$\frac{1}{R_{m}(\Gc)}\leq \Lambda_k(\inflap)\leq \frac{1}{R_k(\Gc)}$$ 
where $m$ can be chosen as:
\begin{enumerate}
    \item 
    $\SNc(f_k)$ for any viscosity eigenfunction $f_k$ corresponding to $\Lambda_k(\inflap)$;
    \item 
    $\limsup
    \limits_{p\to+\infty} \SNc(f_{k,p})$ for any eigenfunction $f_{k,p}$ corresponding to $\lambda_k(\Delta_p)$.
\end{enumerate}  
\end{corollary}
\begin{proof}
    The proof follows from \Cref{Thm:lower_bound_perfect_nodal_count_limit_eigenpairs}, \Cref{thm:relate-3bounds} and \Cref{thm:k-inequality}.
\end{proof}

In the next corollary we observe that if we consider limits of $p$-Laplacian eigenfunctions or viscosity eigenfunctions, beside the perfect nodal domains, also the strong nodal domain count is significative to bound the variational $\infty$-eigenpair in terms of the packing radii. 

\begin{remark}\label{remark:
equality_first_variational_inf_eigenvalues_packing_radii}
We note that from the last \Cref{Cor:packingradii_inequality} it is straightforward to prove the equalities $\Lambda_1(\inflap)=1/R_1(\Gc)$ and $\Lambda_2(\inflap)=1/R_2(\Gc)$ in \Cref{thm:k-inequality}. Indeed from \Cref{Theorem_1st_eigen_characterization} we know that $\SNc(f_{1,p})=1$ and $\SNc(f_{2,p})\geq 2$ for any $p\in(1,\infty)$. 
\end{remark}

\begin{remark}\label{cor:tree}
   Additionally, as in \Cref{remark_variational_eigenvalues_on_trees_cheeger_constants}, it is possible to prove that on \textbf{trees}, where the nodal count is exact for $p\in(2,\infty)$ it holds the equality 
    \begin{equation}
        \Lambda_{k}(\inflap)=\frac{1}{R_k(\Gc)}\qquad \forall k=1,\dots,N.
    \end{equation}
    We refer to \cite{ge2025} for a formal proof.  
\end{remark}

\begin{example}\label{Ex_inf_eigenvalues_complete_graph}
Assume that $\Gc$ is a complete graph on $N$ vertices, without boundary, and whose edges have all length equal to $1$. From \Cref{thm:k-inequality} we have that 
\begin{equation}
   2=\frac{1}{R_N(\Gc)}=\frac{1}{R_2(\Gc)} \leq \Lambda_k(\inflap)\leq \frac{1}{R_N(\Gc)}=2 \qquad \forall k\in \{2,\dots,N\}.
\end{equation}

Thus all the variational eigenvalues, except the first one, are equal to $2$. In particular, since the $\gamma$-multiplicity is always greater than the variational one, we have $\gamma$-mult$(2)\geq N-1$. On the other hand, the eigenspace corresponding to $2$, $X_2$, is such that $X_2\cap S^{N-1}\subseteq S^{N-1}\setminus\{\underline{1},-\underline{1}\}$, where $\underline{1}$ is the constant vector. Note that $S^{N-1}\setminus\{\underline{1},-\underline{1}\}$ is homotopy equivalent to $S^{N-2}$, so its genus is equal to $N-1$. It follows that $\gamma(X_2)\leq N-1$ which proves the equality $\gamma(X_2)=N-1$. We conclude observing that from \Cref{EX:inf_noweights_spectrum} there are no other eigenvalues different from the variational ones.   
\end{example}
From the bounds that relate the $\infty$-eigenvalues to the packing radii of the graph, and the monotonicity of the spectral functions in \Cref{lemma:p-monotonic}, it is possible to establish some weaker bounds relating the $p$-Laplacian eigenvalues (for any $p$) to the packing radii of the graph. 
\begin{proposition}\label{pro:p-Lap-r_k}
For any given $p\ge 1$,  and $k=1,\cdots,N$, 
$$\frac{1}{|\nodeweight|}\Big(\frac{1}{R_m(\Gc)}\Big)^{p}\le \lambda_k(\Delta_p)\le \frac{|E|}{2}\Big(\frac{1}{R_k(\Gc)}\Big)^{p}$$
where the subscript $m$ can be taken as $\SNc(f_{k,p})$ or as $\SNc(f_{k,q})$, for any eigenfunction $f_{k,q}$, corresponding to $\lambda_k(\Delta_q)$ with $q\ge p-\epsilon$ for sufficiently small $\epsilon>0$.
\end{proposition}

\begin{proof}

By Lemma \ref{lemma:p-monotonic} and Theorem \ref{thm:k-inequality}, we have 
\begin{equation}
    \lambda_k(\Delta_p)\le \frac{|E|}{2}\Lambda_k^p(\inflap)\le \frac{|E|}{2}\left(\frac{1}{R_k(\Gc)}\right)^{p}.    
\end{equation}
Next, we move  to the lower bound of $\lambda_k(\Delta_p)$. 
We start from proving that if $\lambda_k(\Delta_p)$ has an eigenfunction with $m$ nodal domains, then 
 \begin{equation}
     \frac{1}{|\nodeweight|}\left(\frac{1}{R_m(\Gc)}\right)^{p}\le \lambda_k(\Delta_p).
 \end{equation}
Indeed, assume that $f$ is an eigenfunction corresponding
to $\lambda_k(\Delta_p)$ with multiplicity $r$, and $f$ has $m$ strong nodal domains 
which are denoted by $V_1,\ldots,V_{m}$. 
Let $\hat{f}=\sum_{i=1}^m\frac{1}{\|f|_{V_i}\|_\infty}f|_{V_i}$. Then, $\|\hat{f}|_{V_i}\|_\infty=1$ for all $i$. Recalling the expression of the $p$-Rayleigh quotient in \cref{eq_p_rayleigh_quotient} it is known that%
\footnote{Here, we use the property that for any  $\hat{f}\in \mathrm{span}(f|_{V_1},\cdots,f|_{V_m})$, if $f$ is an  eigenfunction  corresponding to  $\lambda_k(\Delta_p)$ with $p\in(1,\infty)$, then $\mathcal{R}_p(\hat{f})^{1/p}\le \lambda_k(\Delta_p)^{1/p}$. However,  if $f$ is an eigenfunction  corresponding to  $\Lambda_k
$,  $\mathcal{R}_\infty(\hat{f})\le \Lambda_k
$ may not hold even if we set $\|\hat{f}|_{V_i}\|_\infty=\|f\|_\infty$ $\forall i$ (see \cref{example:nodal_domains_inf_eigenpairs}).}
$\mathcal{R}_p(\hat{f})\le \lambda_k(\Delta_p)$ (see \cite{Tudisco1}), yielding the following bounds
\begin{equation}
    \mathcal{R}_\infty(\hat{f})\le |\nodeweight|^{\frac1p}\mathcal{R}_p(\hat{f})^{\frac1p}\le |\nodeweight|^{\frac1p}\lambda_k(\Delta_p)^{\frac1p}.
\end{equation} 
Next, we estimate $\mathcal{R}_\infty(\hat{f})$. Let $v_i\in V_i$ be such that $|\hat{f}(v_i)|=1$ for any $i=1,\cdots,m$. 
First we claim that given any pair $i,j$ with $i\ne j$, and $\hat{v}_0:=v_i,\hat{v}_1,\cdots,\hat{v}_{l-1},\hat{v}_l:=v_j$ a shortest path connecting $v_i$ and $v_j$, then:
\begin{equation}\label{p-eigenvalues_and_packing_radii_ Claim_1}
\rayl_{\infty}(\hat{f})\geq 2/\mathrm{dist}(v_i,v_j).
\end{equation}
Second we claim that for any $i\in \{1,\dots,m\}$, given $\hat{v}_0:=v_i,\hat{v}_1,\cdots,\hat{v}_{l-1},\hat{v}_l:=v_b\in \boundary$ a shortest path connecting $v_i$ to the boundary,  then:
\begin{equation}\label{p-eigenvalues_and_packing_radii_ Claim_2}
\rayl_{\infty}(\hat{f})\geq 1/\mathrm{dist}_{\boundary}(v_i),
\end{equation}
where $\mathrm{dist}_{\boundary}(v_i)$ denotes the distance of $v_i$ from the boundary.
As the proof of the two claims is very similar we prove only the first one. To prove the first claim we have to consider two possible situations: the first where $\hat{f}(v_i)=-\hat{f}(v_j)$, and the second where $\hat{f}(v_i)=\hat{f}(v_j)$. Start from the first, then,
\begin{equation}
\begin{aligned}
    \dist(v_i,v_j)\mathcal{R}_\infty(\hat{f})&= \dist(v_i,v_j)\sup_{(u,v)\in \edgeset} \edgelength_{uv}|\hat{f}(u)-\hat{f}(v)|
    \\ &=  \sum_{t=0}^{l-1}\frac{\sup_{(u,v)\in \edgeset} \edgelength_{uv}|\hat{f}(u)-\hat{f}(v)|}{\edgelength_{\hat{v}_t \hat{v}_{t+1}} }
\geq \sum_{t=0}^{l-1}\frac{\edgelength_{\hat{v}_t\hat{v}_{t+1}}}{\edgelength_{\hat{v}_t\hat{v}_{t+1}}}|\hat{f}(\hat{v}_t)-\hat{f}(\hat{v}_{t+1})|
\\ &\ge |\hat{f}(\hat{v}_0)-\hat{f}(\hat{v}_l)|=|\hat{f}(v_i)-\hat{f}(v_j)|=2    
\end{aligned}
\end{equation}
The case $\hat{f}(v_i)=\hat{f}(v_j)=1$ can be dealt similarly,  indeed if $v_i$ and $v_j$ lie in different  nodal domains, there exists $m_0\in\{1,\cdots,l-1\}$ such that $\hat{f}(\hat{v}_{m_0})\le0$, and  then, changing only the last step of the last inequality:
\begin{equation}
\sum_{t=0}^{l-1}|\hat{f}(\hat{v}_t)-\hat{f}(\hat{v}_{t+1})|\ge \Big(|\hat{f}(\hat{v}_0)-\hat{f}(\hat{v}_{m_0})|+|\hat{f}(\hat{v}_{m_0})-\hat{f}(\hat{v}_l)|\Big)\ge \big(|\hat{f}(v_i)|+|\hat{f}(v_j)|\big)=2.
\end{equation}
In any case, we obtain $\mathcal{R}_\infty(\hat{f})\ge\frac{2}{\mathrm{dist}(v_i,v_j)}$, for any $i\ne j$. Therefore, 
\begin{equation}
\mathcal{R}_\infty(\hat{f})\ge\frac{2}{\min\limits_{i\ne j}\mathrm{dist}(v_i,v_j)}\ge \frac{1}{R_m(\Gc)}.    
\end{equation}

As a consequence of the two claims we have $\rayl_\infty(\hat{f})\geq 1/R_m(\Gc)$, i.e. the thesis.
Finally, let us focus on the case $q\ge p-\epsilon$ for some sufficiently small $\epsilon>0$. 
By Lemma \ref{lemma:p-monotonic} and the above claim, for any $q\ge p$, and for any eigenfunction $f$ corresponding to $\lambda_k(\Delta_q)$ who has $m$ nodal domains,
\begin{equation}\label{eq:pq}
|\nodeweight|^{\frac1p}\lambda_k(\Delta_p)^{\frac1p}\ge|\nodeweight|^{\frac1q}\lambda_k(\Delta_q)^{\frac1q}\ge \frac{1}{R_m(\Gc)}.    
\end{equation}
Suppose the contrary, that for any $\epsilon>0$, there exist $q\ge p-\epsilon$ and an  eigenfunction $f_{k,q}$ corresponding to $\lambda_k(\Delta_q)$  such that $|\nodeweight|^{\frac1p}\lambda_k(\Delta_p)^{\frac1p}<1/R_m(\Gc)$, where $m=\SNc(f_{k,q})$. Then, such $q$ must lie in $(p-\epsilon,p)$. So, there exist a sequence $q_i<q$ with $q_i\to p$, $i\to +\infty$, and an eigenfunction  $f_{k,q_i}$ such that $|\nodeweight|^{\frac1p}\lambda_k(\Delta_p)^{\frac1p}<1/R_{m_i}(\Gc)$, where $m_i=\SNc(f_{k,q_i})$. Let $m=\inf\limits_{i} m_i=\min\limits_im_i$.  
It is clear that $|\nodeweight|^{\frac1p}\lambda_k(\Delta_p)^{\frac1p}<1/R_{m}(\Gc)$ and thus there is a fixed gap $\delta=1 / R_{m}(\Gc)-|\nodeweight|^{\frac1p}\lambda_k(\Delta_p)^{\frac1p}>0$. Therefore, 
\begin{equation}
|\nodeweight|^{\frac1p}\lambda_k(\Delta_p)^{\frac1p}+\delta \le \frac{1}{R_{m_i}(\Gc)}\le |\nodeweight|^{\frac{1}{q_i}}\lambda_k(\Delta_{q_i})^{\frac{1}{q_i}} \quad \forall i,
\end{equation}
which is a contradiction with $\lim_{i\to+\infty}|\nodeweight|^{\frac{1}{q_i}}\lambda_k(\Delta_{q_i})^{\frac{1}{q_i}}=|\nodeweight|^{\frac1p}\lambda_k(\Delta_p)^{\frac1p}$. 
\end{proof}

\begin{remark}
Note that considering the inequalities between the higher order Cheeger constants and the $p$-Laplacian variational eigenvalues \Cref{ThM:p-eigenpairs_and_cheeger_constants}, and taking $p\to+\infty$, we get $0\le \Lambda_k(\Delta_\infty)
\le 2 \max_{(u,v)\in\edgeset}\edgelength_{uv}$ which is trivial. Thus the estimate in Theorem \ref{thm:k-inequality} becomes more significative than the one in \Cref{ThM:p-eigenpairs_and_cheeger_constants} when $p$ is sufficiently large. In particular, when $p>>2$, the geometric  quantity $R_k(\Gc)$  can be used to obtain a  better  estimate  of the  variational eigenvalues of the graph $p$-Laplacian than the  multi-way Cheeger constant $h_k(\Gc)$. 
As an example, let $\Gc$ be an unweighted cycle graph of order $N\ge 3$ without boundary. Then, it is easy to see  $1/R_k(\Gc)=2h_k(\Gc)$. And, we immediately obtain $1/|\nodeweight|(1/R_m(\Gc))^{p}>2^{p-1}/p^p h_k^p(\Gc)$ when $p^p>n/2$. Similarly,  $|E|(1/R_k(\Gc))^{p}< 2^{p-1}h_k(\Gc)$ if $k\le N/2$ and $p>2+\log_2 N$. This implies  that for $k\le N/2$ and  $p>2+\log_2 N$, the estimate for the variational eigenvalue $\lambda_k(\Delta_p)$ using the quantity $R_k(\Gc)$ is more accurate than the one obtained by means of the multi-way Cheeger constant $h_k(\Gc)$. 
\end{remark}

Next we observe that the $\infty$-variational spectrum usually does not exhaust the spectrum, neither in the cases where it does for any $p$ smaller than $\infty$.

\begin{example}\label{ex:path-nonvariational}
Consider the standard unweighted path graph $P_N$ with $N$ vertices. 
From Corollary \ref{cor:tree}, we have 
$$\Lambda_k(\inflap)=\frac{2}{\lfloor\frac{N-1}{k-1}\rfloor},\;\;k=2,3,\ldots,N,$$ 
where $\lfloor\cdot\rfloor$ is the usual round-down function. 
Note that $N-1-\lfloor (N-1)/2\rfloor\ge2$ for $N\ge4$. Thus, for $N\ge4$, from \Cref{EX:inf_noweights_spectrum}, there exists a non-variational eigenvalue $\Lambda=2/\ell\in(\Lambda_2(\inflap),\Lambda_3(\inflap))$ for some $\ell\in \{1,2,\cdots,N-1\}$. 
In fact, when $N\ge4$, for any $\ell\in \{\lfloor (N-1)/2\rfloor+1,\ldots,N-2\}$, $2/\ell$ is an eigenvalue but not a variational eigenvalue. 

\end{example}

In  \Cref{thm:k-inequality} we have proved that it holds the inequality $\Lambda_k(\inflap)\leq 1/R_k(\Gc)$, with the equality always holding on trees. However, equality does not always hold for general graphs, as we prove in the next example. 

\begin{example}\label{ex:odd-cycle}
Let $C_N$ be the unweighted cycle graph of order $N$, i.e. with $N$ vertices, uniform weights equal to $1$, i.e. $\edgelength_{uv}=1$ for all $(u,v)\in\edgeset$ and $\nodeweight_u=1$ for all $u\in\nodeset$ and no boundary. Then we state that:
\begin{enumerate}
\item[{claim} 1] If $N=2\ell+1$ with  $\ell\ge2$, then $\Lambda_3(\inflap)<1/R_3(\Gc)$.
\item[{claim} 2] If $N=4\ell+1$ with $\ell\ge 1$, then $\Lambda_{2\ell+1}(\inflap)<1/R_{2\ell+1}(\Gc)$.
\end{enumerate}

Note that the above statements proves that, e.g. on the cycle graph with 5 vertices, we have $\Lambda_3(\inflap)<1/R_3(\Gc)$. Similarly, the cycle graph on 9 vertices satisfies $\Lambda_3(\inflap)<1/R_3(\Gc)$ and $\Lambda_5(\inflap)<1/R_5(\Gc)$, while interestingly  $\Lambda_4(\inflap)=1/R_4(\Gc)$.

\begin{proof}[Proof of claim 1 and claim 2]
We start noting that the cycle graph is dual to itself (see \Cref{Sec:Duality} and Remark 5 in \cite{tudisco2022nonlinear} for details), i.e., the hypergraph obtained as the dual of $\mathcal{C}_N$ is $\mathcal{C}_N$ itself. 
Thus the $\infty$-Laplacian dual, on a cycle graph, is the graph 1-Laplacian.  In particular from \Cref{Thm:duality} we have that on a cycle graph:
\begin{equation}
    \Lambda_k(\inflap)=\Lambda_k(\Delta_1) \qquad \forall k
\end{equation}
Moreover it is possible to check the following equality 
\begin{equation}
h_k(\Gc)=\frac{2}{\lfloor N/k\rfloor}=\frac{1}{R_k(\Gc)} \qquad \forall k.
\end{equation}
Indeed for $h_k(\Gc)$, we can suppose that each $A_i$ in the family of subsets  $(A_1,\cdots,A_k)$ realizing $h_k(\Gc)$ induces a connected subgraph of $C_N$, and thus $c(A_i)=2/|A_i|$. Since $\sum_{i=1}^k|A_i|\le N$, it is easy to see $h_k(\Gc)=2/\lfloor N/k\rfloor.$
Similarly for $R_k(\Gc)$, since for any $k$ vertices in $K_N$ $\sum_{i=1}^k\mathrm{dist}(v_i,v_{i+1})\le  N$, where w.l.o.g. we are assuming that they are ordered increasingly, it is easy to see $2R_k(\Gc)= \lfloor N/k\rfloor.$
To complete the proof, it is sufficient to show that $\Lambda_3(\Delta_\infty)=\Lambda_3(\Delta_1)<h_3(\Gc)=1/R_3(\Gc)$ and $\Lambda_{2\ell+1}(\Delta_\infty)=\Lambda_{2\ell+1}(\Delta_1)\leq h_{2\ell+1}(\Gc)=1/R_{2\ell+1}(\Gc)$, respectively.

Start with the cycle graph $C_{2\ell+1}$ where $\ell\ge2$. It is easy to prove that 
\begin{equation}
    h_2(\Gc)=\frac{2}{\ell}<h_3(\Gc)=\frac{2}{\lfloor\frac{2\ell+1}{3}\rfloor}.
\end{equation}
Then for any $i=1,\cdots,2\ell+1$ consider the subset of nodes 
\begin{equation}
U_i=\{i,i+1,\cdots,i+\ell-1\},
\end{equation}
where with a small abuse of notation we identify the node $2k+1+j$ with the node $j$, for any $j=1,\cdots,\ell$.
Denote by $\chi_{U_i}$ the indicator function of the subset $U_i$,  and note that the cardinality of any $U_i$ is equal to $\ell$.
For $j=0,1,\cdots,2(2\ell+1)-1$, let \[\triangle_{\pm j}=\left\{\left.t_1(-1)^j\chi_{U_{1+j\ell}}+t_2(-1)^{j+1}\chi_{U_{1+(j+1)\ell}}\pm t_3\chi_V\right|\,t_1+t_2+t_3=1,t_1,t_2,t_3\ge0\right\}\] where $\chi_{\nodeset}=(1,1,\cdots,1)$ is  the indicator function of $\nodeset$. 

It is possible to verify that for any $j\in\{0,1,\cdots,2(2\ell+1)-1\}$ and any $f\in \triangle_{\pm j}$,  $\rayl_1(f)\le 2/\ell$. Indeed, considering $\Delta_{+j}$ and assuming w.l.o.g. $j$ odd, some easy computations show that:
\begin{equation}
\begin{aligned}
    &\|\grad f\|_1=|t_1-t_2|+t_1+t_2\leq \max\{2t_1,2t_2\}\leq 2 \\
    &\|f\|_1=|t_1-t_3|\ell+t_3+(t_2+t_3)\ell=|t_1-t_3|k+t_3-t_1\ell+\ell
\geq k
\end{aligned}
\end{equation}
Now consider the centrally symmetric compact subsets $\mathcal{S}:=\bigcup\limits_{j=0}^{2(2\ell+1)-1}(\triangle_{+j}\cup \triangle_{-j})$ and 
\[\mathcal{S}_1:=\bigcup_{j=0}^{2(2\ell+1)-1}\left\{\left.t(-1)^j\chi_{U_{1+j\ell}}+(1-t)(-1)^{j+1}\chi_{U_{1+(j+1)\ell}}\right|\,0\le t\le1\right\}.\]
It is easy to see that $\mathcal{S}_1$ is the union of $2(2\ell+1)$ closed segments and it is odd-homeomorphic to the 1-dimensional sphere  $\mathbb{S}^1$. 
Also, note that $\mathcal{S}$ is the geometric/simplicial join of $\mathcal{S}_1$ and the two-point set  $\{\chi_{\nodeset},-\chi_{\nodeset}\}$. 
By the property of geometric join operator $*$, 
$\mathcal{S}=\mathcal{S}_1*\{\chi_{\nodeset},-\chi_{\nodeset}\}\cong \mathbb{S}^1*\mathbb{S}^0\cong\mathbb{S}^2$, meaning that the centrally symmetric compact set $\mathcal{S}$  is actually a topological sphere of dimension 2. 
Hence, $\gamma(\mathcal{S})=3$.

In particular, $\Lambda_3(\Delta_1)\le \sup\limits_{f\in \mathcal{S}}{\mathcal{R}_1}(f)\le 2/\ell$. Together with the fact that $\Lambda_3(\Delta_1)\ge\Lambda_2(\Delta_1)=h_2(\Gc)=2/\ell$, we finally obtain $\Lambda_3(\Delta_1)=2/\ell=\Lambda_2(\Delta_1)$, where we used the well-known equality $\Lambda_2(\Delta_1)=h_2(\Gc)$, \cref{thm:1-eigenpair_and_cheeger_constants}. 
As a consequence, we get $\Lambda_3(\Delta_1)<h_3(\Gc)$.

Second consider the  cycle graph $C_{4\ell+1}$, as before it is very easy to see $h_{2\ell+1}(\Gc)=2$.  
Then, we take 
\begin{equation}
    f_i=\chi_{\{2i-1,2i\}}+\frac12(-1)^{i-1}\chi_{\{2i\}} \qquad\forall i=1,\cdots,2\ell    
\end{equation}
and consider $X=\mathrm{span}(\chi_{\nodeset},f_1,f_2,\cdots,f_{2\ell})\subset \mathbb{R}^{4\ell+1}$ to be the linear subspace of dimension $2\ell+1$. Easily we have $\gamma(X\cap \mathbb{S}^{4\ell})=2\ell+1$ and to conclude it is sufficient to show that, for any $f\in X\setminus\{0\}$, $\rayl_1(f)<2$. 
Suppose the contrary, i.e. that there exists  $f=t_0 \chi_{\nodeset}+\sum_{i=1}^{2\ell}t_if_i$, for some  nonzero vector $(t_0,t_1,\cdots,t_{2\ell})$, such that $\rayl_1(f)=2$. 
Then,
\begin{equation}\label{eq:R1(f)=1}
f(j)f(j+1)\le0 \;\;\text{ for any }   j=1,\cdots,4\ell+1
\end{equation} where we identify $4\ell+2$ with 1. 
If $t_0=0$, then there exists some $i$ such that $t_i\ne0$, but this implies that 
\begin{equation}
    f(2i-1)f(2i)=t_i^2\Big(1+\frac12(-1)^{i-1}\Big)>0,
\end{equation}
which contradicts \eqref{eq:R1(f)=1}. 
So, it must hold $t_0\ne0$. Without loss of generality, we may assume that $t_0=1$. 

Keeping \eqref{eq:R1(f)=1} in mind, note that 
\begin{equation}
    f(4\ell+1)f(1)=f(1)=1+t_1\le0 \quad  and  \quad f(1)f(2)=(1+t_1)\left(1+\frac32t_1\right)\le0
\end{equation}
thus, we have $t_1=-1$ and $f(2)=-1/2$. 
Similarly, by 
\begin{equation}
f(2)f(3)=-1/2(1+t_2)\le0 \quad \text{and} \quad f(3)f(4)=(1+t_2)\Big(1+\frac12t_2\Big)\le0,
\end{equation}
we can conclude $t_2=-1$. 
Repeating the above process, we can get $t_i=-1$ for any $i=1,\cdots,2\ell$, which means $f=\chi_{\nodeset}-\sum_{i=1}^{2\ell}f_i$. 
However, $f(4\ell)f(4\ell+1)=f(4\ell)=1-(1-1/2)=1/2>0$. 
Therefore, the assumption $\rayl_1(f)=2$ does not hold. 
We have proved that \[\Lambda_{2\ell+1}(\Delta_1)\le\sup_{f\in X\setminus\{0\}}{\rayl_1}(f)=\max_{f\in X\cap \mathbb{S}^{4\ell}}\rayl_1(f)<2=h_{2\ell+1}(\Gc).\]
\end{proof}
\end{example}

\subsection{Relations with independence numbers}
In this section, we observe that expressing the $\infty$-eigenvalues in terms of the radii of inscribed balls, together with bounds in terms of the packing radii of the graph, yields information about the independence numbers of the graph. \textbf{
Throughout this section we assume that the graph has no boundary}. Following
\cite{ABIAD2019indep,FIRBY1997indep,fiol1997eigenvalue,MeirMoon75,Topp91}, the $\ell$-independence number (also referred to as the $\ell$-packing number in \cite{MeirMoon75}) of a graph is the size of the largest set of vertices such that any two vertices in the set are at distance at least $\ell$. 
\begin{definition}\label{k_indip_number}
The $\ell$-independence number of a graph $\Gc$ is defined as:
$$
\alpha_\ell(\Gc)=\max  \{k:\exists 
v_1,\cdots,v_k\in \nodeset \text{ with }\mathrm{dist}(v_i,v_j)\ge \ell,\ \forall i\ne j \},
$$
where $0\leq \ell \leq \mathrm{diam}(\Gc)$. A set $\{v_1,\cdots, v_{\alpha_\ell(\Gc)}\}$ as above is called a maximum $\ell$-independent set.
\end{definition}
Clearly, $\alpha_\ell(\Gc)\geq \alpha_{\ell+1}(\Gc)$ for any $\ell$.

We also note that most of the literature on independence numbers considers only unweighted graphs, i.e., $\edgelength_{uv}=1$ for all edges, and uses the term independence number to denote $\alpha_2(\Gc)=:\alpha(\Gc)$. 

Our interest in these quantities stems from the fact that independence numbers can be directly related to the packing radii of the graph, as follows:
\begin{equation}\label{eq_indip_n_and_packing_radii}
R_k(\Gc)=\max\limits_{\alpha_{\ell}(\Gc)\ge k}\ell/2 \quad \forall k\in\N, 
\qquad 
\alpha_\ell(\Gc)=\max\limits_{2R_k(\Gc)\ge \ell}k 
\quad \forall \ell\leq \mathrm{diam}(\Gc).
\end{equation}
In particular, as a direct consequence of \eqref{eq_indip_n_and_packing_radii}, we obtain the following inequalities:
\begin{equation}\label{inequalities_independence_radius}
    R_{\alpha_{\ell}(\Gc)}(\Gc)\ge \ell/2
    \qquad \text{and} \qquad 
    \alpha_{2R_k(\Gc)}(\Gc)\ge k.
\end{equation}
These inequalities, together with the bounds for the variational $\infty$-eigenvalues in terms of the packing radii of the graph (see \Cref{thm:k-inequality}), readily yield the following proposition.
\begin{proposition}\label{pro:independence}
Let $\Gc$ be a graph without boundary, then for any $0<\ell\leq \mathrm{diam}(\Gc)$ it holds 
$$\Lambda_{\alpha_{\ell}(\Gc)}(\inflap)\le\frac2\ell \qquad \text{and}\qquad \alpha_\ell(\Gc) \le \#\Big\{k:\Lambda_k(\inflap)\le \frac{1}{R_{\alpha_\ell(\Gc)}(\Gc)}\Big\}.$$
Moreover, if $\edgelength_m:=\min\{\edgelength_{uv}|\;(u,v)\in\edgeset\}$, then for any $\ell\geq 2\edgelength_m^{-1}$ we have:
$$\alpha_\ell(\Gc)\leq \#\{k:\Lambda_k(\inflap)\ge \edgelength_m\}.$$
\end{proposition}
\begin{proof}
Both the first two inequalities in the thesis follow by the inequalities:
\begin{equation}
\Lambda_{\alpha_{\ell}(\Gc)}(\inflap)\le \frac{1}{R_{\alpha_{\ell}(\Gc)}(\Gc)}\le \frac{2}{\ell}.  
\end{equation}

To prove the second part of the thesis, it suffices to prove that $\Lambda_{N-\alpha_{\ell}(\Gc)+1}(\inflap)\ge \edgelength_m$ for any $\ell\geq 2\edgelength_m^{-1}$. Let  $\{u_1,\cdots,u_{\alpha_\ell}(\Gc)\}$ be a maximum $\ell$-independent set. 
Denote by $\chi_{u_i}$ the indicator function of the vertex $u_i$, $i=1,\cdots,\alpha(\Gc)$ and let $\pi=\mathrm{span}(\chi_{u_1},\cdots,\chi_{u_{\alpha_\ell(\Gc)}})$. Note that, by definition, $\edgelength_m$ is the reciprocal of the maximal distance between two connected nodes; thus, any two vertices $u_i$ and $u_j$ are not directly connected. As a consequence we note that for any $f\in \pi\setminus\{0\}$, $\rayl_\infty(f)\geq \min_i\{\rayl_\infty(\chi_{u_i})\}\geq \edgelength_m$. 

By the intersection property of the Krasnoselskii's $\mathbb{Z}_2$-genus, for any centrally symmetric subset $S\subset\R^N$ with $\gamma(S)\ge N-\alpha_{\ell}(\Gc)+1$,  $S\cap (\pi\setminus0)\ne\varnothing$, see \cref{Lemma_Krasnoselskii_intersection}. Therefore,
\begin{equation}
\begin{aligned}
\Lambda_{N-{\alpha_\ell(\Gc)}+1}(\inflap)
&=\inf\limits_{ S:\gamma(S)\ge N-{\alpha_\ell(\Gc)}+1}\sup\limits_{g\in S}\mathcal{R}_\infty(g) \\
&\ge \inf\limits_{ S:\gamma(S)\ge N-\alpha_{\ell}(\Gc)+1}\inf\limits_{f\in S\cap (\pi\setminus0)}\mathcal{R}_\infty(f) 
\ge \inf\limits_{f\in \pi\setminus0}\mathcal{R}_\infty(f)=\edgelength_m,
\end{aligned}
\end{equation}
completing the proof.
\end{proof}
Note that in the particular case of an unweighted graph, the canonical independence number $\alpha(\Gc)$ satisfies the upper bound
\begin{equation}
    \alpha(\Gc)\leq \min\left\{\#\Big\{k:\Lambda_k\le \frac{1}{R_{\alpha(\Gc)}(\Gc)}\Big\},\; \#\{k:\Lambda_k\ge 1\}\right\}.
\end{equation}
In the next remark, we observe that the independence numbers can provide bounds for the variational multiplicity of the $\infty$-eigenvalues.

\begin{remark}\label{Rmrk_Multiplicity_independence_numbers}
Recall that from Theorem \ref{thm:k-inequality} we have $\Lambda_2(\inflap)=2/\mathrm{diam}(\Gc)$ and $\Lambda_k(\inflap)\le 1/R_k(\Gc)$ for all $k>1$. Let $m$ denote the variational multiplicity of $\Lambda_2$, that is,
\begin{equation}
\Lambda_2(\inflap)=\dots=\Lambda_{m+1}(\inflap)=\frac{1}{R_2(\Gc)}=\frac{2}{\mathrm{diam}(\Gc)}\,.
\end{equation}
Now consider the $\mathrm{diam}(\Gc)$-independence number $\alpha_{\mathrm{diam}(\Gc)}(\Gc)$. By definition, there exists a maximum $\mathrm{diam}(\Gc)$-independent set $\{v_1,\dots,v_{\alpha_{\mathrm{diam}(\Gc)}(\Gc)}\}$ such that $\mathrm{dist}(v_i,v_j)=\mathrm{diam}(\Gc)$ for all $ i\neq j$. This equality immediately implies
\begin{equation}
    R_{\alpha_{\mathrm{diam}(\Gc)}(\Gc)}(\Gc)=\dots=R_2(\Gc)=\frac{\mathrm{diam}(\Gc)}{2}\,,
\end{equation}
which in turn yields
\begin{equation}
\frac{1}{R_2(\Gc)}=\frac{1}{R_{\alpha_{\mathrm{diam}(\Gc)}(\Gc)}(\Gc)}
\ge \Lambda_{\alpha_{\mathrm{diam}(\Gc)}(\Gc)}(\inflap)
\ge \Lambda_2(\inflap)=\frac{1}{R_2(\Gc)}\,.
\end{equation}
Hence, the $\mathrm{diam}(\Gc)$-independence number provides a lower bound for the multiplicity of $\Lambda_2(\inflap)$, namely,
\begin{equation}
    m\geq \alpha_{\mathrm{diam}(\Gc)}(\Gc)-1\,.
\end{equation}
It is not difficult to see that this bound, which can be transferred to the Krasnosel'skii multiplicity of $\Lambda_2(\inflap)$ (recall that $\gamma\text{-}\mathrm{mult}(\Lambda_2(\inflap))\geq \mathrm{mult}(\Lambda_2(\inflap))$), is generally stronger than the bound provided in \Cref{Prop:gamma_mult_inf_eigenvalues}.
Finally, we note that whenever $\Lambda_k(\inflap)=1/R_k(\Gc)$, repeating the discussion above yields
\[
\mathrm{mult}(\Lambda_k(\inflap))\geq \alpha_{2R_k(\Gc)}(\Gc)-k+1.
\]
\end{remark}

Next, we investigate the relationships between the maximum number of perfect nodal domains associated with an $\infty$-eigenvalue, $\Lambda$, and the $2\Lambda^{-1}$-independence number of the graph. Recall that the higher the number of perfect nodal domains, the stronger the relation between $\Lambda$ and the packing radii of the graph, see \Cref{thm:k-inequality}.

\begin{proposition}\label{prop_k-Indipendece-perfect_nodal_domains-inequalities}
Let $\Gc$ be a graph without boundary and assume $\Lambda$ to be an $\infty$-eigenvalue, then 
$$\max_{f\in X_{\Lambda}}\PNc(f)\leq \alpha_{2/\Lambda}(\Gc).$$
\end{proposition}

\begin{proof}
Assume that $f\in X_{\Lambda}$ is the eigenfunction that realizes the maximum number of perfect nodal domains relative to $\Lambda$, and let $A_1,\cdots, A_p$ be its perfect nodal domains. By definition for any $A_i$ there exists some $v_i\in A_i$ such that $|f(v_i)|=\|f\|_{\infty}$. An easy computation shows that $\dist(v_i,v_j)\geq 2/\Lambda$ for any $i, j$ (see the proof of Theorem \ref{thm:k-inequality} for the details). Thus, by definition of independence number we have the thesis:
\begin{equation}
    \max_{f\in X_{\Lambda}}\PNc(f)\leq \alpha_{2/\Lambda}(\Gc)\,.
\end{equation}
\end{proof}

The result above provides an upper bound on the number of perfect nodal domains; however, \Cref{thm:k-inequality} suggests that lower-bounding the maximum number of perfect nodal domains in terms of the independence numbers would be more interesting. Unfortunately, this is not generally possible, as we see in the next example

\begin{example}\label{example_indep_perfect_nodal_dom}
Consider the star graph in \Cref{fig:star_graph}.
\begin{figure}
\centering
    \begin{tikzpicture}[inner sep=0.5mm, scale=.5, thick]
    \node (1) at (0,0) [circle,draw] {$v_1$};
    \node (2) at (-4,2) [circle,draw] {$v_2$};
    \node (3) at (-5,-2) [circle,draw] {$v_3$};
    \node (4) at (5,0) [circle,draw] {$v_4$};
\draw[color=black] (1) --node[left, sloped,xshift=10pt,yshift=4pt]{\tiny $\edgelength_{12}=2$} (2); 
\draw[color=black] (1) --node[left, sloped,xshift=8pt,yshift=4pt]{\tiny $\edgelength_{13}=1$} (3);
\draw[color=black] (1) --node[left, sloped,xshift=12pt,yshift=4pt]{\tiny $\edgelength_{14}=3/2$} (4);
%
\end{tikzpicture}~~
    \caption{Graph discussed in \cref{example_indep_perfect_nodal_dom},  where the inequality     $\max_{f\in X_{\Lambda}}\PNc(f)\leq \alpha_{2/\Lambda}(\Gc)$, between maximum number of perfect nodal domains associated to an $\infty$-eigenvalue $\Lambda$ and the $(2/\Lambda)$-independence number, is strict. }
    \label{fig:star_graph}
\end{figure}
From \Cref{prop:Infinity_eigenfunctions_characterization} we have that $\Lambda=2$ is an $\infty$-eigenvalue with eigenfunction $f$
\begin{equation}
    f(v_1)=1, \quad f(v_2)=0, \quad f(v_3)=-1, \quad f(v_4)=-1/3. 
\end{equation}
Clearly the maximum number of perfect nodal domains corresponding to $\Lambda=2$ is $2$. Indeed, since the path connecting $v_1$ to $v_3$ is the only one having length $\Lambda/2=1$, any eigenfunction $f$ corresponding to $\Lambda=2$ has to satisfy $f(v_1)=-f(v_3)$. Moreover, the edges $(v_1,v_2)$ and $(v_1,v_4)$ are shorter than $(v_1,v_3)$, and both $|\grad f(v_1,v_2)|$ and $|\grad f(v_1, v_4)|$ have to be smaller than $|\grad f(v_1,v_3)|$, so $f(v_2)\geq 0$ and $f(v_4)\geq -1/3$. On the other hand, simply taking the nodes $\{v_2,\:v_3,\: v_4\}$, it is not difficult to observe that $\alpha_{2/\Lambda}(\Gc)=\alpha_1(\Gc)=3\geq 2$.   
\end{example}

Nonetheless, there are particular cases in which it is actually possible to lower bound the maximum number of perfect nodal domains in terms of the independence number. Below we study the case of unweighted graphs. Note that in many cases the independence numbers are used only to derive information about the topology of the graph and thus the weights that represent the lengths can be neglected.

\begin{proposition}\label{prop_k-Indipendece-perfect_nodal_domains-inequalities_unweighted_case}
Assume $\Gc$ to be an unweighted graph without boundary, i.e. $\edgelength_{uv}=1$ for any $(u,v)\in \edgeset$, and let $\Lambda=2/l$ be an $\infty$-eigenvalue, then  
\begin{itemize}
    \item If $l$ is even
$$\max_{f\in X_{\Lambda}}\PNc(f)= \alpha_{l}(\Gc).$$
\item If $l$ is odd
$$\alpha_{l}(\Gc)-\betti\leq\max_{f\in X_{\Lambda}}\PNc(f)\leq \alpha_{l}(\Gc),$$
where $\betti$ denotes the number of independent cycles of $\Gc$\,.
\end{itemize}
\end{proposition}

\begin{proof}
 The upper bounds are a direct consequence of \Cref{prop_k-Indipendece-perfect_nodal_domains-inequalities}.
 Then note that from \Cref{EX:inf_noweights_spectrum}, $l\in\{1,\dots,\mathrm{diam}(\Gc)\}$ and $f$ is an eigenfunction of $\Lambda$ if and only if it satisfies the $\mathrm{SP}_l$ condition  discussed in \Cref{EX:inf_noweights_spectrum}. 
 Thus now we define an appropriate function, $f$, which satisfies the $\mathrm{SP}_l$ condition and which has $\alpha_l(\Gc)$ nodal domains if $l$ is even, or more than $\alpha_{l}(\Gc)-\betti$ perfect nodal domains when $l$ is odd.
Start from \cref{k_indip_number} and let $v_1,\dots,v_{\alpha(\Gc)}$ be $\alpha_{l-1}(\Gc)$ nodes such $\dist(v_i,v_j)\geq l\;\forall i\neq j$\,. Observe that without loss of generality we can assume $\dist(v_1,v_2)=l$\,.
Then, for any $i=1,\dots,\alpha(\Gc)$, consider 
\begin{equation}
A_i=\left\{v\in \nodeset\;\Big|\; \dist(v,v_i)<\frac{l}{2} \right\}\,,
\end{equation}
clearly $A_i\cap A_j=\varnothing$.
Now, first let $l$ be even and define the function
\begin{equation}
\begin{aligned}
f(v)=&\max\Big\{1-\frac{2}{l}\dist(v_1,v),0\Big\}-\max\Big\{1-\frac{2}{l}\dist(v_2,v),0\Big\}\\
&+\sum_{i=3}^{\alpha(\Gc)} \max\Big\{1-\frac{2}{l}\dist(v_i,v),0\Big\}\,,
\end{aligned}
\end{equation}
since $\dist(v_1,v_2)=l$, $f$ is easily proved to satisfy the $\mathrm{SP}_l$ condition. By construction, $f(u)\neq 0$ if and only if $u\in\cup_{i=1}^{\alpha(\Gc)} A_i$, and $f(v_i)=1=\|f\|_{\infty}\;\forall\;i=1,\dots,\alpha(\Gc)$. 
Moreover, since $l$ is even, for any pair $v_i,v_j\;i\neq j=1,\dots,\alpha(\Gc)$, on the shortest path that joins $v_i$ and $v_j$, there is at least one node $v$ such that $\dist(v,v_i)\geq\frac{l}{2}$,  $\dist(v,v_j)\geq\frac{l}{2}$, implying $f(v)=0$.
The last fact implies that $f$ has exactly $\alpha_{l}(\Gc)$ perfect nodal domains. 

Consider now the case $l$ odd.
Define a graph $\Gc'$ with $\alpha_l(\Gc)$ nodes, $v_1',\dots,v_{\alpha_l(\Gc)}'$, such that $v_i'\sim v_j'$ if there exist and edge $(u,w)\in E(\Gc)$ such that $u\in A_i$ and $w\in A_j$. Observe that, by assumption $v_1'\sim v_2'$. Moreover, letting $\beta(\Gc')$ be the number of loops of the new graph $\Gc'$, we get $\beta(\Gc')\leq \betti$. Now, first, remove $\beta(\Gc')$ edges from $\Gc'$ in such a way to reduce the last to a forest, $\Gc''$ where it still holds $v_1'\sim v_2'$. Then, since $\Gc''$ is a bipartite graph, define a sign function $\mathrm{sign}(v_i')=\pm 1\;\forall i=1,\dots,\alpha(\Gc)$, such that, if $\mathrm{sign}(v_i')=+1$, $\mathrm{sign}(v_j')=-1$ for any $v_j'\sim v_i'$ in $\Gc''$.
Observe that adding back the removed edges to $\Gc''$, the $\mathrm{sign}$ function induces, on $\Gc'$, at most $\alpha_l(\Gc)-\beta(\Gc')\geq \alpha_l(\Gc)-\betti$ nodal domains.
Finally, on the graph $\Gc$, define the function 
\begin{equation}
f(v)=\sum_{i=1}^{\alpha(\Gc)} \mathrm{sign}(v_i') \max\Big\{1-\frac{2}{l}\dist(v_i,v),0\Big\}\,.
\end{equation}
Observe that, since $\dist(v_1,v_2)=l$ and $\mathrm{sign}(v_1')\mathrm{sign}(v_2')=-1$, $f$ satisfies the $\mathrm{SP}_l$ condition, being, thus, an eigenfunction associated to $\Lambda$.
Moreover, since  $A_i\sim A_j$ if and only if $v_i'\sim v_j'$, we can state that $\mathrm{sign}(f|_{A_i})=\mathrm{sign}(v_i')$, $|f(v_i)|=1=\|f\|_{\infty}$ and $\PNc(f)=\SNc(f)=\SNc(\mathrm{sign})\geq  \alpha_{l}(\Gc)-\betti$, which concludes the proof.
\end{proof}
The last \Cref{prop_k-Indipendece-perfect_nodal_domains-inequalities_unweighted_case} with the inequalities from \Cref{thm:infty_nodal-count} leads to bound, in the unweighted case, the quantities $\alpha_{l}(\Gc)$ in terms of the position of $\Lambda=2/l$ in the variational $\infty$-spectrum.
\begin{corollary}
In the same hypothesis of \Cref{prop_k-Indipendece-perfect_nodal_domains-inequalities_unweighted_case}, consider any $l\in\{1,\dots,diam(\Gc)\}$ and assume that $2/l=\Lambda<\Lambda_{k+1}(\inflap)$, then 
$$\alpha_{l}(\Gc)\leq k+\betti\,.$$
\end{corollary}
\begin{proof}
It is a direct consequence of \Cref{prop_k-Indipendece-perfect_nodal_domains-inequalities_unweighted_case} and \Cref{thm:infty_nodal-count}.
\end{proof}

\subsection{Relations with matching number}
In this section, we investigate a further connection between topological properties of the graph and the $\infty$-Laplacian spectrum. As in the previous paragraph \textbf{we assume the boundary to be empty and concentrate on the unweighted case, i.e. $\edgelength_{uv}=1$ for all $(u,v)\in\edgeset$} and $\nodeweight_u=1$ for all $u\in\nodeset$. 
We relate the multiplicity of the maximum eigenvalue of the $\infty$-Laplacian to the matching number of the graph. A matching is a set of pairwise non-adjacent edges; that is, no two edges share common vertices. A maximum matching is a matching that contains the largest possible number of edges. Then the matching number is defined as follows \cite{West}:
\begin{definition}\label{Matching_number}
The \textbf{matching number}, $\matchingN$ of a graph $\Gc$ is defined as the size of a maximum matching.    
\end{definition}
\begin{proposition}\label{pro:matching-2}
Let $\Gc$ be an unweighted graph without boundary. Then it holds the equality $\Lambda_{N}(\inflap)=2$ with $N=|\nodeset|$ and
$$\matchingN\le \mathrm{mult}(2)\le  \gamma\text{-}\mathrm{mult}(2) 
\le 2\matchingN.$$
Moreover, if $\Gc$ is a connected bipartite graph with the two parts $\nodeset=A\cup B$, and  $|A|\le|B|$, then $\mathrm{mult}(2)= \gamma\text{-}\mathrm{mult}(2)=|A|$. Hence, the equality $\matchingN= \mathrm{mult}(2)=  \gamma\text{-}\mathrm{mult}(2)$ holds on any connected bipartite graph. 
\end{proposition}
\begin{proof}
Let $v_1v_2$, $v_3v_4$, $\cdots$, $v_{2\matchingN-1}v_{2\matchingN}$ be $\matchingN$ edges that realize the maximum matching number. 
For each $i=1,\cdots,\matchingN$, let $f_i$ be defined as 
\begin{equation}
    f_i(v)=\begin{cases}
1,&\text{ if }v=v_{2i-1},\\
-1,&\text{ if }v=v_{2i},\\
0,&\text{ otherwise}.\\
\end{cases}
\end{equation}
Let $\pi=\mathrm{span}(f_1,\cdots,f_\matchingN)$. Then, $\dim(\pi)=\matchingN$. And every nonzero $f\in \pi$ is an eigenfunction corresponding to the largest eigenvalue, due to the fact  $\mathcal{R}_\infty(f)=2$. 
Note that for any centrally symmetric subset $S\subset\R^{N}$ with $\gamma(S)\ge N-\matchingN+1$, $S\cap (\pi\setminus0)\ne\varnothing$. 
So, 
\begin{equation}
\begin{aligned}
    \Lambda_{N-\matchingN+1}(\inflap)=
\inf\limits_{ \gamma(S)\ge N-\matchingN+1}\sup\limits_{g\in S}\mathcal{R}_\infty(g) \ge
\inf\limits_{f\in \pi\setminus0}\mathcal{R}_\infty(f)=2
\end{aligned}
\end{equation}
which implies $\Lambda_{N-\matchingN+1}(\inflap)=\cdots=\Lambda_N(\inflap)=2$. 
In consequence, the multiplicity of the maximum eigenvalue $2$ is at least $\matchingN$, $\matchingN\le \mathrm{mult}(2)$.

Next we prove the right hand side inequality in the thesis.
We consider $M:=\{f\ne0: \mathcal{R}_\infty(f)=2\}$ and claim that $\gamma(M)\le N-\alpha_2(\Gc)$. Suppose that the claim was untrue, i.e. that $\gamma(M)\ge N-\alpha_2(\Gc)+1$. Then, let  $\{u_1,\cdots,u_{\alpha_2(\Gc)}\}$ be a maximum independent set, denote by $\chi_{u_i}$ the indicator function on any vertex $u_i$, and consider $\pi'=\mathrm{span}(\chi_{u_1},\cdots,\chi_{u_{\alpha(\Gc)}})$. 
By \Cref{Lemma_Krasnoselskii_intersection}, it can be observed that $M\cap \pi'\setminus0\ne \varnothing$, however $\rayl_\infty(g)=1$ for any $g\in \pi'$, yielding a contradiction. We have proved the claim $\gamma(M)\le N-\alpha_2(\Gc)$.

Moreover, it is well known that $N\le \alpha_2(\Gc)+2\matchingN$. Indeed, if $\{v_1,v_2\}$, ..., $\{v_{2\matchingN-1},v_{2\matchingN}\}$ is a maximum matching, $U:=\nodeset\setminus \{v_1,\cdots,v_{2\matchingN}\}$ must be an independent set (otherwise, there is an edge $\{u_1,u_2\}$ in $U$, which contradicts the maximality assumption). 

From the inequalities proved above, we obtain 
\begin{equation}
    \gamma(M)\le N-\alpha_2(\Gc)\le2\matchingN,
\end{equation}
meaning that the $\gamma$-multiplicity of the  maximum eigenvalue $2$ is at most $2\matchingN$. 

Finally, we focus on bipartite graphs. Suppose $A$ and $B$ are two parts of $\nodeset$ such that every edge has one vertex in $A$ and the other vertex in $B$. 
Since $|A|\le |B|$, it is easy to check that $\alpha_2(\Gc)=|B|$ and $\matchingN=|A|$. 
According to the above proof, we have 
\[|A|=\matchingN\le \mathrm{mult}(2)\le \gamma\text{-}\mathrm{mult}(2)\le N-\alpha_2(\Gc)=N-|B|=|A|\]
which yields  $\mathrm{mult}(2)= \gamma\text{-}\mathrm{mult}(2)=|A|$.
\end{proof}

\begin{remark}
  We note that there are also examples where the right-hand-side equality in \cref{pro:matching-2} holds. This is, e.g., the case of the unweighted triangle graph, where we have $\mathrm{mult}(2)=  \gamma\text{-}\mathrm{mult}(2) 
= 2\matchingN$. Indeed, it can be easily observed that $\matchingN=1$ and $2=1/R_2(\Gc)=\Lambda_2(\Delta_\infty)\leq \Lambda_3(\Delta_\infty)=2$, which yields the desired equality.
\end{remark}

We remark here that the multiplicity of the largest eigenvalue of the graph 1-Laplacian is bounded from below and above by $\alpha_2(\Gc)$ and $2\alpha_2(\Gc)$, respectively, see \Cref{Subsec:1_lap_indep_numb}. It is very interesting to note that Proposition \ref{pro:matching-2} analogously states that the multiplicity of the largest eigenvalue of the graph $\infty$-Laplacian is bounded by $\matchingN$ and $2\matchingN$ from below and above, respectively. 

\section{Other relevant topics}

\subsection{Minimal $k$-partition}\label{sec:spectral-minmial-partition}

In this section, by means of the $p$-Laplacian eigenvalue problem with homogeneous Dirichlet boundary conditions, we introduce the notion of spectral minimal $k$-partition. This problem was originally studied in \cite{Conti05,Helffer09} for both linear and nonlinear spectral minimal partition problems, in \cite{BobkovParini18} for the $p$-Laplacian on bounded open domains, and in \cite{miclo2008eigenfunctions} for the Laplacian on metric graphs. In particular, in \cite{BobkovParini18} the authors characterized the higher-order Cheeger constants as the limit, as $p\to 1$, of the spectral minimal partition cost of the $p$-Laplacian.
With the aim of relating the packing radii of the graph to the spectral minimal partition costs of the $\infty$-Laplacian, we introduce the following notion of $p$-Laplacian disjoint spectral minimal $k$-partition cost:
\begin{equation}\label{p_lap_disjoint_spec_min_partition}
    \widehat{\lambda_k}(\Delta_p):=\min\limits_{\text{disjoint }V_1,\ldots,V_k\subset \internalnodes}\;
    \max_i \lambda_1\big(\Delta_p,\Gc(V_i)\big).
\end{equation}
Here, $\lambda_1\big(\Delta_p, \Gc(V_i)\big)$ denotes the smallest eigenvalue of the Dirichlet $p$-Laplacian on the graph $\Gc(V_i)$ induced by $V_i$, that is, the graph having $V_i$ as interior nodes and $\nodeset\setminus V_i$ as boundary nodes.
Analogously, we define the $p$-Laplacian non-adjacent spectral minimal $k$-partition cost as
\begin{equation}\label{p_lap_nonadj_spec_min_partition} 
\widetilde{\lambda_k}(\Delta_p):=\min\limits_{\text{non-adjacent }V_1,\ldots,V_k\subset \internalnodes}\;
\max_i \lambda_1\big(\Delta_p, \Gc(V_i)\big).
\end{equation}
In particular, in the case $p=\infty$, we write
\begin{equation}\label{Inf_lap_spec_min_partition}
\begin{aligned}
\widehat{\Lambda}_k(\Delta_{\infty})&:= \min\limits_{\text{disjoint }V_1,\ldots,V_k}
\max_i \Lambda_1\big(\Delta_\infty,\Gc(V_i)\big),\\
\widetilde{\Lambda}_k(\Delta_{\infty})&:=\min\limits_{\text{non-adjacent }V_1,\ldots,V_k}
\max_i \Lambda_1\big(\Delta_{\infty}, \Gc(V_i)\big),
\end{aligned}
\end{equation}
where $\Lambda_1\big(\Delta_\infty, \Gc(V_i)\big)$ denotes the smallest eigenvalue of the Dirichlet $\infty$-Laplacian on the graph $\Gc(V_i)$ with homogeneous boundary conditions on $\nodeset\setminus V_i$. These quantities represent the disjoint and non-adjacent $\infty$-Laplacian spectral minimal $k$-partition costs, respectively.
First, we show that the $p$-Laplacian spectral minimal partition costs converge to their $\infty$-Laplacian counterparts as $p\to\infty$.
\begin{lemma}
\label{lemma:limit1} Let $\widehat{\lambda_k}(\plap), \widehat{\Lambda_k}(\Delta_\infty)$, $\widetilde{\lambda_k}(\plap)$ and $ \widetilde{\Lambda_k}(\inflap)$ be defined as above. 
Then
$$\lim_{p\rightarrow \infty}\widehat{\lambda_k}(\Delta_p)^{\frac1p}=\widehat{\Lambda_k}(\inflap)
\qquad \text{and}\qquad
\lim_{p\rightarrow \infty}\widetilde{\lambda_k}(\Delta_p)^{\frac1p}=\widetilde{\Lambda_k}(\inflap) \qquad \forall \;k=1,\dots,|\nodeset|$$
\end{lemma}
\begin{proof}
    The proof easily follows from Lemma \ref{lemma:p-monotonic}. Indeed for any $V_i\subset \internalnodes$, 
    \begin{equation}
    \begin{aligned}
    |\nodeweight|^{-\frac1p}\Lambda_1\big(\inflap,\Gc(V_i)\big)&\leq |\nodeweight(V_i)|^{-\frac1p}\Lambda_1\big(\inflap,\Gc(V_i)\big)   \leq\lambda_1^{\frac1p}\big(\Delta_p, \Gc(V_i)\big)\\
&\leq |E(V_i)\big|^{\frac1p}\Lambda_1\big(\inflap,\Gc(V_i)\big)\leq|E\big|^{\frac1p}\Lambda_1\big(\inflap, \Gc(V_i)\big)
    \end{aligned}
    \end{equation}
where $E(V_i)$ is the number of edges incident to some node in $V_i$ and $|\nodeweight(V_i)|=\sum_{v\in V_i}\nodeweight_{v}$.
It follows that for any family of disjoint subsets $V_1,\dots,V_k\subset \internalnodes$,
\begin{equation}
|\nodeweight|^{-\frac1p}\max_{i=1,\dots,k}\Lambda_1\big(\inflap, \Gc(V_i)\big)\leq \max_{i=1,\dots,k}\lambda_1^{\frac1p}\big(\Delta_p, \Gc(V_i)\big) 
\leq|E\big|^{\frac1p}\max_{i=1,\dots,k}\Lambda_1\big(\inflap, \Gc(V_i)\big)
\end{equation}
Taking the minimum over all the possible families of disjoint subsets and than the limit leads to the desired equality

\begin{align}\widehat{\Lambda_k}(\inflap)&=\min_{\substack{\text{disjoint}\\ V_1,\dots,V_k}}\max_{i=1,\dots,k}\Lambda_1\big(\inflap,\Gc(V_i)\big)\leq \liminf_{p\rightarrow\infty} \min_{\substack{\text{disjoint}\\ V_1,\dots,V_k}} \max_{i=1,\dots,k}\lambda_1^{\frac1p}\big(\Delta_p, \Gc(V_i)\big)= \\
&=\liminf_{p\rightarrow\infty} \widehat{\lambda_k}(\Delta_p)^{\frac1p}
\leq \limsup_{p\rightarrow\infty} \widehat{\lambda_k}(\Delta_p)^{\frac1p}=\\
&=  \limsup_{p\rightarrow\infty} \min_{\substack{\text{disjoint} V_1,\dots,V_k}} \max_{i=1,\dots,k}\lambda_1^{\frac1p}\big(\Delta_p, \Gc(V_i)\big)\leq \\
&\leq \min_{\substack{\text{disjoint}\\ V_1,\dots,V_k}}\max_{i=1,\dots,k}\Lambda_1\big(\inflap,\Gc(V_i)\big)=\widehat{\Lambda_k}(\inflap)\,.
\end{align}
The same argument yields the equality $\widetilde{\Lambda_k}(\inflap)=\lim_{p\rightarrow\infty}\widetilde{\lambda_k}(\Delta_p)^{\frac1p}\,.$
\end{proof}

Now we can  prove the following result that relates the spectral minimal partition costs of the $\infty$-Laplacian to the 
packing radii. In particular, as a corollary, we show that the $k$-th packing radius on unweighted graphs can be written as the arithmetic mean of the disjoint and non-adjacent $\infty$-Laplacian spectral $k$-partition costs.
\begin{theorem}\label{Thm_spectral_min_partitions}
Let $\edgelength_M:=\min_{(u,v)\in\edgeset} \edgelength_{(u,v)}$ be the reciprocal of the length of the longest edge of a graph $\Gc$. Then the $\infty$-Laplacian spectral minimal partition costs and the packing radii of $\Gc$ satisfy the following inequalities:  
    $$
    R_k(\Gc) 
 \leq \frac{1}{\widehat{\Lambda}_k(\inflap)}\leq R_k(\Gc)+\frac{1}{2\edgelength_M} 
    \qquad \text{and} \qquad 
 R_k(\Gc)-\frac{1}{2\edgelength_M}\leq\frac{1}{\widetilde{\Lambda}_k(\inflap)}\leq R_k(\Gc).
    $$
\end{theorem}
\begin{proof}
    Start by recalling the definition of $R_k(\Gc)$:
    \begin{equation}
    R_k(\Gc)=\max_{v_1,\dots,v_k\in \internalnodes}\min_{i,j=1,\dots,k} \frac{d(v_i,v_j)}{2}\,.
    \end{equation}
Then let $v_1,\dots,v_k$ be a maximizing set of nodes in the definition of $R_k(\Gc)$, and for any $i=1,\dots,k$, let
\begin{equation}
\widehat{V}_i=\{u\in \nodeset\,|\;d(u,v_i)< R_k(\Gc)\} \quad \text{and} \quad \widetilde{V}_i=\{u\in \nodeset\,|\;d(u,v_i)< R_k(\Gc)-1/(2\edgelength_M)\}\,.
\end{equation}
By construction, the families of sets $\{\widehat{V}_i\}$ and $\{\widetilde{V}_i\}$ are suitable sets in the definitions, respectively, of $\widehat{\Lambda}_k(\inflap)$ and $\widetilde{\Lambda}_k(\inflap)$. Hence, from the characterization of the first $\infty$-eigenvalue, see \Cref{thm:k-inequality}, we observe
\begin{equation}
\frac{1}{\widehat{\Lambda}_k(\inflap)}\geq \min_i \frac{1}{\Lambda_1(\inflap,\widehat{V}_i)}=\min_i \max_{u\in \widehat{V}_i}d_{\nodeset\setminus \widehat{V}_i}(u)\geq R_k(\Gc)\,.
\end{equation}
Analogously we observe that, for any $i$, 
\begin{equation}
\frac{1}{\widetilde{\Lambda}_k(\inflap)}\geq  R_k(\Gc)- \frac{1}{2\edgelength_M}\,.
\end{equation}

We miss to prove the upper bounds, let $V_1,\dots V_k$ be a maximizing family of sets in the definition of $\widehat{\Lambda}_k(\inflap)$\,, i.e.
\begin{equation}
\widehat{\Lambda}_k(\inflap)=\max_i \frac{1}{\max_{u\in V_i}d_{\nodeset\setminus V_i}(u)}\,,
\end{equation}
then define $v_i^{max}=\argmax_{u\in V_i}d_{\nodeset\setminus V_i}(u)$ and $\partial V_i:=\{u\in V_i\;| \;\exists (u,v)\in E\; \text{with}\; v\in \nodeset\setminus V_i\}$\,.
Note that any path connecting $v_i^{max}$ to $v_j^{max}$ has to be incident both to $\partial V_i$ and $\partial V_j$, thus
\begin{align}\label{eq:1_Thm_spectral_min_partitions_a}
d(v_i^{max},v_j^{max})\geq d(v_i^{max},\partial V_j) +d(\partial V_j,v_j^{max})\geq \frac{1}{\widehat{\Lambda}_k(\inflap)}+\frac{1}{\widehat{\Lambda}_k(\inflap)}-\frac{1}{\edgelength_M}
\end{align}
where we have used that for any $i$,
$d(v_i^{max},\nodeset\setminus V_i)\leq d(v_i^{max},\partial V_i)+1/\edgelength_M$ and $d(v_i^{max},\partial V_j)\geq d_{\nodeset\setminus V_i}(v_i^{max}) $.
Clearly \eqref{eq:1_Thm_spectral_min_partitions_a} implies that 
\begin{equation}
 R_k(\Gc)\geq \frac{1}{\widehat{\Lambda}_k(\inflap)}-\frac{1}{2\edgelength_M}.  
\end{equation}
The same proof, assuming $V_1,\dots, V_k$ to be a maximizing family of subsets in the definition of $\widetilde{\Lambda}_k$, leads to  
\begin{equation}
    R_k(\Gc)\geq \frac{1}{\widetilde{\Lambda}_k(\inflap)},
\end{equation}
in this case, in fact, the subsets $\{V_i\}$ are non-adjacent and thus $d(v_i^{max},v_j^{max})\geq d(v_i^{max},\nodeset\setminus V_i)+ d(v_j^{max},\nodeset\setminus V_j)$.
\end{proof}

\begin{corollary}\label{cor_unweighted_spectral_min_partitions}
    In the same hypotheses of \Cref{Thm_spectral_min_partitions}, assuming also that the graph $\Gc$ is unweighted, we have:
     $$R_k(\Gc)=\frac12\Big(\frac{1}{\widehat{\Lambda}_k(\inflap)}+\frac{1}{\widetilde{\Lambda_k}(\inflap)}\Big)\,.$$
\end{corollary}
\begin{proof}
Defining the families $\{\widehat{V_i}\}$ and $\{\widetilde{V_i}\}$ as in the proof of $\Cref{Thm_spectral_min_partitions}$, it is possible to note that for any $i$,
$\max_{u\in \widehat{V}_i}d_{\nodeset\setminus \widehat{V}_i}(u)=\lfloor R_k(\Gc)+1/2\rfloor$. Indeed every edge has length $1$, i.e. the distance is an integer and $R_k(\Gc)\in \{l/2\,|\;l=1,\dots,N\}$.
The last equality yields the lower bound
\begin{equation}
    \frac{1}{\widehat{\Lambda}_k(\inflap)}\geq \Big\lfloor R_k(\Gc)+\frac{1}{2} \Big\rfloor\,.
\end{equation}
Analogously, we observe that, for any $i$, $\max_{u\in \widetilde{V}_i}d_{\nodeset\setminus \widetilde{V}_i}(u)=\lfloor R_k(\Gc)\rfloor$, which leads to 
\begin{equation}
\frac{1}{\widetilde{\Lambda}_k(\inflap)}\geq\lfloor R_k(\Gc)\rfloor\,.
\end{equation}
Finally combining the two lower bounds and upper bounds from \Cref{Thm_spectral_min_partitions} we get 
\begin{equation}
2R_k(\Gc)= \Big\lfloor R_k(\Gc)+\frac{1}{2} \Big\rfloor+\lfloor R_k(\Gc) \rfloor=2R_k(\Gc)     \leq \frac{1}{\widehat{\Lambda}_k(\inflap)}+\frac{1}{\widetilde{\Lambda}_k(\inflap)} \leq  2R_k(\Gc)+\frac{1}{2}\,,
\end{equation}
here, however, both $1/\widehat{\Lambda}_k(\inflap)$ and $1/\widetilde{\Lambda}_k(\inflap)$ and $2R_k(\Gc)$ are distances among nodes and thus they are integers, concluding the proof.
\end{proof}

A similar result relates the $2k$-packing radius of the graph to the $(k,2)$-partitions costs $\widehat{\Lambda}_{k,2}(\inflap)$ and $\tilde{\Lambda}_{k,2}(\inflap)$ defined by:
\begin{equation}
\begin{aligned}
\widehat{\Lambda}_{k,2}(\inflap)&:= \min\limits_{\text{disjoint }V_1,\ldots,V_k}\max_i \Lambda_2\big(\inflap, \Gc(V_i)\big)\\ \widetilde{\Lambda}_{k,2}(\inflap)&:=\min\limits_{\text{non-adjacent }V_1,\ldots,V_k}\max_i \Lambda_2\big(\inflap, \Gc(V_i)\big).
\end{aligned}
\end{equation}

\begin{theorem}\label{Thm_spectral_min_partitions_2}
Let $\edgelength_M:=\min_{(u,v)\in\edgeset} \edgelength_{(u,v)}$ be the reciprocal of the length of the longest edge of $\Gc$. Then the minimal $(k,2)$-partition costs satisfy the following inequalities:  
$$ R_{2k}(\Gc)
 \leq \frac{1}{\widehat{\Lambda}_{k,2}(\inflap)}\leq 
    R_{2k}(\Gc)+\frac{1}{2\edgelength_M} 
    \quad \text{and} \quad 
R_{2k}(\Gc)-\frac{1}{2\edgelength_M}\leq\frac{1}{\widetilde{\Lambda}_{k,2}(\inflap)}\leq R_{2k}(\Gc)\,.$$  
\end{theorem}
\begin{proof}
    As in the proof of the previous Theorem \ref{Thm_spectral_min_partitions}, let $v_1,\dots,v_{2k}$ be a maximizing set of nodes in the definition of $R_{2k}(\Gc)$. Then for any $i=1,\dots,k$, define the following subset of nodes:
\begin{equation}
    \begin{aligned}
    &\widehat{V}_i:=\Big\{u\in \nodeset\,|\;\min\{d(u,v_{2i-1}), d(u,v_{2i})\} < R_{2k}(\Gc)\Big\}\,,\\ 
    &\widetilde{V}_i:=\Big\{u\in \nodeset\,|\;\min\{d(u,v_{2i-1}), d(u,v_{2i})\}< R_{2k}(\Gc)-\frac{1}{2\edgelength_M} \Big\}\,.
    \end{aligned}
\end{equation}
By construction, the families of sets $\{\widehat{V}_i\}$ and $\{\widetilde{V}_i\}$ are suitable sets in the definitions, respectively, of $\widehat{\Lambda}_{k,2}(\inflap)$ and $\widetilde{\Lambda}_{k,2}(\inflap)$. 
Hence, from the characterization of the second $\infty$-eigenvalue \Cref{thm:k-inequality} and the definition of $\widehat{V}_i$, we observe
\begin{equation}
\frac{1}{\widehat{\Lambda}_{k,2}(\inflap)}\geq \min_i \frac{1}{\Lambda_2\big(\inflap,\Gc(\widehat{V}_i)\big)}=\min_i R_2\big(\Gc(\widehat{V_i})\big)\geq R_{2k}(\Gc)\,,
\end{equation}

where we are using $R_2\big(\Gc(\widehat{V_i})\big)$ to denote the second packing radius of the subgraph of $\Gc$ having $\widehat{V_i}$ as internal nodes and $\nodeset\setminus\widehat{V_i}$ as boundary nodes.
Analogously, since for any $i$ we can observe that by definition $R_2\big(\Gc(\widetilde{V}_i)\big)\geq R_{2k}(\Gc)-1/(2\edgelength_M)$, we can conclude  
\begin{equation}
\frac{1}{\widetilde{\Lambda}_k(\inflap)}\geq R_{2k}(\Gc) -\frac{1}{2\edgelength_M}\,.
\end{equation}
To prove the upper bounds, let $\widehat{V}_1,\dots \widehat{V}_k$ be a maximizing family of sets in the definition of $\widehat{\Lambda}_{k,2}(\inflap)$\,, i.e.
\begin{equation}
\widehat{\Lambda}_{k,2}(\inflap)=\max_i \frac{1}{R_2\big(\Gc(\widehat{V}_i)\big)}\,.
\end{equation}
Given the definition of $R_2\big(\Gc(\widehat{V}_i)\big)$, for any $i=1,\dots,k$ we can take two nodes $v_{2i-1},v_{2i}\in \widehat{V}_i$ such that $d(v_{2i-1},v_i)\geq 2/\widehat{\Lambda}_{k,2}(\inflap)$, $d(v_{2i-i}, \nodeset\setminus \widehat{V}_i)\geq 1/\widehat{\Lambda}_{k,2}(\inflap)$, $d(v_{2i}, \nodeset\setminus \widehat{V}_i)\geq 1/\widehat{\Lambda}_{k,2}(\inflap)$\,.
Then, the same discussion presented in the proof of Theorem \ref{Thm_spectral_min_partitions} allows us to observe
\begin{equation}\label{eq:1_Thm_spectral_min_partitions_2}
d(v_i,v_j)\geq \frac{2}{\widehat{\Lambda}_{k,2}(\inflap)}-\frac{1}{\edgelength_M} \quad \forall\: i,j=1,2k \,.
\end{equation}
Finally, the inequality \eqref{eq:1_Thm_spectral_min_partitions_2} yields $R_{2k}(\Gc)\geq 1/\widehat{\Lambda}_{k,2}(\inflap)-1/(2\edgelength_M)$.
The same proof, assuming $V_1,\dots, V_{2k}$ to be a maximizing family of subsets in the definition of $\widetilde{\Lambda}_{2,k}(\inflap)$, yields $R_{2k}(\Gc)\geq 1/\widetilde{\Lambda}_{k,2}(\inflap)$ and concludes the proof.
\end{proof}

\begin{remark}
We can similarly get 1-Laplacian analogs of some of the results in Section \ref{sec:spectral-minmial-partition}. For example, a 1-Laplacian version of Lemma \ref{lemma:limit1} reads as:
$$\lim_{p\rightarrow 1}\widehat{\lambda_k}(\Delta_p)^{\frac1p}=\widehat{\Lambda_k}(\Delta_1)\qquad \text{and}\qquad \lim_{p\rightarrow 1}\widetilde{\lambda_k}(\Delta_p)^{\frac1p}=\widetilde{\Lambda_k}(\Delta_1),$$
and a 1-Laplacian version of Theorem \ref{Thm_spectral_min_partitions} is
\begin{align}
h_k(\Gc):&=\min_{\{A_1,\dots, A_k\}\in\mathcal{D}_k(\Gc)}\;\max_{i} c(A_i)=\min\limits_{\text{disjoint }V_1,\ldots,V_k\text{ in }\internalnodes}\max_i h_1\big(\Gc(V_i)\big)   \\&= \min\limits_{\text{disjoint }V_1,\ldots,V_k\text{ in }\internalnodes}\max_i \Lambda_1\big(\Delta_1,\Gc(V_i)\big)=\widehat{\Lambda_k}(\Delta_1).
\end{align} 
Hence, compared to the variational eigenvalues, the spectral minimal partition costs are more closely related to certain geometric quantities such as higher order Cheeger constants and packing radii.
\end{remark}

\section{Conclusions}
We have devoted this manuscript to the study of the tight relationships between the graph $p$-Laplacian nonlinear eigenproblem and the geometry of the graph. In particular, we have shown that the relationships between eigenpairs and geometric quantities become more informative as one moves away from the linear case $p=2$. Indeed, in the two extreme cases $p=1$ and $p=\infty$, the eigenvalues encode exact geometric quantities, namely the normalized cut of subgraphs and distances between nodes, respectively. 
In more detail, by focusing on the subclass of variational eigenpairs, we have proved that the $p$-Laplacian spectrum provides information about several geometric and combinatorial quantities of the graph, including the Cheeger constants, the packing radii, the independence numbers, and the matching numbers. Moreover, we have established different results that quantify how tight these relationships are. 
Nonetheless, the nonlinear spectral theory of the graph $p$-Laplacian is still far from complete. Throughout the manuscript, and already in the introduction, we have identified several open problems whose resolution could have a significant impact on the further development of the theory.

\section*{Acknowledgments.} 
The authors would like to express their sincere gratitude to the anonymous reviewers for their comprehensive, detailed, thoughtful and insightful comments and suggestions, which have greatly enhanced the quality of this paper.


\bibliographystyle{abbrv}
\bibliography{strings.bib, references2}

\begin{thebibliography}{100}

\bibitem{Elmoataz2}
E.~Abderrahim, D.~Xavier, L.~Zakaria, and L.~Olivier.
\newblock Nonlocal infinity laplacian equation on graphs with applications in
  image processing and machine learning.
\newblock {\em Mathematics and Computers in Simulation}, 102:153--163, 2014.
\newblock 4th International Conference on Approximation Methods and Numerical
  Modeling in Environment and Natural Resources.

\bibitem{ABIAD2019indep}
A.~Abiad, G.~Coutinho, and M.~Fiol.
\newblock On the k-independence number of graphs.
\newblock {\em Discrete Mathematics}, 342(10):2875--2885, 2019.
\newblock Algebraic and Extremal Graph Theory.

\bibitem{Alagmir_2011_p-resistances}
M.~Alamgir and U.~Luxburg.
\newblock Phase transition in the family of p-resistances.
\newblock In J.~Shawe-Taylor, R.~Zemel, P.~Bartlett, F.~Pereira, and
  K.~Weinberger, editors, {\em Advances in Neural Information Processing
  Systems}, volume~24. Curran Associates, Inc., 2011.

\bibitem{Amghibech1}
S.~Amghibech.
\newblock Eigenvalues of the {D}iscrete p-{L}aplacian for {G}raphs.
\newblock {\em Ars Combinatoria}, 67:283–302, 04 2003.

\bibitem{Amghibech2}
S.~Amghibech.
\newblock Bounds for the largest p-{L}aplacian eigenvalue for graphs.
\newblock {\em Discr. Math.}, 306:2762--2771, 2006.

\bibitem{BAI_SCF}
Z.~Bai, D.~Lu, and B.~Vandereycken.
\newblock Robust rayleigh quotient minimization and nonlinear eigenvalue
  problems.
\newblock {\em SIAM Journal on Scientific Computing}, 40(5):A3495--A3522, 2018.

\bibitem{Band}
R.~Band, I.~Oren, and U.~Smilansky.
\newblock Nodal domains on graphs - {H}ow to count them and why?
\newblock In {\em Proc. Sympos. Pure Math., Providence, RI: Amer. Math. Soc},
  volume~77, 11 2007.

\bibitem{Berge}
C.~Berge.
\newblock {\em Hypergraphs}, volume~45 of {\em North-Holland Mathematical
  Library}.
\newblock North-Holland Publishing Co., Amsterdam, 1989.
\newblock Combinatorics of finite sets, Translated from the French.

\bibitem{bergermann2023nonlinear}
K.~Bergermann, M.~Stoll, and F.~Tudisco.
\newblock A nonlinear spectral core-periphery detection method for multiplex
  networks.
\newblock {\em arXiv preprint arXiv:2310.19697}, 2023.

\bibitem{Berkolaiko2}
G.~Berkolaiko.
\newblock A lower bound for nodal count on discrete and metric graphs.
\newblock {\em Communications in mathematical physics}, 278(3):803--819, 2008.

\bibitem{berkolaiko2025eigenvalues}
G.~Berkolaiko and M.~Hofmann.
\newblock Eigenvalues of the discrete p-laplacian via graph surgery.
\newblock {\em arXiv preprint arXiv:2509.06686}, 2025.

\bibitem{Berkolaiko_2017}
G.~Berkolaiko, J.~B. Kennedy, P.~Kurasov, and D.~Mugnolo.
\newblock Edge connectivity and the spectral gap of combinatorial and quantum
  graphs.
\newblock {\em Journal of Physics A: Mathematical and Theoretical},
  50(36):365201, aug 2017.

\bibitem{Berkolaiko1}
G.~Berkolaiko, H.~Raz, and U.~Smilansky.
\newblock Stability of nodal structures in graph eigenfunctions and its
  relation to the nodal domain count.
\newblock {\em Journal of Physics A: Mathematical and Theoretical},
  45(16):165203, 2012.

\bibitem{Biy1}
T.~Biyiko{\u g}lu.
\newblock A discrete nodal domain theorem for trees.
\newblock {\em Linear algebra and its Applications}, 360:197--205, 2003.

\bibitem{Biy2}
T.~Biyiko{\u g}lu, J.~Leydold, and P.~Stadler.
\newblock Nodal {D}omain {T}heorems and {B}ipartite {S}ubgraphs.
\newblock {\em The electronic journal of linear algebra}, 13:344--351, 01 2005.

\bibitem{Stadler_book-Laplacian-eigenvectors}
T.~Biyiko{\u{g}}u, J.~Leydold, and P.~F. Stadler.
\newblock {\em Laplacian Eigenvectors of Graphs Perron-Frobenius and
  Faber-Krahn Type Theorems}.
\newblock Springer, 2007.

\bibitem{Blum}
G.~Blum, S.~Gnutzmann, and U.~Smilansky.
\newblock Nodal domains statistics: {A} criterion for quantum chaos.
\newblock {\em Physical Review Letters}, 88(11):114101, 2002.

\bibitem{BobkovParini18}
V.~Bobkov and E.~Parini.
\newblock On the higher cheeger problem.
\newblock {\em Journal of the London Mathematical Society}, 97(3):575--600,
  2018.

\bibitem{Bozorgnia_power_method}
F.~Bozorgnia.
\newblock Convergence of inverse power method for first eigenvalue of p-laplace
  operator.
\newblock {\em Numerical Functional Analysis and Optimization},
  37(11):1378--1384, 2016.

\bibitem{BOZORGNIA_approximation}
F.~Bozorgnia and A.~Arakelyan.
\newblock Approximation of the second eigenvalue of the p-laplace operator in
  symmetric domains.
\newblock {\em Journal of Computational and Applied Mathematics}, 434:115349,
  2023.

\bibitem{bozorgnia2024infinity}
F.~Bozorgnia, L.~Bungert, and D.~Tenbrinck.
\newblock The infinity laplacian eigenvalue problem: reformulation and a
  numerical scheme.
\newblock {\em Journal of Scientific Computing}, 98(2):40, 2024.

\bibitem{Bhuler}
T.~B\"uhler and M.~Hein.
\newblock Spectral clustering based on the graph $p$-{L}aplacian.
\newblock {\em Proceedings of the 26th International Conference on Machine
  Learning}, 01 2009.

\bibitem{BUNGERT_gradient_flows_and_nonlinear_power_method}
L.~Bungert and M.~Burger.
\newblock Chapter 13 - gradient flows and nonlinear power methods for the
  computation of nonlinear eigenfunctions.
\newblock In E.~Trélat and E.~Zuazua, editors, {\em Numerical Control: Part
  A}, volume~23 of {\em Handbook of Numerical Analysis}, pages 427--465.
  Elsevier, 2022.

\bibitem{BungertNonlineardecomp}
L.~Bungert, M.~Burger, A.~Chambolle, and M.~Novaga.
\newblock Nonlinear spectral decompositions by gradient flows of
  one-homogeneous functionals.
\newblock {\em Anal. PDE}, 14(3):823--860, 2021.

\bibitem{bungert2021eigenvalue}
L.~Bungert and Y.~Korolev.
\newblock Eigenvalue problems in $\mathrm{L}^{\infty}$: Optimality conditions,
  duality, and relations with optimal transport.
\newblock {\em Communications of the American Mathematical Society},
  2:345--373, 2022.

\bibitem{Bungert1}
L.~Bungert, Y.~Korolev, and M.~Burger.
\newblock Structural analysis of an {L}-infinity variational problem and
  relations to distance functions.
\newblock {\em Pure and Applied Analysis}, 2, 01 2020.

\bibitem{Burger23}
M.~Burger.
\newblock Nonlinear eigenvalue problems for seminorms and applications.
\newblock In {\em I{CM}---{I}nternational {C}ongress of {M}athematicians.
  {V}ol. 7. {S}ections 15--20}, pages 5234--5255. EMS Press, Berlin, 2023.

\bibitem{burger2016spectral}
M.~Burger, G.~Gilboa, M.~Moeller, L.~Eckardt, and D.~Cremers.
\newblock Spectral decompositions using one-homogeneous functionals.
\newblock {\em SIAM Journal on Imaging Sciences}, 9(3):1374--1408, 2016.

\bibitem{calder2018game}
J.~Calder.
\newblock The game theoretic $p$-laplacian and semi-supervised learning with
  few labels.
\newblock {\em Nonlinearity}, 32(1):301--330, dec 2018.

\bibitem{calder2019consistency}
J.~Calder.
\newblock Consistency of lipschitz learning with infinite unlabeled data and
  finite labeled data.
\newblock {\em SIAM Journal on Mathematics of Data Science}, 1(4):780--812,
  2019.

\bibitem{chang1981variational}
K.-C. Chang.
\newblock Variational methods for non-differentiable functionals and their
  applications to partial differential equations.
\newblock {\em Journal of Mathematical Analysis and Applications},
  80(1):102--129, 1981.

\bibitem{chang2016spectrum}
K.~C. Chang.
\newblock Spectrum of the 1-laplacian and cheeger's constant on graphs.
\newblock {\em Journal of Graph Theory}, 81(2):167--207, 2016.

\bibitem{Zhang_1-Lap_cheeger_cut}
K.~C. Chang, S.~Shao, and D.~Zhang.
\newblock The 1-laplacian cheeger cut: Theory and algorithms.
\newblock {\em Journal of Computational Mathematics}, 33(5):443--467, 2015.

\bibitem{ZhangNodalDO}
K.~C. Chang, S.~Shao, and D.~Zhang.
\newblock Nodal domains of eigenvectors for 1-{L}aplacian on graphs.
\newblock {\em Advances in Mathematics}, 308:529--574, 2017.

\bibitem{chang2021nonsmooth}
K.-C. Chang, S.~Shao, D.~Zhang, and W.~Zhang.
\newblock Nonsmooth critical point theory and applications to the spectral
  graph theory.
\newblock {\em Science China Mathematics}, 64(1):1--32, 2021.

\bibitem{Cheeger}
J.~Cheeger.
\newblock A lower bound for the smallest eigenvalue of the {L}aplacian.
\newblock In {\em Problems in analysis ({S}ympos. in honor of {S}alomon
  {B}ochner, {P}rinceton {U}niv., {P}rinceton, {N}.{J}., 1969)}, pages
  195--199. Princeton Univ. Press, Princeton, NJ, 1970.

\bibitem{chung1997spectral}
F.~R. Chung.
\newblock {\em Spectral graph theory}, volume~92.
\newblock American Mathematical Soc., 1997.

\bibitem{clarke1990optimization}
F.~H. Clarke.
\newblock {\em Optimization and nonsmooth analysis}.
\newblock SIAM, 1990.

\bibitem{COHEN_introducing_p_lap_spectra}
I.~Cohen and G.~Gilboa.
\newblock Introducing the p-laplacian spectra.
\newblock {\em Signal Processing}, 167:107281, 2020.

\bibitem{cohen2005stability}
D.~Cohen-Steiner, H.~Edelsbrunner, and J.~Harer.
\newblock Stability of persistence diagrams.
\newblock In {\em Proceedings of the twenty-first annual symposium on
  Computational geometry}, pages 263--271, 2005.

\bibitem{Colbrook}
M.~J. Colbrook.
\newblock Finiteness and exponential growth of full graph $p$-spectra, July
  2026.

\bibitem{Conti05}
M.~Conti, S.~Terracini, and G.~Verzini.
\newblock On a class of optimal partition problems related to the {F}u\v c\'ik
  spectrum and to the monotonicity formulae.
\newblock {\em Calc. Var. Partial Differential Equations}, 22(1):45--72, 2005.

\bibitem{Courant}
R.~Courant and D.~Hilbert.
\newblock {\em Methods of Mathematical Physics}, volume~1.
\newblock Interscience, New York, 1953.

\bibitem{cristofari2020total}
A.~Cristofari, F.~Rinaldi, and F.~Tudisco.
\newblock Total variation based community detection using a nonlinear
  optimization approach.
\newblock {\em SIAM Journal on Applied Mathematics}, 80(3):1392--1419, 2020.

\bibitem{Cuesta99}
M.~Cuesta, D.~de~Figueiredo, and J.-P. Gossez.
\newblock The beginning of the {F}u\v cik spectrum for the {$p$}-{L}aplacian.
\newblock {\em J. Differential Equations}, 159(1):212--238, 1999.

\bibitem{daneshgar2010isoperimetric}
A.~Daneshgar, H.~Hajiabolhassan, and R.~Javadi.
\newblock On the isoperimetric spectrum of graphs and its approximations.
\newblock {\em Journal of Combinatorial Theory, Series B}, 100(4):390--412,
  2010.

\bibitem{daneshgar2012nodal}
A.~Daneshgar, R.~Javadi, and L.~Miclo.
\newblock On nodal domains and higher-order {C}heeger inequalities of finite
  reversible {M}arkov processes.
\newblock {\em Stochastic Processes and their Applications}, 122(4):1748--1776,
  2012.

\bibitem{Davies}
E.~B. Davies, G.~Gladwell, J.~Leydold, and P.~Stadler.
\newblock Discrete nodal domain theorems.
\newblock {\em Linear algebra and its Applications}, 336(1):51--60, 2001.

\bibitem{deidda2023PhdThesis}
P.~Deidda.
\newblock {\em The graph p-Laplacian eigenvalue problem}.
\newblock PhD thesis, Universit{\`a} degli studi di Padova, 2023.

\bibitem{deidda2024_inf_eigenproblem}
P.~Deidda, M.~Burger, M.~Putti, and F.~Tudisco.
\newblock The graph \(\boldsymbol{\infty }\)-laplacian eigenvalue problem.
\newblock {\em SIAM Journal on Mathematical Analysis}, 58(1):1--34, 2026.

\bibitem{ge2025}
P.~Deidda, C.~Ge, S.~Liu, F.~Tudisco, and D.~Zhang.
\newblock Perturbation of $p$-laplacian eigenfunctions on trees.
\newblock {\em in progress}, 2025.

\bibitem{DEIDDA2023_nod_dom}
P.~Deidda, M.~Putti, and F.~Tudisco.
\newblock Nodal domain count for the generalized graph p-laplacian.
\newblock {\em Applied and Computational Harmonic Analysis}, 64:1--32, 2023.

\bibitem{deidda2024_spec_energy}
P.~Deidda, N.~Segala, and M.~Putti.
\newblock Graph \({p}\)-laplacian eigenpairs as saddle points of a family of
  spectral energy functions.
\newblock {\em SIAM Journal on Matrix Analysis and Applications},
  46(2):1540--1561, 2025.

\bibitem{drabek2012variational}
P.~Dr{\'a}bek.
\newblock On the variational eigenvalues which are not of
  ljusternik-schnirelmann type.
\newblock In {\em Abstract and Applied Analysis}, volume 2012. Hindawi, 2012.

\bibitem{drabek2010variational}
P.~Dr{\'a}bek and P.~Tak{\'a}{\v{c}}.
\newblock On variational eigenvalues of the p-laplacian which are not of
  ljusternik--schnirelmann type.
\newblock {\em Journal of the London Mathematical Society}, 81(3):625--649,
  2010.

\bibitem{dructu2019random}
C.~Dru{\c{t}}u and J.~M. Mackay.
\newblock Random groups, random graphs and eigenvalues of p-laplacians.
\newblock {\em Advances in Mathematics}, 341:188--254, 2019.

\bibitem{Duval}
A.~Duval and V.~Reiner.
\newblock {P}erron–{F}robenius type results and discrete versions of nodal
  domain theorems.
\newblock {\em Linear algera and its Applications}, 294(1):259--268, 1999.

\bibitem{ekeland1999convex}
I.~Ekeland and R.~Temam.
\newblock {\em Convex analysis and variational problems}.
\newblock SIAM, 1999.

\bibitem{elmoataz2011infinity}
A.~Elmoataz, X.~Desquesnes, Z.~Lakhdari, and O.~Lezoray.
\newblock On the infinity laplacian equation on graph with applications to
  image and manifolds processing.
\newblock In {\em International Conference on Approximation Methods and
  Numerical Modelling in Environment and Natural Resources}, 2011.

\bibitem{Elmoataz1}
A.~Elmoataz, M.~Toutain, and D.~Tenbrinck.
\newblock On the p-{L}aplacian and $\infty$-{L}aplacian on graphs with
  applications in image and data processing.
\newblock {\em SIAM Journal on Imaging Sciences}, 8(4):2412--2451, 2015.

\bibitem{Esposito}
L.~Esposito, B.~Kawohl, C.~Nitsch, and C.~Trombetti.
\newblock The {N}eumann eigenvalue problem for the $\infty$-{L}aplacian.
\newblock {\em Rendiconti Lincei - Matematica e Applicazioni}, 26, 05 2014.

\bibitem{Fazeny}
A.~Fazeny.
\newblock $p$-laplacian operators on hypergraphs.
\newblock Master's thesis, Friedrich Alexander University of Erlangen
  Nur\"nberg, 2023.

\bibitem{Fazeny23}
A.~Fazeny, D.~Tenbrinck, and M.~Burger.
\newblock Hypergraph p-laplacians, scale spaces, and information flow in
  networks.
\newblock In {\em Scale Space and Variational Methods in Computer Vision},
  pages 677--690. SSVM, 2023.

\bibitem{Fazeny24}
A.~Fazeny, D.~Tenbrinck, K.~Lukin, and M.~Burger.
\newblock Hypergraph {$p$}-{L}aplacians and scale spaces.
\newblock {\em J. Math. Imaging Vision}, 66(4):529--549, 2024.

\bibitem{Gilboa_rayleight_quot_minimization}
T.~Feld, J.-F. Aujol, G.~Gilboa, and N.~Papadakis.
\newblock Rayleigh quotient minimization for absolutely one-homogeneous
  functionals.
\newblock {\em Inverse Problems}, 35(6):064003, may 2019.

\bibitem{Fiedler}
M.~Fiedler.
\newblock Eigenvectors of acyclic matrices.
\newblock {\em Czechoslovak Mathematical Journal}, 25(4):607--618, 1975.

\bibitem{fiol1997eigenvalue}
M.~A. Fiol.
\newblock An eigenvalue characterization of antipodal distance-regular graphs.
\newblock {\em the electronic journal of combinatorics}, pages R30--R30, 1997.

\bibitem{FIRBY1997indep}
P.~Firby and J.~Haviland.
\newblock Independence and average distance in graphs.
\newblock {\em Discrete Applied Mathematics}, 75(1):27--37, 1997.

\bibitem{flores2022analysis}
M.~Flores, J.~Calder, and G.~Lerman.
\newblock Analysis and algorithms for $\ell_p$-based semi-supervised learning
  on graphs.
\newblock {\em Applied and Computational Harmonic Analysis}, 60:77--122, 2022.

\bibitem{Fucik}
S.~Fucik, J.~Necas, J.~Soucek, and V.~Soucek.
\newblock {\em Spectral analysis of nonlinear operators}, volume 346.
\newblock Springer, 2006.

\bibitem{trillos2020error}
N.~Garc{\'\i}a~Trillos, M.~Gerlach, M.~Hein, and D.~Slep{\v{c}}ev.
\newblock Error estimates for spectral convergence of the graph laplacian on
  random geometric graphs toward the laplace--beltrami operator.
\newblock {\em Foundations of Computational Mathematics}, 20(4):827--887, 2020.

\bibitem{Garcia20}
N.~Garc\'ia~Trillos, M.~Gerlach, M.~Hein, and D.~Slep\v{c}ev.
\newblock Error estimates for spectral convergence of the graph {L}aplacian on
  random geometric graphs toward the {L}aplace-{B}eltrami operator.
\newblock {\em Found. Comput. Math.}, 20(4):827--887, 2020.

\bibitem{GarciaSlepcev16}
N.~Garc\'{\i}a~Trillos and D.~Slep\v{c}ev.
\newblock Continuum limit of total variation on point clouds.
\newblock {\em Arch. Ration. Mech. Anal.}, 220(1):193--241, 2016.

\bibitem{Garcia16}
N.~Garc\'ia~Trillos, D.~Slep\v{c}ev, J.~von Brecht, T.~Laurent, and X.~Bresson.
\newblock Consistency of {C}heeger and ratio graph cuts.
\newblock {\em J. Mach. Learn. Res.}, 17:Paper No. 181, 46, 2016.

\bibitem{ge2023new}
C.~Ge, S.~Liu, and D.~Zhang.
\newblock New graph invariants based on $ p $-laplacian eigenvalues.
\newblock {\em arXiv preprint arXiv:2310.08189}, 2023.

\bibitem{zhang2023pLap_noddom}
C.~Ge, S.~Liu, and D.~Zhang.
\newblock Nodal domain theorems for $ p $-laplacians on signed graphs.
\newblock {\em Journal of Spectral Theory}, 13(3):937--989, 2023.

\bibitem{ge2025computing}
C.~Ge and O.~Qin.
\newblock Computing the $ p $-laplacian eigenpairs of signed graphs.
\newblock {\em arXiv preprint arXiv:2501.07929}, 2025.

\bibitem{Ghoussoub}
N.~Ghoussoub.
\newblock {\em Duality and Perturbation Methods in Critical Point Theory}.
\newblock Cambridge Tracts in Mathematics. Cambridge University Press, 1993.

\bibitem{Gilboa_total_variation_spectral_framework}
G.~Gilboa.
\newblock A total variation spectral framework for scale and texture analysis.
\newblock {\em SIAM Journal on Imaging Sciences}, 7(4):1937--1961, 2014.

\bibitem{gilboa2018nonlinear}
G.~Gilboa.
\newblock {\em Nonlinear eigenproblems in image processing and computer
  vision}.
\newblock Springer, 2018.

\bibitem{gilboa2021iterative}
G.~Gilboa.
\newblock Iterative methods for computing eigenvectors of nonlinear operators.
\newblock In {\em Handbook of mathematical models and algorithms in computer
  vision and imaging: mathematical imaging and vision}, pages 1--28. Springer,
  2021.

\bibitem{godsil01}
C.~Godsil and G.~Royle.
\newblock {\em Algebraic graph theory}, volume 207 of {\em Graduate Texts in
  Mathematics}.
\newblock Springer-Verlag, New York, 2001.

\bibitem{govc2016definition}
D.~Govc.
\newblock On the definition of the homological critical value.
\newblock {\em Journal of Homotopy and Related Structures}, 11(1):143--151,
  2016.

\bibitem{grove1995new}
K.~Grove and S.~Markvorsen.
\newblock New extremal problems for the riemannian recognition program via
  alexandrov geometry.
\newblock {\em Journal of the American Mathematical Society}, pages 1--28,
  1995.

\bibitem{hein2010inverse}
M.~Hein and T.~B{\"u}hler.
\newblock An inverse power method for nonlinear eigenproblems with applications
  in 1-spectral clustering and sparse pca.
\newblock {\em Advances in Neural Information Processing Systems}, 23, 2010.

\bibitem{oberwolfach15}
M.~Hein, H.~D. Lenz, and D.~Mugnolo.
\newblock Mini-workshop: {D}iscrete {$p$}-{L}aplacians: {S}pectral {T}heory and
  {V}ariational {M}ethods in {M}athematics and {C}omputer {S}cience.
\newblock {\em Oberwolfach Rep.}, 12(1):399--447, 2015.
\newblock Abstracts from the mini-workshop held February 8--14, 2015, Organized
  by Matthias Hein, Daniel Lenz and Delio Mugnolo.

\bibitem{Hein_Beyond}
M.~Hein and S.~Setzer.
\newblock Beyond spectral clustering - tight relaxations of balanced graph
  cuts.
\newblock In J.~Shawe-Taylor, R.~Zemel, P.~Bartlett, F.~Pereira, and
  K.~Weinberger, editors, {\em Advances in Neural Information Processing
  Systems}, volume~24. Curran Associates, Inc., 2011.

\bibitem{Helffer09}
B.~Helffer, T.~Hoffmann-Ostenhof, and S.~Terracini.
\newblock Nodal domains and spectral minimal partitions.
\newblock {\em Ann. Inst. H. Poincar\'e{} C Anal. Non Lin\'eaire},
  26(1):101--138, 2009.

\bibitem{herbster2010triangle}
M.~Herbster.
\newblock A triangle inequality for p-resistance.
\newblock 2010.

\bibitem{herbster2009predicting}
M.~Herbster and G.~Lever.
\newblock Predicting the labelling of a graph via minimum p-seminorm
  interpolation.
\newblock In {\em NIPS Workshop 2010: Networks Across Disciplines: Theory and
  Applications}, 2009.

\bibitem{Hua}
B.~Hua and L.~Wang.
\newblock Dirichlet p-{L}aplacian eigenvalues and {C}heeger constants on
  symmetric graphs.
\newblock {\em Advances in Mathematics}, 364:106997, 04 2020.

\bibitem{HYND_approximation_of_the_least_Rayleight_quotient}
R.~Hynd and E.~Lindgren.
\newblock Approximation of the least rayleigh quotient for degree p homogeneous
  functionals.
\newblock {\em Journal of Functional Analysis}, 272(12):4873--4918, 2017.

\bibitem{Jarlebring_an_inverse}
E.~Jarlebring, S.~Kvaal, and W.~Michiels.
\newblock An inverse iteration method for eigenvalue problems with eigenvector
  nonlinearities.
\newblock {\em SIAM Journal on Scientific Computing}, 36(4):A1978--A2001, 2014.

\bibitem{Mulas2022_plap_hyper}
J.~Jost, R.~Mulas, and D.~Zhang.
\newblock p-laplace operators for oriented hypergraphs.
\newblock {\em Vietnam Journal of Mathematics}, 50(2):323--358, 2022.

\bibitem{zhang2021discrete}
J.~Jost and D.~Zhang.
\newblock Discrete-to-continuous extensions: piecewise multilinear extension,
  min-max theory and spectral theory.
\newblock {\em arXiv preprint arXiv:2106.04116}, 2021.

\bibitem{jostzhang24+}
J.~Jost and D.~Zhang.
\newblock Cheeger inequalities on simplicial complexes.
\newblock {\em Ann. Sc. Norm. Super. Pisa Cl. Sci. (5)}, 2024.

\bibitem{Lind3}
P.~Juutinen and P.~Lindqvist.
\newblock On the higher eigenvalues for the $\infty$-eigenvalue problem.
\newblock {\em Calculus of Variations and Partial Differential Equations},
  23(2):169--192, 2005.

\bibitem{Lind2}
P.~Juutinen, P.~Lindqvist, and J.~Manfredi.
\newblock The $\infty$-eigenvalue problem.
\newblock {\em Archive for Rational Mechanics and Analysis}, 148:89--105, 1999.

\bibitem{Kawohl2003}
B.~Kawohl and V.~Fridman.
\newblock Isoperimetric estimates for the first eigenvalue of the $p$-{L}aplace
  operator and the {C}heeger constant.
\newblock {\em Commentationes Mathematicae Universitatis Carolinae},
  44(4):659--667, 2003.

\bibitem{KELLER201680}
M.~Keller and D.~Mugnolo.
\newblock General cheeger inequalities for p-laplacians on graphs.
\newblock {\em Nonlinear Analysis: Theory, Methods \& Applications},
  147:80--95, 2016.

\bibitem{Lanza24}
A.~Lanza, S.~Morigi, and G.~Recupero.
\newblock Variational graph p-laplacian eigendecomposition under
  p-orthogonality constraints.
\newblock {\em Computational Optimization and Applications}, 2024.

\bibitem{laubmann2025duality}
J.~Laubmann, M.~Friedrich, and D.~Tenbrinck.
\newblock Duality perspective on nonlinear eigenproblems.
\newblock {\em arXiv preprint arXiv:2511.19188}, 2025.

\bibitem{lawler1988bounds}
G.~F. Lawler and A.~D. Sokal.
\newblock Bounds on the {L}$^2$ spectrum for {M}arkov chains and {M}arkov
  processes: a generalization of {C}heeger’s inequality.
\newblock {\em Transactions of the American mathematical society},
  309(2):557--580, 1988.

\bibitem{lee2014multiway}
J.~R. Lee, S.~O. Gharan, and L.~Trevisan.
\newblock Multiway spectral partitioning and higher-order cheeger inequalities.
\newblock {\em Journal of the ACM (JACM)}, 61(6):1--30, 2014.

\bibitem{Li}
P.~Li and O.~Milenkovic.
\newblock Submodular hypergraphs: p-{L}aplacians, {C}heeger inequalities and
  spectral clustering.
\newblock In {\em Proceedings of the 35th International Conference on Machine
  Learning}, volume~80 of {\em Proceedings of Machine Learning Research}, pages
  3014--3023. PMLR, 10--15 Jul 2018.

\bibitem{Lind4}
P.~Lindqvist.
\newblock Note on a nonlinear eigenvalue problem.
\newblock {\em The Rocky Mountain Journal of Mathematics}, pages 281--288,
  1993.

\bibitem{Luo10}
D.~Luo, H.~Huang, C.~Ding, and F.~Nie.
\newblock On the eigenvectors of {$p$}-{L}aplacian.
\newblock {\em Mach. Learn.}, 81(1):37--51, 2010.

\bibitem{manfredi2010asymptotic}
J.~Manfredi, M.~Parviainen, and J.~Rossi.
\newblock An asymptotic mean value characterization for p-harmonic functions.
\newblock {\em Proceedings of the American Mathematical Society},
  138(3):881--889, 2010.

\bibitem{manfredi2015nonlinear}
J.~J. Manfredi, A.~M. Oberman, and A.~P. Sviridov.
\newblock Nonlinear elliptic partial differential equations and p-harmonic
  functions on graphs.
\newblock {\em Differential Integral Equations}, 28(1-2):79--102, 2015.

\bibitem{MeirMoon75}
A.~Meir and J.~W. Moon.
\newblock Relations between packing and covering numbers of a tree.
\newblock {\em Pacific J. Math.}, 61(1):225--233, 1975.

\bibitem{MengZhang25}
Z.~Meng and D.~Zhang.
\newblock Combinatorial courant-fischer-weyl minimax principle on cheeger $ k
  $-constants of weighted forests.
\newblock {\em arXiv preprint arXiv:2510.06301}, 2025.

\bibitem{miclo2008eigenfunctions}
L.~Miclo.
\newblock On eigenfunctions of markov processes on trees.
\newblock {\em Probability Theory and Related Fields}, 142:561--594, 2008.

\bibitem{milnor2016morse}
J.~Milnor.
\newblock Morse theory.(am-51), volume 51.
\newblock In {\em Morse Theory.(AM-51), Volume 51}. Princeton university press,
  2016.

\bibitem{MulasZhang}
R.~Mulas and D.~Zhang.
\newblock Spectral theory of {L}aplace operators on oriented hypergraphs.
\newblock {\em Discrete Math.}, 344(6):Paper No. 112372, 22, 2021.

\bibitem{Oren}
I.~Oren.
\newblock Nodal domain counts and the chromatic number of graphs.
\newblock {\em Journal of Physics A: Mathematical and Theoretical},
  40(32):9825, 2007.

\bibitem{osher2005iterative}
S.~Osher, M.~Burger, D.~Goldfarb, J.~Xu, and W.~Yin.
\newblock An iterative regularization method for total variation-based image
  restoration.
\newblock {\em Multiscale Modeling \& Simulation}, 4(2):460--489, 2005.

\bibitem{palais1970critical}
R.~S. Palais.
\newblock Critical point theory and the minimax principle.
\newblock In {\em Proc. Symp. Pure Math}, volume~15, pages 185--212. American
  Mathematical Society Providence, 1970.

\bibitem{parini2010second}
E.~Parini.
\newblock The second eigenvalue of the p-laplacian as p goes to 1.
\newblock {\em International Journal of Differential Equations},
  2010(1):984671, 2010.

\bibitem{PARK2011}
J.-H. Park.
\newblock On a resonance problem with the discrete p-laplacian on finite
  graphs.
\newblock {\em Nonlinear Analysis: Theory, Methods \& Applications},
  74(17):6662--6675, 2011.

\bibitem{Peres}
Y.~Peres and S.~Sheffield.
\newblock {Tug-of-war with noise: A game-theoretic view of the $p$-Laplacian}.
\newblock {\em Duke Mathematical Journal}, 145(1):91 -- 120, 2008.

\bibitem{Reff}
N.~Reff and L.~J. Rusnak.
\newblock An oriented hypergraphic approach to algebraic graph theory.
\newblock {\em Linear Algebra Appl.}, 437(9):2262--2270, 2012.

\bibitem{Robeva}
E.~Robeva and A.~Seigal.
\newblock Duality of graphical models and tensor networks.
\newblock {\em Inf. Inference}, 8(2):273--288, 2019.

\bibitem{rockafellar2015convex}
R.~T. Rockafellar.
\newblock Convex analysis.
\newblock In {\em Convex analysis}. Princeton university press, 2015.

\bibitem{roith2023continuum}
T.~Roith and L.~Bungert.
\newblock Continuum limit of lipschitz learning on graphs.
\newblock {\em Foundations of Computational Mathematics}, 23(2):393--431, 2023.

\bibitem{rombach2014core}
M.~P. Rombach, M.~A. Porter, J.~H. Fowler, and P.~J. Mucha.
\newblock Core-periphery structure in networks.
\newblock {\em SIAM Journal on Applied Mathematics}, 74(1):167--190, 2014.

\bibitem{rudin1992nonlinear}
L.~I. Rudin, S.~Osher, and E.~Fatemi.
\newblock Nonlinear total variation based noise removal algorithms.
\newblock {\em Physica D: nonlinear phenomena}, 60(1-4):259--268, 1992.

\bibitem{saito23}
S.~Saito and M.~Herbster.
\newblock Generalizing {$p$}-{L}aplacian: spectral hypergraph theory and a
  partitioning algorithm.
\newblock {\em Mach. Learn.}, 112:241--280, 2023.

\bibitem{saito23effective_p_resistance}
S.~Saito and M.~Herbster.
\newblock Multi-class graph clustering via approximated effective
  $p$-resistance.
\newblock In {\em Proceedings of the 40th International Conference on Machine
  Learning}, volume 202 of {\em Proceedings of Machine Learning Research},
  pages 29697--29733. PMLR, 23--29 Jul 2023.

\bibitem{saito2018hypergraph}
S.~Saito, D.~P. Mandic, and H.~Suzuki.
\newblock Hypergraph p-{L}aplacian: A differential geometry view.
\newblock In {\em Thirty-Second AAAI Conference on Artificial Intelligence},
  2018.

\bibitem{Schaub20}
M.~T. Schaub, A.~R. Benson, P.~Horn, G.~Lippner, and A.~Jadbabaie.
\newblock Random walks on simplicial complexes and the normalized {H}odge
  1-{L}aplacian.
\newblock {\em SIAM Rev.}, 62(2):353--391, 2020.

\bibitem{chang2016signless}
S.~Shao, C.~Yang, and D.~Zhang.
\newblock Dual {C}heeger constants, signless 1-{L}aplacians and maxcut.
\newblock {\em Sci. China Math.}, 68(11):2773--2790, 2025.

\bibitem{shi1992signed}
C.-J. Shi.
\newblock A signed hypergraph model of the constrained via minimization
  problem.
\newblock {\em Microelectronics journal}, 23(7):533--542, 1992.

\bibitem{slepcev2019analysis}
D.~Slep\v{c}ev and M.~Thorpe.
\newblock Analysis of $p$-laplacian regularization in semisupervised learning.
\newblock {\em SIAM J. Math. Anal.}, 51(3):2085--2120, 2019.

\bibitem{Solimini}
S.~Solimini.
\newblock Chapter 7 multiplicity techniques for problems without compactness.
\newblock {\em Handbook of Differential Equations: Stationary Partial
  Differential Equations}, 2:519--599, 12 2005.

\bibitem{spielman2019spectralandalgebraic}
D.~Spielman.
\newblock Spectral and algebraic graph theory.
\newblock {\em Yale lecture notes, draft of December}, 4:47, 2019.

\bibitem{struwe}
M.~Struwe.
\newblock {\em Variational methods}, volume 991.
\newblock Springer, 2000.

\bibitem{Topp91}
J.~Topp and L.~Volkmann.
\newblock On packing and covering numbers of graphs.
\newblock {\em Discrete Math.}, 96(3):229--238, 1991.

\bibitem{trevisan2017lecture}
L.~Trevisan.
\newblock Lecture notes on graph partitioning, expanders and spectral methods.
\newblock {\em University of California, Berkeley, https://people. eecs.
  berkeley. edu/luca/books/expanders-2016. pdf}, 2017.

\bibitem{trillos2018variational}
N.~G. Trillos and D.~Slep{\v{c}}ev.
\newblock A variational approach to the consistency of spectral clustering.
\newblock {\em Applied and Computational Harmonic Analysis}, 45(2):239--281,
  2018.

\bibitem{tudisco2021nonlinear}
F.~Tudisco, A.~R. Benson, and K.~Prokopchik.
\newblock Nonlinear higher-order label spreading.
\newblock In {\em Proceedings of the Web Conference 2021}, pages 2402--2413,
  2021.

\bibitem{Tudisco1}
F.~Tudisco and M.~Hein.
\newblock A nodal domain theorem and a higher-order {C}heeger inequality for
  the graph $ p $-{L}aplacian.
\newblock {\em Journal of Spectral Theory}, 8(3):883--908, 2018.

\bibitem{Tudisco_core_periphery}
F.~Tudisco and D.~J. Higham.
\newblock A nonlinear spectral method for core--periphery detection in
  networks.
\newblock {\em SIAM Journal on Mathematics of Data Science}, 1(2):269--292,
  2019.

\bibitem{Tudisco_core_periphery_hypergraphs}
F.~Tudisco and D.~J. Higham.
\newblock Core-periphery detection in hypergraphs.
\newblock {\em SIAM Journal on Mathematics of Data Science}, 5(1):1--21, 2023.

\bibitem{tudisco2018community}
F.~Tudisco, P.~Mercado, and M.~Hein.
\newblock Community detection in networks via nonlinear modularity
  eigenvectors.
\newblock {\em SIAM J. Applied Mathematics}, 78:2393--2419, 2018.

\bibitem{tudisco2022nonlinear}
F.~Tudisco and D.~Zhang.
\newblock Nonlinear spectral duality.
\newblock {\em Annali della Scuola Normale Superiore di Pisa, Classe di
  Scienze}, 2026.

\bibitem{upadhyaya2021self}
P.~Upadhyaya, E.~Jarlebring, and F.~Tudisco.
\newblock The self-consistent field iteration for p-spectral clustering.
\newblock {\em arXiv preprint arXiv:2111.09750}, 2021.

\bibitem{West}
D.~B. West.
\newblock {\em Introduction to graph theory}.
\newblock Prentice Hall, Inc., Upper Saddle River, NJ, 1996.

\bibitem{Xu}
H.~Xu and S.-T. Yau.
\newblock Nodal domain and eigenvalue multiplicity of graphs.
\newblock {\em Journal of Combinatorics}, 3(4):609--622, 2012.

\bibitem{yang1954theorems}
C.-T. Yang.
\newblock On theorems of borsuk-ulam, kakutani-yamabe-yujob{\^o} and dyson, i.
\newblock {\em Annals of Mathematics}, 60(2):262--282, 1954.

\bibitem{Yao07}
X.~Yao and J.~Zhou.
\newblock Numerical methods for computing nonlinear eigenpairs. {I}.
  {I}so-homogeneous cases.
\newblock {\em SIAM J. Sci. Comput.}, 29(4):1355--1374, 2007.

\bibitem{Zelazo2007}
D.~Zelazo, A.~Rahmani, and M.~Mesbahi.
\newblock Agreement via the edge laplacian.
\newblock In {\em 2007 46th IEEE Conference on Decision and Control}, pages
  2309--2314, 2007.

\bibitem{ZHANG_top_mult}
D.~Zhang.
\newblock Topological multiplicity of the maximum eigenvalue of graph
  1-laplacian.
\newblock {\em Discrete Mathematics}, 341(1):25--32, 2018.

\bibitem{zhang2021homological}
D.~Zhang.
\newblock Homological eigenvalues of graph {$p$}-{L}aplacians.
\newblock {\em J. Topol. Anal.}, 17(2):555--606, 2025.

\end{thebibliography}

\appendix

\section{Other classes of variational eigenvalues}\label{appendix_variational_spectrum}

This section is devoted to discuss other classes of variational eigenvalues that have been introduced in the literature besides the Kransoselskii spectrum.

\subsection{Drábek Spectrum}
Drábek, in \cite{drabek2012variational, drabek2010variational}, introduced a family of variational eigenvalues, denoted by $\{\lambda_k^D\}_k$, defined as follows:
\begin{equation}\label{Drabek_variational_eigenvalues}
\lambda_k^D := \underset{A \in \Fc_k^D}{\inf} \underset{f \in A}{\sup} \mathcal{R}_{p}^p(f),
\end{equation}
where 
\begin{equation}
\Fc_k^D := \{ A \subset S_p \;|\; \exists h \in C(S^{k-1}_2, A), \; h \text{ is odd and surjective} \}.
\end{equation}
Since $\Fc_k^D \subset \Fc_k$, meaning that any subset $A$ in $\Fc_k^D$ has genus at least $k$, it follows that $\lambda_k \leq \lambda_k^D$. 
In \cite{drabek2012variational, drabek2010variational}, it is also proved that $\lambda_1 = \lambda_1^D$ and $\lambda_2 = \lambda_2^D$, with the equality being based on the characterization of $\lambda_1$ as a minimum and the characterization of $\lambda_2$ given in \eqref{Eq_lambda_2_on_curves}. However, the question of whether the higher variational eigenvalues satisfy $\lambda_k = \lambda_k^D$ for $k \geq 3$ remains an open problem.

\subsection{Yang spectrum}
A different class of variational eigenvalues is given by the Yang spectrum introduced in \cite{zhang2021homological} and based on the notion of the Yang index \cite{yang1954theorems}.

First, we give an introduction of the Yang index of a set \cite{yang1954theorems}. 
For any given  centrally symmetric set $A\subset \R^N\setminus\{\vec0\}$, the antipodal map $-$ is a  continuous involution with no fixed point, and $A$ is invariant under $-$, i.e., $-A=A$. Let $C_*(A)$ be the singular chain complex with $\mathbb{Z}_2$-coefficients, and denote by $-_\#$ the chain map of $C_*(A)$ induced
by the antipodal map $-$. 
Consider the symmetric $k$-chains that form a subgroup $C_k(A,-)$ of $C_k(A)$,  and
the boundary operator $\partial_k$ that maps  $C_k(A,-)$ to  $C_{k-1}(A,-)$. Then, these
subgroups form a subcomplex $C_*(A,-)$, and we can then define the corresponding cycles $Z_k(A,-)$, boundaries $B_k(A,-)$, and homology groups $H_k(A,-)$ in a standard way, respectively. 
Let $\tau:Z_k(A,-)\to\mathbb{Z}_2$ be homomorphisms inductively defined by
\[\tau(z)=\begin{cases}
\mathrm{Ind}(c),&\text{ if }k=0,
\\
\tau(\partial_k z),&\text{ if }k\ge1
\end{cases}
\]
if
$z=-_\#(c)+c$, where the index of a 0-chain $c=\sum n_i\sigma_i$ is  defined by $\mathrm{Ind}(c):=\sum n_i$. It is known that $\tau$ is well-defined and $\tau B_k(A,-)=0$, and thus it induces the index homomorphism $\tau_*:H_q(A,-)\to\mathbb{Z}_2$ by $\tau_*([z])=\tau(z)$ (see \cite{yang1954theorems}).

\begin{definition}[\cite{yang1954theorems}]\label{def:Yang}
The Yang index of a centrally symmetric compact set $A$ in $\mathbb{R}^N\setminus\{\vec0\}$ is defined as   
\[
\gamma_Y(A) :=
\begin{cases}
\min\limits\{k\in\mathbb{Z}^+: \tau_*H_k(A,-)=0\} & \text{if}\; A\ne\emptyset,\\
0 & \text{if}\; A=\emptyset.
\end{cases}
\]

\end{definition}

With this notion, we can now introduce the family of Yang variational eigenvalues \cite{zhang2021homological}, denoted by $\{\lambda_k^Y\}_k$, defined as follows
\begin{equation}\label{Yang_variational_eigenvalues}
\lambda_k^Y:=\underset{A\in \La_k^Y}{\inf}\underset{f\in A}{\sup}\mathcal{R}_{p}^p(f)\,,
\end{equation}
where 
\begin{equation}
\Fc_k^Y:=\{A\subset S_p|\; \gamma_Y(A)\ge k\}\,.
\end{equation}
Since $\Fc_k^D\subset \Fc_k^Y\subset \Fc_k$,  we have $\lambda_k\leq \lambda_k^Y\leq \lambda_k^D$.

\subsection{$C^2$-``spheres'' spectrum ($p>2$)}

In \cite{deidda2024_spec_energy,deidda2023PhdThesis} the authors observed that for $p>2$ the graph $p$-Laplacian eigenvalue equation can be written as a constrained linear eigenvalue problem 
\begin{equation}\label{p-lap_eq_reformulation}
\frac{1}{2}K^T(|\grad f|^{p-2}\odot\grad f)=\lambda \nodeweight \odot |f|^{p-2}\odot f
\quad \Leftrightarrow \;
 \begin{cases}
     \frac{1}{2}K^T(\edgeweight^*\odot \grad f)=\lambda \nodeweight\odot \nodeweight^*\odot f\\
\edgeweight^*_{uv}=|\grad f(u,v)|^{p-2}\quad \forall (u,v)\in\edgeset\\
\nodeweight^*_u=|f(u)|^{p-2}\qquad\forall u\in\internalnodes.
 \end{cases}   
\end{equation}
In particular, since the cardinality of the spectrum of the linear eigenvalue problem weighted in $\edgeweight^*$ and $\nodeweight^*$ is finite, the eigenvalues can be counted and assigned an index accordingly, we denote them by $\{\lambda_{k}(\edgeweight^*,\nodeweight^*)\}_{k}$. The index assigned to $\lambda$ in this way has been proved to provide information about the behaviour of the $p$-Rayleigh quotient in a neighbourhood of the eigenfunction $f$. Indeed, locally around $f$, the behaviour of $\rayl_p$ and the Rayleigh quotient, $\rayl_{2,\edgeweight^*,\nodeweight^*}$ of the linear eigenvalue problem weighted in $\edgeweight^*$ and $\nodeweight^*$ is the same, where
\begin{equation}\label{weighted_linear_rayleigh}
    \rayl_{2,[\edgeweight^*,\nodeweight^*]}(f)=\frac{\|\grad f\|_{\edgeweight^*}}{\|f\|_{\nodeweight^*}}=\frac{\Big(\frac{1}{2} \sum_{(u, v) \in \edgeset} \edgeweight^*_{uv}|\grad f(u, v)|^2\Big)^{\frac{1}{2}}}{\Big(\sum_{u \in \internalnodes}\nodeweight_u\nodeweight_u^* |f(u)|^2\Big)^{\frac{1}{2}}}.
\end{equation}
In particular, denoting  by $\morse[f](\rayl_p)$ (the Morse index of $\rayl_p$ in $f$) the dimension of the largest subspace in which the Hessian matrix of $\rayl_p$ at $f$ is negative definite~\cite{milnor2016morse}, i.e. the number of local decreasing directions of $\rayl_p$ in $f$, we have the following proposition  

\begin{proposition}[Prop 3.2 \cite{deidda2024_spec_energy}]\label{increasing_directions}
  Let $(f,\lambda)$ be a $p$-Laplacian eigenpair, with $p>2$. Consider the  weights $\nodeweight^*=|\eigenfunction|^{p-2}$, $\edgeweight^*=|\grad\eigenfunction|^{p-2}$ and assume that $(f, \lambda)=\big(f_{k}(\edgeweight^*, \nodeweight^*), \lambda_{k}(\edgeweight^*, \nodeweight^*)\big)$, with multiplicity $m$, i.e. $\lambda_{k}(\edgeweight^*, \nodeweight^*)=\dots=\lambda_{k+m-1}(\edgeweight^*, \nodeweight^*)$. Then:
  \begin{align}
    \morse[f](\rayl_p)
    &=k-1\,,\\
    \morse[f](-\rayl_p)
    &=N-k-m+1\,.
  \end{align}
  where $N$ is the dimension of the space, i.e. $N=|\internalnodes|$.
   \end{proposition}

We point out that the norm $\|f\|_{\nodeweight^*}$ in \eqref{weighted_linear_rayleigh} should be actually denoted by $\|f\|_{\nodeweight\nodeweight^*}$, however, as with the $p$-norm, we omit to report the dependence on $\nodeweight$. We also note that $\nodeweight^*$ could be zero on some nodes, meaning that, generally $\|f\|_{\nodeweight^*}$ is only a seminorm. In particular the authors of \cite{deidda2024_spec_energy} note that $\rayl_{2,[\edgeweight^*,\nodeweight^*]}$ is well defined on $(\mathrm{Ker}(\diag(\nodeweight^*))\cap \mathrm{Ker}(K^T\diag(\edgeweight^*)K))^{\perp}$ admitting values in $[0,\infty]$. In particular, given an eigenvalue $\lambda$ of the generalized linear eigenvalue problem
 \begin{equation}
     \frac{1}{2}K^T(\edgeweight^*\odot \grad f)=\lambda \nodeweight\odot \nodeweight^*\odot f,
 \end{equation}
 the authors of \cite{deidda2024_spec_energy} define the multiplicity $m$ of $\lambda$, e.g. in \cref{increasing_directions}, as $\mathrm{mult}(\lambda)=\mathrm{dim}\{f\;|\;\frac{1}{2}K^T(\edgeweight^*\odot \grad f)=\lambda \nodeweight\odot \nodeweight^*\odot f\}$. Finally we recall that $\infty$ is considered a well defined eigenvalue whenever there exists some $f\in \mathrm{Ker}(\diag(\nodeweight^*))$ but $f\not \in \mathrm{Ker}(K^T\diag(\edgeweight^*)K)$, in which case $(f,0)$ satisfies the reversed eigenvalue equation
  \begin{equation}
      \nodeweight\odot \nodeweight^*\odot f=\lambda\frac{1}{2}K^T(\edgeweight^*\odot \grad f),
 \end{equation}

In the remainder of this section we show that \cref{increasing_directions} can be used to derive new insights about the variational eigenvalues. Indeed, considering the Krasnoselskii spectrum, the $k$-th variational eigenvalue is typically defined as the $\min\max$ of $\rayl_p$ over the closed and symmetric subsets of $S_p$ that can be mapped surjectively onto a sphere of dimension at least $k-1$. 
Heuristically, we can interpret the $k$-th variational eigenvalue as a maximum point of $\rayl_p$ over a subset whose local ``dimension'' is greater than $k-1$. This suggests a comparison between the number of increasing and decreasing directions of $\rayl_p$ in the neighborhood of a variational eigenfunction and the index of the eigenvalue in the linear case (i.e., its linear index). As we show in \cref{Comparison_theorem}, this approach provides bounds for the variational index in terms of the linear index.
However, the standard definition of Krasnoselskii variational eigenpairs does not provide sufficient regularity to obtain precise information about the behavior of the $p$-Rayleigh quotient in a neighborhood of an eigenfunction. To address this limitation, we introduce a new class of variational eigenvalues defined via a min-max principle over a class of sufficiently smooth submanifolds of $S_p$. 

Observe that when $p>2$ in \Cref{Deformation_lemma_Intro},  we can assume the deformation to be a $C^2$ homeomorphisms from the sphere to the sphere, yielding the following lemma whose proof can be found in \cref{Appendix:a}

\begin{lemma}[Deformation Lemma]\label{Deformation_lemma_smooth}
Let $p>2$ and assume $c$ to be a regular value of $\rayl_{p}$. Then there exists $\epsilon>0$ and a family of $\phi\in C^2( [0,1]\times S_p,S_p)$ such that
\begin{enumerate}
    \item $\phi(t,\cdot)$ is a $C^2$ odd homeomorphism for any $t\in[0,1]$,
    \item $\phi(1,\rayl_p^{c+\epsilon})\subset \rayl_p^{c-\epsilon}$\,.
\end{enumerate} 
\end{lemma}

Given the above deformation Lemma, in the inf-sup definition of the variational eigenvalues, we can consider families of embedded submanifolds of the sphere that have some regularity. In particular we consider the family of embedded $C^2$ $k$-dimensional spheres in $S_p$ given by:
\begin{equation}
\mathcal{F}_k^S=\{A\subset S_p|\, \exists \varphi\in C^2(S^{k-1}_2,A), \varphi^{-1}\in C^2(A,S^{k-1}_2)\,\}\,,
\end{equation}
where $S^{k-1}_2$ is the $(k-1)$-dimensional sphere. Observe that the family $\mathcal{F}_k^S$ is preserved by the deformations in \Cref{Deformation_lemma_smooth} for any $k$ and thus we can define the following family of variational eigenvalues
\begin{equation}\label{submanifolds_variational_eigenvalues}
\lambda^S_k=\underset{A\in\Fc_k^S}{\inf}\underset{f\in A}{\max}\rayl_{p}^p(f).
\end{equation}
Observe that since $\Fc^S_k\subset  \Fc_k$ the eigenvalues $\lambda_k^S$ satisfy the following inequalities with respect to the Krasnoselskii variational eigenvalues:
\begin{equation}\label{variational_inequalities_1}
\lambda_k\leq \lambda_k^S\,.
\end{equation}

Now, thanks to Lemma \ref{increasing_directions}, given a variational $p$-Laplacian eigenvalue, $\lambda$, we are able to compare its position in the variational spectrum with its index as eigenvalue of the corresponding generalized weighted Laplacian eigenvalue problem \eqref{p-lap_eq_reformulation}.
\begin{theorem}\label{Comparison_theorem}
Let $(\eigenfunction,\lambda)$ be a $p$-Laplacian eigenpair such that $$\lambda=\lambda_{h}^S(\plap).$$ 
Define $\edgeweight^*=|\grad\eigenfunction|^{p-2}$, and $\nodeweight^*=|\eigenfunction|^{p-2}$ and assume $\lambda=\lambda_{h-m+1}(\edgeweight^*,\nodeweight^*)=\dots=\lambda_{h}(\edgeweight^*,\nodeweight^*)$, where $\lambda_{j}(\edgeweight^*,\nodeweight^*)$ are the generalized eigenvalues of the weighted Laplacian as given in Lemma \eqref{p-lap_eq_reformulation} and $m$ is their multiplicity.
Then 
$$k\leq h.$$
\end{theorem}
\begin{proof}
We first assume $k>h$. Then by characterization \eqref{submanifolds_variational_eigenvalues}, for any $\epsilon>0$ there exists $f_{\epsilon}\in S_p$ and a subset $A_{k-1}^{\epsilon}\subset S_p$, $C^2$-diffeomorphic to a $(k-1)$-dimensional sphere such that $$\lambda+\epsilon=\mathcal{R}_{p}(f_{\epsilon})=\max_{f\in A_{k-1}^{\epsilon}}\mathcal{R}_{p}(f)\,.$$
Then, denoting by $T_g(X)$ the tangent space to $X$ in $g$, we have that for any tangent vector $\xi\in T_{\eigenfunction_{\epsilon}}(A_{k-1}^{\epsilon})\subset T_{\eigenfunction_{\epsilon}}(S_p) $ and any curve, $\gamma(t)$, in $A_{k-1}^{\epsilon}$ such that $\gamma(0)=f_{\epsilon}$ and $\gamma'(0)=\xi$, we have 
\begin{equation}\label{Comparison_1}
\frac{d^2}{dt^2}\mathcal{R}_{p}\big(\gamma(t)\big)\Big|_{t=0}=\left\langle\xi,\frac{\partial^2}{\partial f^2}\mathcal{R}_{p}(\eigenfunction_{\epsilon})\xi\right\rangle+\left\langle \frac{\partial}{\partial f}\mathcal{R}_{p}(\eigenfunction_{\epsilon}),\gamma''(t)\right\rangle\leq 0\,.
\end{equation}
In particular, up to subsequences, we can assume that $f=\lim_{\epsilon\rightarrow 0}f_{\epsilon}$ and $\pi:=\lim_{\epsilon\rightarrow 0}T_{\eigenfunction_{\epsilon}}(A_{k-1}^{\epsilon})\subset T_{f}(S_p)$. 
Moreover, since $f$ is a $p$-Laplacian eigenfunction it satisfies $\partial/\partial f \big( \rayl_{p}(f)\big)=0$. So, for any $\xi\in \pi$, from \eqref{Comparison_1},  we have 
\begin{equation}\label{tg_derivative}
\frac{\partial^{2}}{\partial \epsilon^{2}}\bigg(\frac{\|\grad(f+\epsilon \xi)\|_p^p}{\|f+\epsilon \xi\|_p^p}\bigg)\bigg|_{\epsilon=0}=\frac{p(p-1)}{2}\frac{\partial^{2}}{\partial \epsilon^{2}}\bigg(\frac{\|\grad(f+\epsilon \xi)\|_{\edgeweight^*}^2}{\|f+\epsilon \xi\|_{\nodeweight^*}^2}\bigg)\bigg|_{\epsilon=0}\leq 0.
\end{equation}
where the equality between the second derivatives of the $p$-Rayleigh quotient and the $(\edgeweight^*,\nodeweight^*)$-Rayleigh quotient has been proved in Prop 3.3 in \cite{deidda2024_spec_energy}.
Moreover, if $\lambda=\lambda_{h}(\edgeweight^*,\nodeweight^*)$, from Lemma 3.2 in \cite{deidda2024_spec_energy} we know that for any $\xi\in T_{f}(S_{\nodeweight^*})=T_{f}(S_p)$ with 
$\xi \in \,\mathrm{span} \{f_{j}(\edgeweight^*,\nodeweight^*)\}_{j>h}$ then
   \begin{equation}\label{span_derivative}
   \frac{\partial^{2}}{\partial \epsilon^{2}}\bigg(\frac{\|\grad(f+\epsilon \xi)\|_p^p}{\|f+\epsilon \xi\|_p^p}\bigg)\bigg|_{\epsilon=0}=\frac{p(p-1)}{2}\frac{\partial^{2}}{\partial \epsilon^{2}}\bigg(\frac{\|\grad(f+\epsilon \xi)\|_{\edgeweight^*}^2}{\|f+\epsilon \xi\|_{\nodeweight^*}^2}\bigg)\bigg|_{\epsilon=0}>0
   \end{equation}
where  $S_{\nodeweight^*}=\{f\;|\;\|f\|_{\nodeweight^*}=1\}$ and $\{f_{j}(\edgeweight^*,\nodeweight^*)\}_{j>h}$ are the linear eigenvectors relative to the eigenvalues  $\{\lambda_{j}(\edgeweight^*,\nodeweight^*)\}_{j>h}$ .
Now observe that $\mathrm{dim}\big(\pi\,\cap\, T_{\eigenfunction}(S_p)\big)=k-1$ and
$\mathrm{dim}\big(\mathrm{span} \{\eigenfunction_j(\edgeweight^*,\nodeweight^*)\}_{j>h}\big)=N-h>N-k$. Thus 
$$T_{\eigenfunction}(S_p\cap A_k)\cap \mathrm{span} \{\eigenfunction_j(\edgeweight^*,\nodeweight^*)\}_{j>h}\neq \emptyset$$
which is a contradiction by \eqref{span_derivative} and \eqref{tg_derivative}.
\end{proof}

Now, we present an example of application of the above Theorem to the study of the variational spectrum when there are more eigenvalues than the dimension of the space.

\begin{remark}
   Consider the graph in Figure~\ref{Fig1.2} from \cite{DEIDDA2023_nod_dom}, with the corresponding eigenpairs of the corresponding $p$-Laplacian. 
    Since the number of eigenpairs is greater than the dimension of the space, at least one of the eigenvalues is a non-variational one. However, the 
    definition of the variational eigenvalues does not help to identify which eigenvalue is variational and which one is not. Differently, \Cref{Comparison_theorem} allows us to conclude that the eigenvalue $\lambda=2+2^{p-1}$ is not a variational eigenvalue for $p>2$. Indeed, it is easy to show that $(f,\lambda)=\Big(\big(0,1,0,-1\big),2+2^{p-1}\Big)$ is the second eigenpair of the $(\edgeweight^*,\nodeweight^*)$-weighted linear eigenvalue problem \eqref{p-lap_eq_reformulation}: 
    \begin{equation}
        \grad^T(\edgeweight^* \odot \grad f)=\lambda \nodeweight^*\odot f\,,
    \end{equation} 
    and that $\lambda$ is a simple eigenvalue, i.e.:
 \begin{equation}
     \lambda_1(\edgeweight^*,\nodeweight^*)<\lambda=\lambda_2(\edgeweight^*,\nodeweight^*)<\lambda_3(\edgeweight^*,\nodeweight^*)=\lambda_4(\edgeweight^*,\nodeweight^*)=\infty\,.
 \end{equation}

 Indeed, it is easy to check that $\lambda_1(\edgeweight^*,\nodeweight^*)=0$ is the first eigenvalue with eigenvector $f_1=(1,1,1,1)$ while the $\lambda_3(\edgeweight^*,\nodeweight^*)=\lambda_4(\edgeweight^*,\nodeweight^*)=\infty$ corresponds to the eigenspace generated by the vectors $f_3=(1,0,0,0)$ and $f_4=(0,0,1,0)$ that are in the Kernel of the matrix $\diag(\nodeweight^*)$. Here note that saying that $(f,\infty)$ is an eigenpair we mean that $(f,0)$ is a well defined eigenpair of the opposite eigenvalue problem $ 2\nodeweight^*\odot f=\lambda\grad^T(\edgeweight^* \odot \grad f)$.  
 
Thus, if $\lambda$ was a variational eigenvalue $\lambda=\lambda_k^S(\plap)$, \Cref{Comparison_theorem} would yield, $k\leq 2$. However, observe that the level set $A_\epsilon:=\{f\in S_p:\rayl_p(f)< \lambda_2+\epsilon\}$ is a connected open set and thus a path-connected open set. Then, for any $x\in A$, there exists a simple odd continuous map $\ell:[-1,1]\to A_\epsilon$ such that $\ell(-1)=x=\ell(1)$ and $\ell(0)=-x$, and thus $\ell([-1,1)])$ gives a symmetric  closed curve which can be approximated by $C^2$-closed curve. In particular, passing to the limit for $\epsilon$ that goes to zero, it follows that $\lambda_2= \lambda_2^S$ and from the characterization of the second eigenvalue in \Cref{characterization_2nd_variational_eigenvalue}, we know that $\lambda$ can not be the second variational eigenvalue. In fact, if that was the case, $\lambda$ should be the smallest eigenvalue among the eigenvalues that are strictly larger than the eigenvalue $0$. Nonetheless, the eigenvalue $2$ is always between $0$ and $\lambda$, see Fig.~\ref{Fig1.2}. Hence we can conclude that $\lambda=2+2^{p-1}$ is not a variational eigenvalue.
\end{remark}

\begin{figure}
  \centering
  \begin{tikzpicture}[inner sep=1.5mm, scale=.5, thick]

    \node (1) at (0,2) [circle,draw] {1};
    \node (2) at (3,0) [circle,draw] {2};
    \node (4) at (3,4) [circle,draw] {4};
    \node (3) at (6,2) [circle,draw] {3};

    \draw [-] (1.south) -- (2.west);
    \draw [-] (1.north) -- (4.west);
    \draw [-] (4.east) -- (3.north);
    \draw [-] (4.south) -- (2.north);
    \draw [-] (2.east) -- (3.south);
    \node at (15,2) 
    {\quad \begin{minipage}{.65\textwidth}
        \begin{enumerate} 
        \item $f=(1,1,1,1),\quad \lambda=0$
        \item $f=(1,0,-1,0),\quad \lambda=2$
        \item $f=(0,1,0,-1),\quad \lambda=2+2^{p-1}$
        \item $f=(1,0,1,-2^{\frac{1}{p-1}}),\;(1,-2^{\frac{1}{p-1}},1,0)\\[.4em]\lambda=1+\big(1+2^{\frac{1}{p-1}}\big)^{p-1}$
        \item $f=(1,-1,1,-1), \quad \lambda=2^{p}$
        \end{enumerate}
      \end{minipage}};
  \end{tikzpicture}
  \caption{ Left: Example graph in which the corresponding 
    $p$-Laplacian $\plap$ with $\edgelength_{uv}=1 \;\forall (u,v)\in\edgeset$, has more eigenvalues than the dimension of the space. Right: Set of five eigenvalues and corresponding eigenfunctions.}\label{Fig1.2}
\end{figure}

\subsection{Max-min Eigenvalues}

Next, we note that the Deformation \Cref{Deformation_lemma_Intro} can also be used to construct deformations that map $\rayl_p^{c-\epsilon}$ into $\rayl_p^{c+\epsilon}$ for any regular value $c$ and sufficiently small $\epsilon$. By reversing the roles of the two sublevel sets $\rayl_p^{c-\epsilon}$ and $\rayl_p^{c+\epsilon}$ in \Cref{Critical_point_theorem}, we see that:
\begin{equation}
    \Lambda := \sup_{A} \inf_{f\in A} \rayl_p(f)
\end{equation}
also produces critical values of $\rayl_p$.

In particular, for any family of sets $\mathcal{F}_k^*$ considered in previous sections, it is possible to define a sup-inf class of variational eigenvalues, analogous to the linear case:
\begin{equation}\label{def_max_min_eigenvalues}
      \tilde{\lambda}_k^{*} := \sup_{B \in \Fc_{N-k+1}^*} \inf_{f \in B} \rayl_{p}^p(f).
\end{equation}

We highlight that, unlike in the linear case---where all these definitions coincide---it remains an open problem whether the various families of min-max and max-min variational eigenvalues in the nonlinear setting are equal or distinct.

\begin{openproblem}\label{eq:conj-minmax}
    Let $\lambda_k(\plap)$, $\lambda_k^D(\plap)$, $\lambda_k^Y(\plap)$, $\lambda_k^S(\plap)$ (for $p>2$) be the Krasnoselskii, Drábek, Yang and $C^2$-spheres min-max variational eigenvalues defined in \eqref{Krasnoselskii_variationa_eigenvalues_Intro}, \eqref{Drabek_variational_eigenvalues}, \eqref{Yang_variational_eigenvalues}, \eqref{submanifolds_variational_eigenvalues} and the corresponding max-min eigenvalues $\tilde{\lambda}_k(\plap)$, $\tilde{\lambda}_k^D(\plap)$, $\tilde{\lambda}_k^Y(\plap)$, $\tilde{\lambda}_k^S(\plap)$ (for $p>2$) defined accordingly to \eqref{def_max_min_eigenvalues}. Are they all the same or is there any example where they are different?
\end{openproblem}

\section{Proofs of \Cref{sec:p_Lap_eigenvalue}}
\label{Appendix:a}

\subsection{Proofs of \cref{subsec:variational_spectrum}}

\begin{proof}[Proof of \cref{Deformation_lemma_smooth}]
Let $B_{\eta}(\rayl_{\plap}^{-1}(c))$ and $B_{2\eta}(\rayl_{\plap}^{-1}(c))$ be two symmetric open neighborhoods of $\rayl_{\plap}^{-1}(c)$ in $S_p$ such that $B_{\eta}(\rayl_{\plap}^{-1}(c))\subset B_{2\eta}(\rayl_{\plap}^{-1}(c))\subset S_p\setminus K$. Consider $\xi\in C^{1}(S_p,\R^+)$ a symmetric cutoff function such that 
$$\xi|_{B_\eta}=1\,,\qquad \xi|_{S_p\setminus B_{2\eta}}=0\,.$$
Then introduce the deformation function $\phi:\R\times S_p \rightarrow S_p$ that solves the following Cauchy problem
\begin{equation}
\begin{cases}
\partial t\phi(t,f)=-\xi(\phi(t,f))\partial f \rayl_{p}(\phi(t,f))\\
\phi(0,f)=f
\end{cases}
\end{equation}
Since $\partial f\big(\rayl_{p}(f)\big)$ is an odd $C^1$ function, it follows that $\partial f\big(\rayl_{p}(f)\big)\in C^1(S_p,\R^M)$ and $\partial f\big(\rayl_{p}(f)\big)=-\partial f\big(\rayl_{p}(-f)\big)$, for any $t\in[-1,1]$. Hence $\phi(t,f)$ is a $C^2$ odd function, $\phi(t,\cdot)\in C^2(S_p\cap \Dc_0,S_p\cap \Dc_0)$ and $\phi(t,f)=-\phi(t,-f)$.
Finally observe that
 $$\frac{\partial}{\partial t}(\rayl_{p}(\phi(f,t))=-\xi(\phi(f,t))\Big\|\frac{\partial}{\partial f} \rayl_{p}(\phi(f,t)\Big\|^2\leq 0\,.$$
In addition, since $c$ is a regular value and $\rayl_{p}^{-1}(c)$ is a compact set, $\|\partial/\partial f\big(\rayl_{p}(f)\big)\|$ admits a minimum greater than zero on $\rayl_{p}^{-1}(c)$, thus for any $t>0$ there exists an $\epsilon>0$ such that 
$$\rayl_{p}(\phi(f,t))<c-\epsilon \quad \forall f \in \rayl_{p}^{-1}(c)\,.$$
To conclude observe that by the continuity of $\phi$, there exists a neighborhood $U_{\epsilon}$ (that we can assume w.l.o.g to be $\rayl_p^{c+\epsilon}$ of $\rayl_{p}^{-1}(c)$ such that 
\begin{equation}
\rayl_{p}(\phi(f,t))<c-\epsilon \quad \forall f \in U_{\epsilon}\,.
\end{equation}
\end{proof}

\subsection{Proofs of \cref{subsec:infinity_1_Notation}}

\begin{proof}[Proof of \Cref{Thm_subgradient_chain_rule}]
The corollary is equivalent to state the following equality between subsets of $\R^N$:
$$\partial \|f\|_{\edgeweight,p}=\nodeweight^{-1}\odot\Big( A^T \Big(\edgeweight\odot \partial_{Af}(\|Af\|_{\edgeweight,p})\Big)\Big)\,.$$

To prove it,  note that if $\Xi\in \partial_{G=Af} \|G\|_{\edgeweight,p}$, then by definition    
    \begin{equation}
        \|Ag\|_{\edgeweight,p}-\|Af\|_{\edgeweight,p}\geq \langle \Xi, Ag-Af\rangle_{\edgeweight}=\big\langle \nodeweight^{-1}\odot \big(A^T \big(\edgeweight\odot \Xi\big)\big), g-f\big\rangle_{\nodeweight} \qquad \forall g.
    \end{equation}
    thus $\nodeweight^{-1}\odot\big( A^T\big( \edgelength\odot \Xi\big)\big)\in \partial \Phi(f)$. 
    On the other hand, by convex duality we can write:
    \begin{equation}
        \Phi(f)=\|Af\|_{\edgeweight,p}=\sup_\xi\big(\langle f, \xi\rangle_{\nodeweight}-\chi_C(\xi)\big)=\chi_C^*(f)
    \end{equation}
    where $C=\{\nodeweight^{-1}\odot\big(A^T\big(\edgeweight\odot\Xi)\big)|\;\|\Xi\|_{\edgeweight,q}\leq 1\}$ and $\chi_c$ is the convex characteristic function which takes values $0$ and $\infty$. 
    
    In particular by convex duality we have 
        \begin{equation}
       \sup_g\langle \xi, g\rangle_{\nodeweight}-\Phi(g) = \Phi^*(\xi)=\chi_c^{**}(\xi)=\chi_c(\xi)=\begin{cases}
           0 \qquad &\text{if } \xi\in C\\
           \infty \qquad &\text{otherwise}
       \end{cases}
    \end{equation}
This, as a corollary of \Cref{Lemma_norm_subgradient}, proves that $\xi\in \partial\Phi(f)$ if and only if $\xi \in C$. In particular, since for $\xi\in \partial\Phi(f)$ we have $\xi=\nodeweight^{-1}\odot\big(A^T \big(\edgeweight\odot\Xi\big)\big)$ and $\langle \xi, f\rangle_{\nodeweight}=\Phi(f)$, necessarily $\|\Xi\|_{\edgeweight,q}=1$ and  $\langle \Xi, G\rangle_{\edgeweight}\leq \|G\|_{\edgeweight,p} $ for any $G\in \R^M$. Moreover $\langle  \Xi, Af\rangle_{\edgeweight}=\|Af\|_{\edgeweight,p}$ and thus from \Cref{Lemma_norm_subgradient} we have $\Xi \in \partial_{G=Af}\|G\|_{\edgeweight,p}$, which concludes the proof.

\end{proof}


\section{Proofs of \Cref{Sec:Duality}}\label{appendix_duality}

\begin{proof}[Proof of \cref{Lemma_duality}]
    As in the proof of Lemma 2.6 of \cite{zhang2021discrete} we observe \textbf{1)} that if $\xi \in \partial \Psi(f)$ with $f\neq 0$, then $\Psi^*(\xi)=1$; \textbf{2)} that if $\Psi^*(\xi)=\Psi(x)=1$ then $\xi\in \partial\Psi(f)$ if and only if $f\in \partial\Psi^*(\xi)$. 
    Now if $(f,\Lambda)$ is a critical pair of $\Phi(Af)/\Psi(f)$, there exist $\Xi\in \partial_{G=Af}\Phi(G)$ and $\xi\in \partial_f\Psi(f)$ such that 
    \begin{equation}
        \nodeweight^{-1}\odot A^T (\edgeweight\odot \Xi)=\Lambda \xi
    \end{equation}
    where we have used the characterization of the subgradients of $f\rightarrow \Phi(Af)$ \Cref{Thm_subgradient_chain_rule}. From 1) we have $\Psi^*(\xi)=1$ and $\Phi^*(\Xi)=1$. Moreover assuming w.l.o.g $\Psi(f)=1$ then  $\Phi(Af)=\Lambda$. In particular since the subgradient is scale invariant we have $\Xi\in \partial_{G=Af/\Lambda}\Phi(G)$ which by 2) yields $Af/\Lambda\in \partial_{\Theta=\Xi}\Phi^*(\Theta)$ and similarly $f\in \partial_{\xi}\Psi^*(\xi)$. In particular, from the generalized eigenvalue equation above we have  
    \begin{equation}
        f\in \partial_{\eta=\nodeweight^{-1}\odot (A^T(\edgeweight\odot \Xi))}\Psi^*(\eta).
    \end{equation}
     Now note that by \Cref{Thm_subgradient_chain_rule}
     \begin{equation}
     \begin{aligned}
     \partial_{\Theta=\Xi}\Psi^*(\nodeweight^{-1}\odot( A^T(\edgeweight\odot \Theta)))
     &=\edgeweight^{-1} \edgeweight A \nodeweight^{-1} \nodeweight \partial_{\eta=\nodeweight^{-1}\odot( A^T(\edgeweight\odot \Xi
     ))}\Psi^*(\eta)=\\
     &= A \partial_{\eta=\nodeweight^{-1}\odot( A^T(\edgeweight\odot \Xi
     ))}\Psi^*(\eta).
     \end{aligned}
     \end{equation}
     we have proved that 
     \begin{equation}
     Af \in \partial_{\Theta=\Xi}\Psi^*(\nodeweight^{-1}\odot( A^T(\edgeweight\odot \Theta))) \qquad Af/\Lambda \in \partial_{\Theta=\Xi}\Phi^*(\Theta),
     \end{equation}
    thus $(\Xi,\Lambda)$ is a generalized critical pair of $\Psi^*(\nodeweight^{-1}\odot(A^T(\edgeweight\odot \Theta)))/\Phi^*(\Theta)$.
\end{proof}


\section{Proofs of \Cref{Sec:1-Lapl_spectrum}}

\begin{proof}[Proof of \Cref{thm:1-eigenpair_and_cheeger_constants}]
   The inequality $h_{\SNc(f_k)}(\Gc)\leq\Lambda_k(\Delta_1)$ is a consequence of \Cref{lemma:1-eigenvalues_are_isoperimetric_constants} i.e. $\Lambda_k(\Delta_1)=c(A_i)$ for any $A_i$ nodal domain of $f_k$. Indeed, the family of the nodal domains, $A_1,\dots,A_{S(f_k)}$, induced by $f_k$ is a suitable family in the definition of $h_{S(f_K)}(\Gc)$ and hence we can derive the lower bound in the thesis 
    \begin{equation}
     h_{\SNc(f_k)}(\Gc)\leq \max_{i=1,\dots,\SNc(f_k)} c(A_i)=\Lambda_k(\Delta_1)\,.
    \end{equation}
    Next we prove the upper bound of the thesis. Consider a maximizing family, $\{A_i\}_{i=1}^k$, in the definition of $h_k(\Gc)$ and, for any subset $A_i$, consider the corresponding characteristic function 
    \begin{equation}
    \rchi_{A_i}(u)=\begin{cases}
        1\quad if \;u\in A_i\\
        0\quad otherwise
    \end{cases}\,.
    \end{equation}
Then, we define $\pi=\mathrm{span}\{\rchi_{A_i}\}_{i=1}^k$ and, since the subsets $A_i$ are pairwise disjoint, we note that the Krasnoselskii index of $\pi$, which matches the dimension of $\pi$, is equal to $k$. 
    Thus, we can state that $\Lambda_k(\Delta_1)\leq \max_{f\in \pi}\rayl_1(f)$\, and, to conclude, it is sufficient to prove $\rayl_1(f)\leq h_k(\Gc)$\,:
    \begin{equation}
        \begin{aligned}
\rayl_1(\sum_{i}\alpha_i\rchi_{A_i})=\frac{\|\sum_i \alpha_i \grad \rchi_{A_i}\|_1}{\|\sum_i \alpha_i \rchi_{A_i}\|_1}\leq \frac{\sum_i |\alpha_i|\|\grad \rchi_{A_i}\|_1}{\sum_i |\alpha_i|\|\rchi_{A_i}\|_1}\leq \max_i \frac{\|\grad \rchi_{A_i}\|}{\|\rchi_{A_i}\|_1}=h_k(\Gc)\,,
        \end{aligned}
    \end{equation}
where, since the $A_i$ are pairwise disjoint, we have used the relations $\|\sum_i \alpha_i \grad \rchi_{A_i}\|_1 \leq \sum_i |\alpha_i|\|\grad \rchi_{A_i}\|_1$,  $\|\sum_i \alpha_i \rchi_{A_i}\|_1=\sum_i |\alpha_i|\|\rchi_{A_i}\|_1$, the inequality $\big(\sum_i \alpha_i\big)/\big(\sum_i \beta_i\big) \leq \max_i \big(\alpha_i/\beta_i\big)$ and the equality \eqref{eq-cheeger_constant_and_characteristic_functn}.   
To conclude, we only miss to prove the equalities $h_1(\Gc)=\Lambda_1(\Delta_1)$ and
$h_2(\Gc)=\Lambda_2(\Delta_1).$ The first equality follows from the upper and lower bounds, observing that $\SNc(f_1)\geq 1$\,
To prove the second equality we need some extra work. Let $\Gamma\subset S_1$ be a subset of genus greater than $2$, such that 
\begin{equation}
    \Lambda_2(\Delta_1)=\max_{f\in\Gamma} \rayl_1(f)\,.
\end{equation}
Since $\Gamma$ is a closed and symmetric subset of $S_1$ with genus greater then $2$ it has to be connected and thus, necessarily, contains a closed and connected curve. Next, on $\Gamma$ we consider the following continuous function :
\begin{equation}
    \begin{aligned}
    \Psi:\Gamma&\longrightarrow \R\\
    f&\mapsto \rayl_1(f^+)-\rayl_1(f^-)\,,
    \end{aligned}
\end{equation}
where $f^+$ and $f^-$ are, respectively, the positive and negative part of $f$. Since $\Gamma$ contains a symmetric closed curve, the function $\Psi$ admits at least one zero $f_{\Gamma}$\,. Moreover, since 
\begin{equation}
    \min \{\rayl_1(f^+),\rayl_1(f^-)\}\leq \rayl_1(f)\leq \max \{\rayl_1(f^+),\rayl_1(f^-)\}\,,
\end{equation}
$f_\Gamma$, necessarily, satisfies the following equality:
\begin{equation}
\rayl_1(f_{\Gamma})=\rayl_1(f_{\Gamma}^-)=\rayl_1(f_{\Gamma}^+)\,.
\end{equation}
Next, we can use an argument from \cite{Hua, Tudisco1}, to claim that there exist two subsets of nodes $A^+\subset \{f_{\Gamma}^+>0\}$
and $A^-\subset \{f_{\Gamma}^->0\}$ such that 
$\rayl_1(f^+_{\Gamma})\geq c(A^+)$ and $\rayl_1(f^-_{\Gamma})\geq c(A^-)$\,.
Note that the claim implies that 
\begin{equation}
\begin{aligned}
    \Lambda_2(\Delta_1)=\max_{f\in \Gamma}\rayl_1(f)\geq \rayl_1(f_{\Gamma})&=\{\rayl_1(f_{\Gamma}^+),\rayl_1(f_{\Gamma}^-)\}
    \\&\geq \max\{c(A^+),c(A^-)\}\geq h_2(\Gc),
\end{aligned}
\end{equation}
i.e. the missing lower bound.
In particular the claim is proved in \Cref{remark_cheeger_constants_and_support} as a consequence of the first part of the theorem.

\end{proof}

\section{$p$-Laplacian spectrum on complete graphs}

We shall compute the $p$-Laplacian eigenvalues of unweighted complete graphs without boundary. This was done in Amghibech's work  \cite{Amghibech1}, but the original computation  contains a small mistake. Here we present a corrected proof. 

In this Appendix, we use $\Delta_p$ to denote the unnormalized $p$-Laplacian (i.e., $w_{uv}=1$ and $\mu_v=1$ for all edges and vertices) on a complete unweighted graph. 

\begin{theorem}\label{thm:complete}
Let $\Gc$ be a complete unweighted graph with $\nodeset=\{1,\cdots,N\}$. Then the nonzero eigenvalues of $\Delta_p$ are $N-\alpha-\beta+(\alpha^{\frac{1}{p-1}}+\beta^{\frac{1}{p-1}})^{p-1}$, where $\alpha,\beta\in\mathbb{Z}_+$ with $\alpha+\beta\le N$.
\end{theorem}

\begin{statement}\label{st:a1A-b1B}
Under the above setting, if $A$ and $B$ are two disjoint nonempty subsets of $\nodeset$, then there exist $a>0$ and $b>0$ such that $a\vec\chi_A-b\vec\chi_B$ is an eigenfunction corresponding to the  eigenvalue $N-|A|-|B|+(|A|^{\frac{1}{p-1}}+|B|^{\frac{1}{p-1}})^{p-1}$ of $\Delta_p$, where $\vec\chi_A$ stands for the characteristic function of $A$.
\end{statement}

\begin{proof}
Note that if $a\vec\chi_A-b\vec\chi_B$ is an eigenfunction of $\Delta_p$, the corresponding $p$-Laplacian eigenequation 
can be written in component-wise form as
$$\begin{cases}
|B|(a+b)^{p-1}+(N-|A|-|B|)a^{p-1}=\lambda a^{p-1}\\
|B|b^{p-1}-|A|a^{p-1}=0\\
|A|(a+b)^{p-1}+(N-|A|-|B|)b^{p-1}=\lambda b^{p-1}.
\end{cases}$$
By solving the equations above, we determine $a,b,\lambda$ as follows:  $a=c|B|^{\frac{1}{p-1}}$, $b=c|A|^{\frac{1}{p-1}}$ and $\lambda=N-|A|-|B|+(|A|^{\frac{1}{p-1}}+|B|^{\frac{1}{p-1}})^{p-1}$, where $c\ne0$. 
\end{proof}

\begin{statement}\label{st:eigen-aAbB}
If $f$ is an  
eigenfunction corresponding to a positive eigenvalue of $\Delta_p$, then there exist disjoint nonempty subsets $A,B\subset \nodeset$, and $a,b>0$, such that $f=a\vec\chi_A-b\vec\chi_B$.
\end{statement}

\begin{proof}
Let $(f,\lambda)$ be an eigenpair with $\lambda>0$. Without out loss of generality, we may assume  $f(1)\le\cdots\le f(N)$ due to the complete symmetry of the graph $\Gc$.

Clearly, there exist $1\le\alpha<\beta\le N$ such that
$$f(1)\le\cdots\le f(\alpha)<0=f(\alpha+1)=\cdots=f(\beta-1)=0<f(\beta)\le\cdots\le f(N).$$
Then, the eigenequation of $\Delta_p$ becomes
\begin{equation}\label{eq:complete:eigen-system}
\begin{cases}
\lambda=-\displaystyle{\sum_{k=1}^i\left(-1+\frac{f(k)}{f(i)}\right)^{p-1}+\sum_{k=i}^N\left(1-\frac{f(k)}{f(i)}\right)^{p-1}}&\text{ for }i=1,\cdots,\alpha\\[1em]
\lambda=\displaystyle{\sum_{k=1}^i\left(1-\frac{f(k)}{f(i)}\right)^{p-1}-\sum_{k=i}^N\left(-1+\frac{f(k)}{f(i)}\right)^{p-1}}&\text{ for }i=\beta,\cdots,N\\[1.5em]
0=\sum_{k=1}^\alpha\left(-f(k)\right)^{p-1}-\sum_{k=\alpha+1}^N f(k)^{p-1}& 
\end{cases}    
\end{equation}
We shall prove that $f(1)=\cdots=f(\alpha)$, and we split  the rest of the 
proof into two cases:
\begin{itemize}
 \item $\mathbf{p>2}$:
 In this case, we use 
 the above equations to express $\lambda$ at  $i=\alpha$ and $i=1$. 
 %
     \begin{align}
\lambda&\xlongequal{i=1\text{ in }\eqref{eq:complete:eigen-system}}\sum_{k=1}^N\left(1-\frac{f(k)}{f(1)}\right)^{p-1} \notag
\\
&\le \sum_{k=1}^\alpha\left(\frac{f(k)}{f(\alpha)}-\frac{f(k)}{f(1)}\right)^{p-1}+\sum_{k=\alpha+1}^N\left(1-\frac{f(k)}{f(1)}\right)^{p-1} \label{eq:alpha-1}\\
&= \sum_{k=\alpha+1}^N\left(\frac{f(k)}{f(1)}-\frac{f(k)}{f(\alpha)}\right)^{p-1}+\sum_{k=\alpha+1}^N\left(1-\frac{f(k)}{f(1)}\right)^{p-1} \label{eq:alpha-2}
\\
&\le \sum_{k=\alpha+1}^N\left(1-\frac{f(k)}{f(\alpha)}\right)^{p-1}\label{eq:alpha-3}
\\
&\le -\sum_{k=1}^\alpha\left(-1+\frac{f(k)}{f(\alpha)}\right)^{p-1}+\sum_{k=\alpha+1}^N\left(1-\frac{f(k)}{f(\alpha)}\right)^{p-1}\xlongequal{i=\alpha\text{ in }\eqref{eq:complete:eigen-system}}\lambda\notag
     \end{align}

 %
 where the inequality  \eqref{eq:alpha-1} is due to the fact that $1\le f(k)/f(\alpha)$, $k=1,\cdots,\alpha$, the equality \eqref{eq:alpha-2} follows from $\sum_{k=1}^\alpha(-f(k))^{p-1}=\sum_{k=\alpha+1}^N f(k)^{p-1}$, the inequality  \eqref{eq:alpha-3} is based on the elementary inequality $a^t+b^t\le (a+b)^t$ whenever $a,b>0$, $t>1$, in which we take   $a=f(k)/f(1)-f(k)/f(\alpha)$, $b=1-f(k)/f(1)$ and $t=p-1$. 
 So, all the inequalities are in fact equalities, which yield 
 \begin{equation}
     f(1)=\cdots=f(\alpha).
 \end{equation}
 \item $\mathbf{1<p<2}$:
  In this case, we  use the  equations in \eqref{eq:complete:eigen-system} to represent $\lambda$ when $i=\alpha$ and $i=N$. We shall use the elementary inequality $a^t\ge (a+b)^t-b^t$ (or equivalently, $a^t+b^t\ge (a+b)^t$) whenever $a,b>0$, $0<t<1$, in which we shall take $t=p-1$.
  %
%
%
\begin{align}
\lambda&\xlongequal{i=N\text{ in }\eqref{eq:complete:eigen-system}}\sum_{k=1}^N \left(1-\frac{f(k)}{f(N)}\right)^{p-1}\notag
\\
& \begin{aligned}
    \ge&\sum_{k=1}^\alpha\left(\frac{f(k)}{f(\alpha)}-\frac{f(k)}{f(N)}\right)^{p-1}-\sum_{k=1}^\alpha\left(-1+\frac{f(k)}{f(\alpha)}\right)^{p-1}\\
    &+\sum_{k=\alpha+1}^N\left(1-\frac{f(k)}{f(N)}\right)^{p-1}\end{aligned}\label{eq:p<2:alpha}
\\
&   \begin{aligned}
=&\sum_{k=\alpha+1}^N\!\!\left(\frac{f(k)}{f(N)}-\frac{f(k)}{f(\alpha)}\right)^{p-1}+\sum_{k=\alpha+1}^N \left(1-\frac{f(k)}{f(N)}\right)^{p-1}\\
&-\sum_{k=1}^\alpha\left(-1+\frac{f(k)}{f(\alpha)}\right)^{p-1}
\end{aligned}\label{eq:p<2:alpha2}
\\
&\ge \sum_{k=\alpha+1}^N\left(1-\frac{f(k)}{f(\alpha)}\right)^{p-1}-\sum_{k=1}^\alpha\left(-1+\frac{f(k)}{f(\alpha)}\right)^{p-1}\xlongequal{i=\alpha\text{ in }\eqref{eq:complete:eigen-system}}\lambda\label{eq:p<2:alpha3}
\end{align}
%
%
%
where the inequality \eqref{eq:p<2:alpha} comes from 
\begin{equation}
    (1-f(k)/f(N))^{p-1}\ge (f(k)/f(\alpha)-f(k)/f(N))^{p-1}-(-1+f(k)/f(\alpha))^{p-1},
\end{equation} 
which holds for any $k=1,\cdots,\alpha$, the equality \eqref{eq:p<2:alpha2} is due to the equality $\sum_{k=1}^\alpha(-f(k))^{p-1}=\sum_{k=\alpha+1}^N f(k)^{p-1}$ and the inequality \eqref{eq:p<2:alpha3} follows from 
\begin{equation}
(f(k)/f(N)-f(k)/f(\alpha))^{p-1}+(1-f(k)/f(N))^{p-1}\ge (1-f(k)/f(\alpha))^{p-1}.
\end{equation}
Again, all the inequalities are actually 
 equalities, which imply  $f(1)=\cdots=f(\alpha).$
\end{itemize}

If we change $f$ to $-f$, we obtain that $f$ is constant on $\{\beta,\beta+1,\cdots,N\}$; this achieves the proof of the remainder, that is, we have proved that $f=a\vec\chi_A-b\vec\chi_B$ for some $a,b>0$, where  $A=\{1,\cdots,\alpha\}$ and $B=\{\beta,\beta+1,\cdots,N\}$. 
\end{proof}

Combining Statements \ref{st:a1A-b1B} and \ref{st:eigen-aAbB}, we complete the proof of Theorem \ref{thm:complete}.

\end{document}